\documentclass[a4paper,12pt,twoside]{article}
\usepackage[latin1]{inputenc}                 
\usepackage{amsfonts,amssymb,amsmath,exscale} 
\usepackage[dvips]{graphicx,color,psfrag}     
\usepackage{fancyhdr}                         
\usepackage{caption}                          
\usepackage{rotating}
\usepackage{defcml}                           
\usepackage[numbers]{natbib}                  
\Floatboxname{Box}
\usepackage{verbatim}                         

\usepackage{amsmath,amsthm,amsfonts,amssymb,amscd}
\usepackage{subfig}
\usepackage{graphicx,bm,color}
\usepackage{algorithm}
        \usepackage{setspace}

\usepackage{algpseudocode}
\usepackage{stmaryrd}
\usepackage[framed,numbered,autolinebreaks,useliterate]{mcode}
\usepackage{xargs}                      
\usepackage[color=gray!20, backgroundcolor=yellow, textwidth=1.85cm]{todonotes}
\newcommandx{\ju}[2][1=]{\todo[linecolor=blue,backgroundcolor=blue!25,bordercolor=blue,#1]{#2}}

\usepackage{mathtools}

\usepackage{booktabs,ctable,multirow,longtable} 
\usepackage[latin1]{inputenc}                 
\usepackage{amsfonts,amssymb,amsmath,exscale} 
\usepackage[dvips]{graphicx,color,psfrag}     
\usepackage{fancyhdr}                         
\usepackage{caption}                          
\usepackage{rotating}
\usepackage{defcml}                           
\usepackage[numbers]{natbib}                  
\usepackage{url}
\Floatboxname{Box}
\newcolumntype{\xx}{\bm{x}}

\newcommand{\xx}{\boldsymbol{x}}
\newcommand{\R}{\mathbb{R}}
\newcommand{\sm}{\mathcal{S}\mathcal{S}_m}
\newcommand{\sx}{\mathcal{S}\mathcal{S}_\xi}

\newcommand{\cs}{\chi^\star}
\newcommand{\cn}{\chi^{j-1}}
\newcommand{\C}{\mathcal{C}}

\newcommand{\K}{\mathcal{K}}
\newcommand{\N}{\mathcal{N}}
\newcommand{\var}{\texttt}
\definecolor{gray}{gray}{0.6}

\newtheorem{Remark}{Remark}[section]
\newtheorem{theorem}{Proposition}

\newcommand*\widefbox[1]{\fbox{\hspace{2em}#1\hspace{2em}}}

\def\R{\mathbb R}
\usepackage{empheq}

\fancypagestyle{plain}{%
	\fancyhead{}%
	\fancyfoot[c]{\sffamily\thepage}%
}
\makeatletter
\def\cleardoublepage{\clearpage\if@twoside \ifodd\c@page\else
	\hbox{}
	\vspace*{\fill}
	\thispagestyle{empty}
	\newpage
	\if@twocolumn\hbox{}\newpage\fi\fi\fi}

\usepackage{scalerel,stackengine}
\stackMath
\newcommand\reallywidecheck[1]{%
	\savestack{\tmpbox}{\stretchto{%
			\scaleto{%
				\scalerel*[\widthof{\ensuremath{#1}}]{\kern-.6pt\bigwedge\kern-.6pt}%
				{\rule[-\textheight/2]{1ex}{\textheight}}
			}{\textheight}%
		}{0.5ex}}%
	\stackon[1pt]{#1}{\scalebox{-1}{\tmpbox}}%
}

\begin{document}
\Titel{
Global-local techniques for hydraulic fracture
      }
\Autor{F. Aldakheel, N. Noii, T. Wick, M. Wheeler, P. Wriggers}
\Report{02--I--17}
\Journal{

}
%



\thispagestyle{empty}
\vspace{-1cm}
\ce{\bf\large Bayesian inversion for unified ductile phase-field fracture \\[3mm]}
\vskip .35in

\ce{
	Nima Noii\(^{a,b}\), Amirreza Khodadadian\(^{b}\),  Jacinto Ulloa\(^{c}\), Fadi Aldakheel\(^{a,}\)\footnote{Corresponding author (Fadi Aldakheel).\\[3mm] 
		E-mail addresses: noii@ifam.uni-hannover.de (N. Noii); khodadadian@ifam.uni-hannover.de (A. Khodadadian);
		jacintoisrael.ulloa@kuleuven.be (J. Ulloa); 	
		aldakheel@ikm.uni-hannover.de (F. Aldakheel);  thomas.wick@ifam.uni-hannover.de (T. Wick); stijn.francois@kuleuven.be (S. Fran\c{c}ois); wriggers@ikm.uni-hannover.de (P. Wriggers).
}}\ce{
	Thomas Wick\(^{b,d}\), Stijn Fran\c{c}ois\(^{c}\), Peter Wriggers\(^{a,d}\)} \vskip .25in
\vspace{-0.2cm}
\ce{\(^a\) Institute of Continuum Mechanics} \ce{Leibniz Universit\"at Hannover, An der Universit\"at 1, 30823 Garbsen, Germany}\vskip .25in
\vspace{-0.2cm}
\ce{\(^b\) Institute of Applied Mathematics} \ce{Leibniz Universit\"at Hannover, Welfengarten 1, 30167 Hannover, Germany} \vskip .22in
\vspace{-0.2cm}
\ce{\(^c\) Department of Civil Engineering} \ce{KU Leuven, Kasteelpark Arenberg 40, 3001 Leuven, Belgium} \vskip .25in
\vspace{-0.2cm}
\ce{\(^d\) Cluster of Excellence PhoenixD (Photonics, Optics, and
	Engineering - Innovation} \ce{Across Disciplines), Leibniz Universit\"at Hannover, Germany}\vskip .25in

\begin{Abstract}
The prediction of crack initiation and propagation in ductile failure processes are challenging tasks for the design and fabrication of metallic materials and structures on a large scale. Numerical aspects of ductile failure dictate a sub-optimal calibration of  plasticity- and fracture-related parameters for a large number of material properties. These parameters enter the system of partial differential equations  as a forward model. Thus, an accurate estimation of the material parameters enables the precise determination of the material response in different stages, particularly for the post-yielding regime, where crack initiation and propagation take place. In this work, we develop a Bayesian inversion framework for ductile fracture to provide accurate knowledge regarding the \textit{effective} mechanical parameters. To this end, synthetic and experimental observations are used to estimate the posterior density of the unknowns. To model the ductile failure behavior of solid materials, we rely on the phase-field approach to fracture, for which we present a \textit{unified} formulation that allows recovering different models on a variational basis. In the variational framework, incremental minimization principles for a class of gradient-type dissipative materials are used to derive the governing equations. The overall formulation is revisited and extended to the case of anisotropic ductile fracture. Three different models are subsequently recovered by certain choices of parameters and constitutive functions, which are later assessed through Bayesian inversion techniques. A \textit{step-wise} Bayesian inversion method is proposed to determine the posterior density of the material unknowns for a ductile phase-field fracture process. To estimate the posterior density function of ductile material parameters, three common Markov chain Monte Carlo (MCMC) techniques are employed: (i) the Metropolis-Hastings algorithm, (ii) delayed-rejection adaptive Metropolis, and (iii) ensemble Kalman filter combined with MCMC. To examine the computational efficiency of the MCMC methods, we employ the $\hat{R}-convergence$ tool. The resulting framework is algorithmically described in detail and substantiated with numerical examples.
\\[2mm]
\end{Abstract} 
\vspace{-0.5cm}
{\small\tableofcontents} 
\sectpa[Section1]{Introduction} 

Fracture in the form of evolving crack surfaces in ductile solid materials exhibits dominant plastic deformation. In comparison to brittle materials, the crack evolves at a slow rate and is accompanied by a huge plastic distortion. The prediction of such failure mechanisms due to crack initiation and growth coupled with elastic-plastic deformations is an intriguingly challenging task and plays an extremely important role in various engineering applications. 
 

Recently, in the setting of continuum mechanics, a new perspective was proposed for embedding microscopic mechanisms into the macromechanical continuum formulation, based on a multi-field incremental variational framework for gradient-extended standard dissipative solids~\cite{miehe2011,miehe2013mixed}. Typical examples are theories of gradient-enhanced damage \cite{peerlings1998gradient,kiefer2018gradient,junker2021efficient,barfusz2021single}, phase-field models \cite{miehe+welschinger+hofacker10a,BourFraMar08,kuhn2010continuum}, and strain
gradient plasticity \cite{de1996some,polizzotto1998thermodynamics,liebe2001theory}. Such models incorporate non-local effects based on length scales,
which reflect properties of the material \textit{micro-structure size} with respect to the \textit{macro-structure size}. In this context, the term \textit{size effects} is used to describe the influence of the macro-structure size on the mechanical response during inelastic deformations. Thus, micro-structure interaction effects are introduced through the so-called \textit{local length-scale}, which describes the gradient information of the quantity of interest within neighboring material points (e.g., the damage or ductility zones). From a mathematical point of view, local length-scales regularize both the plastic response as well as the crack discontinuities. Hence, it resolves the loss of ellipticity of the governing equations and avoids pathological mesh-dependence in post-critical ranges, as well documented in \cite{miehe2014variational,aldakheel16,NoiiGL18}. Within the variational framework for gradient-extended dissipative phenomena, the modeling challenge is two-fold. 
\begin{itemize}
	\item First, the derivation of well-posed theoretical formulations for describing the \textit{forward model}. Hereby, variational phase-field modeling is considered, which is a regularized approach to fracture with a strong capability to simulate complicated failure processes. This includes crack initiation (also in the absence of a crack tip singularity)  \cite{KUMAR2020104027,TaTiBouMaMau17,GoeNov10}, propagation, coalescence, and branching, without additional ad-hoc criteria \cite{BourFraMar08,NoiiWick2019}. A summary of multiphysics phase-field fracture models is outlined in \cite{Wi20_book}.
	\item The second challenge is to elucidate the \textit{backward model} to estimate the model parameters and other univariate quantities of interest. A Bayesian estimation model (as an \textit{inverse model}) is here used for the ductile fracture problem to provide accurate knowledge regarding the effective mechanical parameters.	
\end{itemize}

\sectpb[Section13]{Ductile phase-field fracture as a forward model} 

A variety of studies have recently extended the phase-field approach to fracture towards the ductile case. The essential idea is to couple the evolution of the crack phase-field to an elasto-plasticity model. Initial works on this topic include~\citep{alessi2014,du2015,ambati2015,aldakheel16,alessi2017,borden2016,kuhn2016,ulloa2016} (see~\cite{alessi2018comparison}  for an overview). Phase-field models for ductile fracture were subsequently developed in the context of cohesive-frictional materials~\citep{choo2018,kienle2019}, porous plasticity~\citep{aldakheel+wriggers+miehe18} including thermal effects~\citep{dittmann2020}, the virtual element method (VEM) \cite{aldakheel2019virtual}, fiber pullout behavior \cite{storm2021comparative}, hydraulic fracture~\citep{aldakheel2020microscale,heider2020phase,aldakheel2020global}, degradation of the fracture toughness \cite{yin2020ductile}, multi-surface plasticity~\citep{fang2019} and fatigue~\citep{ulloa2021}, among others.

The majority of the ductile phase-field models found in the literature are based on local plasticity. In this setting, a strong localization of plastic strains may occur during the post-critical regime, while the damage gradient, as well as the displacement field, suffer jumps~\citep{alessi2015,tanne2017}. These occurrences are particularly relevant in the case of perfect plasticity due to the absence of a plastic regularization mechanism. Thus, from a numerical perspective, the use of local plasticity in phase-field models does not ensure mesh-objective simulations in the post-critical regime and may lead to non-realistic localized responses, such as ductile fracture with damage growth in non-plasticized regions~\citep{miehe2017phase}. To address these problems, phase-field models coupled to gradient-extended plasticity have been proposed in the literature, which incorporate a \emph{plastic internal length scale}, in the spirit of~\cite{muhlhaus1991variational}. The resulting formulation allows for a physically meaningful description of the coupled plasticity-damage evolution and mesh-objective finite element simulations. Models of this class were considered in~\citep{rodriguez2018,ulloa2021}, where a variationally consistent energetic formulation was adopted to derive the coupled system of partial differential equations (PDEs) that governs the gradient-extended elastic-plastic damage response. This approach is consistent with the models proposed in~\citep{dittmann2018} in a finite-strain setting and the extensions to micromorphic regularization~\citep{aldakheel16,miehe2016bphase,miehe2017phase}, where the governing equations were derived from rate-type variational principles, namely, the principle of virtual power. 

In this study, we present a unified formulation for ductile phase-field fracture based on variational principles, rooted in incremental energy minimization, for gradient-extended dissipative solids \cite{miehe2013mixed,miehe+hofacker+schaenzel+aldakheel15}. The coupling of plasticity to the crack phase-field is achieved by a constitutive work density function, which is characterized by a degraded stored elastic energy and the accumulated dissipated energy due to plasticity and damage. Three different models are subsequently recovered by certain choices of parameters and constitutive functions. Specifically, two phase-field models coupled to local plasticity are derived, followed by a model that considers gradient extended plasticity. The overall formulation is revisited and extended to the case of anisotropic ductile fracture. Thereby, at a specific material point, the stress state relates to the given direction (resembling solids enhanced with stiff fibers), which entails a \textit{deformation-direction-dependent} solid material. Hence, similar to \cite{noii2020adaptive}, a stiffness parameter is introduced to enforce the crack phase-field evolution according to the preferred fiber orientation.

\sectpb[Section13]{Bayesian inversion as a backward model} 

Providing reliable mechanical parameters is essential in computational mechanics to construct models with accurate predictive ability. Many of these a priori unknown parameters cannot be estimated directly through experimental procedures, and a significant effort is often needed to obtain reliable values. Furthermore, material parameters fluctuate randomly in space, giving rise to spatial uncertainty through the geometry. Consequently, the development of a sound statistical framework emerges as an interesting approach to reliably estimate mechanical properties. 

Bayesian inversion is a probabilistic technique used to identify unknown parameters. Hereby, a forward model (e.g., a system of PDEs) obtains a set of given data according to the prior density (the initial/prior information of the parameters) and gives a response related to the given unknowns. The output of the inverse problem is the posterior density, which is related to a reference observation (e.g., from experimental or synthetic measurements). This distribution provides very useful information concerning the parameter range, its standard deviation, and expectation. 

Markov Chain Monte Carlo (MCMC) methods are frequently employed to extract the posterior distribution of a parameter of interest \cite{smith2013uncertainty}. We propose several candidates according to the prior density (e.g., uniform/Gaussian) and determine whether the proposed one is rejected or accepted. The Metropolis-Hastings algorithm is one of the most common MCMC techniques due to its efficiency and easiness. In \cite{khodadadian2020bayesian},  a Bayesian framework according to this algorithm was developed to identify the material parameters in brittle fracture. The method suffers from {\it slow convergence}, since most of the candidates are rejected. To enhance its efficiency, several techniques have been recently developed, such as delayed rejection adaptive Metropolis (DRAM) \cite{haario2006dram},  differential evolution adaptive Metropolis (DREAM) \cite{laloy2012high}, and ensemble Kalman filter with MCMC (EnKF-MCMC)
\cite{emerick2012combining}. The reader is referred to \cite{khodadadian2020bayesian,adeli2020effect,adeli2020comparison,mirsian2019new,khodadadian2020bayesiann} for the application of Bayesian inversion in applied sciences.

A Bayesian approach (as a backward model) to estimate model parameters for brittle fracture in elastic solids has recently been proposed in \cite{khodadadian2020bayesian}. Specifically, the Metropolis-Hastings algorithm is devised therein to approximate model parameters based on synthetic measurements which are obtained through a sufficiently refined discretization space as the replacement of experimental observations. Thereafter, a Bayesian inversion framework towards hydraulic phase-field fracture was designed for transversely isotropic and layered orthotropic poroelastic materials \cite{noii2020bayesian}. Specifically, the DRAM algorithm was extended for parameter identification. 


In this study, we develop a Bayesian inversion framework to identify the effective mechanical parameters in ductile fracture. To this end, we first introduce a methodology to estimate the unknowns in different stages, i.e., elastic, plastic, and fracturing responses for isotropic and anisotropic materials. Three specific MCMC techniques are used to estimate the material parameters using synthetic and experimental measured data. Afterwards, a fair comparison is drawn between two improved models to determine the better convergence rate. As previously mentioned, having accurate information regarding the material parameters will enhance the accuracy of the PDE-based model. For instance, in the present case of ductile fracture, the goal is to predict the dissipative response in different stages, anticipating crack initiation and its propagation during time. We will thus employ the inferred parameters in the model equations and compare the response with the initial knowledge, highlighting the role of the probabilistic approach in improving the model's performance. 

\begin{figure}[t!]
	\subfloat{\includegraphics[clip,trim=11cm 0cm 3cm 3cm, width=8.3cm]{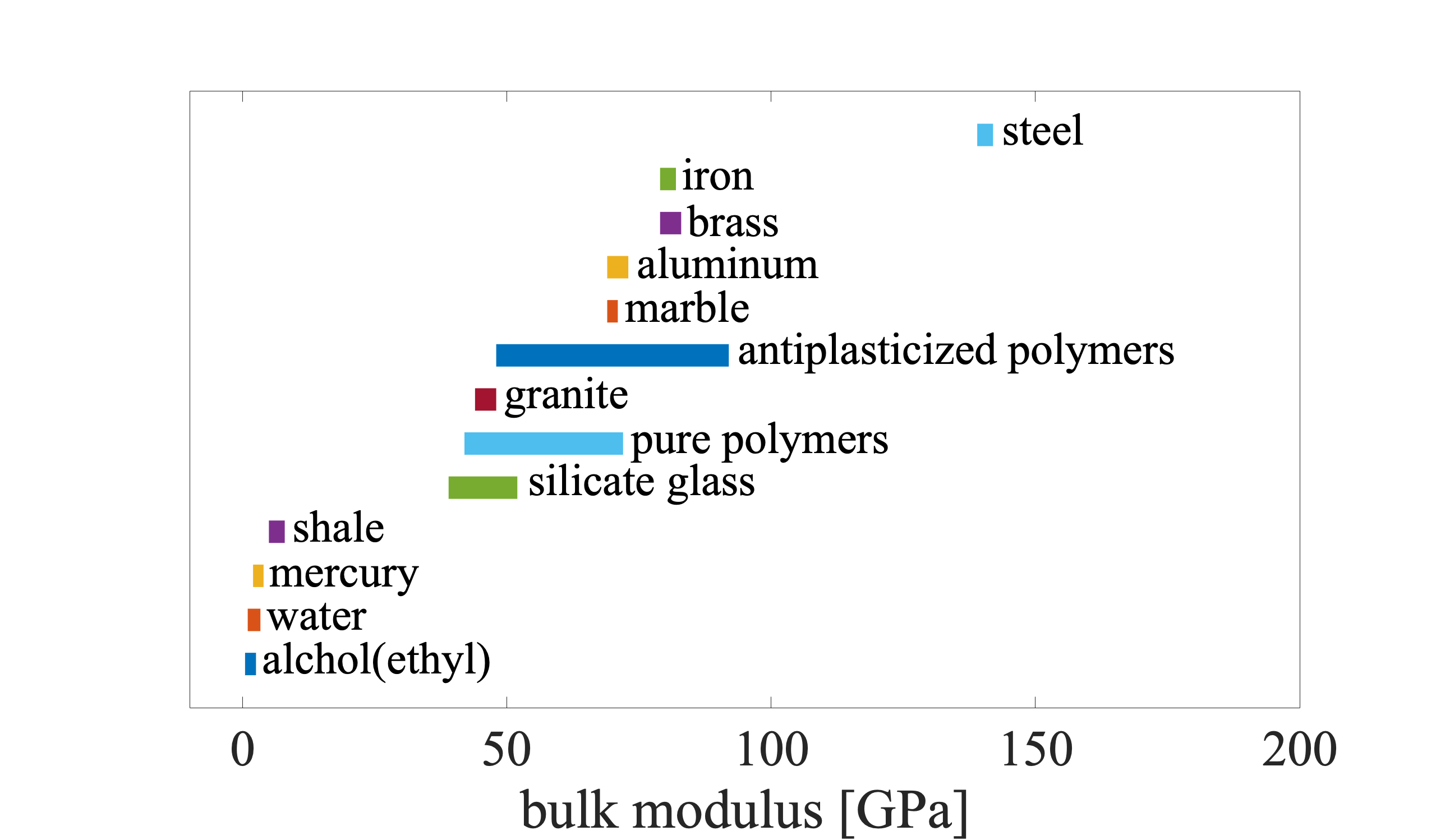}} 	\subfloat{\includegraphics[clip,trim=11cm 0cm 3cm 3cm, width=8.1cm]{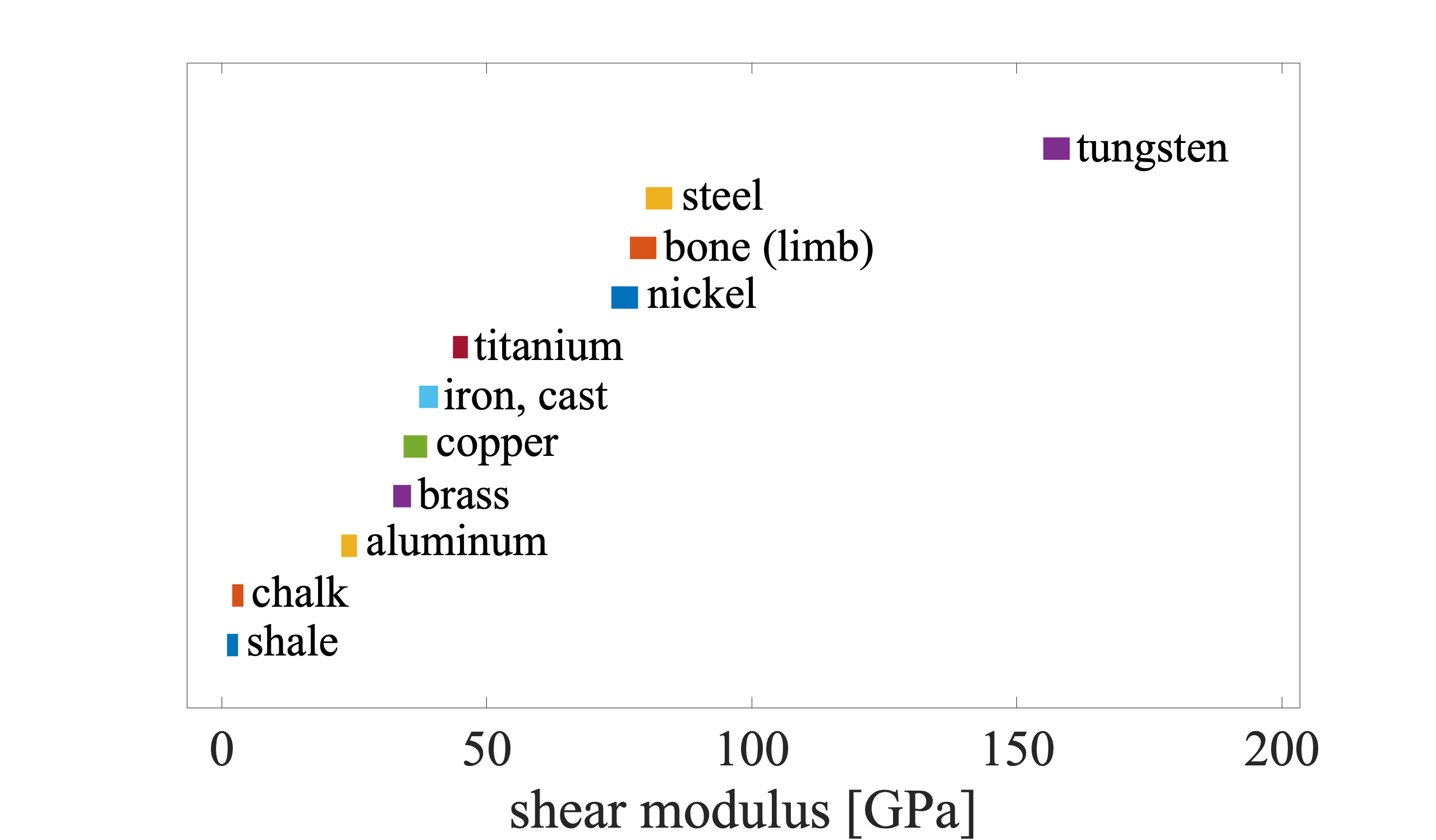}}  
	\\[5mm]
	\subfloat{\includegraphics[clip,trim=11cm 0cm 3cm 3cm, width=8.3cm]{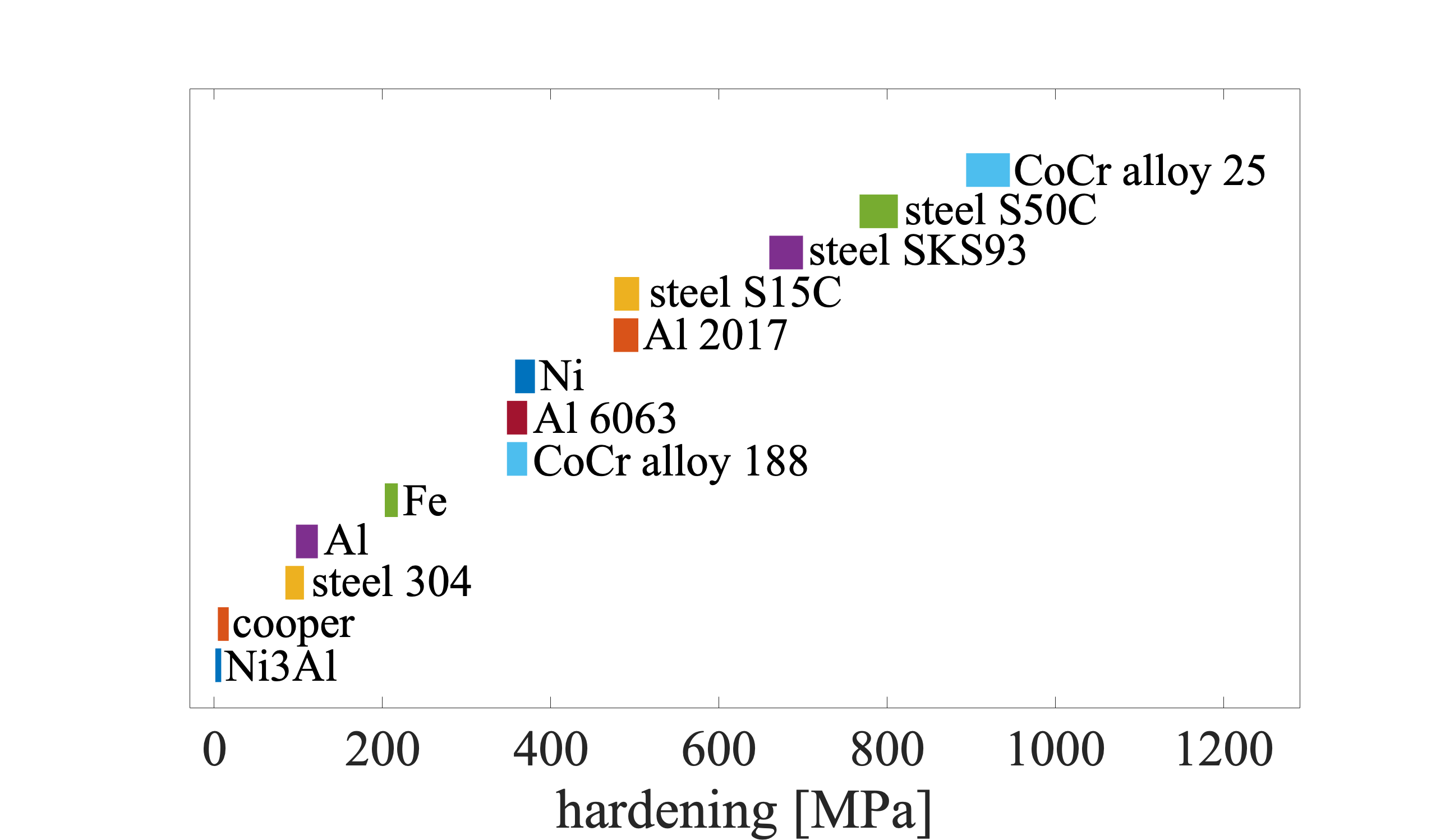}} 	
	\subfloat{\includegraphics[clip,trim=11cm 0cm 3cm 3cm, width=8.1cm]{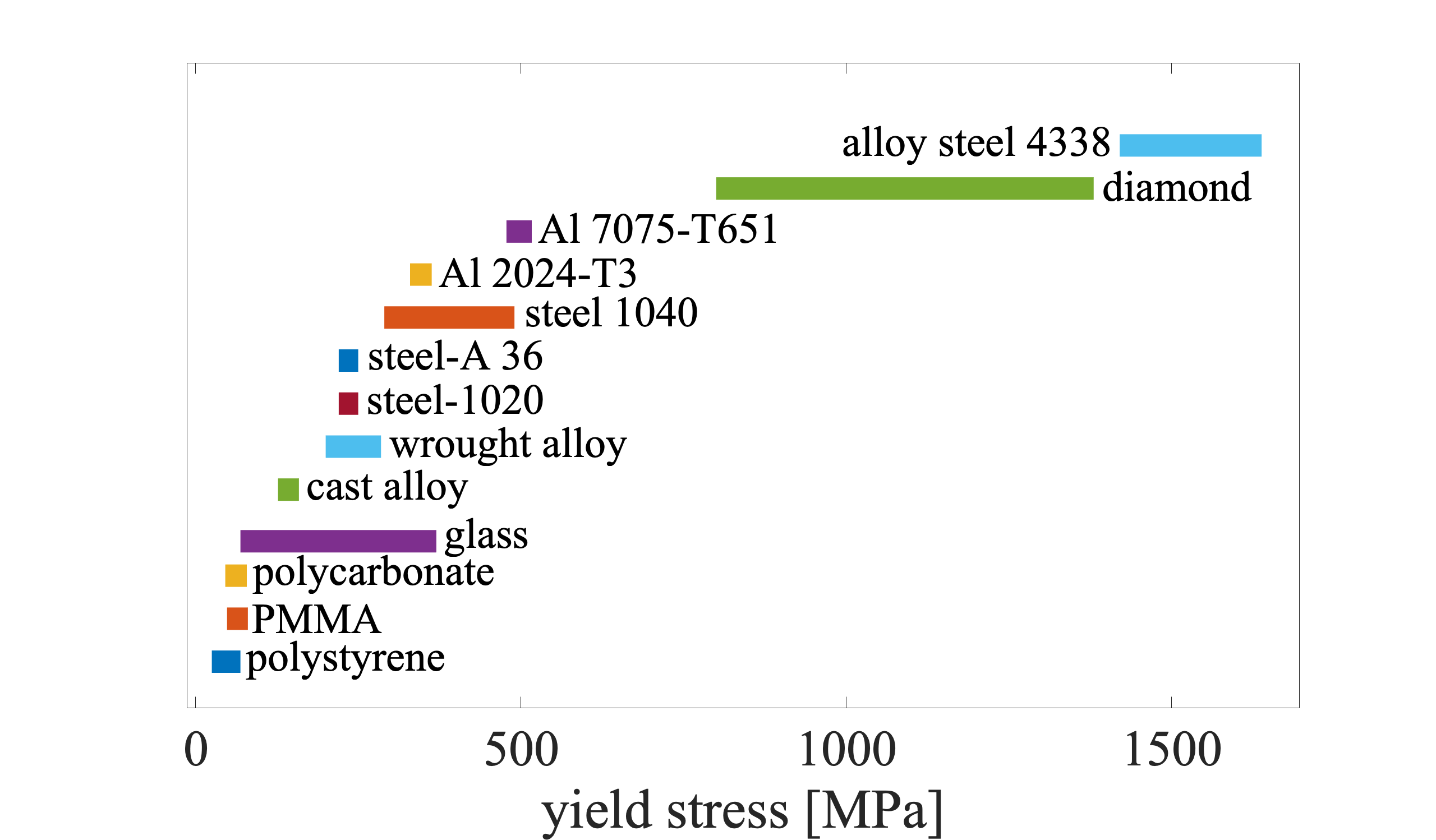}}\\[5mm]
	\subfloat{\includegraphics[clip,trim=11cm 0cm 3cm 3cm, width=8.3cm]{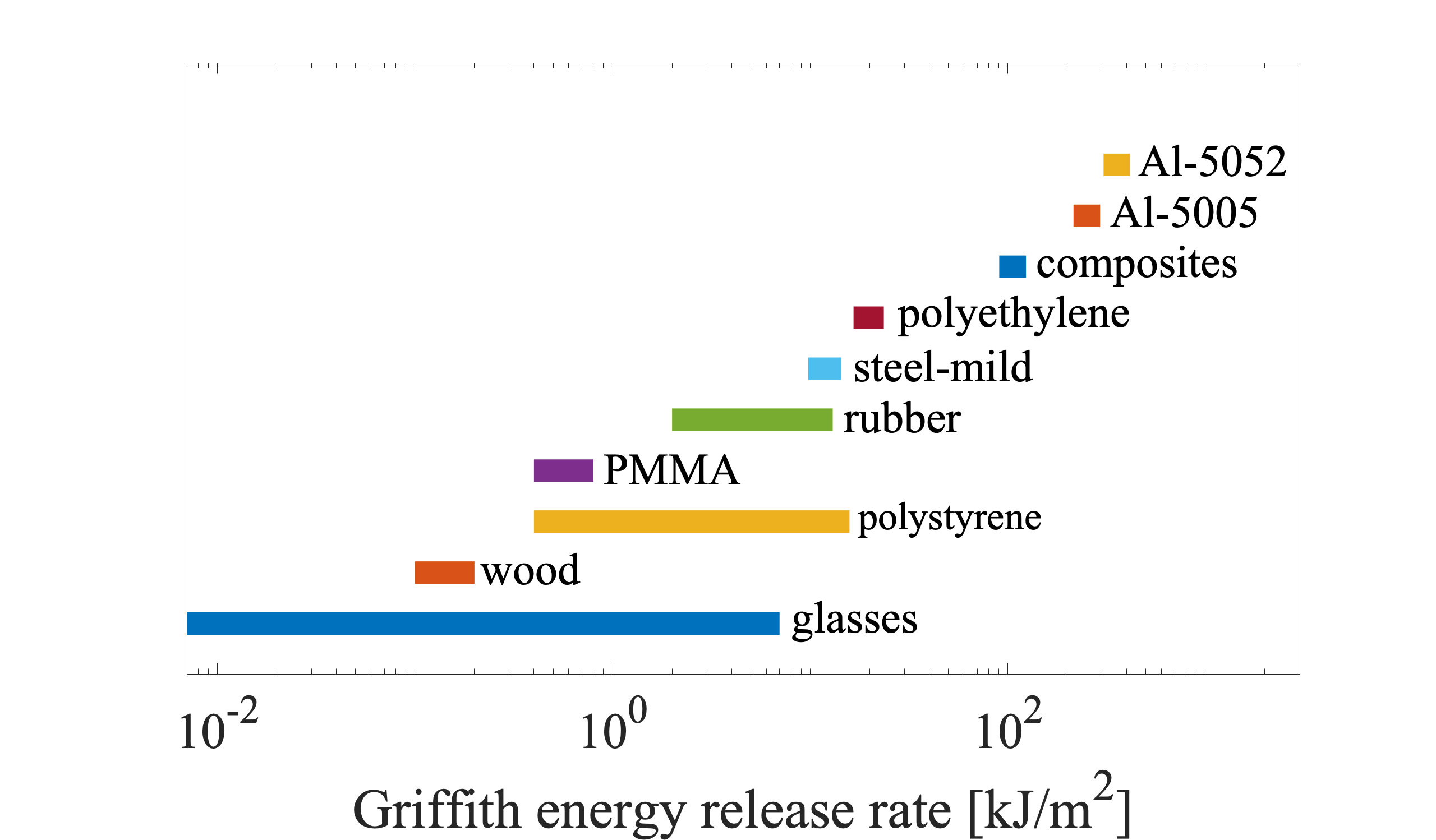}} 
	\subfloat{\includegraphics[clip,trim=11cm 0cm 3cm 3cm, width=8.4cm,height=5.425cm]{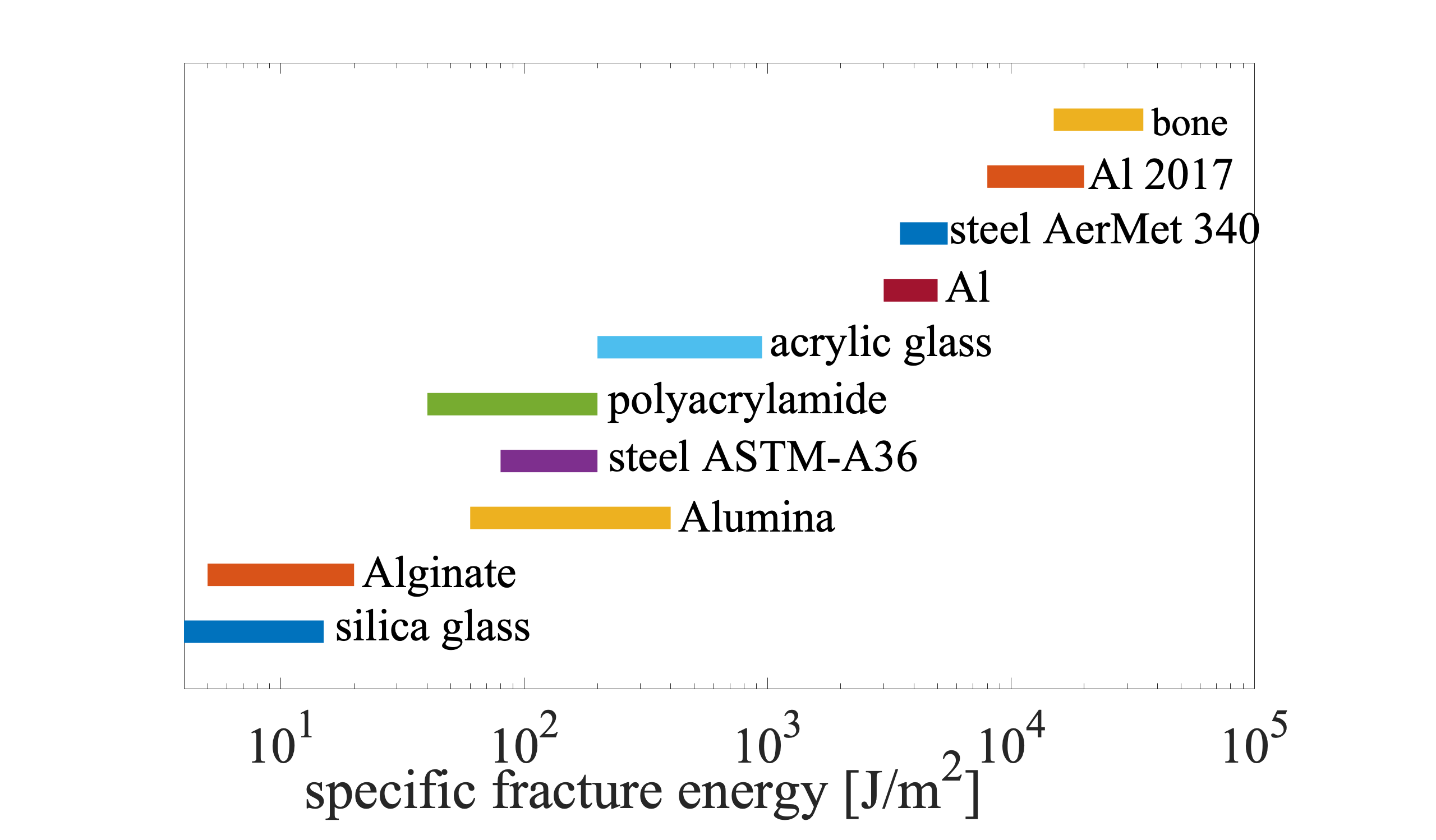}}   
	\caption{The range of different effective parameters for different materials (source: \cite{kuryaeva1997influence}--\cite{chen2017flaw}).}
	\label{material}
\end{figure}

\sectpb[Section13]{Physical interpretation of the ductile parameters} 
In ductile fracture, crack propagation is affected by several material properties. Figure \ref{material} shows the range of the different parameters for a wide variety of materials. The effective parameters required in the models are introduced below.

\begin{itemize}
	\item 
	The \textbf{bulk modulus~$K$} indicates how much the solid will compress as a result of an applied external pressure, and denotes the relation between a change in pressure and the resulting decrease/increase in fractional volume compression. See, e.g.,  \cite{kuryaeva1997influence,pariseau2017design,riggleman2010antiplasticization,giancoli2016physics}.
	
	\item The \textbf{shear modulus $\mu$} is a positive constant, smaller than $K$, 
	which indicates the response of the solid to shear stress (the ratio of shear stress to shear strain). Large shearing stresses give rise to flow and permanent deformation or fracture. See, e.g., \cite{giancoli2016physics,callister2018materials,pariseau2017design}.
	
	The elastic properties of the solid can be alternatively described in terms of the Young's modulus and the Poisson's ratio, or any other pair of Lam{\'e}'s parameters. In our previous work \cite{khodadadian2020bayesian}, it was reported that due to the boundness of the Poisson's ratio ($-1<\nu<\frac{1}{2}$) and Lam{\'e}'s first parameter ($\lambda>\frac{2\mu}{3}$), the Poisson's ratio and  Lam{\'e}'s first parameter are not appropriate for Bayesian inference. Hence, for the elasticity identification, $K$ and $\mu$ are selected, where $K>0$ and $\mu>0$. 
	
	\item The \textbf{Griffith's energy release rate~$G_c$} indicates the necessary energy (absolutely positive) to drive crack growth in elastic media. It measures the amount of energy dissipated in a localized fractured state, and therefore has units of energy per unit area. The energy rate is directly related to the toughness, indicating that in tougher materials, more energy is required to initiate fracture \cite{guo2013experimental,grif1,ambati2016}.
	
	\item The \textbf{hardening modulus $H$} characterizes the resistance of the material during plastic deformation. Hardening has an essential impact on crack initiation, specifically in phase transition zones \cite{eller2014plasticity}. See, e.g., \cite{li2017high,gomatam2006comprehensive}.
	
	\item The \textbf{yield stress $\sigma_Y$} represents the maximum stress that can be applied without exceeding a specified value of permanent strain. For a solid, it denotes the stress related to the yield point (the starting point of plasticity) where the material starts to deform in the plastic regime \cite{callister2018materials}.
	
	\item The \textbf{critical value $\alpha_{\mathrm{crit}}$} stems from a physical assumption that fracture evolution is promoted once a threshold value for the accumulated plastic strain has been reached~\cite{ambati2016}.
	
	\item The \textbf{specific fracture energy $\psi_c$} characterizes the dissipated energy during a complete damage process in a homogeneous volume element \cite{marigo2016}. This property can be interpreted as the amount of strain energy density (strain on a unit volume of material) that a given material can absorb before it fractures \cite{chen2017flaw}. A study for different materials can be found in \cite{chen2017flaw}.
\end{itemize}
In Section \ref{bayesian}, the role of these quantities in different stages of the deformation process will be clarified. 

The paper is structured as follows.
 In Section 2, a unified version for ductile phase-field fracture models is presented in a variational setting, making use of incremental energy minimization. In Section 3, we first introduce the MCMC techniques and explain the three specific Bayesian estimation methods. Afterwards, a parameter identification setting for ductile fracture is presented based on MCMC. Thereafter a review of different possibilities (specific observations) for the implementation of the Bayesian inversion framework is introduced. In Section 4, we employ the presented techniques to precisely estimate the effective parameters in different stages of the deformation process. This information allows us to enhance the accuracy of the models, as evidenced by very good agreements between simulations and experiments, highlighting the noticeable efficiency of the MCMC techniques. Furthermore, a fair comparison between the proposed models is outlined. Finally, the conclusions are drawn in Section 5.  
\sectpa[Section2]{Phase-field modeling of ductile fracture in anisotropic elastic-plastic materials}

In this section, we summarize the material models considered for the phase-field approach to ductile fracture. Three models found in the literature are revisited and extended to anisotropic fracture, considering the case of transversely isotropic  materials. To this end, a unified formulation is first provided in Sections~\ref{sec:basic_def}--\ref{sec:eq_general}. The three examined models are then recovered in Sections~\ref{sec:model1}--\ref{sec:model3}. These models will be analyzed in subsequent sections using Bayesian inversion techniques, aiming for parameter identification in anisotropic elastic-plastic fracturing materials.


\sectpb[Section21]{Basic continuum mechanics}
\label{sec:basic_def}

\begin{figure}[!t]
	\centering
	{\includegraphics[clip,trim=0.5cm 12cm 0cm 8.5cm, width=16cm]{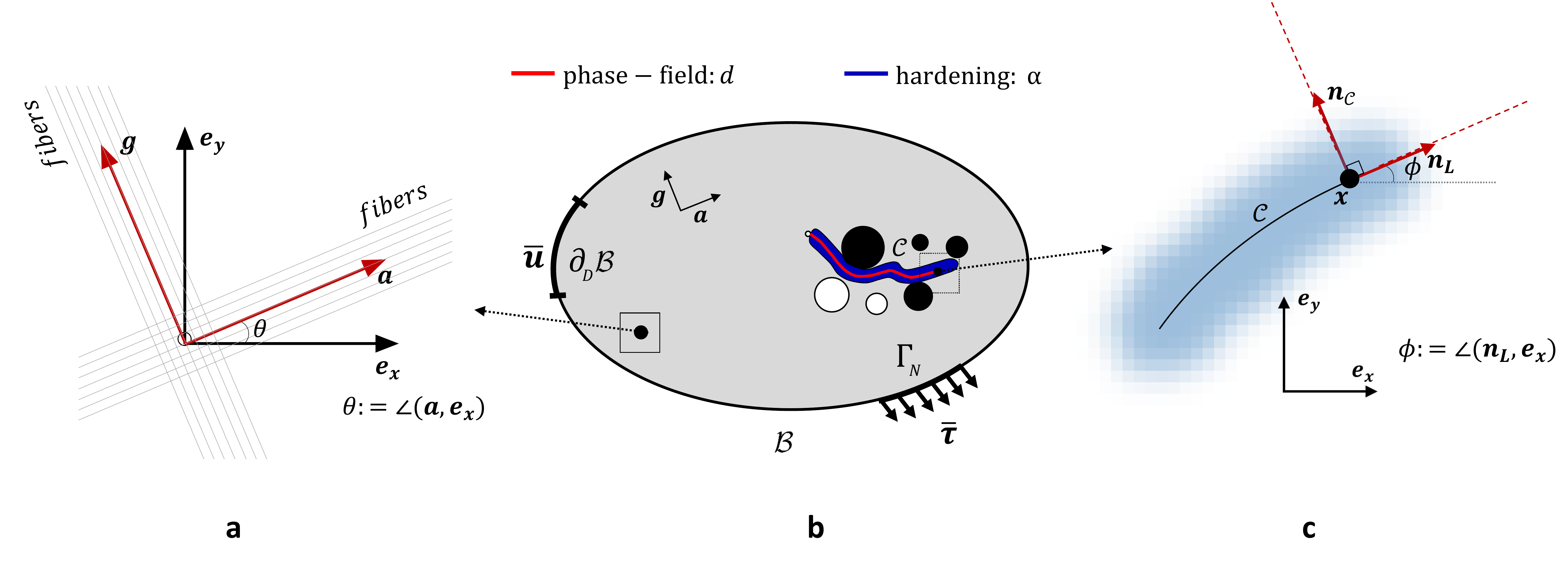}}  
	\caption{ \textcolor{black}{Problem outline and setup of the notation. (a) Global Cartesian coordinate system with unit vectors
		$(\bm{e}_{\bm{x}},\bm{e}_{\bm{y}})$ and local orthogonal principal
		material coordinates corresponding to the first and second families of fibers $(\bm{a},\bm{g})$, (b) solid with a crack inside of a plastic zone and boundary conditions, and (c) tangential and unit vectors denoted by $(\bm{n}_L, \bm{n}_\mathcal{C})$ at the crack tip point $\bm{x}$.}}
	\label{Figure1}
\end{figure}

Let $\calB\subset{\mathbb{R}}^{\delta}$ be an arbitrary solid domain, $\delta=\{2,3\}$ with a smooth  boundary $\partial\calB$ (Figure \ref{Figure1}). We assume Dirichlet boundary conditions on $\partial_D\calB $ and Neumann boundary conditions on $\partial_N \calB := \Gamma_N \cup \mathcal{C}$, where $\Gamma_N$  denotes the outer domain boundary and {$\calC\in \mathbb{R}^{\delta-1}$} is the crack boundary, as illustrated in Figure \ref{Figure1}b. 

The response of the fracturing solid at material points $\Bx\in\calB$ and time $t\in \calT = [0,T]$ is described by the displacement field $\Bu(\Bx,t)$ and the crack phase-field $d(\Bx,t)$ as
\begin{equation}
\Bu: 
\left\{
\begin{array}{ll}
\calB \times \calT \rightarrow \mathbb{R}^\delta \\[2mm]
(\Bx, t)  \mapsto \Bu(\Bx,t)
\end{array}
\right.
\AND
d: 
\left\{
\begin{array}{ll}
\calB \times \calT \rightarrow [0,1] \\[2mm]
(\Bx, t)  \mapsto d(\Bx,t)
\end{array}
\right.
\WITH
\dot{d} \ge 0 .
\label{s2-fields}
\end{equation}
Intact and fully fractured states of the material are characterized by $d(\Bx,t)=0$ and $d(\Bx,t)=1$, respectively. In order to derive the variational formulation, the following space is first defined. For an arbitrary $A\subset\mathbb{R}^\delta$, we set
\begin{align}
\mathrm{H}^1(\calB,A):=\{v:\calB\times\calT\rightarrow A\quad:\quad v\in \mathrm{H}^1(A)\}.
\end{align}
We also denote the vector valued space $\mathbf{H}^1(\calB,A):=\left[\mathrm{H}^1(\calB,A)\right]^\delta$ and define
\begin{equation}
\calW_{\overline{\Bu}}^{\Bu}:=\{\Bu\in\mathbf{H}^1(\calB,\mathbb{R}^\delta) \quad\ \colon\quad \ \Bu=\overline{\Bu} \ \text{on} \ \partial_D\calB\}.
\label{eq:spaces_u}
\end{equation}
Concerning the crack phase-field, we set
\begin{equation}
\calW^{d}:=\mathrm{H}^1(\calB) \AND  \calW^{d}_{d_n}:=\{d\in \mathrm{H}^1(\calB,{\color{black}[0,1]}) \quad \colon\quad \ d \geq d_n \},
\label{eq:spaces_d}
\end{equation}
where $d_n$ is the damage value in a previous time instant. Note that $\calW^{d}_{d_n}$ is a non-empty, closed and convex subset of $\calW^{d}$, and introduces the evolutionary character of the phase-field, incorporating an irreversibility condition in incremental form. 

Focusing on the isochoric setting of von Mises plasticity theory, we define the plastic strain tensor $\Bve^p(\Bx,t)$ and the hardening variable $\alpha(\Bx,t)$ as
\begin{equation}
\Bve^p: 
\left\{
\begin{array}{ll}
\calB \times \calT \rightarrow \mathbb{R}^{\delta\times\delta}_\mathrm{dev} \\[2mm]
(\Bx, t)  \mapsto \Bve^p(\Bx,t)
\end{array}
\right.
\AND
\alpha: 
\left\{
\begin{array}{ll}
\calB \times \calT \rightarrow \mathbb{R}_+ \\[2mm]
(\Bx, t)  \mapsto \alpha(\Bx,t)
\end{array}
\right.
\WITH
\dot{\alpha} \ge 0 ,
\label{s2-fields}
\end{equation}
where $\mathbb{R}^{\delta\times\delta}_\mathrm{dev}:=\{\Be\in\mathbb{R}^{\delta\times\delta} \ \colon \ \Be^T=\Be,\ \tr{[\Be]}=0\}$ is the set of symmetric second-order tensors with vanishing trace. The plastic strain tensor is considered as a local internal variable, while the hardening variable is a possibly non-local internal variable. In particular, $\alpha$ may be introduced to incorporate phenomenological hardening responses and/or non-local effects, for which the evolution equation 
\begin{equation}
\dot\alpha = \sqrt{\frac{2}{3}}\,\vert\dot{\Bve}^p\vert,
\label{s2-evol-alpha}
\end{equation}
is considered. As such, $\alpha$ can be viewed as the equivalent plastic strain, which starts to evolve from the initial condition $\alpha(\Bx,0) = \mathit{0}$. Concerning function spaces, we assume sufficiently regularized plastic responses, i.e., endowed with hardening and/or non-local effects, for which we assume $\Bve^{p}\in\mathbf{Q}:=\mathrm{L}^2(\calB;\mathbb{R}^{\delta\times\delta}_\mathrm{dev})$. Moreover, in view of~\eqref{s2-evol-alpha}, it follows that $\alpha$ is irreversible. Assuming in this section the setting of gradient-extended plasticity, we define the function spaces
%
\begin{equation}
\calW^\alpha_{\alpha_n,\,\Bq}:=\{\alpha\in\calW^\alpha \quad \colon \quad \alpha = \alpha_n + \sqrt{2/3}\,\vert \Bq \vert, \ \Bq \in \mathbf{Q} \},
\label{Walpha}
\end{equation}
{\color{black} where $\calW^\alpha=\mathrm{L}^2(\calB)$ for local plasticity, while $\calW^\alpha=\mathrm{H}^1(\calB)$ for gradient plasticity.} The hardening law~\eqref{s2-evol-alpha} is thus enforced in incremental form by setting $\alpha\in \calW^\alpha_{\alpha_n,\,\Bve^p-\Bve^p_n}$.

The gradient of the displacement field defines the symmetric strain tensor of the geometrically linear theory as
\begin{equation}
\Bve = \nabla_s \Bu = \sym[ \nabla \Bu ] := \frac{1}{2} [\nabla\Bu + \nabla\Bu^T]
.
\label{s2-disp-grad}
\end{equation}
In view of the small strain hypothesis and the isochoric nature of the plastic strains, the strain tensor is additively decomposed into an elastic part $\Bve^e$ and a plastic part $\Bve^p$ as
\begin{equation}
\Bve = \Bve^e + \Bve^p 
\WITH
\tr{[\Bve]} = \tr{[\Bve^e]}
.
\label{s2-strain-e-p}
\end{equation}


For simplicity, the anisotropic material is assumed to be strengthened by a single family of fibers, whose direction is described by a unit vector field $\Ba$ (Figure~\ref{Figure1}). Consequently, the direction-dependent response is characterized by the second-order structural tensor
\begin{equation}
\BM:=\Ba\otimes\Ba .
\end{equation}
The introduction of additional preferred directions can be easily incorporated in future work, following, e.g., \cite{reese2021using}.
The solid $\calB$ is loaded by prescribed deformations and external tractions on the boundary, defined by time-dependent Dirichlet conditions and Neumann conditions
\begin{equation}
\Bu = \overline{\Bu} \ \textrm{on}\ \partial_D\calB
\AND
\Bsigma \cdot \Bn 
= \overline{\Btau} \ \textrm{on}\ \partial_N\calB ,
\label{s2-bcs}
\end{equation}
where $\Bn$ is the outward unit normal vector on the surface $\partial \calB$. The stress tensor $\Bsigma$ is the thermodynamic dual to $\Bve$ and $\bar{\Btau}$ is the prescribed traction vector. Finally, the stress equilibrium is defined as the quasi-static form of the balance of linear momentum
\begin{equation}
\div\,\Bsigma + \overline\Bf = \Bzero,
\label{s2-equil:defo}
\end{equation}
where dynamic effects are neglected and $\overline\Bf$ is a given body force. 


\sectpb[Section22]{Energy quantities and variational principles}
\label{sec:energy_functions}

Let $\BfrakC$ denote the set of constitutive state variables. In the most general setting considered in this study, one has
\begin{equation}
\BfrakC := \{ \Bve, \Bve^p, \alpha, d, \nabla \alpha, \nabla d \}
.
\label{state}
\end{equation}
A pseudo-energy density per unit volume is then defined as $W:=W(\BfrakC)$, which is additively decomposed into an elastic contribution $W_{elas}$,  a plastic contribution $W_{plas}$, and a (regularized) fracture contribution~$W_{frac}\;$:
\begin{equation}
\boxed{W(\BfrakC;\BM):= {W}_{elas}(\Bve,\Bve^p ,d ,\alpha;\BM) + {W}_{plas}(\alpha, d, \nabla \alpha) +
{W}_{frac}(d, \nabla d ; \BM) . } 
\label{psuedo-energy}
\end{equation}
%
We note that $W$ is a state function that contains both energetic and dissipative contributions. With this function at hand, a pseudo potential energy functional can be written~as
\begin{equation}
	\calE(\Bu,\Bve^p,\alpha,d;\BM) := \int_{\calB} {W}(\BfrakC;\BM) \, dv 
	 \; - \; 
	\vphantom{\frac{d}{dt}}
	\mathcal{E}_{ext} (\Bu) , 
\label{potential-functional}
\end{equation}
where $\mathcal{E}_{ext}$ denotes the work of external loads: %
\begin{equation}
	\mathcal{E}_{ext} (\Bu) := 
	\int_{\calB} \overline\Bf \cdot \Bu\, dv  +
	\int_{\partial_N\calB} \overline{\Btau} \cdot \Bu\, da.
\end{equation}

In \textit{variationally consistent models}, the governing equations of the fracturing elasto-plastic solid can be derived from knowledge of the energy functional~\eqref{potential-functional} by invoking rate-type variational principles~\citep{miehe2011,miehe2016bphase} in agreement with the principle of virtual power~\citep{maugin1990,fremond1996}. In such cases, a global rate potential of the form 
\begin{equation}
	\Pi(\dot\Bu,\dot\Bve^p,\dot\alpha,\dot{d};\BM) := \frac{d}{dt}\calE(\Bu,\Bve^p,\alpha,d;\BM) \; + \; \int_\calB \Phi_{vis} (\dot{d},\dot{\alpha}) \, dv,   
\label{potential-functional-cont}
\end{equation}
is defined, where ${\Phi}_{vis}$ denotes the dissipative power density due to viscous resistance forces. In line with previous works~\citep{miehe+welschinger+hofacker10a}, the function
\begin{equation}
{\Phi}_{vis} (\dot{d},\dot{\alpha}) := \frac{\eta_f}{2}\dot{d}^{\,2} + \frac{\eta_p}{2}\dot{\alpha}^2 ,
\label{Dvis}
\end{equation}
is considered, where $\eta_f$ and $\eta_p$ are material parameters that characterize the viscous response of the fracture and plasticity evolutions, respectively. Then, minimization of~\eqref{potential-functional-cont} with respect to $\dot\Bu$, the plasticity variables $(\dot\Bve^p,\dot\alpha)$ subject to the hardening law~\eqref{s2-evol-alpha}, and the crack phase-field $\dot d$ subject to the irreversibility condition $\dot d\geq0$ provide the governing equations for the elasticity problem, the plasticity problem, and the fracture problem, respectively. Such a variational structure results in a convenient numerical implementation based on incremental energy minimization, for which an algorithmic representation of the energy functional~\eqref{potential-functional} is defined as
\begin{equation}
	\Pi^\tau(\Bu,\Bve^p,\alpha,d;\BM) := \calE(\Bu,\Bve^p,\alpha,d;\BM) \; - \; \calE_n \; + \; \Delta  t\int_\calB  \Phi_{vis} \big([d-d_n]/\Delta t,[\alpha-\alpha_n)/\Delta t\big) \, dv  ,  
\label{potential-functional-disc}
\end{equation}
where $\Delta t:=t-t_n$ denotes the time step. The coupled evolution problem then follows as the incremental minimization principle
\begin{equation}
\boxed{ \{\Bu,d,\Bve^p,\alpha\} = \arg\big\{ \min_{\Bu\in\calW_{\overline{\Bu}}^\Bu} \ \min_{d\in\calW^{d}_{d_n}} \ \min_{\{\Bve^p,\alpha\}\in\mathbf{Q}\times\calW^{\alpha}_{\alpha_n,\,\Bve^p-\Bve^p_n}}  \, \Pi^\tau(\Bu,\Bve^p,\alpha,d;\BM)  \big\} . }
\label{inc_min}
\end{equation}
%

\begin{Remark} From equation~\eqref{Dvis}, it is clear that the rate-independent case is recovered by letting $\eta_f\to0$ and $\eta_p\to0$.  In this case, the coupled evolution problem can be equivalently derived in variational form using the energetic formulation for rate-independent systems~\citep{mielke2006,mielke2015}, based on notions of energy balance and stability. This path is followed, for instance, in references~\cite{alessi2015,alessi2017,marigo2016,pham2010,rodriguez2018,ulloa2021}. Moreover, the fact that \eqref{psuedo-energy} is a state function implies that the incremental rate-independent problem exactly recovers the continuous counterpart.
\end{Remark}

For the variational formulation setting, it suffices to define the constitutive energy density functions $W_{elas}$, $W_{plas}$, and $W_{frac}$ to establish the multi-field evolution problem in terms of~\eqref{inc_min}. As we shall recall in the sequel, such a variational structure is not always present in phase-field models for ductile fracture, resulting in greater flexibility at the cost of  a convenient mathematical structure.

\sectpc{Elastic contribution}
The elastic energy density $W_{elas}$ in \eqref{psuedo-energy}
is expressed in terms of the effective strain energy density $\psi_e$. For transversely isotropic materials, $\psi_e$ is defined in terms of the elastic strain tensor $\Bve^e$ and the structural tensor $\BM$. In our formulation, in order to preclude fracture in compression, a decomposition of the effective strain energy density into \textit{damageable} and \textit{undamageable} parts is employed. Thus, we perform
additive decomposition of the strain tensor into \textit{volume-changing}
(volumetric) and \textit{volume-preserving} (deviatoric) counterparts:
\[
\bm\varepsilon^e(\bm u)=\bm\varepsilon^{e,vol}(\bm u)+\bm\varepsilon^{e,dev}(\bm
u),
\]
where the volumetric strain is denoted as $\bm\varepsilon^{e,vol}(\bm
u):=\frac{1}{3}(\bm \varepsilon^e(\bm u):\bm I)\bm I$ and the deviatoric strain is
denoted as $\bm\varepsilon^{e,dev}(\bm u):=\mathbb{P}:\bm \varepsilon^e$. The
fourth-order projection tensor
$\mathbb{P}:=\mathbb{I}-\frac{1}{3}\bm I\otimes\bm I$ 
is introduced to map the full strain tensor onto its deviatoric component.
Therein,
$\mathbb{I}_{ijkl}:=\frac{1}{2}\big(\delta_{ik}\delta_{jl}+\delta_{il}\delta_{jk}\big)$
is the fourth-order symmetric identity tensor. 

The effective strain energy density $\psi_e$ admits the following additive decomposition 
\begin{equation}
	\psi_e(\Bve^e;\BM)=\psi_e^{iso}\big(I_1(\Bve^e),I_2(\Bve^e)\big) + \psi_e^{aniso}\big(I_4(\Bve^e;\BM)\big) .
	\label{eq:strainenergy}
\end{equation}
\textbf{The isotropic strain energy function.} 
The isotropic counterpart admits following additive split:
\begin{equation}
	\psi_e^{iso}\big(I_1,I_2\big):=
	\psi_e^{iso,vol}(I_1)+\psi_e^{iso,dev}(I_1,I_2),
	\label{eq:psi_iso}
\end{equation}
here
\begin{align}
	&\psi_e^{iso,vol}\big(I_1\big)
	=\frac{K}{2}I^2_1
	=\frac{K}{2}\Big(\bm\varepsilon^{e,vol}:\bm I\Big)^2,\notag\\
	&\psi_e^{iso,dev}\big(I_1,I_2\big)
	=\mu\Big(\frac{I_1^2}{3}-I_2\Big)
	=\mu {\bm \varepsilon}^{e,dev}:{\bm \varepsilon}^{e,dev}\WITH \mu>0,
	\label{eq:psi_iso1}
\end{align}
where $I_1$ and $I_2$ denote the invariants through
\begin{equation}
	I_1:=I_1(\Bve^e)=\text{tr}[\Bve^e] \AND I_2:=I_2(\Bve^e)=\tr[(\Bve^e)^2] .
	\label{eq:prin_invar}
\end{equation}
Accordingly, the isotropic strain energy density function given in \req{eq:psi_iso} is additively decomposed into  damageable and undamageable contributions:
\begin{equation}
	\psi_e^{iso}\big(I_1,I_2\big)=\psi_e^{iso,+}(I_1,I_2)+\psi_e^{iso,-}(I_1,I_2),
\end{equation}
where
\begin{align}
	&{\psi_e^{iso,+}}(I_1,I_2)={H{^+}[I_1]}\psi_e^{iso,vol}\big(I_1\big)
	+\psi_e^{iso,dev}\big(I_1,I_2\big)~\AND\notag\\
	&{\psi_e^{iso,-}}(I_1,I_2)=\big(1-{H{^+}[I_1]}\big)\psi_e^{iso,vol}\big(I_1\big)~.	
	\label{eq232_5}
\end{align}
{Therein, $H{^+}[I_1(\Bve^e)]$ is a \textit{positive Heaviside function} which returns one and zero for $I_1(\Bve^e)>0$ and $I_1(\Bve^e)\leq0$, respectively, as shown in \cite{NoiiWick2019}.}

%
\textbf{The anisotropic strain energy function.} 
To complete the formulation,  the anisotropic strain energy function reads
\begin{equation}
	\psi_e^{aniso}\big(I_4(\Bve^e;\BM)\big):=\frac{\chi_a}{2} I_4^2(\Bve^e;\BM)  ,
	\label{eq:psi_aniso}
\end{equation}
where the stiffness parameter $\chi_a$ characterizes the anisotropic deformation response with preferred direction $\Ba$. The \emph{pseudo-invariant} $I_4$ is defined as 
\begin{equation}
	I_4(\Bve^e;\BM)=\text{tr}[\Bve^e\cdot\BM]  .
	\label{eq:pseudo_invar}
\end{equation}

Extending the anisotropic energy into damageable and undamageable parts (see also~\cite{noii2020bayesian}), admits the following splits
\begin{equation}
\psi_e^{aniso}\big(I_4;\BM\big)=\psi_e^{aniso,+}(I_4^+)+\psi_e^{aniso,-}(I_4^-) ,
\end{equation}
where
\begin{equation}
I_4^{\pm}:=\langle I_4(\Bve^e;\BM)\rangle_\pm ,
\label{split_invaniso}
\end{equation}
with the Macaulay bracket $\langle x \rangle_\pm := (x \pm \vert x\vert)/2$. 

%
\textbf{The total elastic strain energy function.} 
%
The elastic contribution to the pseudo-energy density~\eqref{psuedo-energy} finally reads
\begin{equation}
\begin{aligned}
{W}_{elas}(\Bve,\Bve^p,d,\alpha;\BM):= g_e(d,\alpha) \; \big[ \psi_e^{iso,+}(I_1,I_2)&+\psi_e^{aniso,+}(I_4^+) \big ] \\ &+ \psi_e^{iso,-}(I_1,I_2)+\psi_e^{aniso,-}(I_4^-) ,
\end{aligned}
\label{elas-part}
\end{equation}
where $g_e(d,\alpha)$ is the \textit{elastic degradation function}. 

{Following the Coleman-Noll procedure}, the stress tensor is obtained from the potential ${W}_{elas}$ in \req{elas-part} as
\begin{equation}
\begin{aligned}
&\Bsigma = \frac{\partial {W}_{elas}}{\partial\Bve^e}  = \Bsigma^{iso} + \Bsigma^{aniso} \WITH \\
&\Bsigma^{iso}=g_e(d,\alpha)\widetilde\Bsigma^{iso}_{+}+\widetilde\Bsigma^{iso}_{-} \AND  \Bsigma^{aniso}=g_e(d,\alpha)\widetilde\Bsigma^{aniso}_{+}+\widetilde\Bsigma^{aniso}_{-} ,
\end{aligned}
\label{stress}
\end{equation}
where $\widetilde\Bsigma^{iso}$ and $\widetilde\Bsigma^{aniso}$ are the effective stress tensors, given by
\begin{equation}
\begin{aligned}
&\widetilde\Bsigma^{iso}_{+} := 
\frac{\partial \psi_e^{iso,+} }{\partial \Bve^e} = 
K{H{^+}[I_1]}(\bm \varepsilon^e:\bm I)\bm I+2\mu{\bm \varepsilon}^{e,dev}, \\ 
&\widetilde\Bsigma^{iso}_{-} := 
\frac{\partial \psi_e^{iso,-} }{\partial \Bve^e} =K\big(1-{H{^+}[I_1]}\big)(\bm \varepsilon^e:\bm I)\bm I, \AND \\ 
&\widetilde\Bsigma^{aniso}_{\pm} := 
\frac{\partial \psi_e^{aniso,\pm} }{\partial \Bve^e} = \chi_a I_4^\pm\BM. 
\end{aligned}
\end{equation}
%

 
\sectpc{Fracture contribution} 
The phase-field contribution $W_{frac}$ is expressed in terms of the crack surface energy density $\gamma_l$ and the fracture length-scale parameter $l_f$ that governs the regularization. In particular, the sharp-crack surface topology $\calC$ is regularized by a functional $\calC_l$, as outlined in~\cite{aldakheel16} and~\cite{miehe+hofacker+schaenzel+aldakheel15}. This geometrical perspective is in agreement with the framework of~\cite{bourdin2000numerical}, which was conceived as a $\Gamma$-convergence regularization of the variational approach to Griffith fracture~\cite{francfort1998revisiting}. For the case of isotropic materials, the regularized functional~reads
\begin{equation}
\calC_l(d) = \int_{\calB} \gamma_l(d, \nabla d) \, dv .
\label{s2-gamma_l}
\end{equation}
In this work, following~\cite{noii2020adaptive,noii2020bayesian,teichtmeister2017phase}, anisotropic effects are introduced by means of the structural tensor $\BM$. In particular, we assume that $\gamma_l$ admits the additive decomposition 
\begin{equation}
\gamma_l(d, \nabla d;\BM) = \gamma_l^{iso}(d, \nabla d)+\gamma_l^{aniso}(d, \nabla d;\BM) .
\end{equation}
In line with standard phase-field models, a general surface density function for the isotropic part $\gamma_l^{iso}$ is defined as 
\begin{equation}
\gamma_l^{iso}(d, \nabla d):=\frac{1}{c_f}\, \bigg(\frac{\omega(d)}{l_f} + l_f \nabla d \cdot \nabla d \bigg) \WITH c_f:=4\int_0^1\sqrt{\omega(b)}\,db ,
\end{equation}
where $\omega(d)$ is a monotonic and continuous \emph{local fracture energy function} such that $\omega(0)=0$ and  $\omega(1)=1$. A variety of suitable choices for $\omega(d)$ are available in the literature~\citep{kuhn2015,wu2017,wu2018}. Here, the widely adopted linear and quadratic formulations are considered, which yield, respectively, models with and without an elastic stage. Specifically, we define 
\begin{equation}
\omega(d) : = \begin{cases}  d\phantom{^2} \implies c_f=8/3 \quad &\text{model with an elastic stage}, \\ d^2 \implies c_f=2 \quad &\text{model without an elastic stage}.  \end{cases}
\label{wd}
\end{equation}
On the other hand, the anisotropic part $\gamma_l^{aniso}$ reads
\begin{equation}
\gamma_l^{aniso}(d, \nabla d;\bm{M}):=\chi_a \dfrac{l_f}{c_f} \nabla d \cdot \bm{M} \cdot \nabla d  .
\end{equation}
Finally, the fracture contribution to the pseudo-energy density~\eqref{psuedo-energy} reads
\begin{equation}
\begin{aligned}
{W}_{frac}(d, \nabla d ; \BM):= g_f \gamma_l(d, \nabla d;\BM),
\end{aligned}
\label{frac-part}
\end{equation}
where $g_f$ is a parameter that allows to recover different models found in the literature, as will become apparent in the sequel.

\sectpc{Plastic contribution}
The plastic contribution $W_{plas}$ is expressed in terms of an effective plastic energy density $\psi_p$, whose form will depend on the adopted phenomenological model. In line with previous works~\citep{miehe2017phase,rodriguez2018,ulloa2021} let us consider a function in the context of gradient-extended von Mises plasticity:   
\begin{equation}
{\psi}_{p}(\alpha,\nabla\alpha) := \sigma_Y\; \alpha + \frac{H}{2} \alpha^2 + \frac{{\sigma_Y}}{2}\,l_p^2\nabla\alpha\cdot\nabla\alpha ,
\label{psi_p}
\end{equation}
with the initial yield stress $\sigma_Y$, the isotropic hardening modulus $H\ge 0$ and the \emph{plastic length-scale} $l_p$.
The plastic contribution to the pseudo-energy density~\eqref{psuedo-energy} then reads
\begin{equation}
\begin{aligned}
{W}_{plas}(\alpha,d ,\nabla \alpha):= g_p(d) {\psi}_{p}(\alpha,\nabla\alpha),
\end{aligned}
\label{plas-part}
\end{equation}
where $g_p(d)$ is the \textit{plastic degradation function}. The models presented in Sections~\ref{sec:model1} and~\ref{sec:model2} are restricted to local plasticity, for which $l_p=0$, while the model presented in Section~\ref{sec:model3} will include non-local effects, with $l_p>0$. For a variational treatment, it is convenient to invoke the energetic-dissipative decomposition of the plastic energy \req{plas-part}. Thus, the plastic energy density can be further decomposed as
{
\begin{equation}
\begin{aligned}
&{W}_{plas}(\alpha,d ,\nabla \alpha)={W}_{plas}^{ener}(\alpha ,d,\nabla \alpha) + {W}_{plas}^{diss}(\alpha,d),    \WITH \\
&{W}_{plas}^{ener}(\alpha,\nabla \alpha,d):=g_p(d) \frac{1}{2}\big(H\alpha^2 + {\sigma_Y}\,l_p^2\nabla\alpha\cdot\nabla\alpha\big) \AND {W}_{plas}^{diss}(\alpha,d):=g_p(d)\sigma_Y\alpha .
\end{aligned}
\label{plas-part_ener_diss}
\end{equation}
}
%


\sectpb{Stationarity conditions and governing equations}
\label{sec:eq_general}

Let us now derive the \emph{variationally consistent} equations for the multi-field coupled problem. To this end, we seek to find the stationarity conditions for the minimization problem~\eqref{inc_min}. The models presented in Sections~\ref{sec:model1}--\ref{sec:model3} shall take the developments below as canonical forms, and will then deviate from the variationally consistent expressions in favor of greater flexibility.

\sectpc{Elasticity}

The minimization with respect to the displacement field in the variational principle~\eqref{inc_min} yields 
\begin{equation}
\calE_\Bu(\Bu,\Bve^p,\alpha,d;\delta\Bu)=\int_\calB\big[\Bsigma:\Bve(\delta\Bu)-\overline{\Bf}\cdot\delta\Bu\big]\,dv - \int_{\partial\calB_N}\overline{\Btau}\cdot\delta\Bu\,da =0 \quad \forall \,\delta\Bu\in\calW_0^{\Bu}, 
\label{equil_weak}
\end{equation}
which corresponds to the weak form of the mechanical balance equations~\eqref{s2-bcs}, and $\calW_0^{\Bu}$ denotes the function space for the virtual displacement fields, i.e., with homogeneous kinematic boundary conditions.
\sectpc{Fracture}

The directional derivative of~\eqref{potential-functional} with respect to the crack phase-field can be written as
\begin{equation}
\begin{aligned}
\calE_d(\Bu,\Bve^p,\alpha,d;\delta d)=\int_\calB \bigg[ \bigg( & \textcolor{black}{\frac{\partial g_e}{\partial d} (d,\alpha)}\big[\psi_e^{iso,+} + \psi_e^{aniso,+}\big]   \\ & + \frac{g_f}{c_f l_f} \omega'(d) + \textcolor{black}{\partial I_+(d-d_n)} +  \textcolor{black}{g'_p(d)}\psi_p + \frac{\eta_f}{\Delta t} (d-d_n)  \bigg)\delta d  \\ & \hspace*{1cm} +  2\frac{g_f}{c_f} l_f \big( \nabla d \cdot \nabla (\delta d) + \chi_a\nabla d \cdot \BM \cdot \nabla (\delta d) \big) \bigg]\,dv \ni 0 \\ &  \hspace*{7.5cm} \forall \,\delta d\in \calW^d ,
\end{aligned}
\label{phf_weak}
\end{equation}
where the indicator function $I_+\colon \mathbb{R} \to \mathbb{R} \cup \{+\infty\}$ has been introduced to impose the irreversibility condition embedded in $d\in\calW^d_{d_n}$. Let us now define the fracture yield function
\begin{equation}
f_d:= -\textcolor{black}{\frac{\partial g_e}{\partial d} (d,\alpha)}\big[\psi_e^{iso,+} + \psi_e^{aniso,+}\big] - \textcolor{black}{g'_p(d)}\psi_p   - g_f\delta_d\gamma_l  . 
\label{phf_yield}
\end{equation}
The strong form of~\eqref{phf_weak} can then be written as
\begin{equation}
\begin{aligned}
-f_d + \frac{\eta_f}{\Delta t} (d-d_n)+ \partial I_+(d-d_n) \ni 0  \WITH  ({\BI}+\chi_a\BM)\nabla d \cdot \Bn = 0
\quad \mbox{on} \quad \partial\calB . 
\end{aligned}
\label{phf_strong}
\end{equation}
Recalling that  
\begin{equation}
\partial I_+(d-d_n)=\begin{dcases}\{0\} & \text{if}  \quad d>d_n  , \\
\mathbb{R}_-  & \text{if}  \quad d=d_n , \\
\varnothing &  \text{otherwise} , 
\end{dcases}
\label{subdiff_ind}
\end{equation}
the strong form yields, for the rate-dependent case, the evolution equation
\begin{equation}
\frac{\eta_f(d-d_n)}{\Delta t} = f_d \geq 0 \WITH  ({\BI}+\chi_a\BM)\nabla d \cdot \Bn = 0
\quad \mbox{on} \quad \partial\calB .
\label{phf_strong_rdep}
\end{equation}
On the other hand, for the rate-independent case, we obtain the KKT conditions 
\begin{equation}
\begin{aligned}
f_d \leq 0 , \quad (d-d_n)f_d=0   \AND d-d_n\geq0 , \WITH
  ({\BI}+\chi_a\BM)\nabla d \cdot \Bn = 0
\quad \mbox{on} \quad \partial\calB . 
\end{aligned}
\label{phf_strong_rindep}
\end{equation}

The main challenge in solving this evolution problem lies on imposing the irreversibility condition $d\geq d_n$, which allows to replace the set-valued expressions~\eqref{phf_weak} or~\eqref{phf_strong} by equalities. Several alternatives are available in the literature to tackle this problem, including simple penalization methods~\cite{gerasimov2019}, augmented Lagrangian penalization~\cite{wheeler2014}, the primal-dual active set method~\cite{heister2015}, interior point methods~\cite{wambacq2021}, and the complementary system with Lagrange multipliers \cite{mang2020phase}. In this work, we employ the \textit{maximum crack-driving state function} method based on the \textit{history field}, as outlined in~\cite{miehe+hofacker+schaenzel+aldakheel15,MieWelHof10b} and related works. 
\sectpc{Plasticity}
\label{plast_var}
The three models considered in this study employ von Mises plasticity in the rate-independent case, such that $\eta_p=0$ in~\eqref{potential-functional-cont}. Moreover, as will become clear, the plasticity problem does not follow from the incremental minimization problem~\eqref{inc_min} in all three models. In particular, a variational formulation for the plasticity problem that is consistent with the governing equations is only possible if the elastic degradation function introduced in~\eqref{elas-part} does not depend on the hardening variable $\alpha$, that is, {\color{black}$g_{e}:=g_e(d)$}, such that ${W}_{elas}:={W}_{elas}(\Bve,\Bve^p,d;\BM)$. Let us now summarize the variational formulation of such a model in the general gradient-extended case.
{A \textit{free energy density} function for ductile phase-field fracture can be defined~as }
\begin{equation}
	{W}_{free}(\Bve,\Bve^p,d,\alpha;\BM):={W}_{elas}(\Bve,\Bve^p,d;\BM) + {W}_{plas}^{ener}(\alpha ,d,\nabla \alpha) .
	\label{eq:Wfree} 
\end{equation}
Applying the Coleman-Noll procedure to the free energy density function \req{eq:Wfree} yields the following thermodynamic conjugate variables:
\begin{align}
&\Bs^p := -\partial_{\Bve^p} {{W}_{free}}
= \Bsigma \AND\notag\\
& h^p:=\delta_\alpha{{W}_{free}}
=g_p(d)\,H\alpha - \sigma_Y \, l_p^2\div[g_p(d)\,\nabla\alpha].
\label{dual}
\end{align}
In agreement with the classical setting of elasto-plasticity, the yield function is defined as
\begin{equation}
\beta(\Bs^p,h^p;d):= \hbox{$\sqrt{3/2}$}\;\vert \BF^p \vert - h^p- g_p(d)\sigma_Y 
\WITH
\BF^p:= \dev[\Bs^p] = \Bs^p - \frac{1}{3} \mbox{tr} [\Bs^p] \BI .
\label{yield-fcn}
\end{equation}
With the yield function at hand, the strong form of the evolution problem follows from the principle of maximum plastic dissipation
\textcolor{black}{
\begin{equation}
\Phi^p(\dot\Bve^p,\dot\alpha;d) = \sup_{\{\Bs^p,h^p\}} \{ \Bs^p:\dot{\Bve}^p - h^p\dot\alpha  \ \colon \ \beta(\Bs^p,h^p;d)\leq 0 \} ,
\label{pmd}
\end{equation}}
where $\Phi^p$ is the \textit{plastic dissipation potential}. The Euler equations of the maximization principle~\eqref{pmd} follow as the flow rule and hardening law
\begin{equation}
\dot\Bve^p = \lambda^p \frac{\partial \beta}{\partial \Bs^p}  
\AND
\dot\alpha = \lambda^p\frac{\partial \beta}{\partial h^p} ,
\label{plast-evol-eqs}
\end{equation}
together with the KKT conditions conditions
\begin{equation}
\beta \le 0, \quad\quad  \quad\quad \lambda^p \ge 0, 
\quad\quad \mbox{and} \quad\quad
\beta\; \lambda^p = 0 .
\label{palst-kkt}
\end{equation}
Equations~\eqref{plast-evol-eqs} and~\eqref{palst-kkt} constitute the so-called \textit{dual form} of the elasto-plastic problem in strong form. To arrive at a primal formulation~\citep{han1999,ulloa2020}, the dissipation potential is evaluated from~\eqref{pmd} as 
%
\begin{equation}
\Phi_p(\dot\Bve^p; d)=g_p(d)\;\widetilde{\Phi}_p(\dot\Bve^p) \WITH \widetilde{\Phi}_p(\dot\Bve^p)=\sqrt{\frac{2}{3}}\, \sigma_Y  \vert \dot\Bve^p \vert  .
\label{eq:pdiss}
\end{equation}
The Legendre transformation of $\Phi_p$ then reads
\textcolor{black}{
\begin{equation}
\Phi_p^*(\Bs^p,h^p ; d )= \sup_{\{\dot\Bve^p,\dot\alpha\}} \{ \Bs^p:\dot{\Bve}^p - h^p\dot\alpha  - \Phi_p(\dot\Bve^p; d) \} ,
\label{eq:pdiss_trans}
\end{equation}
}
which yields, as a necessary condition, the primal representation of the plasticity evolution problem in the form of a Biot-type equation:
\begin{equation}
\{\Bs^p , h^p\}\in\partial_{\{\dot\Bve^p,\dot\alpha\}}\,\Phi_p(\dot\Bve^p;d) .
\label{eq:pbiot}
\end{equation}
From standard arguments of convex analysis~\cite{rockafellar1970,han1999}, this expression implies the associative flow relations~\eqref{plast-evol-eqs} as well as the loading/unloading conditions~\eqref{palst-kkt}.

To derive the above governing equations from the incremental minimization problem~\eqref{inc_min}, we make use of equations~\eqref{plas-part_ener_diss} and~\eqref{eq:pdiss}, such that the functional derivative of~\eqref{potential-functional-disc} with respect to $\{\Bve^p,\alpha\}$ can be written as 
\begin{equation}
\begin{aligned}
\calE_\alpha(\Bu,\Bve^p,\alpha,d;\delta \Bve^p, \delta \alpha)&=\int_\calB \bigg[\partial_{\Bve^p} W:\delta\Bve^p + \delta_\alpha W\delta\alpha \bigg]\,dv 
\\  &=   \int_\calB \bigg[\partial_{\Bve^p}\big(W_{elas} +  \Phi_p(\Bve^p-\Bve^p_n;d)\big):\delta\Bve^p   + \delta_\alpha{W}_{plas}^{ener}\delta\alpha\bigg]\,dv \ni 0  \\   & \hspace*{5.75cm} \forall \,\delta \Bve^p \in \mathbf{Q} ,  \ \delta \alpha \in \calW^\alpha_{0,\,\delta\Bve^p} .
\end{aligned}
\label{palst_weak}
\end{equation}
In view of equations~\eqref{eq:Wfree} and~\eqref{dual}, the strong form of~\eqref{palst_weak} can be written as
\begin{equation}
\{\Bs^p , h^p\}\in\partial_{\{\Bve^p,\alpha\}}\,\Phi_p(\Bve^p-\Bve^p_n;d) , \WITH
 \nabla \alpha \cdot \Bn = 0
\quad \mbox{on} \quad \partial\calB . 
\label{eq:pbiot_inc}
\end{equation}
Equation~\eqref{eq:pbiot_inc} represents an incremental version of the compact evolution equation~\eqref{eq:pbiot} for the plasticity model, recovered in a variationally consistent manner from the incremental minimization principle~\eqref{inc_min}. Recall that in the present formulation, this was achieved by assuming an elastic degradation function that does not depend on the plastic variables.

\begin{table}[!]
	\caption{Functions and parameters for the three examined models.}
	\vspace{1mm}
	\centering
	\begin{tabular}{lccc}
		Model property                 & $\calM_1$ & $\calM_2$ & $\calM_3$ \\[2mm]\hline\\
		Elastic degradation  $g_e$          & $(1-d)^{2{\alpha}/{\alpha_{\mathrm{crit}}}}$   &  $(1-d)^2$  &  $(1-d)^2$  \\[4mm]
		Fracture constant $g_f$              & $G_c$   &    $2\,l_f c_f{\psi_c}$  & $l_f c_f{w_0}$   \\[4mm]
		Plastic degradation $g_p$           & $1$   &  $(1-d)^2$   &  $(1-d)^2$  \\[4mm]
		Local fracture energy $\omega$         & $d^2$   &  $d$  &  $d$  \\[4mm]
		Crack viscosity $\eta_f$               & $0$   &  $\geq0$  &  $0$  \\[4mm]
		Plastic length-scale $\l_p$           & $0$   & $0$   & $\geq0$  \\[4mm]			Driving scaling factor $\zeta$           & $1$   & $\geq 0$   & $\geq 0$  	\\[4mm]    \hline
		\label{models}
	\end{tabular}
\end{table}

\sectpb{Specific models revisited}
\label{sec:models}

In this section, three benchmark phase-field models  for ductile fracture, hereinafter labeled $\calM_1$, $\calM_2$, and $\calM_3$, are revisited within the framework elaborated in the previous sections. The material parameters and constitutive functions that allow to recover each model from the general formulation are presented in Table~\ref{models}.

\sectpc{Local plasticity with $G_c$ based fracture criteria: Model 1 ($\mathcal{M}_1$)}
\label{sec:model1}

The first model considered in this study takes the work from~\cite{ambati2015} as a point of departure. Therein, an extension of the model proposed in~\cite{du2015} was considered   by further coupling the fracture process to plasticity through dependence of the elastic degradation function on the hardening variable $\alpha$. The model was subsequently extended to finite strains and presented with experimental verification in~\cite{ambati2016}. As discussed in Section~\ref{plast_var}, dependence of the elastic degradation function on $\alpha$  results in lack of variational  consistency for the plasticity evolution problem, in favor of greater flexibility. In this case, consider the 
\begin{equation}
\mbox{Global Primary Fields}: \ 
\BfrakU := \{ \Bu, d \} , 
\label{primary_m1}
\end{equation}
and the
\begin{equation}
\mbox{Constitutive State Variables}: \ 
\BfrakC := \{ \Bve, \Bve^p, \alpha, d, \nabla d \}
.
\label{state_m1}
\end{equation}
With the constitutive choices shown in Table~\ref{models} for $\calM_1$, the following forms of the governing equations presented in Section~\ref{sec:eq_general} are obtained.


The strong form of the crack phase-field evolution~\eqref{phf_strong} takes the form
\begin{equation}
\begin{aligned}
-2\frac{\alpha}{\alpha_{\mathrm{crit}}}(1-d)^{2\frac{\alpha}{\alpha_{\mathrm{crit}}}-1}\frac{l_f}{G_c}\big(\psi_e^{iso,+} &+ \psi_e^{aniso,+}\big) + \big( {d} - l_f^2  \div [ \nabla d ] - l_f^2 \chi_a \div [\nabla d \cdot \BM] \big)+ \\  &\frac{l_f}{G_c}\partial I_+(d-d_n) \ni 0 \quad \quad \quad \;\mbox{in} \quad \calB,  \\ 
 \mbox{with} \qquad &(\BI+\chi_a\BM)\nabla d \cdot \Bn = 0
\quad \quad \mbox{on} \quad  \partial\calB.
\end{aligned}
\label{phf_strong_m1}
\end{equation}
Herein, $\alpha_{\mathrm{crit}}$ is a threshold material parameter introduced in~\cite{ambati2015} to calibrate the softening response. To enforce the crack irreversibility condition, and, therefore, to cast this {\color{black}inequality constrained boundary value problem (BVP) as an equality constrained BVP}, the history field 
\begin{equation}
\calH(\Bx,t):=\max_{s\in[0,t]}\widetilde{D}\big(\BfrakC(\Bx,s)\big) \WITH  \widetilde{D}:=\zeta\frac{2\,l_f}{G_c}\big(\psi_e^{iso,+} + \psi_e^{aniso,+}\big), 
\end{equation}
is introduced. Here, $\zeta\geq0$ is a scaling parameter that introduces further flexibility in the formulation, allowing to tune the post-critical range (cf.~\cite{aldakheel16}). Equation~\eqref{phf_strong_m1} is then restated as
\begin{equation}
\begin{aligned}
\frac{\alpha}{\alpha_{\mathrm{crit}}}(1-d)^{2\frac{\alpha}{\alpha_{\mathrm{crit}}}-1}\calH - \big( {d} - l_f^2  \div [ \nabla d ] - l_f^2 \chi_a \div [\nabla d \cdot\BM] \big) = 0  \quad &\mbox{in} \quad \calB, \\ 
 \mbox{with} \qquad (\BI+\chi_a\BM)\nabla d \cdot \Bn = 0
\quad &\mbox{on} \quad \partial\calB .
\end{aligned}
\label{phf_strong_m1_hist}
\end{equation}
With the last expression, and in view of~\eqref{equil_weak}, the global primary fields are found as the solution of the following coupled problem: find $\Bu\in\calW^{\Bu}_{\overline{\Bu}}$ and {\color{black}$d\in\calW^{d}$} such that
\begin{align}
\tag{$\calM_1$}
\left\{
\begin{aligned}
&\int_\calB\big[\Bsigma(\Bve,\Bve^p,d,\alpha;\BM):\Bve(\delta\Bu)-\overline{\Bf}\cdot\delta\Bu\big]\,dv - \int_{\partial\calB_N}\overline{\Btau}\cdot\delta\Bu\,da =0 \quad  \forall \,\delta\Bu\in\calW^{\Bu}_0\  ,  \\
&\int_\calB \bigg[ \bigg(\frac{\alpha}{\alpha_{\mathrm{crit}}}(1-d)^{2\frac{\alpha}{\alpha_{\mathrm{crit}}}-1}\calH -  {d}\bigg)\delta d - l_f^2  \nabla d \cdot \nabla (\delta d) - l_f^2 \chi_a  \nabla d \cdot \BM \cdot \nabla (\delta d) \bigg]\,dv = 0    \\ 
&  && \hspace*{-2.75cm}\forall \,{\color{black}\delta d\in \calW^{d}}. 
\end{aligned} 
\right.
\label{phf_weak_m1_hist}
\end{align}
Thus, 
\begin{equation*}
	\fterm{ 
	\calM_1:=\calM_1(\Bve,\Bve^p,d,\alpha;\BM)
	=\calM_1(\BfrakU, \delta \Bu)+\calM_1(\BfrakU, \delta d)=0 \quad\forall\;\; (\delta \Bu,\delta d)\in \Big(\calW^{\Bu}_0,\calW^{d}\Big).}
\end{equation*}

\begin{Remark}
\label{variational_incons}
The introduction of the history field in the displacement Euler-Lagrange equation~\eqref{phf_strong_m1}, finally yielding~\eqref{phf_weak_m1_hist}$_1$, results in a loss of variational consistency with respect to the energy functional~\eqref{potential-functional-disc} due to  the filtering of the maximum history value of $\widetilde{D}$ and the scaling factor $\zeta$ for $\zeta\neq 1$. The upside of this choice is a convenient numerical strategy for solving the original inequality-constrained PDE, and greater flexibility in the model. 
\end{Remark}

\begin{Remark}
\label{zeta_M1}
The role of the parameter $\zeta$, i.e., tuning the post-critical range by scaling the driving force, is already achieved in the present model by means of $\alpha_{\mathrm{crit}}$. Consequently, $\zeta=1$ is assumed hereafter for $\calM_1$.
\end{Remark}

\begin{Remark}
	\label{Cdriv_m1}
	At this point, it is worth noting that the crack driving force in \eqref{phf_weak_m1_hist}$_2$, i.e., $\calH$, is scaled by the hardening variable $\alpha$, such that the crack driving force vanishes for $\alpha\rightarrow0$. As a consequence, fracture cannot occur outside the ductility zone, and a response corresponding to elastic damage followed by plastic damage is not possible in this model due to the strong coupling between damage and plasticity. For a detailed discussion of different possible elastic-plastic-damage evolution response, see \cite{alessi2017}.
\end{Remark}


Concerning the plasticity evolution problem, a variational derivation in the sense of~\eqref{palst_weak} is not possible in the present model due to the dependence of the elastic degradation function $g_p$ on $\alpha$. In this case, the \emph{local} evolution of the plasticity variables $\{\Bve^p,\alpha\}$ according to equations~\eqref{plast-evol-eqs} and~\eqref{palst-kkt} (alternatively, \eqref{eq:pbiot} or~\eqref{eq:pbiot_inc}  in the incremental form) is postulated in a non-variational context. For the present model, the yield function~\eqref{yield-fcn} takes the form
\begin{equation}
\beta= \hbox{$\sqrt{3/2}$}\;\vert \BF^p(\Bve,\Bve^p,d,\alpha;\BM) \vert - (\sigma_Y + H\alpha) .
\label{yield-fcn_m1}
\end{equation}

\sectpc{Local plasticity with $\psi_c$ based fracture criteria: Model 2 $(\calM_2)$} 
\label{sec:model2}

The second model is based on the geometrically conceived approach to the phase-field modeling of ductile fracture, conceptually based on the local plasticity theory described in~\cite{miehe+hofacker+schaenzel+aldakheel15} and considered in subsequent works~\cite{aldakheel2020microscale,storm2021comparative}. The original model is constructed within a variationally consistent framework, in agreement with the incremental energy minimization principle~\eqref{inc_min}. In this case, consider the
\begin{equation}
\mbox{Global Primary Fields}: \ 
\BfrakU := \{ \Bu, d \},
\label{primary_m2}
\end{equation}
and the
\begin{equation}
\mbox{Constitutive State Variables}: \
\BfrakC := \{ \Bve, \Bve^p, \alpha, d, \nabla d \}
.
\label{state_m2}
\end{equation}
With the constitutive choices shown in Table~\ref{models} for $\calM_2$, the following forms of the governing equations described in Section~\ref{sec:eq_general} are obtained.

Letting $l_d:=\sqrt{2}\,l_f$ (cf.~\cite{marigo2016}), and after simple manipulations, the strong form of the crack phase-field evolution~\eqref{phf_strong} can be written as (cf.~\cite{miehe+hofacker+schaenzel+aldakheel15,aldakheel16}):
\begin{equation}
\begin{aligned}
-2(1-d)^2\psi_c\bigg(\frac{\psi_e^{iso,+} + \psi_e^{aniso,+}+ \psi_p}{\psi_c} - 1 \bigg) + 2\psi_c\big( d - l_d^2  \div [ \nabla d ] - l_d^2   \chi_a \div [\nabla d &\cdot\BM] \big) \\ + \frac{\eta_f}{\Delta t} (d-d_n) + \partial_d I_+(d-d_n) \ni 0 \quad &\mbox{in} \quad \calB,  \\ 
 \mbox{with} \qquad (\BI+\chi_a\BM)\nabla d \cdot \Bn = 0
\quad &\mbox{on} \quad \partial\calB , 
\end{aligned}
\label{phf_strong_m2}
\end{equation}
where the role of $\psi_c$ as a specific critical fracture energy density is clearly reflected. To enforce the crack irreversibility condition, and, therefore, to cast this {\color{black}inequality constrained BVP as an equality constrained BVP}, the history field 
\begin{equation}
\calH(\Bx,t):=\max_{s\in[0,t]}\widetilde{D}\big(\BfrakC(\Bx,s)\big) \WITH  \widetilde{D}:=\zeta\bigg(\frac{\psi_e^{iso,+} + \psi_e^{aniso,+}+ \psi_p}{\psi_c} - 1 \bigg),
\label{histfield_m2}
\end{equation}
is introduced. Letting $\eta_d:=\eta_f/(2\psi_c)$, \eqref{phf_strong_m2} is restated as
\begin{equation}
\begin{aligned}
(1-d)\calH - \big( {d} - l_d^2  \div [ \nabla d ] - l_d^2 \chi_a \div [\nabla d \cdot\BM] \big) = \frac{\eta_d}{\Delta t}(d-d_n) \quad &\mbox{in} \quad \calB,  \\ 
 \mbox{with} \qquad (\BI+\chi_a\BM)\nabla d \cdot \Bn = 0
\quad &\mbox{on} \quad \partial\calB .
\end{aligned}
\label{phf_strong_m2_hist}
\end{equation}
%
With the last expression, and in view of~\eqref{equil_weak}, the global primary fields are found as the solution of the following coupled problem:  find $\Bu\in\calW^{\Bu}_{\overline{\Bu}}$ and {\color{black}$d\in\calW^{d}$}, such that
\begin{align}
\tag{$\calM_2$}
\left\{
\begin{aligned}
&\int_\calB\big[\Bsigma(\Bve,\Bve^p,d;\BM):\Bve(\delta\Bu)-\overline{\Bf}\cdot\delta\Bu\big]\,dv - \int_{\partial\calB_N}\overline{\Btau}\cdot\delta\Bu\,da =0 \quad  \forall \,\delta\Bu\in\calW^{\Bu}_0\  ,  \\
&\int_\calB \bigg[ \bigg((1-d)\calH -  {d} + \frac{\eta_d}{\Delta t}(d-d_n)\bigg)\delta d - l_d^2  \nabla d \cdot \nabla (\delta d) - l_d^2 \chi_a  \nabla d \cdot \BM \cdot \nabla (\delta d) \bigg]\,dv = 0    \\ 
&  && \hspace*{-2.75cm}\forall \,{\color{black}\delta d\in\calW^{d}}  . 
\end{aligned}
\right.
\label{phf_weak_m2_hist}
\end{align}
Thus, 
\begin{equation*}
	\fterm{ 
		\calM_2:=\calM_2(\Bve,\Bve^p,d,\alpha;\BM)
		=\calM_2(\BfrakU, \delta \Bu)+\calM_2(\BfrakU, \delta d)=0 \quad\forall\;\; (\delta \Bu,\delta d)\in \Big(\calW^{\Bu}_0,\calW^{d}\Big).}
\end{equation*}
Note that, in light of Remark~\ref{variational_incons}, the introduction of the history field and the scaling parameter $\zeta$ in~\eqref{histfield_m2} results in a loss of variational consistency with respect to the energy functional~\eqref{potential-functional-disc} for the fracture problem.

As opposed to $\calM_1$, the local plasticity evolution problem in the present model is variationally consistent (see Section~\ref{plast_var}). The \emph{local} evolution of the plasticity variables $\{\Bve^p,\alpha\}$ according to the evolution equation~\eqref{eq:pbiot_inc}, which represents an incremental, \emph{primal} version of equations~\eqref{plast-evol-eqs} and~\eqref{palst-kkt}, is then a necessary condition of the minimization principle~\eqref{inc_min}. For the present model, the yield function~\eqref{yield-fcn} takes the form
\begin{equation}
\beta= \hbox{$\sqrt{3/2}$}\;\vert \BF^p(\Bve,\Bve^p,d;\BM) \vert - (1-d)^2(\sigma_Y + H\alpha) .
\label{yield-fcn_m2}
\end{equation}

\sectpc{Non-local plasticity with $w_0$ based fracture criteria: Model 3 ($\calM_3$)} 
\label{sec:model3}
The third model considered in this study is inspired by the variational phase-field models coupled to gradient plasticity proposed in~\cite{miehe+hofacker+schaenzel+aldakheel15,aldakheel16}. The modeling framework adopted therein and in subsequent studies~\citep{miehe2016bphase,dittmann2018,aldakheel+wriggers+miehe18} is consistent with the rate-type variational framework of~\cite{miehe2011}. In the small-strain rate-independent case, similar models were proposed in~\cite{ulloa2016,rodriguez2018,ulloa2021}, where a variationally consistent energetic formulation was adopted to derive the governing equations.  Consider, in this case, the 
\begin{equation}
\mbox{Global Primary Fields}: \  
\BfrakU := \{ \Bu, d ,\alpha \},
\label{primary_m3}
\end{equation}
and the
\begin{equation}
\mbox{Constitutive State Variables}: \  
\BfrakC := \{ \Bve, \Bve^p, \alpha, d ,\nabla\alpha, \nabla d \}
,
\label{state_m3}
\end{equation}
representing a combination of a  first-order gradient plasticity model and a first-order gradient damage model. With the constitutive choices shown in Table~\ref{models} for $\calM_3$, the following forms of the governing equations described in Section~\ref{sec:eq_general} are obtained.

With a slight change of parameters, the fracture problem in the present model admits the same formulation of $\calM_2$ in Section~\ref{sec:model2}. In this case, according to Table~\ref{models}, the strong form of the crack phase-field evolution~\eqref{phf_strong} can be written as:
\begin{equation}
\begin{aligned}
-(1-d)^2w_0\bigg(\frac{\psi_e^{iso,+} + \psi_e^{aniso,+}+ \psi_p}{w_0/2} - 1 \bigg) + w_0\big( d - l_d^2  \div [ \nabla d ] - l_d^2   \chi_a \div [\nabla d \cdot&\BM] \big) \\  + \partial_d I_+(d-d_n) \ni 0  \quad &\mbox{in} \quad \calB, \\ 
 \mbox{with} \qquad (\BI+\chi_a\BM)\nabla d \cdot \Bn = 0
\quad &\mbox{on} \quad \partial\calB , 
\end{aligned}
\label{phf_strong_m3}
\end{equation}
where $w_0$ is a critical fracture energy density. One can show the identity of $w_0=2\psi_c$ holds for brittle fracture, but in the present gradient plasticity model, $w_0\neq2\psi_c$ due to the non-local term in $\psi_p$. Indeed, the main difference of the present model with respect to $\calM_1$ and $\calM_2$ is that the plastic free energy $\psi_p$, defined in~\eqref{psi_p}, is considered here with $l_p>0$, and thus introduces non-local effects in the fracture driving force. To enforce the crack irreversibility condition, we define the history field 
\begin{equation}
\calH(\Bx,t):=\max_{s\in[0,t]}\widetilde{D}\big(\BfrakC(\Bx,s)\big) \WITH  \widetilde{D}:=\zeta\bigg(\frac{\psi_e^{iso,+} + \psi_e^{aniso,+}+ \psi_p}{w_0/2} - 1 \bigg) .
\label{histfield_m3}
\end{equation}
Thus,~\eqref{phf_strong_m3} is restated as
\begin{equation}
\begin{aligned}
(1-d)\calH - \big( {d} - l_d^2  \div [ \nabla d ] - l_d^2 \chi_a \div [\nabla d \cdot\BM] \big) = 0  \quad &\mbox{in} \quad \calB, \\ 
 \mbox{with} \qquad (\BI+\chi_a\BM)\nabla d \cdot \Bn = 0
\quad &\mbox{on} \quad \partial\calB .
\end{aligned}
\label{phf_strong_m3_hist}
\end{equation}
As before, in light of Remark~\ref{variational_incons}, the introduction of the history field and the scaling parameter $\zeta$ in~\eqref{histfield_m2} results in a loss of variational consistency with respect to the energy functional~\eqref{potential-functional-disc} for the fracture problem.

As in $\calM_2$, the plasticity evolution problem for the present model is variationally consistent (see Section~\ref{plast_var}). Moreover, the problem now includes non-local effects modulated by the plastic length-scale $l_p>0$, where the yield function~\eqref{yield-fcn} reads
\begin{equation}
\beta= \hbox{$\sqrt{3/2}$}\;\vert \BF^p(\Bve,\Bve^p,d;\BM) \vert - (1-d)^2(\sigma_Y + H\alpha) + \sigma_Y \,l_p^2 \div[ (1-d)^2\nabla\alpha  ] .
\label{yield-fcn_m3}
\end{equation}
To derive the global PDE governing the evolution of the non-local field $\alpha$, we  take the weak form~\eqref{palst_weak} as a point of departure, such that, for $\vert\Bve^p-\Bve^p_n\vert>0$:
\begin{equation}
\begin{aligned}
\int_\calB \bigg[\partial_{\Bve^p}\big(W_{elas} &+  \Phi_p(\Bve^p-\Bve^p_n;d)\big):\delta\Bve^p   + \delta_\alpha{W}_{plas}^{ener}\delta\alpha\bigg]\,dv \\ &=  \int_\calB \bigg[-\Bsigma(\Bve,\Bve^p,d;\BM):\delta\Bve^p  +  (1-d)^2\sigma_Y\hat\Bn:\delta\Bve^p   \\ & \hspace*{2cm}+(1-d^2)\,H\alpha\delta\alpha - \sigma_Y \, l_p^2\div[(1-d)^2\,\nabla\alpha]\delta\alpha\bigg]\,dv \\ &=  \int_\calB \bigg[-\sqrt\frac{3}{2}\vert\BF^p(\Bve,\Bve^p,d;\BM)\vert  +  (1-d)^2\sigma_Y   \\ & \hspace*{2cm}+(1-d^2)\,H\alpha - \sigma_Y \, l_p^2\div[(1-d)^2\,\nabla\alpha]\bigg]\delta\alpha\,dv = 0  ,
\end{aligned}
\label{palst_weak_m3_1}
\end{equation}
where $\hat\Bn:=(\Bve^p-\Bve^p_n)/\vert \Bve^p-\Bve^p_n \vert$ is the direction of the plastic flow. Note that in~\eqref{palst_weak_m3_1}, we have considered virtual fields 
\begin{equation}
\delta\Bve^p=\hat\Bn\vert\delta\Bve^p\vert=\sqrt{\frac{3}{2}}\hat\Bn\delta\alpha ,
\end{equation}
such that the direction of the plastic flow is fixed, while the virtual equivalent plastic strain $\delta\alpha\in\calW^\alpha_{0,\,\delta\Bve^p}$ is allowed to vary. To solve~\eqref{palst_weak_m3_1}, we must enforce the constraint embedded in $\alpha\in \calW^\alpha_{\alpha_n,\,\Bve^p-\Bve^p_n}$ (Equation~\eqref{Walpha}), such that
\begin{equation}\label{plastic_law}
\alpha=\alpha_n+\sqrt{\frac{3}{2}}\vert\Bve^p-\Bve^p_n\vert .
\end{equation}
Recalling that~\eqref{palst_weak_m3_1} is a weak representation of~\eqref{eq:pbiot_inc} for $\vert\Bve^p-\Bve^p_n\vert>0$, and that~\eqref{eq:pbiot_inc}  implies the incremental version of the plastic flow rule~\eqref{plast-evol-eqs}, a possible way to proceed is to replace the local field $\Bve^p$ by setting, in agreement with~\eqref{plast-evol-eqs}, 
\begin{equation}
\Bve^p=\Bve^p_n+\sqrt{\frac{3}{2}}(\alpha-\alpha_n)\hat\Bn^{trial} \WITH \hat\Bn^{trial}=\frac{\BF^{p,trial}}{\vert \BF^{p,trial} \vert}, 
\label{eq:ep_inc}
\end{equation}
where, from standard arguments of von Mises plasticity, $\BF^{p,trial}:=\BF^{p}(\Bve,\Bve^p_n,d;\BM)$. With this expression at hand, \eqref{palst_weak_m3_1} may be left as a function of the non-local field $\alpha$, subject to the irreversibility condition $\alpha\geq\alpha_n$. 

Finally, the global primary fields are found as the solution of the following coupled problem:  find $\Bu\in\calW^{\Bu}_{\overline{\Bu}}$, {\color{black}$\alpha\in \calW^{\alpha}$}, and {\color{black}$d\in\calW^{d}$}, such that
\begin{align}
\tag{$\calM_3$}
\left\{
\begin{aligned}
&\int_\calB\big[\Bsigma(\Bve,\Bve^p,d;\BM):\Bve(\delta\Bu)-\overline{\Bf}\cdot\delta\Bu\big]\,dv - \int_{\partial\calB_N}\overline{\Btau}\cdot\delta\Bu\,da =0 \quad  \forall \,\delta\Bu\in\calW^{\Bu}_0 \  ,  \\
&\int_\calB \bigg[-\sqrt\frac{3}{2}\vert\BF^p(\Bve,\Bve^p,d;\BM)\vert  +  (1-d)^2\sigma_Y  + \partial_\alpha I_+(\alpha-\alpha_n)   \\ & \hspace*{2cm}+(1-d^2)\,H\alpha + \sigma_Y \, l_p^2(1-d)^2\nabla\alpha\cdot\nabla(\delta\alpha)\bigg]\delta\alpha\,dv \ni 0  \quad \forall \,{\color{black} \delta \alpha\in\calW^{\alpha}}, \\
&\int_\calB \bigg[ \big(\,(1-d)\calH -  {d}\, \big)\delta d - l_d^2  \nabla d \cdot \nabla (\delta d) - l_d^2 \chi_a  \nabla d \cdot \BM \cdot \nabla (\delta d) \bigg]\,dv = 0   \quad  \,{\color{black}\delta d\in\calW^{d}}  .
\end{aligned}
\right.
\label{phf_weak_m3_hist}
\end{align}
Thus, 
\begin{empheq}[box=\widefbox]{align*}
\calM_3:=\calM_3(\Bve,\Bve^p,d,\alpha;\BM)
=\calM_3(&\BfrakU, \delta \Bu) +\calM_3(\BfrakU, \delta \alpha)+\calM_3(\BfrakU, \delta d)=0  \\
	&\forall\;\; (\delta \Bu,\delta \alpha, \delta d)\in \Big(\calW^{\Bu}_0,\calW^\alpha,\calW^{d}\Big).
\end{empheq}

Note that in~\eqref{phf_weak_m3_hist}$_2$, the indication function has been introduced to impose irreversibility of $\alpha$, resulting in a multivalued expression. To show the consistency of this equation with the governing equations derived in Section~\ref{plast_var}, we first recall that the incremental flow rule has been enforced by means of~\eqref{eq:ep_inc}. Then, we note that the strong form of~\eqref{phf_weak_m3_hist}$_2$ yields $\beta \in \partial_d I_+(\alpha-\alpha_n)$, with $\beta$ given in~\eqref{yield-fcn_m3}. In view of~\eqref{subdiff_ind}, it is easy to see that this expression represents the incremental version of the KKT conditions~\eqref{palst-kkt}. To handle the inequality constraint and eliminate the multivalued term $\partial_d I_+(\alpha-\alpha_n)$, a constrained optimization technique is required. For instance, an interior-point method has been recently proposed in~\cite{wambacq2021}.

In the sequel, the phase-field models for ductile fracture formulated in $\calM_1$, $\calM_2$, and $\calM_3$ will be taken as inputs for the Bayesian inversion framework described in a detail in Section 4.  
\sectpa[Section3]{Parameter estimation based on Bayesian inference }
\label{Bayesian}
 
In this section, we review different parameter estimation techniques based on MCMC to identify the mechanical parameters involved in ductile fracture. First, some basic statistical principals are briefly recalled.

In Bayesian estimation, a parametric forward model (e.g., a PDE-based model or a coupled variational inequality system) is used to update the available data (considered as random variables) based on the available information (denoting the prior knowledge). The posterior information is then provided as output \cite{rappel2020tutorial,smith2013uncertainty,wang2020uncertainty}. 


Bayes' formula prescribes the probability of an event according to related prior information and is given by
\begin{align}
P(A|B)=\frac{P(B|A)P(A)}{P(B)},
\label{Bay}
\end{align} 
where $P(A|B)$ denotes the conditional probability of event $A$ happening when $B$ has happened (likewise for $P(A|B)$), and $P(\cdot)$ is the probability of observations $A$ and $B$. Using a probability density function $\pi$, we can rewrite \eqref{Bay} as
\begin{align}
\pi(\chi|m)=\frac{\pi(m|\chi) \pi_0(\chi)}{\pi(m)}.
\end{align}
Here, $\pi_0(\chi)$ is the prior distribution which indicates the available information regarding the parameter $\chi$. For the ductile fracture case, the set of parameters $\chi$ is indicated in~\eqref{all}.
Moreover, $\pi(\chi|m)$ denotes the posterior density, i.e., the probability density of the parameter $\chi$ considering the measurement $m$.  The probability of the parameter $\chi$ with respect to the observation/measurement is described by the likelihood function $\pi(m|\chi)$. The denominator $\pi(m)$ is a constant normalization factor, such that
\begin{align}
\pi(\chi|m)\propto\pi(m|\chi)\pi_0(\chi).
\end{align}
In Bayesian inversion,  the solution of the inverse problem is the posterior density giving the distribution of the unknown parameter values based on the sampled observations.
MCMC is a popular method to calculate this distribution, where a Markov chain is constructed whose stationary distribution is the sought posterior distribution in Bayes' theorem.

In order to identify the unknown parameters, we introduce the following statistical model: 
\begin{align}
\label{81}
\mathbb{M}=f(\xx,\chi)+\varepsilon,
\end{align}
where $\mathbb{M}$ is an $n$-dimensional vector that indicates the measurement, $f$ denotes the PDE-based model, and $\chi=\{\chi_{_1},\chi_{_2}\ldots,\chi_{_k}\}$ is a $k$-dimensional vector denoting the model parameters. The model output $f(\xx,\chi)$ is the response quantity of interest, collected in an $n$-dimensional vector, where $n=n_Tn_C$, with $n_C$ denoting the number of components of the response variable and $n_T$ denoting the number of time steps. In the present work, the force-displacement curve is taken as the response variable $f(\chi)$, such that $n_C=1$. For  the measurement error $\varepsilon$, we employ a Gaussian independent and identically distributed error $\varepsilon\sim\mathcal{N}(0,\sigma^2\,I)$, where $\sigma^2$ is a fidelity parameter.

Given a measurement or observation $m=\mathtt{obs}$, the conditional density reads
\begin{align}
\pi(\mathtt{obs})= \displaystyle\int_{\R^n}\pi(\mathtt{obs}|\chi)\pi_0(\chi)\,d\chi\neq 0.
\end{align}
The inverse problem in the Bayesian framework can thus be stated as follows:
given a measurement $m$, find the posterior density $\pi(\chi|m)$. 

To this end, one makes use of \textit{Bayes'~theorem~of~inverse~problems} \cite{smith2013uncertainty}, which can be stated as follows:
\begin{theorem}[\textit{Bayes'~theorem~for~parameter~estimation}]
	We consider random parameter variables $\chi$ and a specific prior distribution $\pi_0(\chi)$,  and we consider $m$ to be a realization of the random observation variable (denoting the measurement or virtual observation). The posterior distribution considering the measurement $m$ follows as
	\begin{align}\label{82}
	\pi(\chi|m)=\frac{\pi(m|\chi) \pi_0(m)}{\pi(m)}=\frac{\pi(m|\chi) \pi_0(\chi)}{\displaystyle\int_{\mathbb{R}^n}\pi(m|\chi)\pi_0(m)d\chi}.
	\end{align}
\end{theorem}
When using the above relation, one implicitly assumes that observed data is used to construct the posterior density. We should note that in our problem of interest, in case that an experimental measurement is not available, a virtual observation $\mathtt{obs}$ resulting from a fine spatial discretization is alternatively employed. Obviously, the observation is more valuable when a real experiment exists.

If we employ the statistical model \eqref{81} with the assumption that errors are Gaussian independent and identically distributed and $\varepsilon_i\sim N(0,\sigma^2)$, where $\sigma^2$ is fixed, then the likelihood function is
\begin{align}\label{83}
\pi(m|\chi)=L(\chi,\sigma^2|m)=\frac{1}{(2\pi\sigma^2)^{n/2}}\exp\left(-\sm/2\sigma^2\right),
\end{align} 
where
\begin{align}\label{84}
\sm=\displaystyle\sum_{j=1}^{n}[m_j-f_j(\xx,\chi)]^2,
\end{align}
is the sum of square errors.

Considering the given likelihood function (\ref{83}), the posterior distribution has the following form:
\begin{align}
\pi(\chi|m)=\frac{\exp\left(-\sm/2\sigma^2_0\right)}{\displaystyle\int_{0}^{\infty} \exp\left(-\sx/2\sigma^2_0\right)d\xi}\,=\frac{1}{\displaystyle\int_{0}^{\infty}\exp\left(-(\sx-\sm)\right)/2\sigma^2_0d\xi},
\end{align}
where $\sx$ is the sum of squares defined by the integration variable; see~\eqref{84}. From a numerical point of view, we can approximate the integral as 
\begin{align}
\pi(\chi|m)\approx\frac{1}{\displaystyle\sum_{i=1}^{n} \exp\left(-(\sx-\sm)\right)/2\sigma^2_0w^i},
\end{align}
where the quadrature points and weights are denoted, respectively, by $\xi^i~\text{and}~w^i$.

In statistics, MCMC methods comprise a class of algorithms for sampling from a probability distribution. By constructing a Markov chain with the desired distribution as its equilibrium distribution, one can obtain a sample of the desired distribution by observing the chain after a number of steps. The more steps there are, the more closely the distribution of the sample matches the actual desired distribution. 

In Bayesian statistics, the recent development of MCMC methods has been a key step in making it possible to compute large hierarchical models that require integration over hundreds or even thousands of unknown parameters. 
In rare event sampling, they are also used for generating samples that gradually populate the rare failure region. 

Below, different popular MCMC methods are reviewed. These methods will be used to identify the parameters in ductile fracture in Section 4. A detailed comparison between the performance of the methods will be given to clarify their efficiency.

\sectpb[MH]{Metropolis and Metropolis-Hasting Algorithms}\label{MH1}
The Metropolis-Hastings (MH) algorithm   is one of the most common techniques among the MCMC methods due to its simplicity for implementation and also its ability to handle different scientific/engineering problems (specifically when the parameters are not strongly correlated) \cite{smith2013uncertainty}. In order to estimate the posterior distribution, in each iteration, a new candidate  parameter value is proposed based on the current sample value according to a proposal distribution. Then, the acceptance ratio is calculated to decide whether the candidate value is accepted or rejected. The acceptance ratio points out how probable the new candidate value is with respect to the current sample. 

The method was first introduced by Metropolis \cite{metropolis1953equation} based on a random walk. The algorithm starts from the initial guess (the prior value) $\chi^0$. Afterwards, according to the chosen proposal distribution a new candidate $\cn$ is proposed, which possibly depends on the previous candidates. Having the new candidate $\cs$, the acceptance rate is calculated~as 
\begin{align}
\lambda(\cn,\cs)=\min\left(1,\frac{\pi(\cs)}{\pi(\cn)}\right).
\end{align} 
As the next step, a random variable $\mathcal{R}\sim \text{Uniform}\,(0,1)$ is produced. If $\mathcal{R}<\lambda$ the candidate is accepted; otherwise, we reject the new proposal and keep the previous candidate in the chain. We follow this procedure for a sufficiently high number of replications. As seen, the algorithm is simple and efficient, specifically when a suitable proposal density is chosen and a large sampling is used. However, an inappropriate proposal results in a significant decrease in performance. If the proposal is very large, many of the candidates will be rejected; therefore, a good convergence to the target density (posterior distribution) will not be achieved. In contrast, if the proposal is too narrow, although many of the candidates are accepted, the chain movement is very slow, and many of the targets will not be captured.  

In the Metropolis algorithm, a symmetric proposal density $\phi(\cs|\cn)=\phi(\cn|\cs)$ is assumed. According to this condition, a movement towards the proposed candidate from the current point is equal to a backward movement (from the current candidate to the proposed point). The use of a non-symmetric proposal distribution was proposed as an efficient improvement by Hastings \cite{hastings1970monte}. Considering $N$ number of samples, the algorithm is summarized in Algorithm 1.

\begin{algorithm}[ht!]
	\label{algorithm}
	\vspace{0.2cm}
	\textbf{Initialization ($j=0$)}: Generate the initial parameter $\chi^0\sim\pi(\chi^0|\,m)$.\\  

	\textbf{while}~{$j<N$}{\\

		\quad 1. Propose the new candidate $\chi^{*}\sim\phi(\chi^j|\,\chi^{j-1})$ ($\phi$ is the proposal distribution).    \\

		\quad 2. Compute the proposal correction parameter~$\beta(\chi^*|\,\chi^{j-1})=\cfrac{\phi(\chi^{j-1}|\,\chi^*)}{\phi(\chi^*|\,\chi^{\ell-1}))}$.
		\\

		\quad  3. Calculate the acceptance/rejection probability $\lambda (\chi^*|\,\chi^{j-1})=\min\left(1, \cfrac{\pi(\chi^*|\,m)}{\pi(\chi^{j-1}|\,m)}\,\beta\right)$.\\
		
		\quad 4. Draw a random number $\mathcal{R}\sim \text{Uniform}\,(0,1)$.\\
		
		\quad 5.$~\textbf{if}~$ $\mathcal{R}<\lambda~$ \textbf{then}\\ 
		
		\qquad\quad accept the candidate $\chi^*$ and set $\chi^j=\chi^*$\\
		
		\qquad \textbf{else}\\
		
		\qquad\quad reject the candidate $\chi^*$ and set $\chi^{j}=\chi^{j-1}$\\
		
		\qquad \textbf{end if}\\
		
	\quad 6. Set $j=j+1$.\\		
		
	}
	\caption{The MH algorithm.}
\end{algorithm}

We can draw the following conclusions: 

\begin{itemize}
	\item A proposal $\cs$ that results in $\pi(m|\cs)>\pi(m|\cn)$ entails a small sum of squared error and thus leads to candidate acceptance.
	\item A proposal $\cs$ that leads to $\pi(m|\cs)<\pi(m|\cn)$  entails a higher 
	sum of squared error and the proposal may be rejected.
\end{itemize}

Regarding the proposal functions and how they affect the posterior distribution, if the variance is too large, a large percentage of the candidates will be rejected, since they will have smaller likelihoods, and hence the chain will stagnate for long periods. The acceptance ratio will be high if the variance is small, but the algorithm will be slow to explore the parameter space.   

There are different measures to determine if the Markov chain is efficiently sampling from the posteriori density. A good criterion is the acceptance rate (the percentage of accepted candidates). The ratio can be used to tune the proposal density, i.e., reduce its variance. Another efficiency test is the autocorrelation function.  The lag-$\tau$ autocorrelation function  $ACF\colon\mathbb{N}\rightarrow [-1,1]$ is estimated by
\begin{equation}
ACF(\tau)=\frac{\displaystyle\sum_{j=1}^{N-\tau} (\chi_j -\bar{\chi})(\chi_{j+\tau} -\bar{\chi}) }{ \displaystyle\sum_{j=1}^{N}\left( \chi_j-\bar{\chi}\right)^2}=\frac{\text{cov}(\chi_j,\chi_{j+\tau})}{\text{var}(\chi_j)}\ge 0.
\end{equation}
Here, $\chi_j$ denotes the $j$-th element of the Markov chain and
$\bar{\chi}$ is the mean value. Note that $ACF(\tau)$ is positive and monotonically decreasing. The interested readers can refer to \cite{khodadadian2020bayesian}, where the authors studied the effect of $ACF$ on different parameters in phase-field modeling of brittle fracture. A more advanced convergence analysis such as $\hat{R}$-statistics can be implemented when multiple MCMCs with different initial values are used. In Section \ref{Section5}, we will use such a diagnostic tool to compare the performance of the Bayesian techniques.

\sectpb[DRAM]{Delayed Rejection Adaptive Metropolis (DRAM)}\label{DRAM}

At this point, it is worth discussing some improvements in the MH algorithm based on the proposal distribution. The main disadvantage of the model is that the covariance of the proposal should be tuned manually. To improve the efficiency, an alternative to using a fixed proposal distribution in each iteration is to update the distribution according to the available samples (adaptive Metropolis). This approach is useful since the posterior distribution is not sensitive to the proposal distribution.

To adapt the proposal function according to the obtained information, Haario \cite{haario2001adaptive} proposed a technique where the current point is chosen as the proposal center and the covariance function is updated using the estimated data. To this end, one can use the following proposal estimation
\begin{align}
V_j=S_p\text{Cov}\left(\chi^0,\chi^1\ldots,\chi^{j-1}\right)+\epsilon I_j,
\end{align}
where the parameter $\epsilon$ is chosen very small (close to zero) and $S_p=\frac{2.38^2}{j}$ (as the scaling parameter). The covariance function is calculated by
\begin{align}
\mathcal{COV}_j=\text{Cov}(\chi^0,\chi^1,\ldots,\chi^{j})=\frac{1}{j}\left(\displaystyle\sum_{i=0}^{j}
\chi^i\left(\chi^i\right)^T-(n+1)~\hat{\chi}^j\left(\hat{\chi}^j\right)^T\right),
\label{cov1}
\end{align}   
where $\hat{\chi}^j=\frac{1}{j+1}\displaystyle\sum_{i=0}^{j}\chi^j$ \cite{smith2013uncertainty}. The proposal adaptation can be done after a specific number of steps (e.g., 1000) instead of all steps.
The efficiency of the algorithm can be further enhanced by adding a delayed rejection step. Mira \cite{green2001delayed} proposed that instead of a rejected candidate, a second stage is used to propose it from another proposal density. As the first step we propose the new candidate using the Cholesky
decomposition of the covariance function \eqref{cov1}:
$$\chi^{*}=\chi^{j-1}+\mathcal{COV}_j\mathcal{U},$$ 
where
$\mathcal{U}\sim \text{Uniform}~(0,I_j)$ and $I_j$ is the $j$-dimensional identity matrix.
The alternative proposal $\chi^{**}$ is chosen using the proposal function
\begin{align}
\phi(\chi^{**}|\chi^{j-1},\chi^*)=\mathcal{N}(\chi^{j-1},\gamma_2^2 V_j),
\end{align}
where $V_j$ is the covariance matrix estimated by the adaptive algorithm \cite{green2001delayed}. The essential parameter is $\gamma_2$, which will be chosen less than one so that the next stage has a narrower proposal function (normally, $\gamma_2=1/5$ is chosen). We use the following acceptance ratio:
\begin{align}\label{821}\nonumber
\lambda_2(\chi^{**}|\chi^{j-1},\chi^*)&:=\min\left(1,\frac{\pi(\chi^{**}|m) \phi(\chi^*|\chi^{**} ) \phi_2(\chi^{j-1}|\chi^{**},\chi^*) [1-\lambda(\chi^*|\chi^{**})]}{\pi(\chi^{j-1}|m) \phi(\chi^*|\chi^{j-1} ) \phi_2(\chi^{**}|\chi^{j-1},\chi^*) [1-\lambda(\chi^*|\chi^{j-1})]}\right)\\
&=\min\left(1,\frac{\pi(\chi^{**}|m) \phi(\chi^*|\chi^{**} )  [1-\lambda(\chi^*|\chi^{**})]}{\pi(\chi^{j-1}|m) \phi(\chi^*|\chi^{j-1} )  [1-\lambda(\chi^*|\chi^{j-1})]}\right).
\end{align}
Then, similar to the MH algorithm, we follow the Markov chain to accept/reject the candidate. A summary of the process is given in Algorithm 2.

\begin{algorithm}[ht!]
	\vspace{0.2cm}
	\label{algorithm1}
	\textbf{Initialization ($j=0$)}: Generate the initial parameter $\chi^0\sim\pi(\chi^0|\,m)$.\\

	\textbf{while~}{$j<N$}{
		\vspace{0.3cm}

		\quad 1. Propose a new candidate $\chi^{*}=\chi^{j-1}+\mathcal{R}_j\mathcal{Z}_j$  where $\mathcal{R}_j$ is the Cholesky decomposition \\
		
		\vspace{-0.3cm} ~~~~~~~of $\mathcal{V}_j$ and $\mathcal{Z}_j\sim$ Uniform (0, $I_j$) where $I_j$ denotes the identity matrix. \\

		\quad  2. Calculate the acceptance/rejection probability $$\lambda_1 (\chi^*|\,\chi^{j-1})=\min\left(1, \cfrac{\pi(\chi^*|\,m)~\phi(\chi^{j-1}|~\chi^*)}{\pi(\chi^{j-1}|\,m)~\phi(\chi^*|\,\chi^{j-1}))}\right)$$.
		
		\quad 3. Draw a random number $\mathcal{R}\sim \text{Uniform}\,(0,1)$.\\
		
		\quad 4.$~\textbf{if}~$ $\mathcal{R}<\lambda_1~$ \textbf{then}\\ 
		
		\qquad\quad accept the candidate $\chi^*$ and set $\chi^j=\chi^*$\\
		
		\qquad \textbf{else}\\
		
		\qquad\quad(i)  Calculate the alternative candidate
		\quad$\chi^{**}=\chi^{j-1}+\sigma^2\mathcal{R}_j\mathcal{Z}_j$.\\
		
		\qquad\quad (ii)   Calculate the acceptance/rejection probability\\ $$\qquad\qquad\lambda_2 (\chi^{**}|\,\chi^{j-1},\chi^*)=\min\left(1, \cfrac{\pi(\chi^{**}|\,m)~\phi(\chi^*|~\chi^{**})\left(1-\lambda_1(\chi^*|\chi^{**})\right)}{\pi(\chi^{j-1}|\,m)~\phi(\chi^*|~\chi^{j-1})\left(1-\lambda_1(\chi^*|\chi^{j-1})\right)}\right)$$.\\
		
		\vspace{-0.3cm}\qquad\quad (iii)  \textbf{if} ~$\mathcal{R}<\lambda_2$~~\textbf{then}\\
		
		\qquad\qquad\quad~~accept the candidate $\chi^{**}$ and set $\chi^{j}=\chi^{**}$\\
		
		\qquad\quad~~~   \textbf{else}\\
		
		\qquad\qquad\quad~~  reject the candidate $\chi^{**}$ and set $\chi^{j}=\chi^{j-1}$
		\qquad\quad\\
		
		\qquad \textbf{end if}\\
		
		\quad 5. Update the covariance matrix as $\mathcal{V}_j=\text{Cov}(\chi^0,\chi^1,\ldots,\chi^{j})$.\\
		
		\quad 6. Update $\mathcal{R}_j$.\\

				\quad 7. Set $j=j+1$.\\
	}
	\caption{The DRAM algorithm.}
\end{algorithm}

\sectpb[kalman]{MCMC with ensemble-Kalman filter}\label{kalman}
As previously mentioned, a good detection of the proposal density will enhance the Markov chain movement to the target density. Here, we introduce another proposal distribution detection technique using an {\color{black}ensemble-Kalman filter} to obtain
\begin{align}
	\cs=\cn+\Delta \chi,
	\label{prop1}
\end{align}
where $\Delta \chi$ denotes the jump of Kalman-inspired proposal. In order to update the proposal, we can separate \eqref{prop1} into
\begin{align}
\Delta \chi=\mathcal{K}\left( y^{j-1}+s^{j-1}\right).
	\label{prop2}
\end{align} 
The first term indicates the so-called Kalman gain, i.e., 
\begin{align}
\mathcal{K}=\C_{\chi M}\left(\C_{MM}+\C_{M}\right)^{-1},
\end{align}
where $\C_{\chi M}$ is the covariance matrix between the inferred parameters and the PDE-based model, $\C_{MM}$ is the covariance matrix of the PDE response, and $\C_{M}$ denotes the measurement noise covariance matrix \cite{zhang2020improving}. In \eqref{prop2}, $y^{j-1}$ is the residual of candidates with respect to the model. In other words, considering $\bar{m}$ an observation/measurement, $y^{j-1}=\bar{m}-f(\chi^{j-1})$ and  $s^{j-1}\sim\N(0,\mathcal{R})$, related to the density of measurement. 

Considering the ductile fracture process, crack propagation is modeled until full fracture has occurred. To highlight this procedure, we have added criterion (iv) to Algorithm 3. For the two other algorithms (MH and DRAM techniques), the same condition can be considered.

\begin{algorithm}[ht!]
	\vspace{0.2cm}
	\textbf{Initialization ($j=0$)}:  initiate the samples according to the prior density $\chi^0$\\[2mm]
	\textbf{while}~$j<N$ 	\\[3mm]
	\vspace{0.1cm}
	\quad 1. $\bullet$~set \var{FLAG}=true \qquad $\bullet$~set $ n=0$\\[3mm]
	\vspace{0.1cm}
	\textbf{While}~\var{FLAG}~\textbf{do}\\[3mm]
	\vspace{0.1cm}
	\qquad\quad(i) Solve the model equations $\mathcal{M}_1$, $\mathcal{M}_2$, and $\mathcal{M}_3$ considering   $\texttt{TOL}_\mathrm{Stag}$ and\\[3mm]
			\vspace{0.1cm}
	\hspace{1.7cm}  the proposed candidate $\cn$ an then obtain $f(\cn)$\\[3mm]
		\vspace{0.1cm}
		\qquad\quad(ii) update the Kalman gain			
	\begin{align*}
\mathcal{K}=\C_{\chi M} \left(\C_{MM}+\C_{M}\right)^{-1}
	\end{align*} 
	\vspace{0.1cm}
	\qquad \qquad \quad$\bullet$ $\C_{M}$ is the measurement noise covariance matrix\\[3mm]
		\vspace{0.1cm}
	\qquad \quad \hspace{-0.0cm}(iii) shift the ensemble 
\begin{align*}
	\cs=\cn+\K\left( y^{j-1}+s^{j-1}\right) 
\end{align*}
	\qquad \qquad ~$\bullet$ $y^{j-1}\sim\N(\varepsilon^{j-1},\mathcal{R})$ is the residual of the proposed parameters\\[3mm]
	\vspace{0.1cm}
	\qquad \qquad \,\,$\bullet$ $s^{j-1}\sim\N(0,\mathcal{R})$ relates to the measurment error\\[3mm]
		\vspace{0.1cm}
	\qquad \quad \hspace{0.0cm}(iv) \textbf{if} full fracture is occurred \textbf{then}\\[3mm]
		\vspace{0.1cm}
	\qquad \qquad \qquad ~~~~$\bullet$ set \var{FLAG}=false\\[3mm]
	\qquad \qquad~	
	\vspace{0.1cm}
	\qquad \quad \hspace{0.7cm} \textbf{else} \\[3mm]
	\vspace{0.1cm}
	\qquad \qquad \qquad ~~~ $\bullet$ set~ $n=n+1$\\[3mm]
	\qquad \qquad~	
			\vspace{0.1cm}
\qquad \quad \hspace{0.7cm} 	\textbf{end if}\\[3mm]			
	\quad 2. Accept/reject the material approximation.\\[3mm]
	\quad 3. Set $j=j+1$.\\
	\caption{Bayesian inversion with ensemble-Kalman filter.} 
	\label{algorithm3}	
\end{algorithm}
\sectpb[bayesian]{Bayesian inversion for ductile phase-field fracture}
\label{bayesian}
A proper knowledge about the mechanical parameters that influence the behavior of fracturing solids is crucial to observe the model response and precisely predict crack initiation and propagation during different stages of the deformation process. Bayesian inversion techniques are convenient tools to monitor the crack behavior using observations (e.g., the measured data) and solving the inverse problem considering the forward model (here $\calM_1$, $\calM_2$, and $\calM_3$). Below, we review different possibilities for the implementation of Bayesian inversion in the context of elastic-plastic fracturing solids governed by phase-field models.
\begin{itemize}
	\item \textbf{Based on the load-displacement curve}: this approach allows us to observe the crack behavior in all time steps up to complete failure. At time-step $n$, the load-displacement curve can be computed as
	\begin{align}
	F_n=\int_{\partial_D\mathcal{B}} \boldsymbol{n}\cdot\boldsymbol{\sigma}\cdot \boldsymbol{n}\,da,
	\end{align}
where $\Bn$ is the outward unit normal on the surface, defined in \eqref{s2-bcs}. The main advantage of working with this curve is its easiness, since it involves a one-dimensional parameter. However, it is sensitive to the mesh size and the length scale; therefore, a sufficiently small (and thus more computationally expensive) mesh size is needed. 
	\item \textbf{Based on the point-wise primary fields}: this approach  monitors the crack behavior, the displacement, and the equivalent plastic strain in the entire geometry. Here, the Bayesian setting strives to find the inferred parameters $\chi^\star$ which minimize
	\begin{align*}
	\|\bar{\boldsymbol{u}}(\boldsymbol{x})-\boldsymbol{u}(\boldsymbol{x},\chi^\star)-\varepsilon_1\mathcal{I}\|^2+\|\bar{d}(\boldsymbol{x})-d(\boldsymbol{x},\chi^\star)-\varepsilon_2\mathcal{I}\|^2+\|\bar{\alpha}(\boldsymbol{x})-\alpha(\boldsymbol{x},\chi^\star)-\varepsilon_3\mathcal{I}\|^2,
	\end{align*}
where $L^2$-norm can be used, and $\bar{\boldsymbol{u}}$, $\bar{d}$, and $\bar{\alpha}$ are the experimental data throughout geometry with respective measurement errors $\varepsilon_1$,  $\varepsilon_2$, and $\varepsilon_3$. This method is informative and provides precise information since the displacement and phase-filed in the entire geometry are considered. However, it is difficult and perhaps even impossible to obtain the measured data from actual experiments in a point-wise manner. Furthermore, a  small mesh size must be chosen in numerical simulations to guarantee accurate estimations in the whole geometry, further rendering the method computationally prohibitive. 
		
\item \textbf{Based on the point-wise phase-field propagation}: this approach is less complex than the previous method. Reliable experimental values of the crack path can be obtained using X-ray or $\mu$-CT scan in two- or three-dimensional problems; see, e.g., \cite{buljac2018calibration}. However, the computational issue regarding mesh sensitivity persists.
	
	\item \textbf{Based on snapshots of proper orthogonal decomposition (POD)}: this approach employs a reduced order method (ROM) to reduce the computational complexity. Here, the snapshots of the solution (using measurements of the crack phase-field) are used to construct the POD basis \cite{abbaszadeh2021reduced}. If an efficient ROM is used, the computational complexity can be reduced significantly. Similarly, a Global-Local approach \cite{noii2021glductile} can be employed to reduce the computational complexity of the forward model. 
		
	\item \textbf{Based on the effective stress-strain response}: this approach  explains the relation between $\boldsymbol{\varepsilon}$ and $\boldsymbol{\sigma}$. It provides useful information regarding different material properties such as bulk modulus, hardening, and yield strength. Therefore, considering the availability of measurements, it entails an instructive procedure. Nevertheless, the computational costs, i.e., the effect of the mesh size and length scale, must be taken into account.
\end{itemize}   
In this work, we choose the load-displacement curve as the observation, and sufficiently small mesh sizes are used to model the crack propagation. Nevertheless, it is worth noting that POD-ROM approaches have significant simulation advantages (noticeable computational cost reduction), while methods based on the stress-strain response are very informative. These procedures will be addressed in future works.   

\begin{figure}[!]
	\centering
	{\includegraphics[clip,trim=0cm 21cm 0cm 8cm, width=16cm]{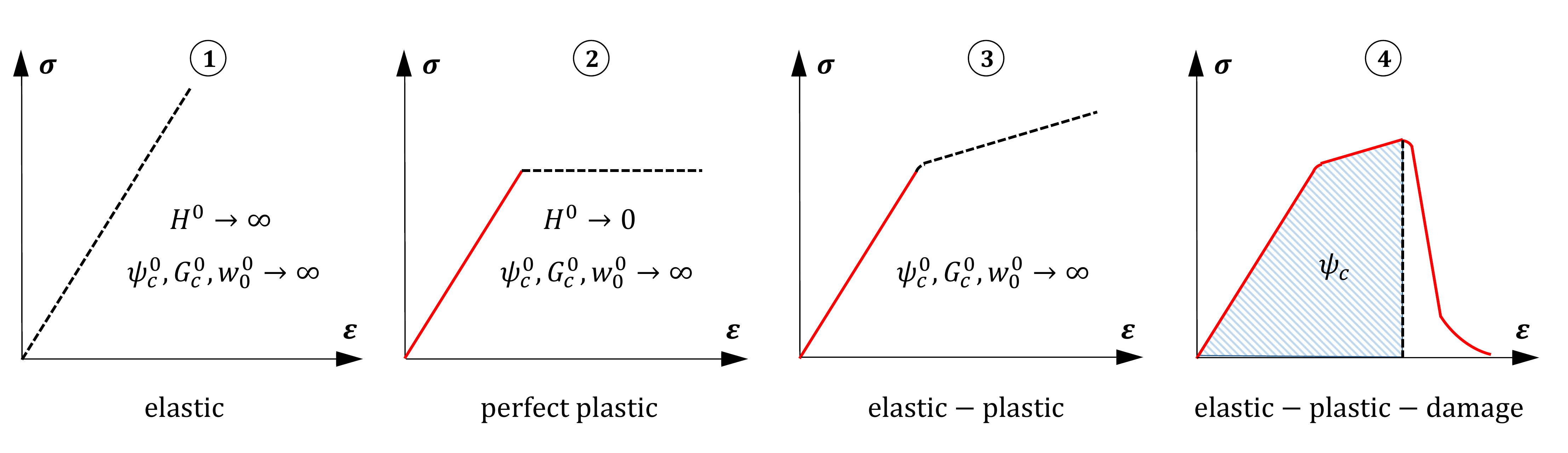}}  
	\caption{Step-wise Bayesian inversion method to determine the posterior density of the material unknowns for ductile phase-field fracture models.}
	\label{FigureBA_ductile}
\end{figure}

Now, we proceed to establish a Bayesian inversion ($\texttt{BI}$) setting to identify the different parameters in ductile fracture. Let us assume that the response of ductile phase-field fracture is either elastic, followed by elastic-plastic, followed by elastic-damage (hereafter $\texttt{E-P-D}$); or elastic, followed by elastic-plastic, followed elastic-plastic-damage (hereafter $\texttt{E-P-DP}$). Next, we aim to determine the candidate $\chi\in(\mu,K,H,\sigma_Y,\psi_c,G_c,w_0,l_p)$ as follows:

\begin{itemize}
	\item[(1)] To find $\tilde{\mu}$ and $\tilde{K}$, we set $H^0\rightarrow\infty$, $l^0_p\rightarrow0$ (in case of $\mathcal{M}_3$), $\chi^0_a\rightarrow0$ (in the anisotropic case), and $G^0_c\rightarrow\infty$ in $\mathcal{M}_1$, $\psi^0_c\rightarrow\infty$ in $\mathcal{M}_2$, and $w^0_0\rightarrow\infty$ in $\mathcal{M}_3$, thus reflecting an elastic response. We then have
	\begin{equation}
	(\tilde{\mu},\tilde{K})=\texttt{BI}\,(\mu^0,K^0,\sigma_Y^0,H^0,l^0_p,G^0_c,\psi^0_c,w^0_0,\chi^0_a).
	\end{equation}
	\item[(2)] To find $\tilde{\sigma_Y}$, we set $H^0\rightarrow 0$, $l^0_p\rightarrow0$ (in case of $\mathcal{M}_3$), $\chi^0_a\rightarrow0$ (in the anisotropic case),  and $G^0_c\rightarrow\infty$ in $\mathcal{M}_1$, $\psi^0_c\rightarrow\infty$ in $\mathcal{M}_2$, and $w^0_0\rightarrow\infty$ in $\mathcal{M}_3$, thus reflecting an ideal plastic response. We then have
	\begin{equation}
	\tilde{\sigma_Y}=\texttt{BI}\,(\tilde{\mu},\tilde{K},Y^0_0,H^0,l^0_p,G^0_c,\psi^0_c,w^0_0,\chi^0_a).
	\end{equation}
	\item[(3)] To find $\tilde{l_p}$ (in case of $\mathcal{M}_3$) and $\tilde{H}$, we set $\chi^0_a\rightarrow0$ (in the anisotropic case), $G^0_c\rightarrow\infty$ in $\mathcal{M}_1$, $\psi^0_c\rightarrow\infty$ in $\mathcal{M}_2$, and $w^0_0\rightarrow\infty$ in $\mathcal{M}_3$, thus reflecting an elastic-plastic response prior to fracture. We then have
	\begin{equation}
	(\tilde{H},\tilde{l_p})=\texttt{BI}\,(\tilde{\mu},\tilde{K},\tilde{\sigma_Y},H^0,l^0_p,G^0_c,\psi^0_c,w^0_0,\chi^0_a).
	\end{equation}
	\item[(4)] To find $\tilde{\chi_a}$ (in the anisotropic case), $\tilde{G_c}$ in $\mathcal{M}_1$, $\tilde{\psi_c}$ in $\mathcal{M}_2$, and $\tilde{w_0}$ in $\mathcal{M}_3$, which reflect a ductile anisotropic fracture response, we have
	\begin{equation}
	(\tilde{G_c},\tilde{\psi_c},\tilde{w_0},\tilde{\chi_a})=\texttt{BI}\, \left(\tilde{\mu},\tilde{K},\tilde{\sigma_Y},\tilde{H},\tilde{l_p},G^0_c,\psi^0_c,w^0_0,\chi^0_a\right).
	\end{equation}
	\item[(5)] Finally, we obtain the following parameter estimation:
	\begin{equation}
	\boxed{
	\label{all}
	(\mu,K,\sigma_Y,H,l_p,G_c,\psi_c,w_0,\chi_a)=\texttt{BI}\, \left(\tilde{\mu},\tilde{K},\tilde{\sigma_Y},\tilde{H},\tilde{l_p},\tilde{G_c},\tilde{\psi_c},\tilde{w_0},\tilde{\chi_a}\right).}
	\end{equation}
\end{itemize}
Figure \ref{FigureBA_ductile} shows the overall procedure, indicating all stages of the deformation process. Note that, from the implementation point of view, we set the limit $\infty$ as $10^8\times E$, where $E$ refers to Young's modulus, while the lower limit is set to 0. From the statistical point of view, we employ the MCMC techniques (Algorithms 1--\ref{algorithm3}) to identify the parameters in Step (1); then, we follow Step (2) to estimate $\sigma_Y$ and pursue the parameter identification procedure until we determine the whole set of material parameters in Step (4). 

In a similar manner, if the response of ductile phase-field fracture is elastic, followed by elastic-damage, followed by elastic-plastic-damage, i.e., $\texttt{E-D-PD}$ (see \cite{alessi2017} for a detailed discussion on different possible evolutions), we first determine the elastic moduli, followed by the the anisotropic fracture properties by assuming that the response is brittle, and finally, we determine the plastic proprieties in the final dissipative stage.

%
%
%
%
\sectpa[Section5]{Numerical examples}
\label{Section5}
This section demonstrates the performance of the proposed Bayesian inversion approaches for parameter estimation within the ductile phase-field fracture models presented earlier. We investigate four numerical examples. To validate the numerical method, the last two examples are concerned with experimental observations, in which the posterior responses are compared with experimental load-displacement curves. The material parameters listed in Table~\ref{material-parameters} are considered, which are initialized based on \cite{ambati2016,ambati2015}.  In the MCMC  the observational noise $\sigma^2=10^{-3}$ is used.

\textbf{Space discretization}. In the numerical simulations, the global primary variables are discretized using finite element basis functions, with bilinear quadrilateral $Q_1$ elements for the two-dimensional problems and trilinear hexahedral $H_1$ elements for the three-dimensional problems.

\textbf{Solution of the nonlinear problems}.  A staggered scheme is used for solving the variational equations resulting from the ductile phase-field fracture models (Section \ref{sec:models}). For model 1 $(\calM_1)$, and model 2 $(\calM_2)$, we alternately {solve} for $d/{\bm u}$ by fixing ${\bm u}/d$ until convergence is reached. Accordingly, for model 3 $(\calM_3)$, we alternately {solve} for ${\bm u}$ by fixing $(\alpha,d)$, and then solve for $\alpha$ by fixing $(\bm u,d)$. Next, we obtain the plastic strain tensor $\Bve^p$ though the incremental plastic evolution equation \req{plastic_law}, and lastly, we find $d$ by fixing $(\bm u,\alpha)$, repeating the procedure until convergence is reached. 

At this point, it is necessary to remark on the convergence criteria for the staggered scheme. Let $n$ and $k$ represent the loading time step and iteration counter of the alternate minimization scheme, respectively. At the fixed loading time step $n$, we obtain a converged state if the following holds: 
\begin{equation}
	\|\calM_\bullet(\bm u^{n,k},\alpha^{n,k},d^{n,k})\|\leq \texttt{Tol}_\texttt{stag.} \WITH \bullet \in\{1,2,3\} \AND
	\texttt{Tol}_\texttt{stag.}\approx10^{-3}.
\end{equation}
Additionally, an iterative Newton solver is used in which the nonlinear equation systems are solved.
The {stopping criterion of the single scale
	and local Newton methods is $\texttt{Tol}_\texttt{N-R}=10^{-10}$. Specifically,
	the relative residual norm is given by $\texttt{Residual}: \| \bm F(\bm x_{k+1}) \|
	\leq \texttt{Tol}_\texttt{N-R} \| \bm F(\bm x_{k}) \|$.}
Here, $\bm F$ refers to the residual of the equilibrium equation of {the
	nonlinear} single scale and local BVPs. The interested reader can refer to \cite{ZiShen16,FaMau17,HeiWi18_pamm,KoKr20,JoLaWi20,JoLaWi20_parallel,graeser2021truncated}
 for the developed linear/nonlinear solvers for phase-field fracture.

\setcounter{table}{0}
\begin{table}[!ht]
	\caption{Material parameters used in the numerical experiments. The fixed values in all examples are mentioned. Other parameters are inferred with Bayesian inversion.}	\vspace{1mm}
	\centering
	\begin{tabular}{cllcccccc}
	  Parameter & Name                  & Unit       & Value     \\[2mm]\hline 
	   $\mu$        & shear modulus    &     $\mathrm{MPa}$ &$\mathrm{BI}$ \\[2mm]
	  $K$      & bulk modulus    &   $\mathrm{MPa}$     &$\mathrm{BI}$            \\[2mm]
	 $H$        & hardening modulus         & $\mathrm{MPa}$ &$\mathrm{BI}$\\[2mm]
	  $\sigma_Y$        & yield stress      & $\mathrm{MPa}$&$\mathrm{BI}$ \\[2mm]
	 	 	 	 $\alpha_{\text{crit}}$        & hardening critical value      & --&$\mathrm{BI}$ \\[2mm]
		 $\psi_c$     & specific fracture energy    & $\mathrm{MPa}$ &$\mathrm{BI}$\\[2mm]
	 $G_c$     & Griffith's energy release rate     & $\mathrm{MPa}$ &$\mathrm{BI}$\\[2mm]
	 $w_0$     & specific fracture toughness    & $\mathrm{MPa}$ &$\mathrm{BI}$\\[2mm]
	 		 		 $\zeta$    &driving scaling factor&--&$\mathrm{BI}$\\[2mm]
	 		 		 		 $\chi_a$ & stiffness parameter &--&$\mathrm{BI}$\\[2mm]
		 $\eta_f$    & crack viscosity       & $\mathrm{N/m^{2}s}$&$10^{-9}$\\[2mm] 
		 $\eta_p$    & plasticity viscosity       & $\mathrm{N/m^{2}s}$&$10^{-9}$\\[2mm] 
	 $\kappa$    & stabilization parameter   & -- & $10^{-8}$& & \\[2mm]
		 $l_d$    & fracture length-scale  &  $\mathrm{mm}$  & $10^{-8}$ & \\[2mm]
		 $l_p$    & plastic length-scale     &$\mathrm{mm}$& $10^{-8}$ \\[2mm]
		 $\phi$       & fiber rotation & degree & $[30^{\circ}, 45^{\circ}, 60^{\circ}]$\\[2mm]
		\hline
		\label{material-parameters}
	\end{tabular}
\end{table}

\begin{figure}[!]
	\centering
	{\includegraphics[clip,trim=5cm 15cm 18cm 7.5cm, width=10cm]{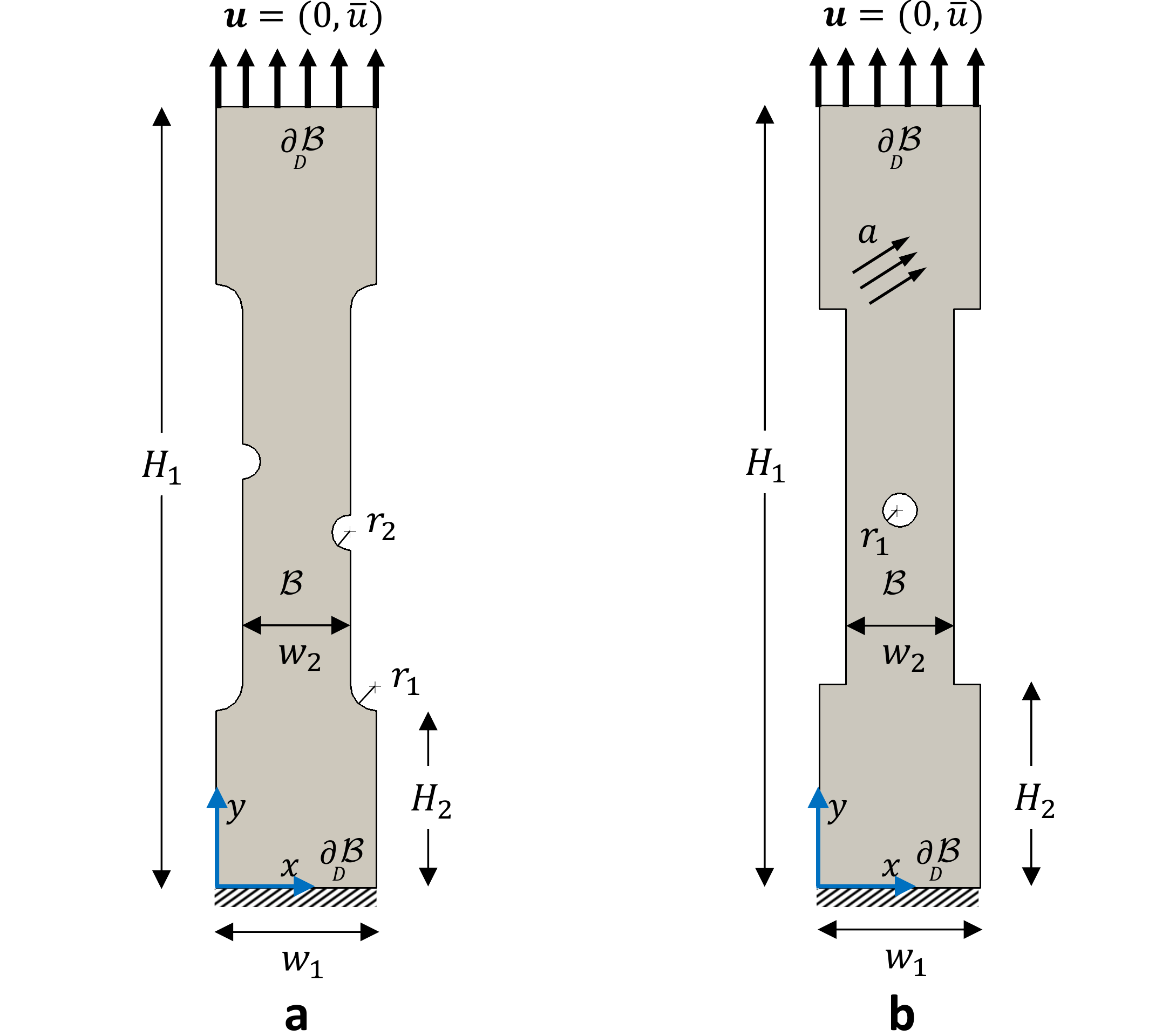}}  
	\caption{ Geometry and boundary conditions of the I-shaped tensile specimen: (a) Example 1, and (b) Example 2.}
	\label{Figure1bvp}
\end{figure}

\begin{table}[!]
	\caption{Example 1: The uniform prior distribution of the inferred parameters.}
	\vspace{1mm}
	\centering
	\begin{tabular}{llcccccccccccccc}
		&Parameter  &$H$        & $\mu$    &   $K$   &   $\sigma_Y$   &$G_c$&$\alpha_{\text{crit}}$     &$\psi_c$   &$w_0$ &$l_p$ &$\zeta$ \\[2mm]\hline\\
		&   min   &150   & 20\,000  & 40\,000 &275 & 5 & 0.05 &20&20&0.5& 0.25& \\[4mm]
		& max &   375  &  40\,000  &  100\,000 & 400 &15 &0.2 & 60 &60&2.5& 10&\\[4mm]
		&   initial  &220  & 25\,000   & 80\,000 &350 & 12 & 0.12 &30& 25 & 1.2&2  \\[4mm]
		\hline
		\label{range1}
	\end{tabular}
\end{table}
\sectpb[Section51]{Example 1: Asymmetrically I-shaped specimen under tensile loading}
To gain a first insight into the performance of the Bayesian inversion approach, the following numerical example is concerned with the asymmetrically notched I-shaped specimen under tension.  The configuration is shown in Figure \ref{Figure1bvp}a. The geometrical dimensions are set as $H_1=110$ mm, $H_2=25$ mm, $r_1=3.625$ mm, $w_1=22$ mm, and $w_2=14.8$ mm, with half-circular notches of radius $r_2=2.5$ mm. The two notches are placed at a vertical distance from the center of $10$ mm.

A monotonic displacement increment  ${\Delta \bar{u}}_y=2\times10^{-3}$ mm is applied in a vertical direction at the top boundary of the specimen. The minimum finite element size in the domain is $0.3$ mm, for which the heuristic requirement $h<l/2$ inside the localization zone is fulfilled. Consequently, the I-shaped domain partition contains 21598 elements. The material and numerical parameters are given in Table \ref{material-parameters} and Table \ref{range1}, respectively. 

\begin{figure}[t!]
	\vspace{0.1cm}
	\hspace{-1cm}
	\subfloat{\includegraphics[width=0.37\textwidth]{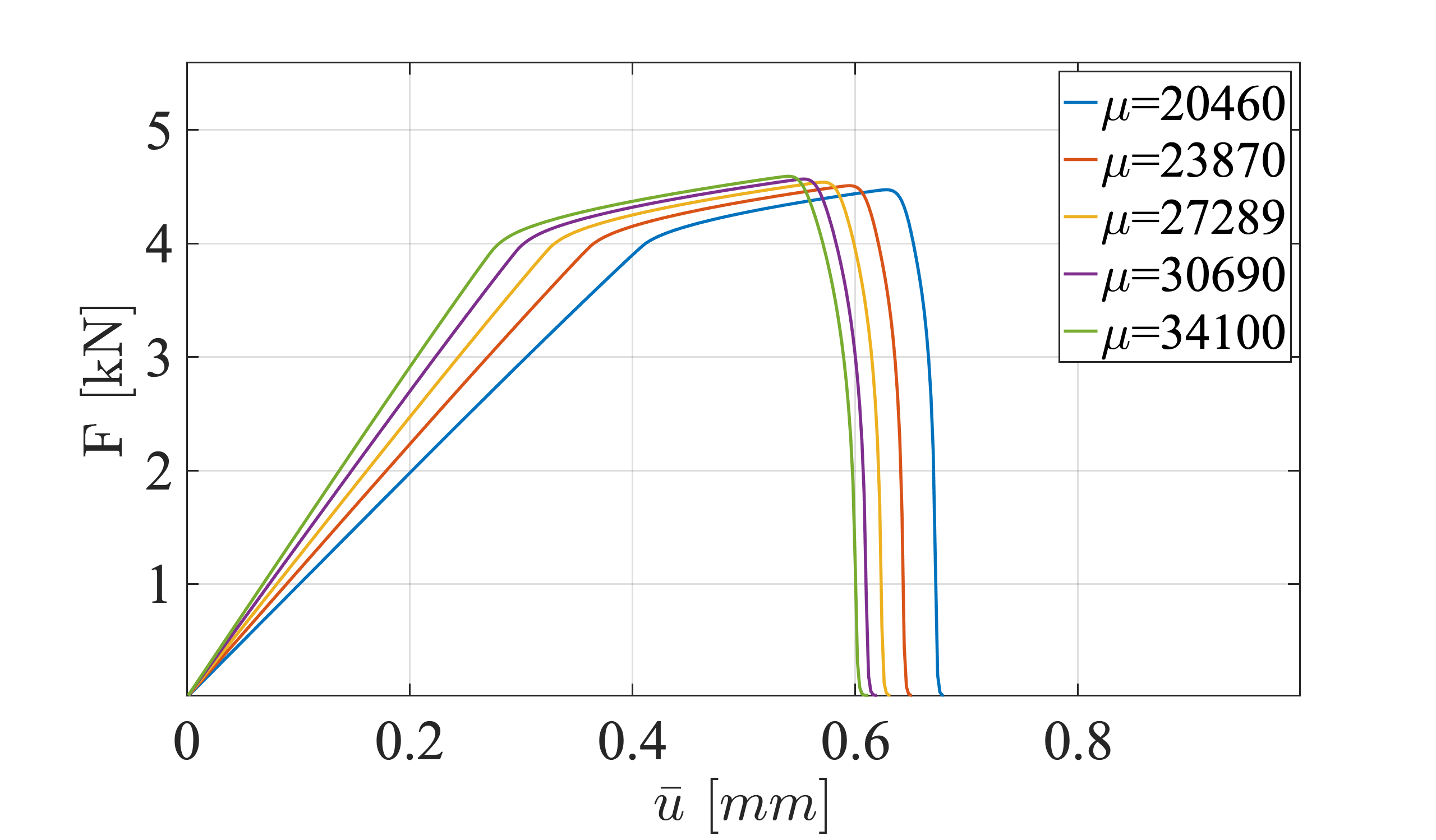}}   \subfloat{\includegraphics[width=0.37\textwidth]{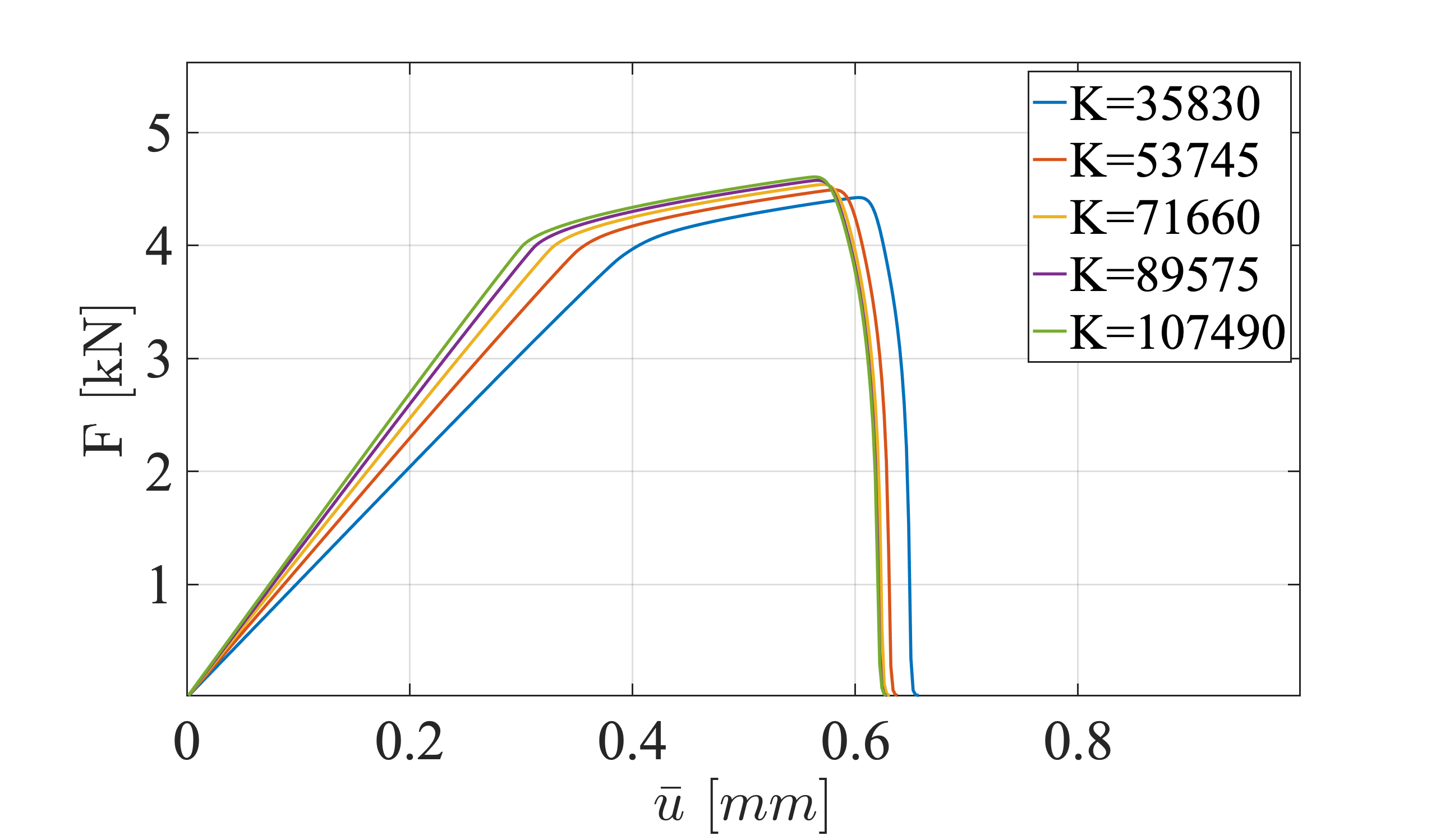}}  	\subfloat{\includegraphics[width=0.37\textwidth]{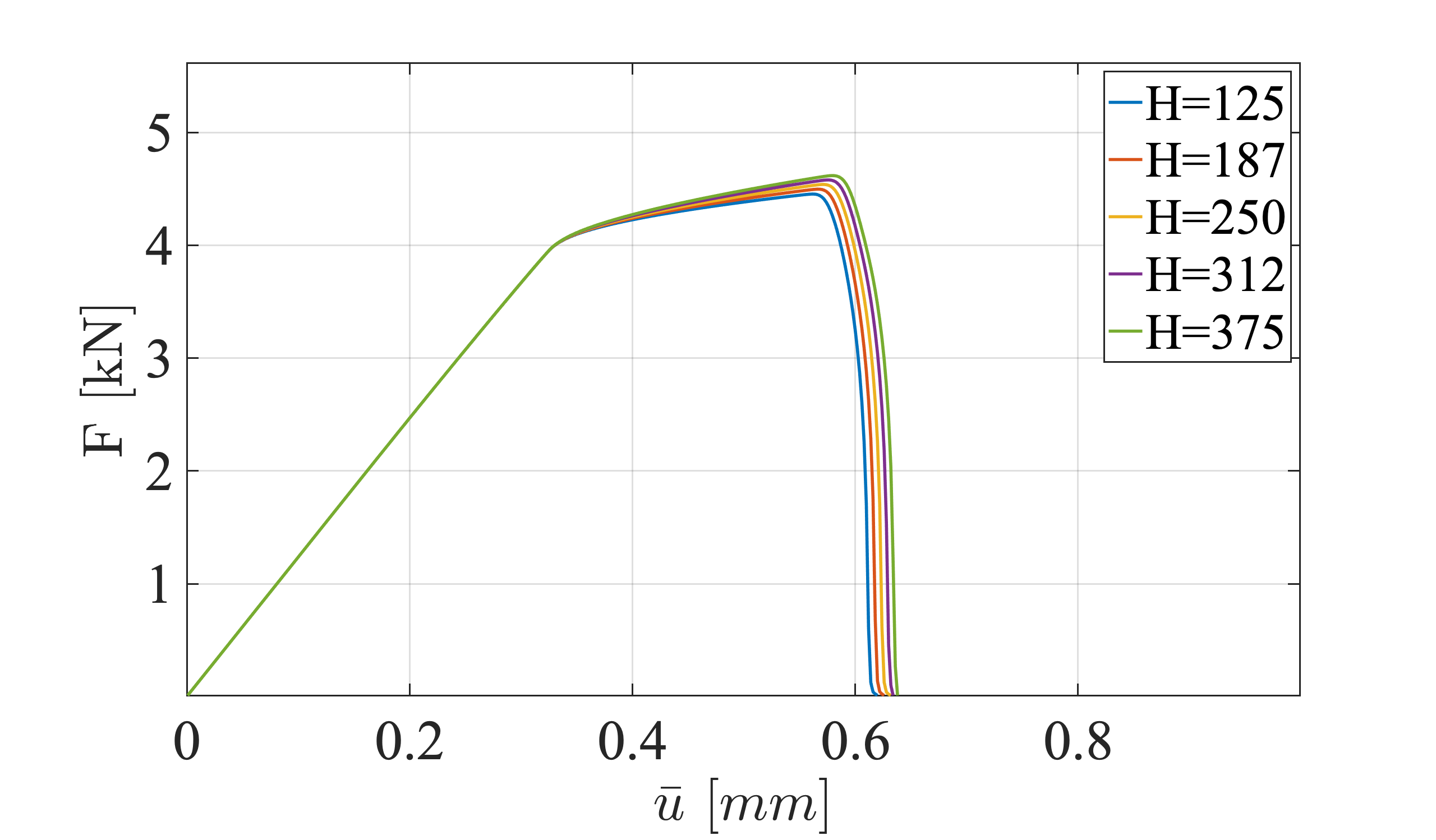}}
	\hfill 
	\vspace{0.1cm}
	\hspace{-1.1cm}
	\subfloat{\includegraphics[width=0.37\textwidth]{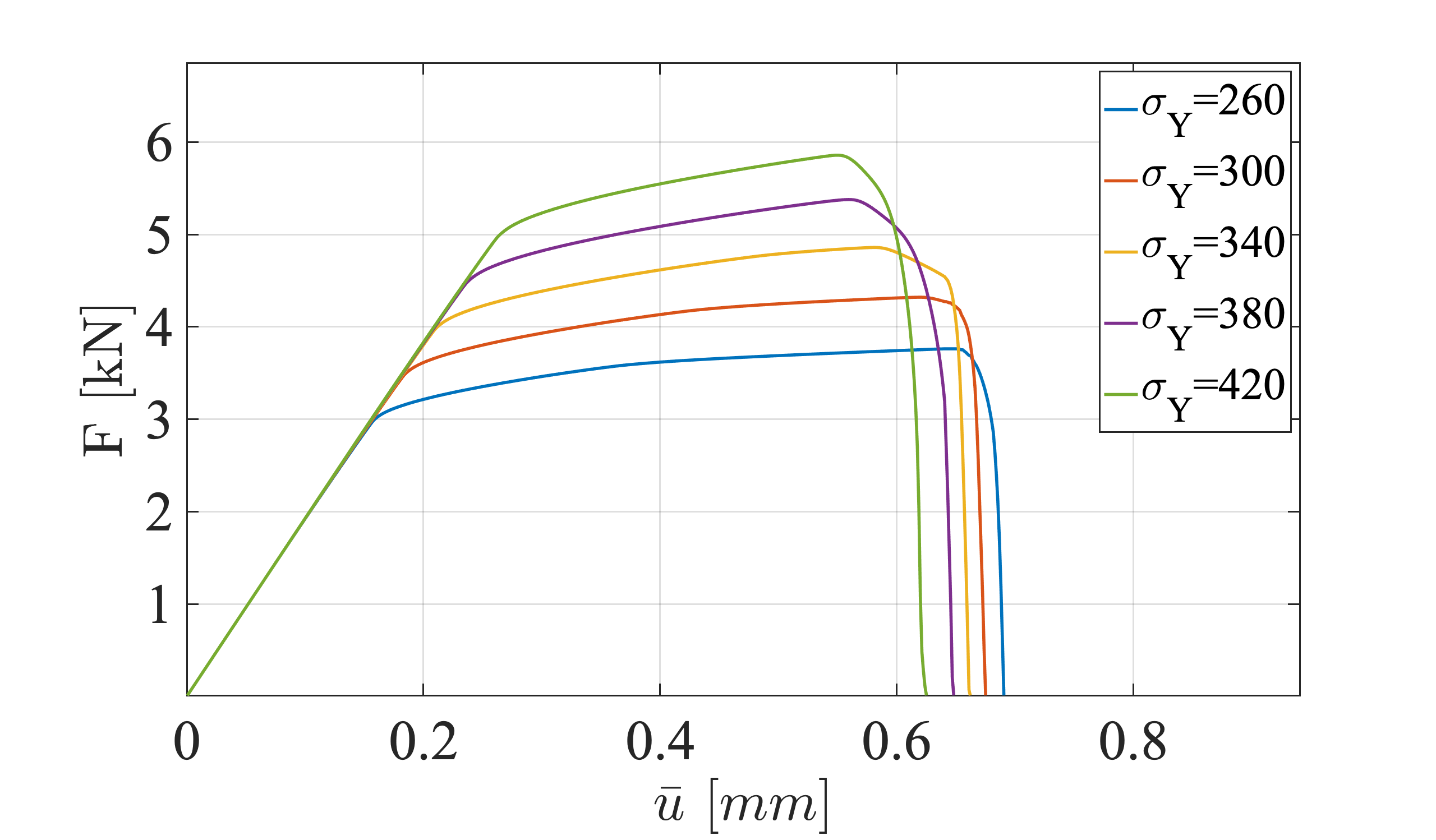}}   \subfloat{\includegraphics[width=0.37\textwidth]{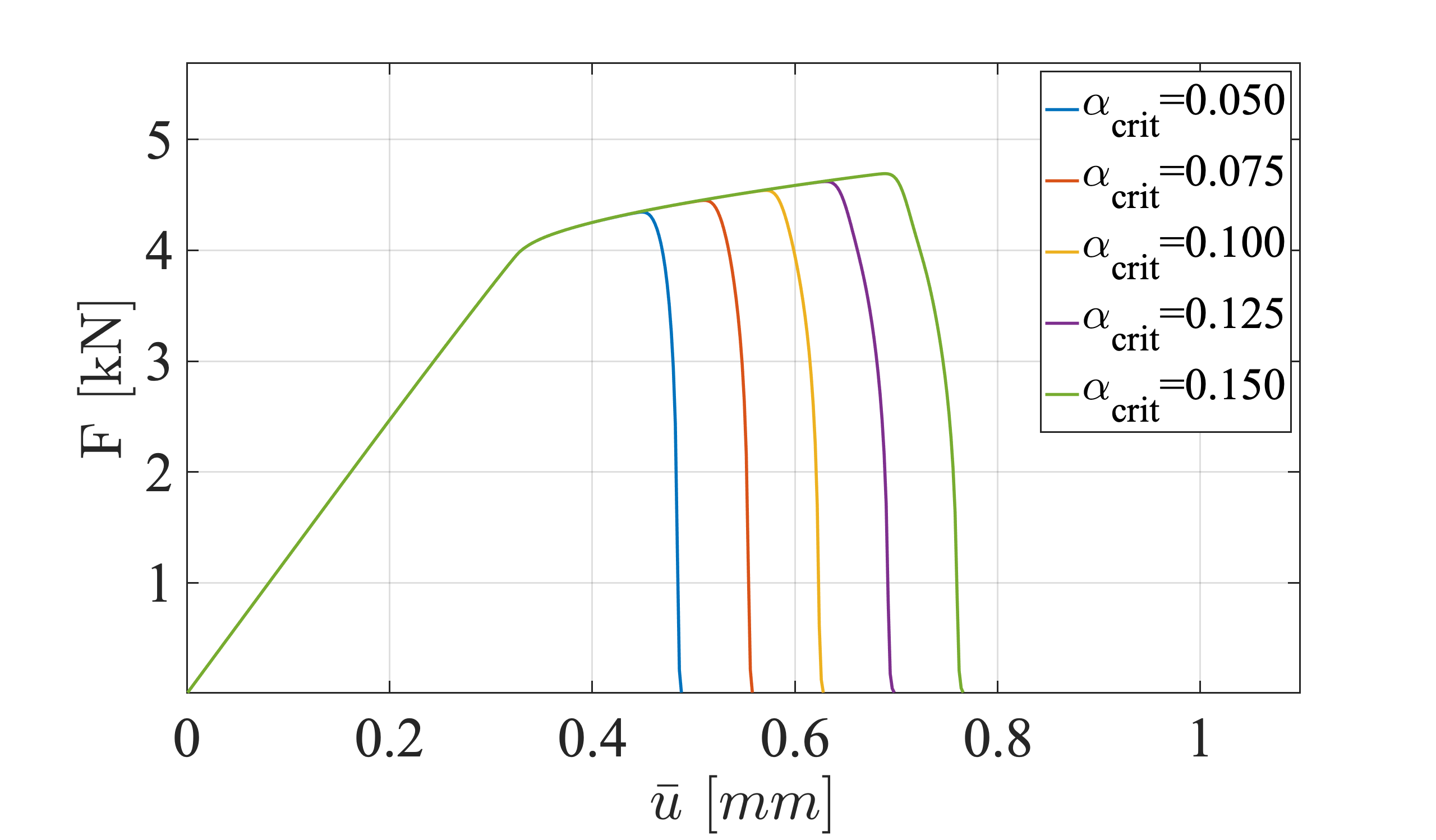}}  	\subfloat{\includegraphics[width=0.37\textwidth]{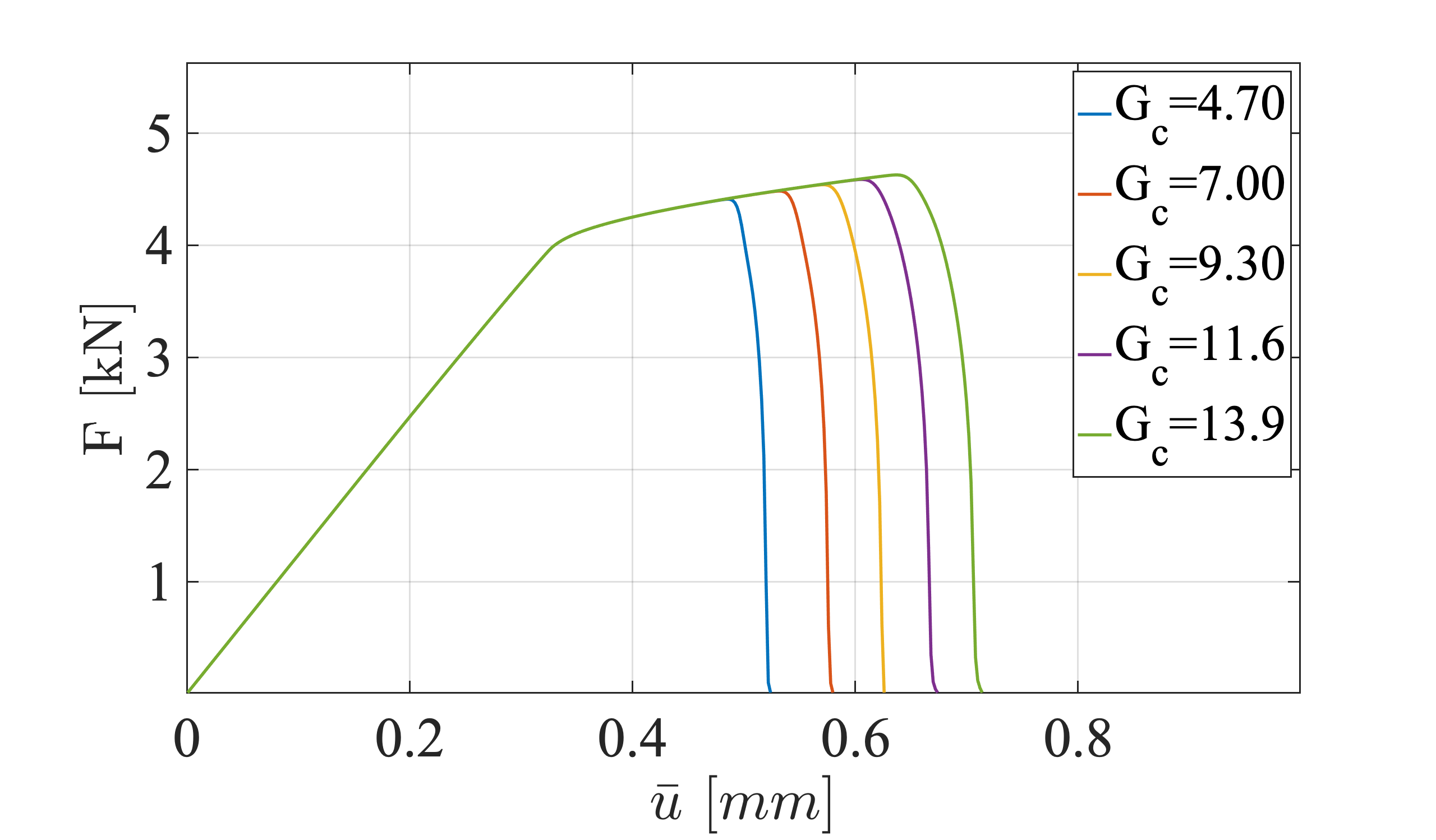}}
	\hfill 
	\vspace{0.1cm}
	\hspace{-1.1cm}
	\subfloat{\includegraphics[width=0.37\textwidth]{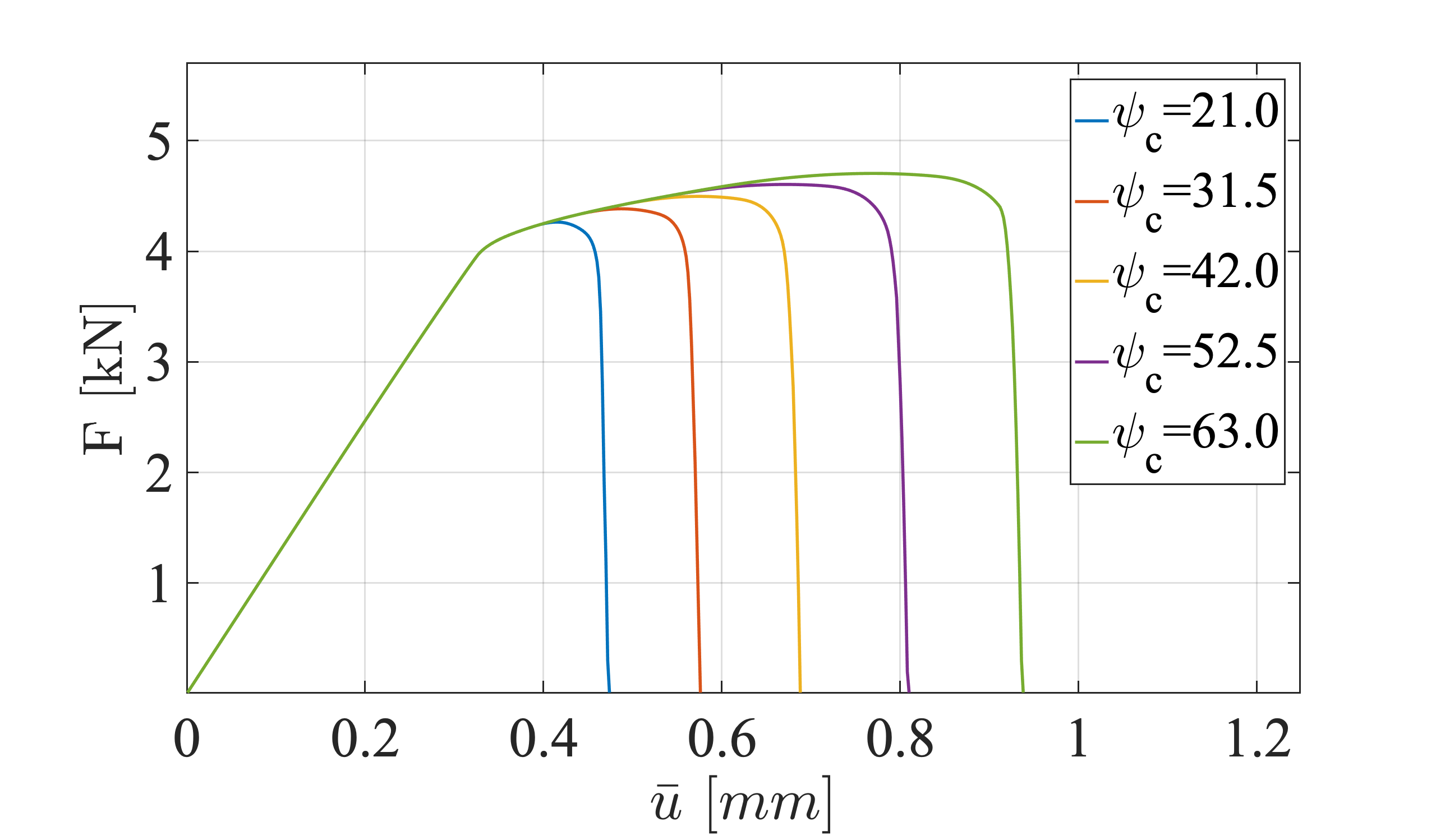}}   \subfloat{\includegraphics[width=0.37\textwidth]{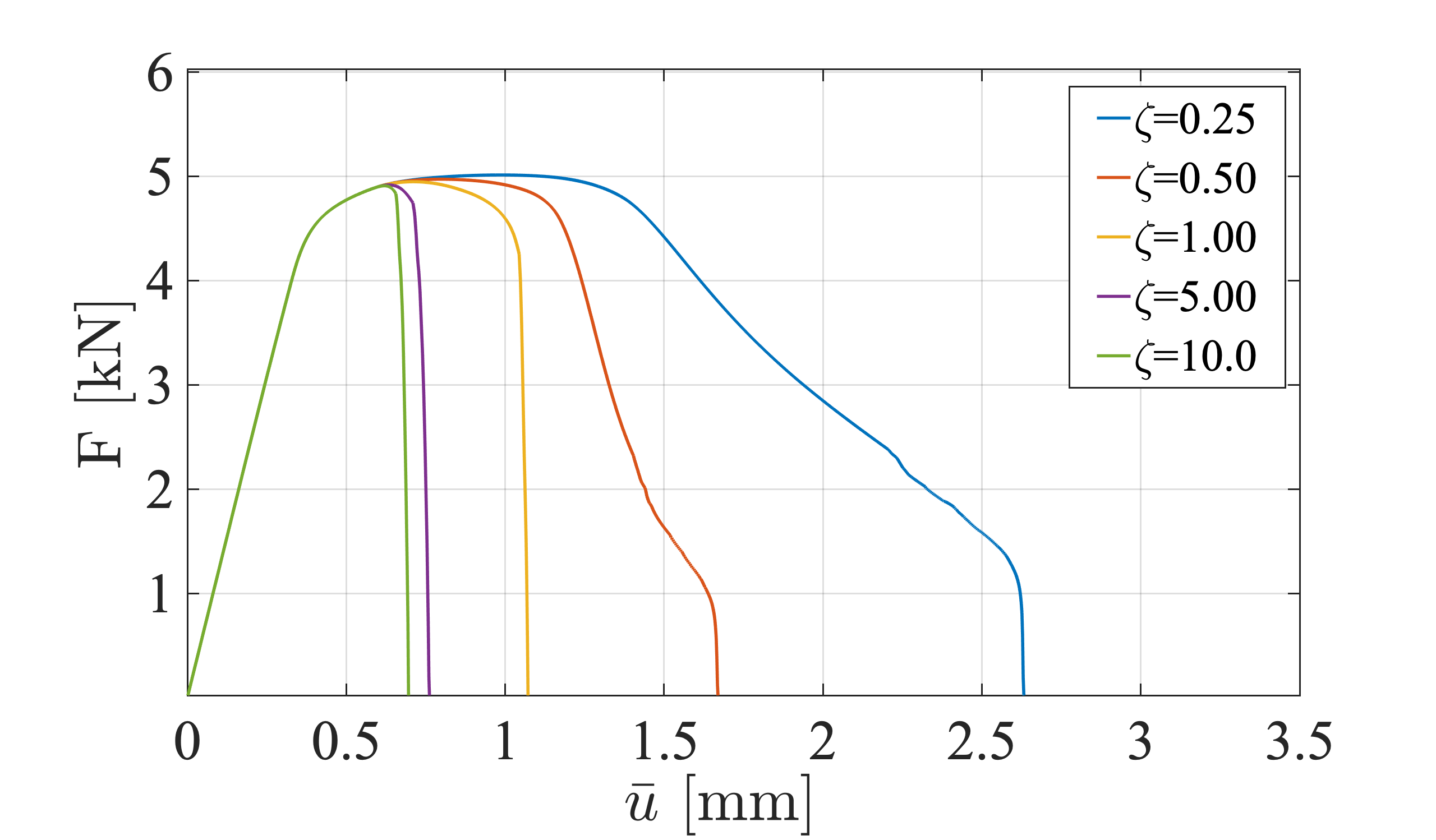}}  	\subfloat{\includegraphics[width=0.37\textwidth]{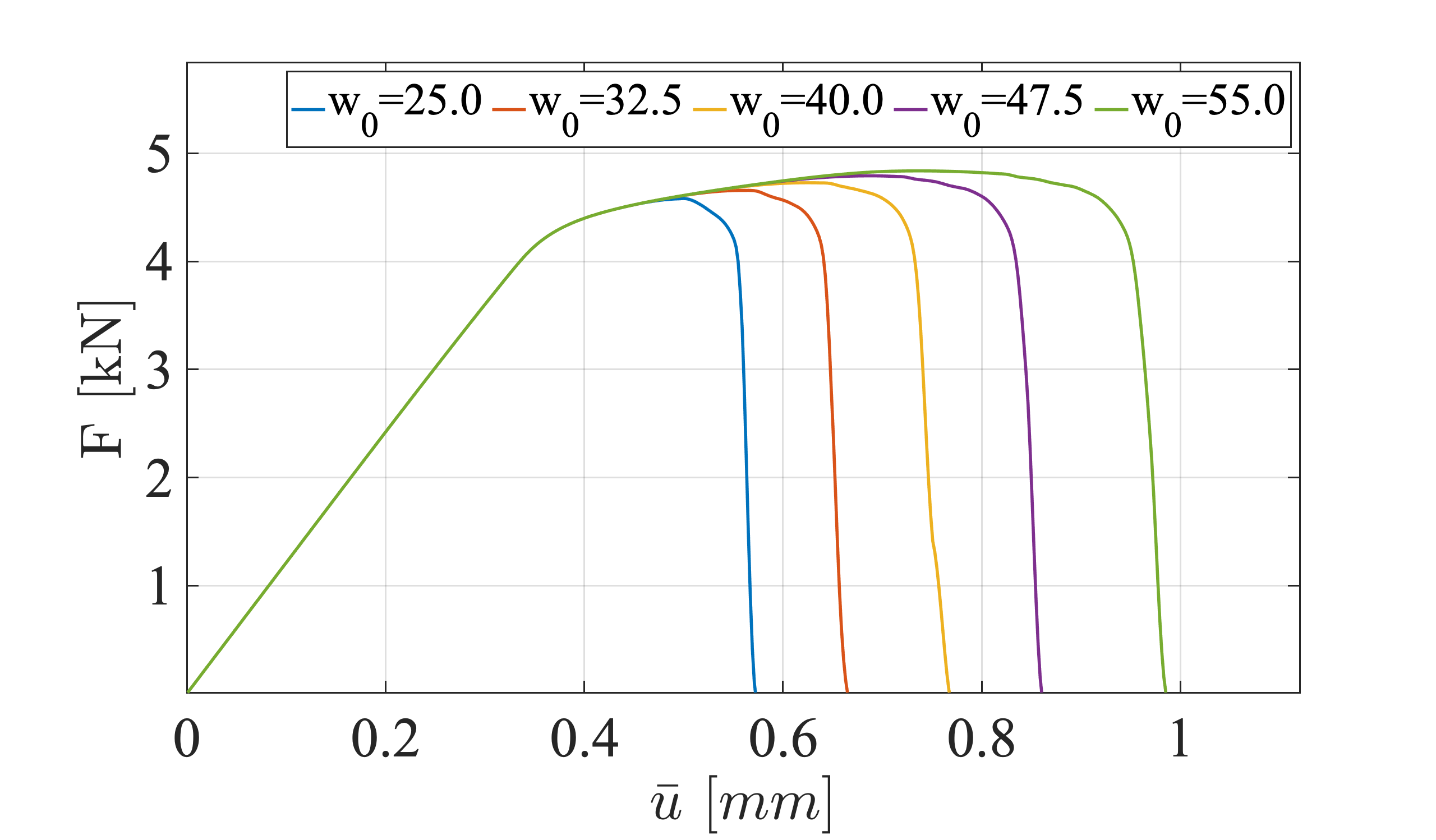}}
	\caption{Example 1: the effect of all studied parameters on  the load-displacement curve. Here, we depict the curves obtained by $\mathcal{M}_1$  for $\mu$, $K$, $H$, $\sigma_Y$, $\alpha_{\text{crit}}$, $G_c$ (the first and the second row). In the third row, we depict the effect of $\psi_c$ ($\mathcal{M}_2$) as well as $\zeta$ and $w_0$ ($\mathcal{M}_3$).}
	\label{Example1}
\end{figure}
\begin{table}[!]
	\caption{Example 1: the mean value of posterior density of  the model parameters for the three models.}
	\vspace{1mm}
	\centering
	\begin{tabular}{ccccccccccccccc}
		&Model  &$H$        & $\mu$    &   $K$   &   $\sigma_Y$   &$G_c$&$\alpha_{\text{crit}}$     &$\psi_c$   &$w_0$ &$l_p$&$\zeta$  \\[2mm]\hline\\
		&    $\mathcal{M}_1$   &240   & 27\,212 & 71\,527 &335 & 10.5 & 0.11 &--&--&--&-- \\[4mm]
		& $\mathcal{M}_2$ &   241  & 27\,120  &  71\,245 & 330 &-- &-- & 41& &--&0.98&\\[4mm]
		&   $\mathcal{M}_3$  &248  & 26\,699   & 74\,500 &328 & -- & -- &--& 38 & 1.25&1.02  \\[4mm]
		\hline
		\label{Example1_posterior}
	\end{tabular}
\end{table}
\begin{table}[!]
	\caption{Example 1: using Bayesian inversion to estimate the equivalence of $\psi_c$  ($\mathcal{M}_2$) with $\alpha_{\text{crit}}$ and $G_c$ ($\mathcal{M}_1$) as well as $w_0$ and $l_p$ ($\mathcal{M}_3$).}
	\vspace{1mm}
	\centering
	\begin{tabular}{cccccccccccc}
		&$\psi_c$  &25        & 35   &   45   &   55   &65&      \\\hline 
		&    $\alpha_{\text{crit}}$   &0.065  & 0.092 & 0.12&0.136& 0.15 &   \\[2mm]
		& $G_c$&   8.25  & 10.1  &  12 & 15 &18   \\[2mm]
		& $w_0$&  20.8  & 29.1   &  38 & 48  &55 &   \\[2mm]
		& $l_p$& 1.05   &  1.11  & 1.25  & 1.4  &1.6 &   \\[2mm]
		\hline
		\label{Example1_equivalence}
	\end{tabular}
\end{table}

For all three models, the shear modulus $\mu$, the bulk modulus $K$, the hardening modulus $H$, and the yield stress $\sigma_Y$ are common. First, we study the effect of the common parameters on the load-displacement curve.  Figure \ref{Example1} shows the diagrams where $\mathcal{M}_1$ is used to obtain the solutions. To monitor different critical values $\alpha_{\text{crit}}$ and energy release rates $G_c$, as well as fracture energies $\psi_c$, $\mathcal{M}_1$ and $\mathcal{M}_2$ are respectively  employed. Moreover, we used $\mathcal{M}_3$ to observe how the curve is affected by specific values of $w_0$ and the parameter $\zeta$, as shown in Figure \ref{Example1}.

\begin{figure}[t!]
 \vspace{0.1cm}
  \hspace{-1cm}
 	\subfloat{\includegraphics[width=0.37\textwidth]{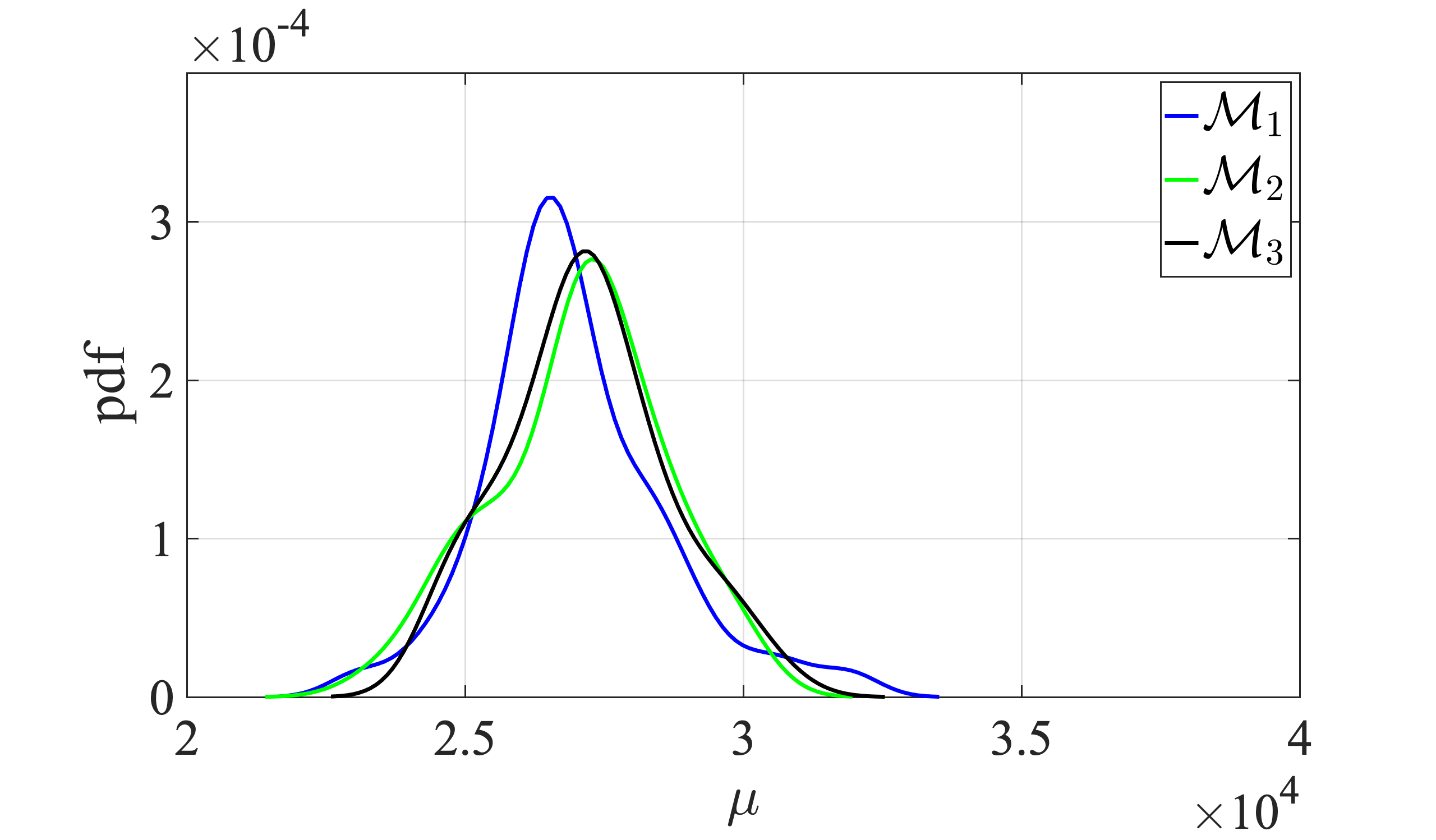}}   \subfloat{\includegraphics[width=0.37\textwidth]{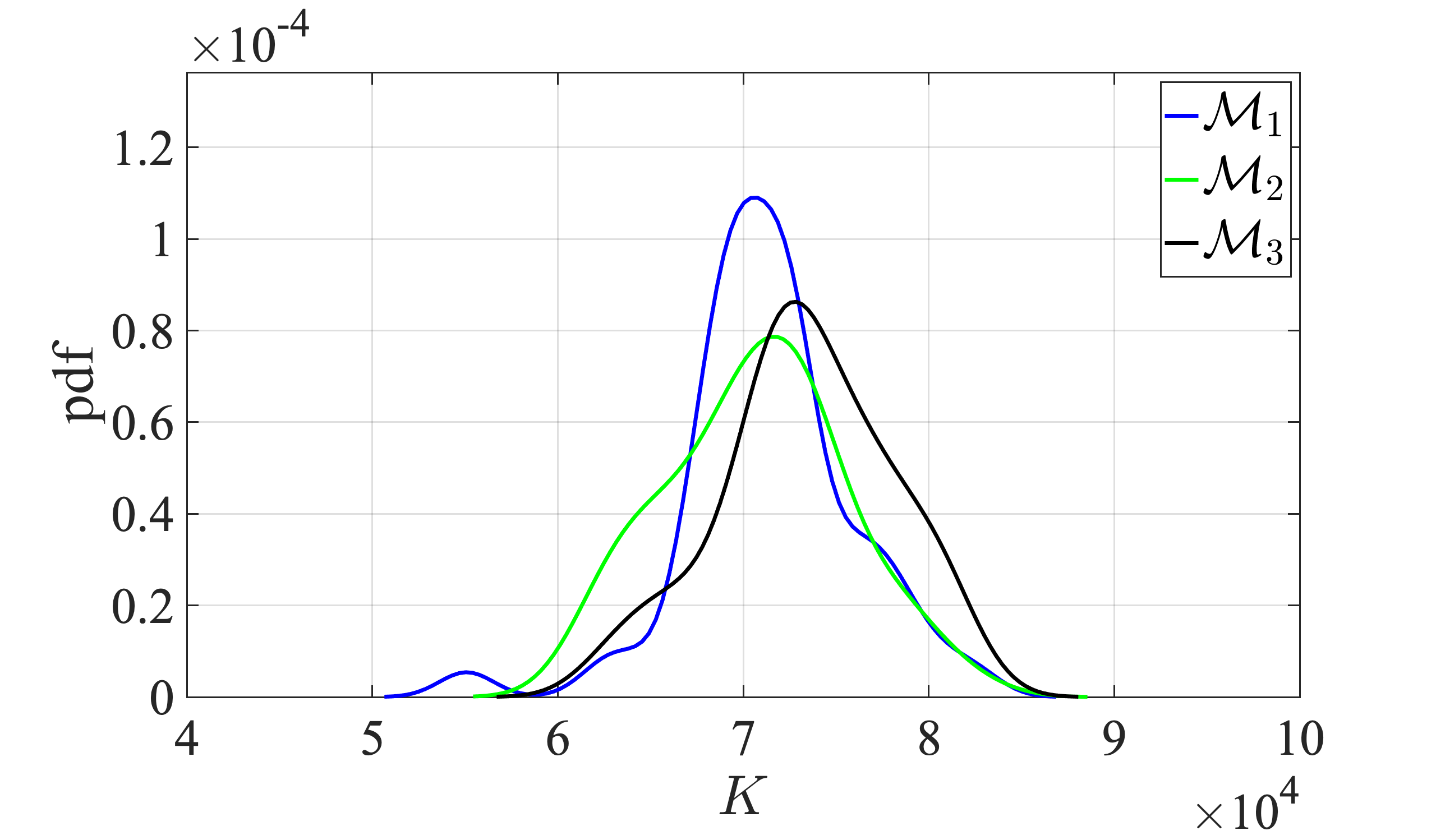}}  	\subfloat{\includegraphics[width=0.37\textwidth]{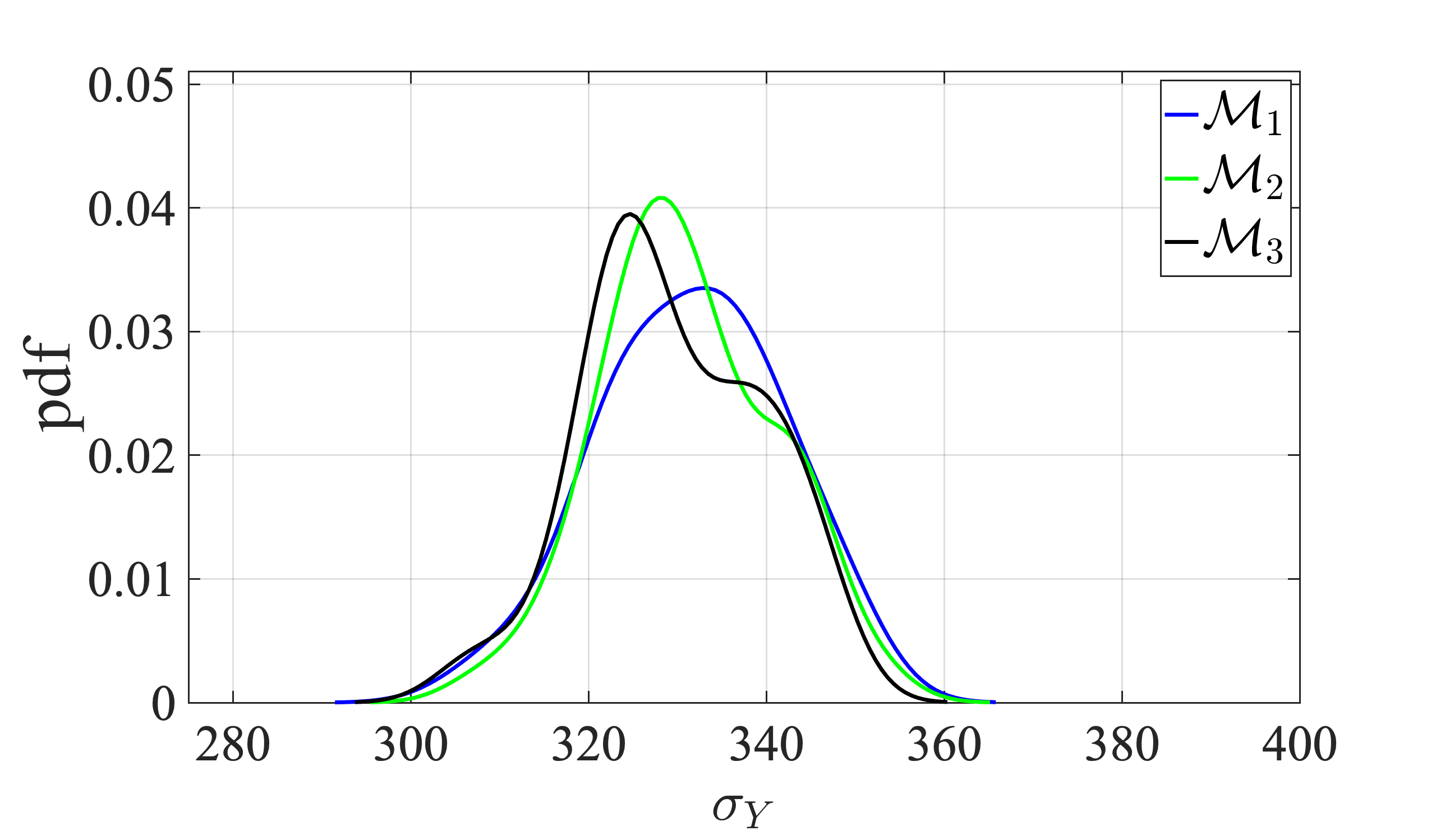}}
 	 \hfill 
 	   \vspace{0.1cm}
  \hspace{-1.1cm}
 	 		 \subfloat{\includegraphics[width=0.37\textwidth]{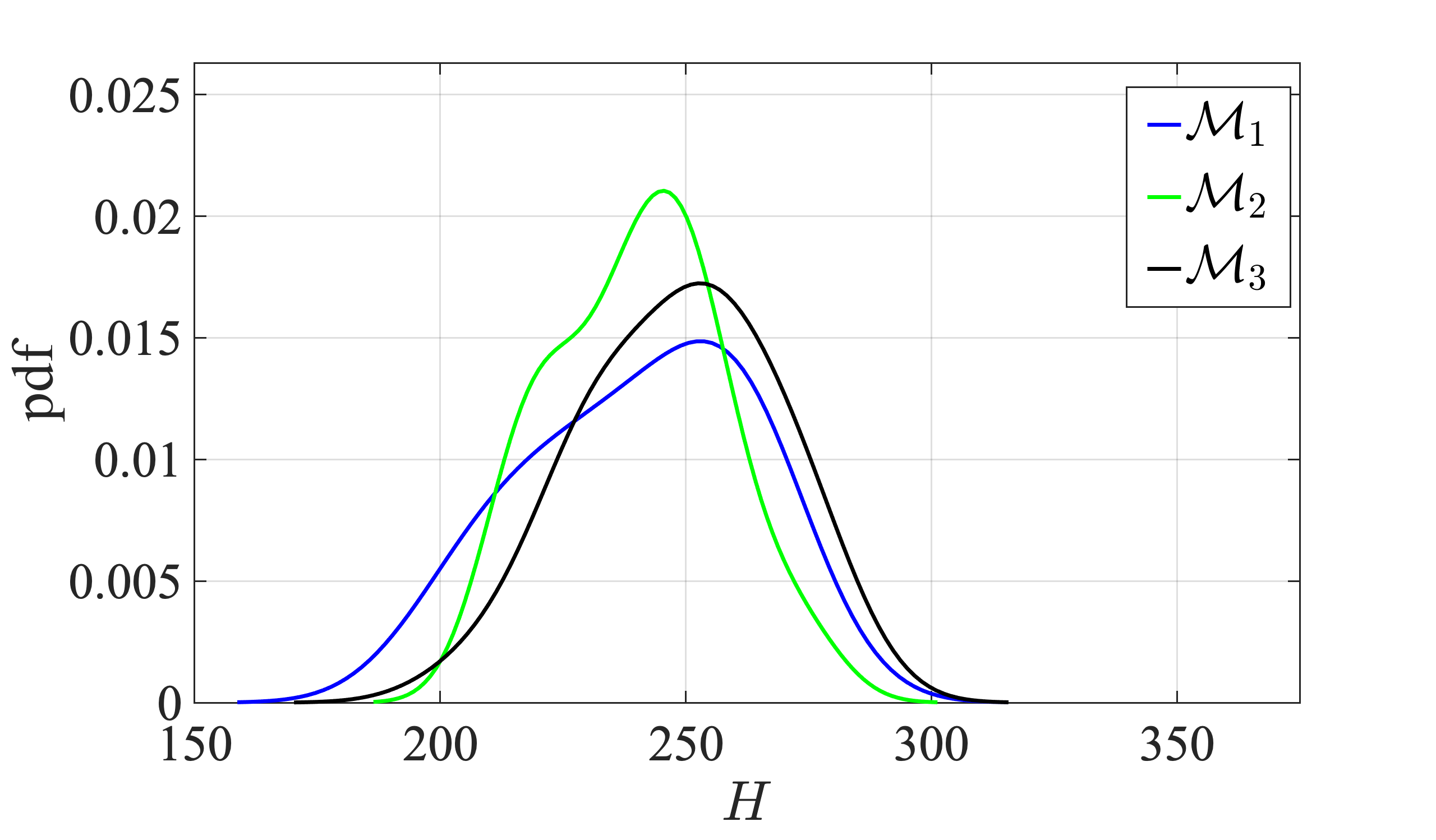}}   \subfloat{\includegraphics[width=0.37\textwidth]{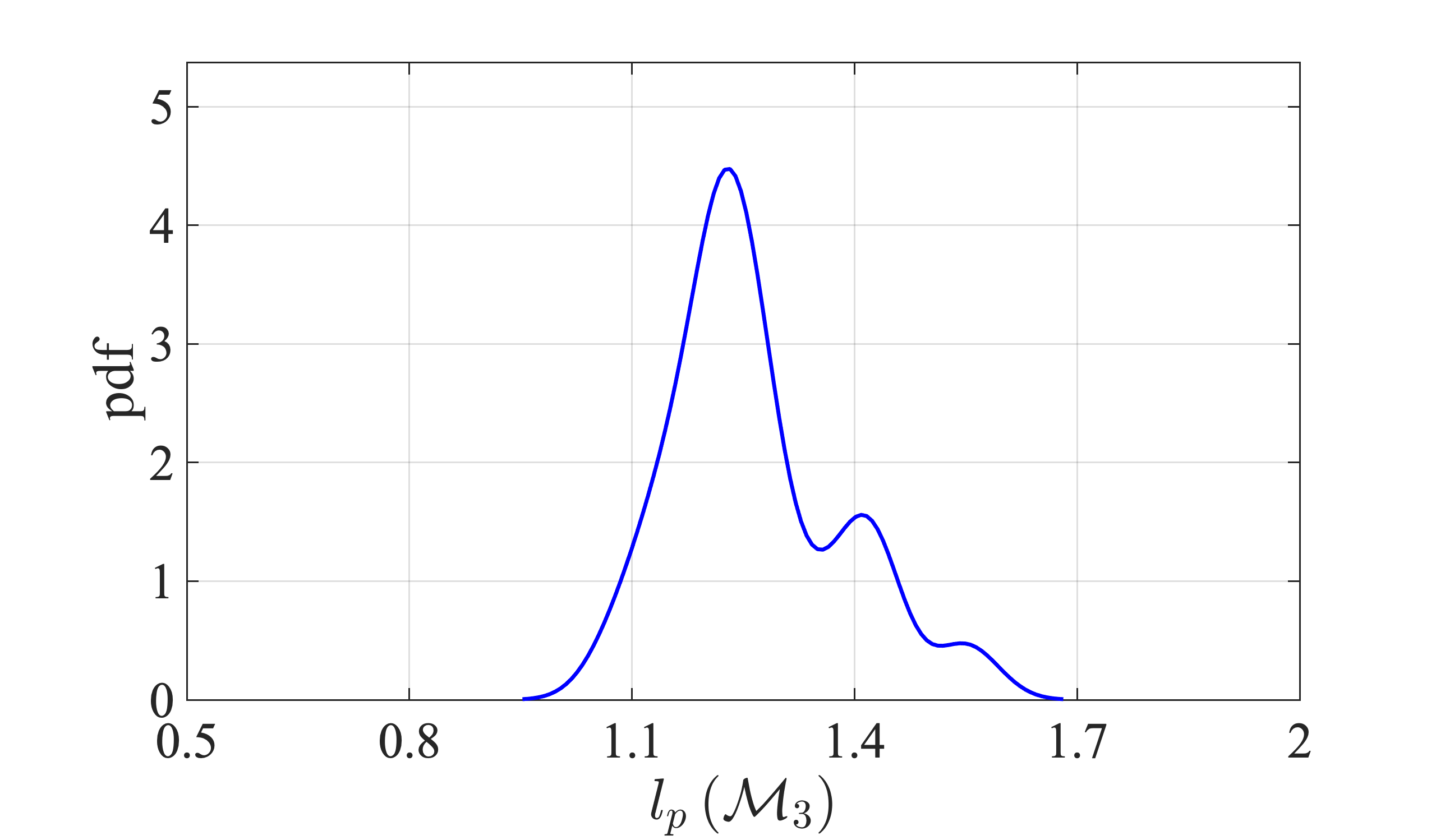}}  	\subfloat{\includegraphics[width=0.37\textwidth]{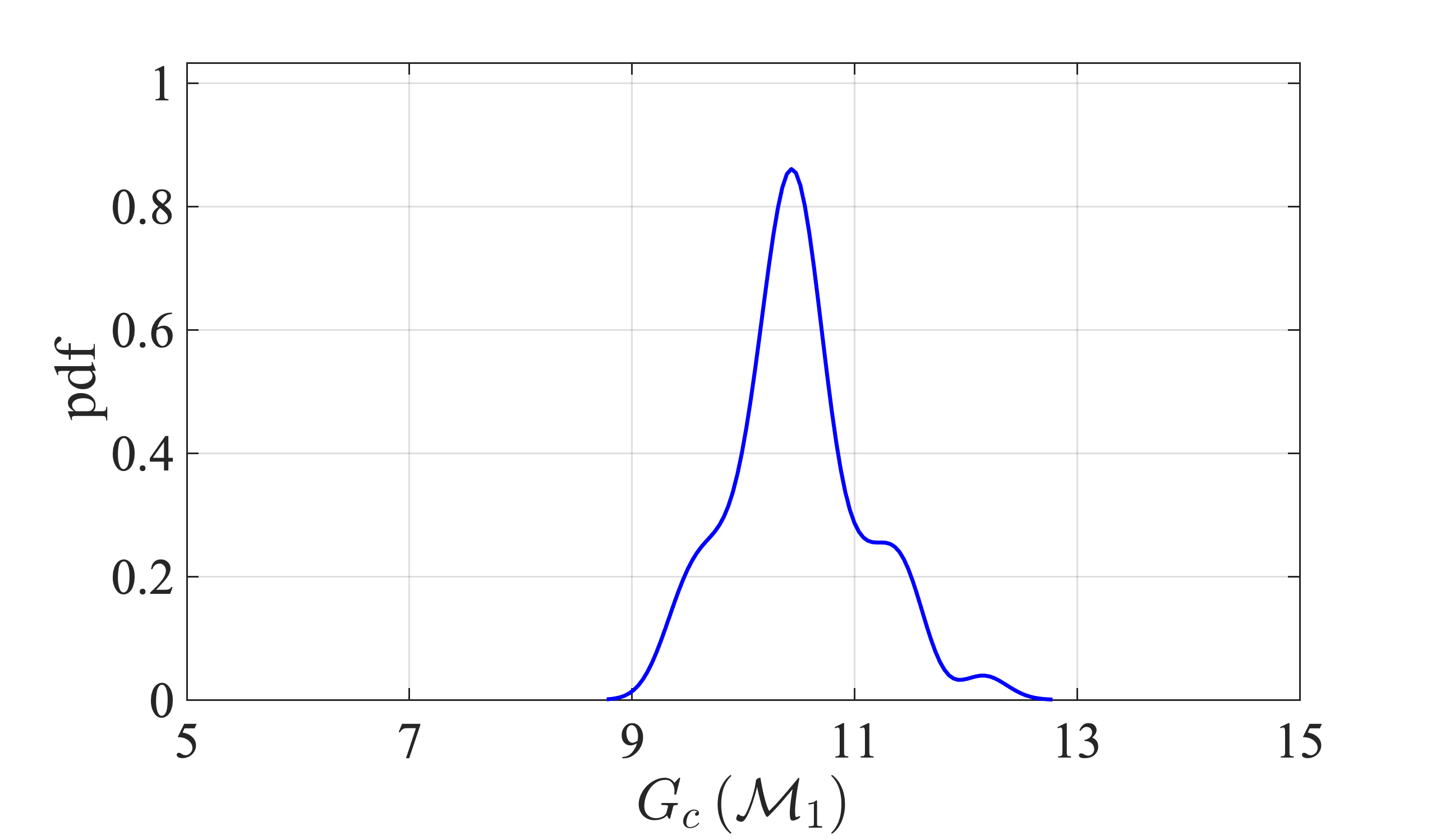}}
 	 		  	 \hfill 
 	 		  	 \vspace{0.1cm}
 	 		  	 \hspace{-1.1cm}
 	 		  	  \subfloat{\includegraphics[width=0.37\textwidth]{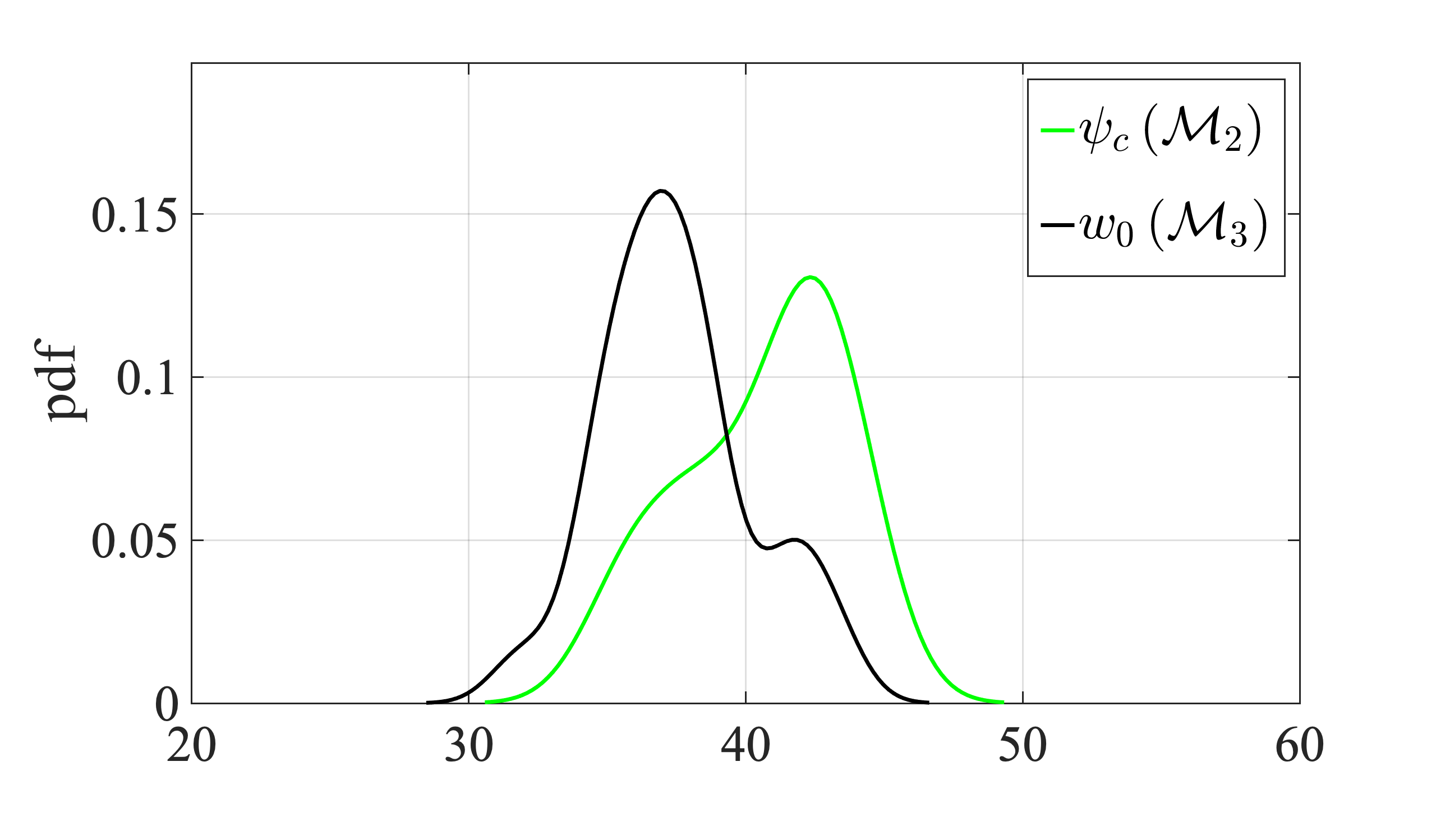}}   \subfloat{\includegraphics[width=0.37\textwidth]{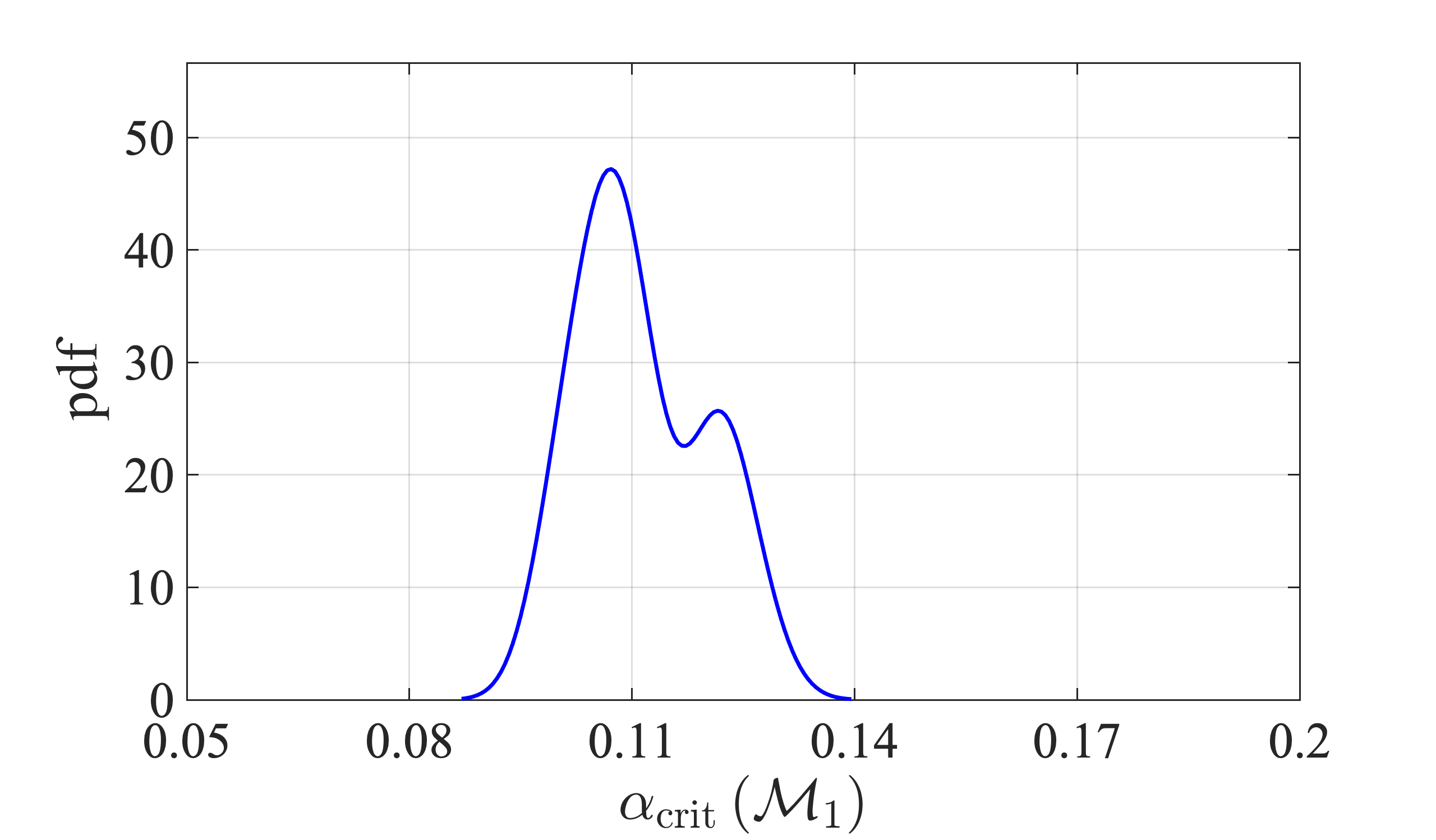}}  	\subfloat{\includegraphics[width=0.37\textwidth]{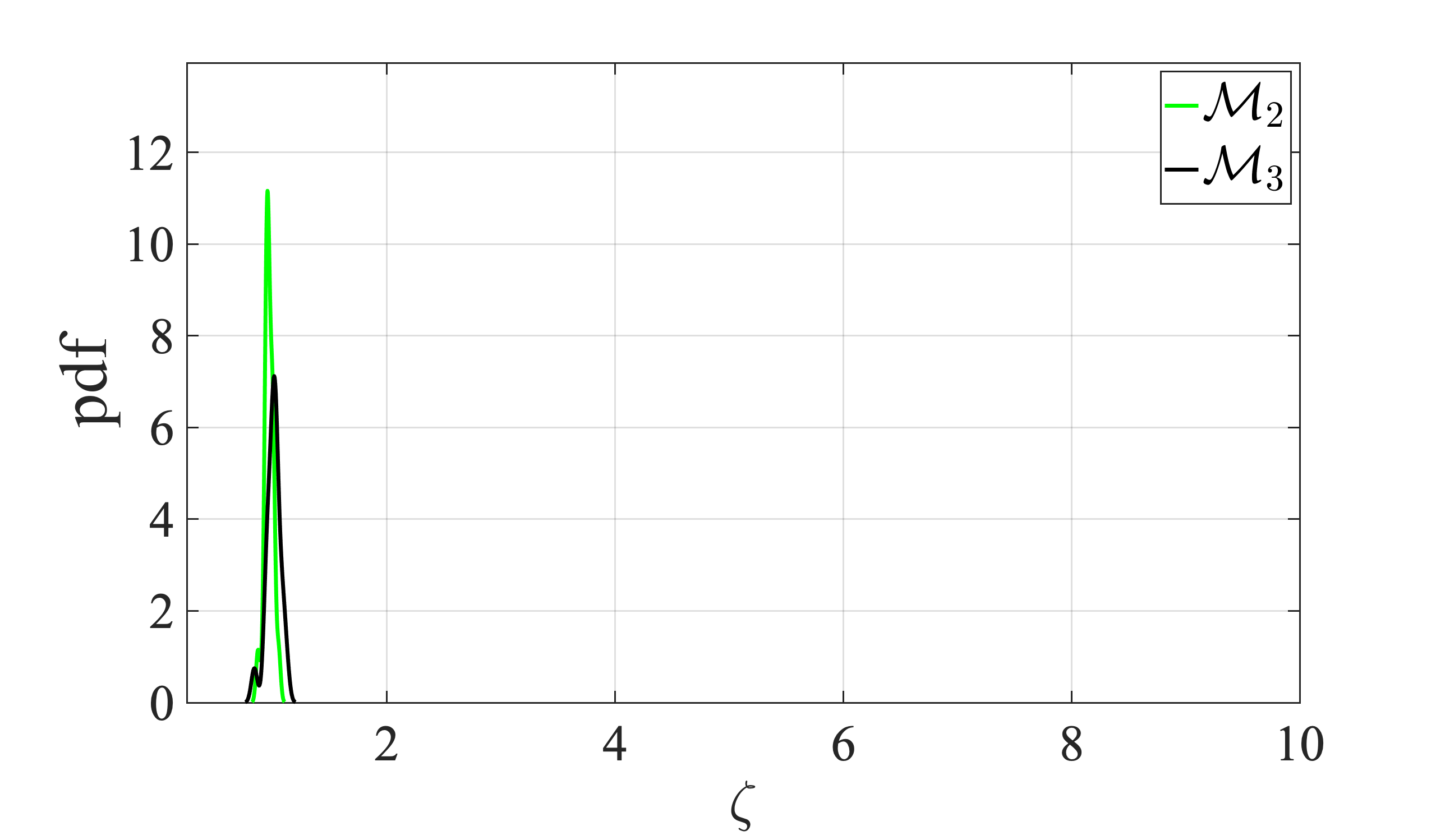}}
 	\caption{Example 1: the posterior distribution of the effective parameters using $\mathcal{M}_1$, $\mathcal{M}_2$, and $\mathcal{M}_3$.}
 	\label{Example1_models}
 \end{figure}

\begin{figure}[t!]
 	\subfloat{\includegraphics[width=0.5\textwidth]{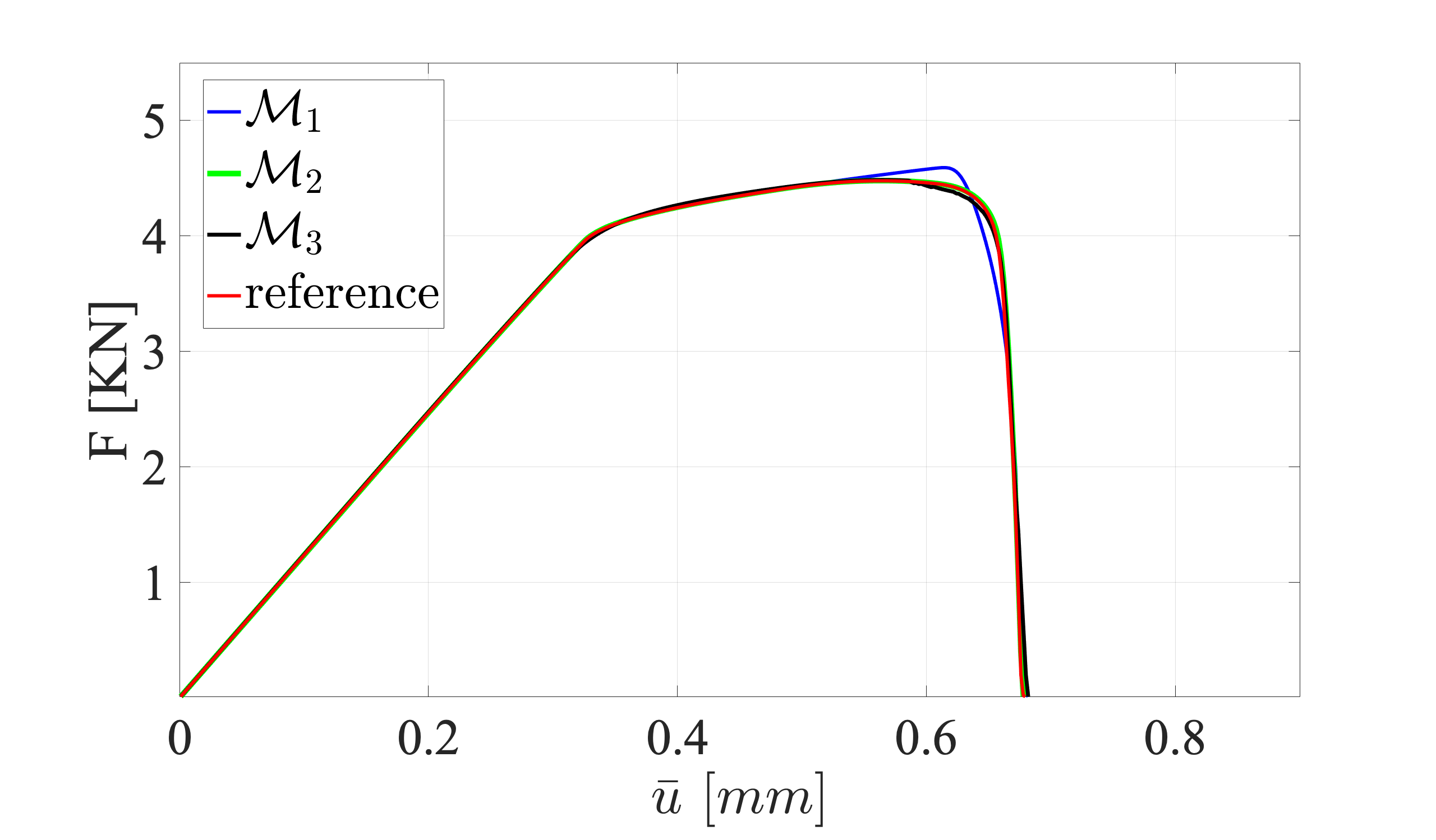}} 	\subfloat{\includegraphics[width=0.5\textwidth]{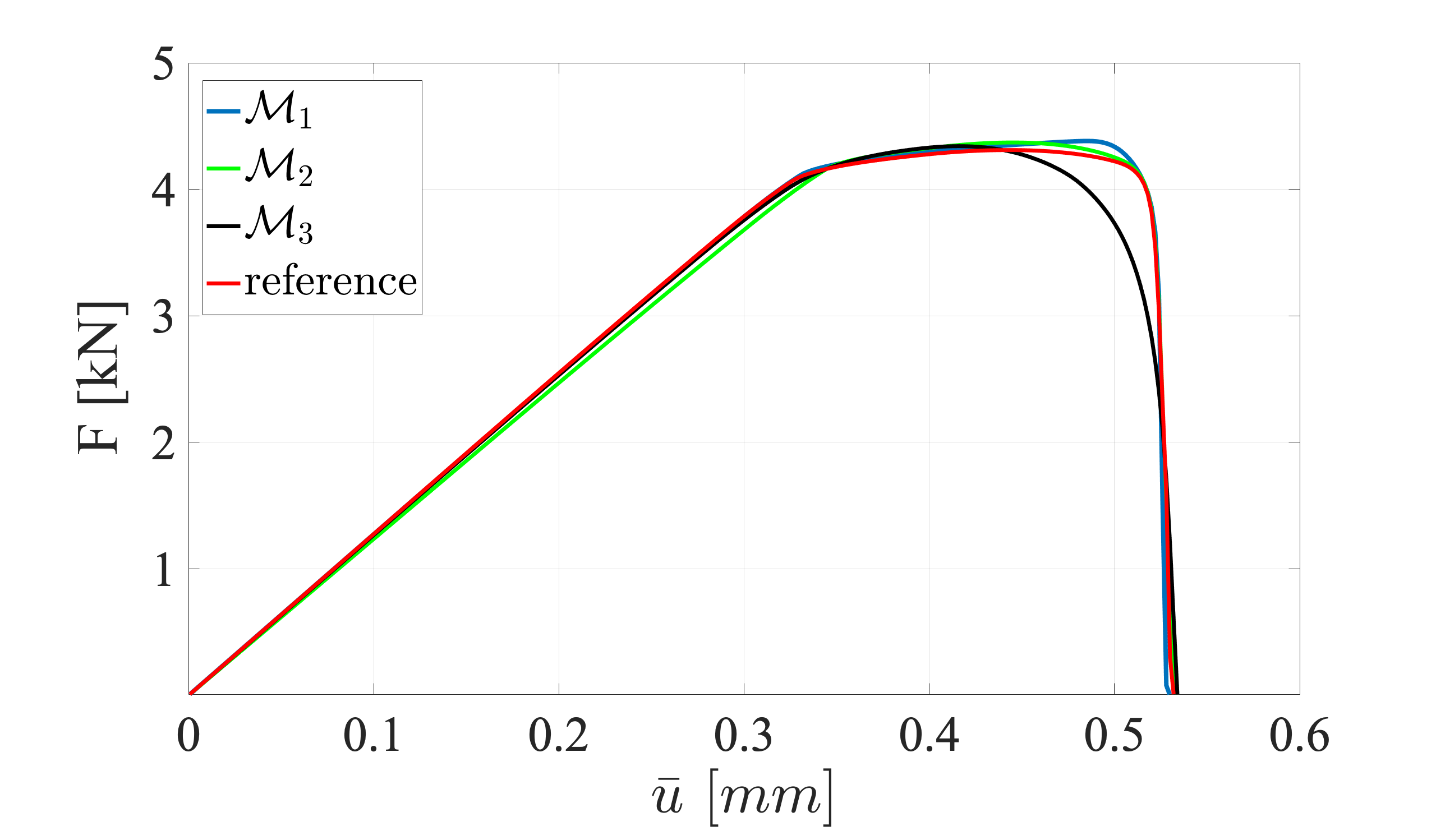}}  
 	\caption{Load-displacement curve computed using the inferred values, employing $\mathcal{M}_1$, $\mathcal{M}_2$, and $\mathcal{M}_3$. The reference values are depicted as well. Here, Table \ref{Example1_posterior} and Table \ref{Example2_posterior} are used for Example 1 (left) and Example 2 (right).}
 	\label{example12}
 \end{figure}

Next, we proceed to identify the effective parameters in the ductile fracture process using the Metropolis-Hastings algorithm introduced in Section \ref{MH1}. The Bayesian framework for ductile fracture is presented in Section \ref{bayesian}.
We employ a uniform distribution to estimate the parameters more accurately, as listed in Table \ref{range1} (the prior densities and the initial values) and use $N=10\,000$ number of candidates.  Regarding the reference values, a synthetic measurement is used, using a total number of degrees of freedom {\color{black}$\mathcal{N}_{\text{dof}}=28\,380$}; the rest of the initial values are summarized in Table \ref{range1}.
The posterior density of the parameters using the three models is shown in Figure \ref{Example1_models}. The mean values of the posterior distributions are used to verify the parameter estimation. The inferred information is listed in Table \ref{Example1_posterior}. To verify the accuracy of the data, we employ the parameters in all three models and compute  the load-displacement curve until the fracture point. Figure \ref{example12} shows the curves resulting from the different models and the reference observation. An excellent agreement indicates that the Bayesian inversion framework identified the parameters correctly, showing a consistent behavior for all three models in all stages.

The accuracy of the Bayesian inversion for all models enables us to provide equivalence for the model parameters. As previously mentioned, all models have four common parameters, but each model is also characterized by its own features. Here, we strive to find an equivalent value for different fracture energies. This allows us to use $\mathcal{M}_1$ and $\mathcal{M}_3$ and derive similar quantities in $\mathcal{M}_2$, and vice versa. To that end, we select the diagram estimated by $\mathcal{M}_2$ as the reference observation, where all parameters are chosen according to the estimated values (see Table \ref{Example1_posterior}). However, $\psi_c$ varies between 25 and 65. We again use the MH algorithm to identify the equivalence of $\psi_c$ in $\mathcal{M}_1$ (i.e, $G_c$ and $\alpha_{\text{crit}}$) and $\mathcal{M}_3$ (i.e., $l_p$ and $w_0$). The estimated quantities are summarized in Table \ref{Example1_equivalence}. Figure \ref{example1_psi} presents the results obtained by the inferred values, where, again, Bayesian inversion provides a very good agreement.


\begin{figure}[!]
	\centering
	\subfloat{{\includegraphics[clip,trim=0cm 21cm 16cm 4.5cm, width=12cm]{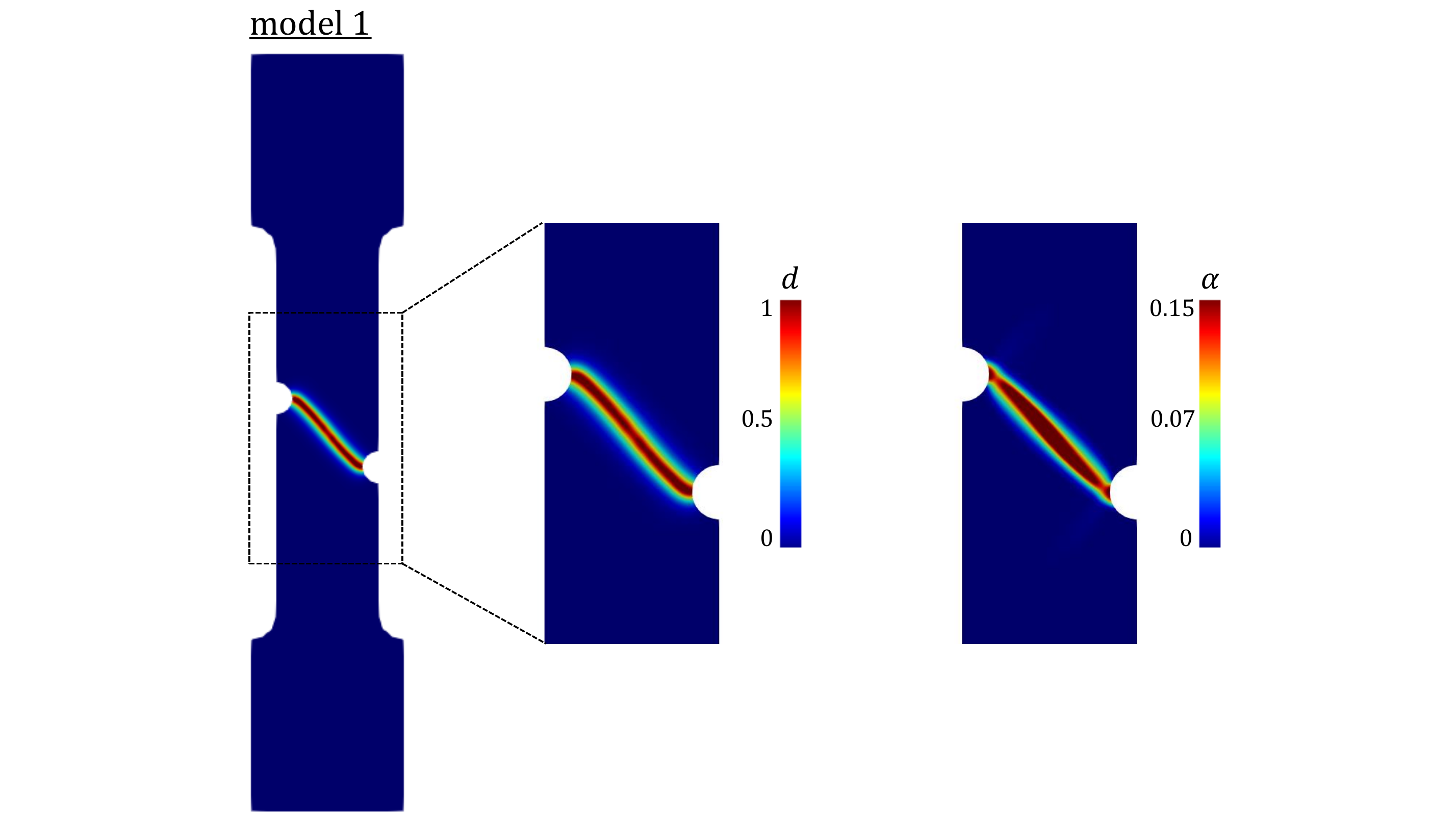}}}\\  
	\subfloat{{\includegraphics[clip,trim=0cm 21cm 16cm 4.5cm, width=12cm]{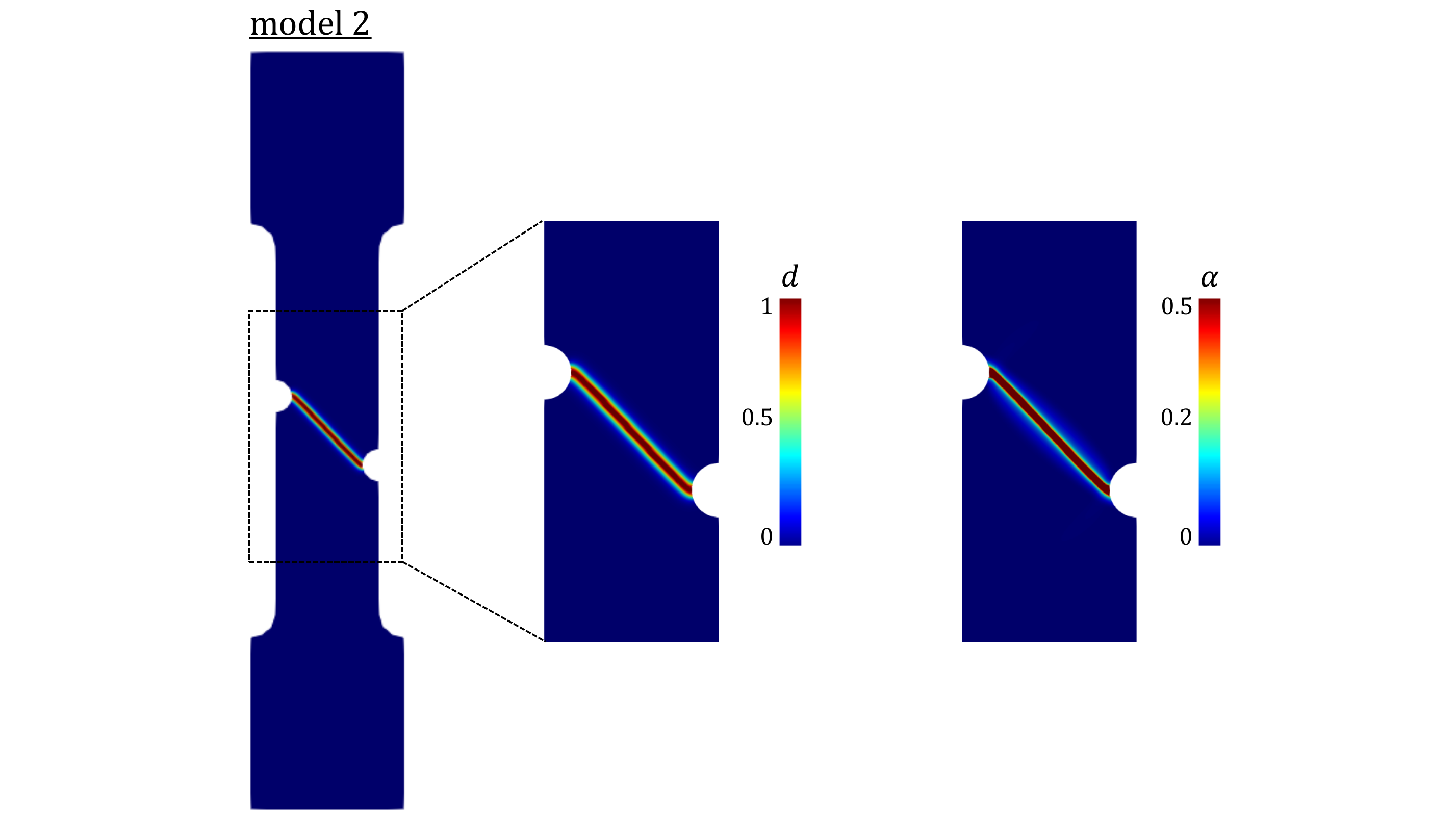}}}\\  
	\subfloat{{\includegraphics[clip,trim=0cm 20cm 16cm 4.5cm, width=12cm]{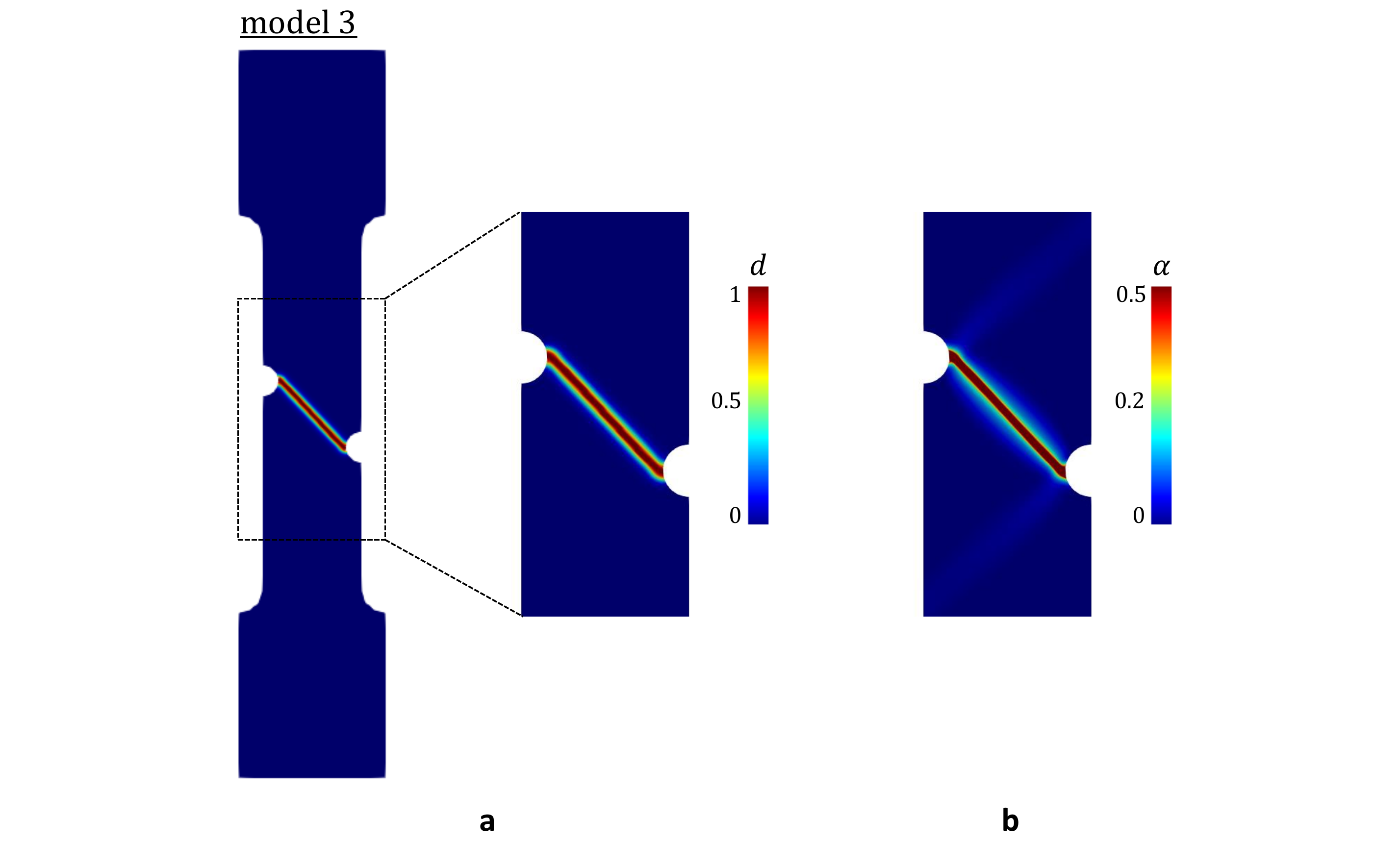}}}  
	\caption{ Example 1: approximated solution obtained through the posterior density of the material parameters at complete failure. The hardening value $\alpha$ and the crack phase-field $d$ are shown for different models.}
	\label{Figure13}
\end{figure}

\begin{figure}[!t]
	\centering
   	{\includegraphics[width=0.85\textwidth]{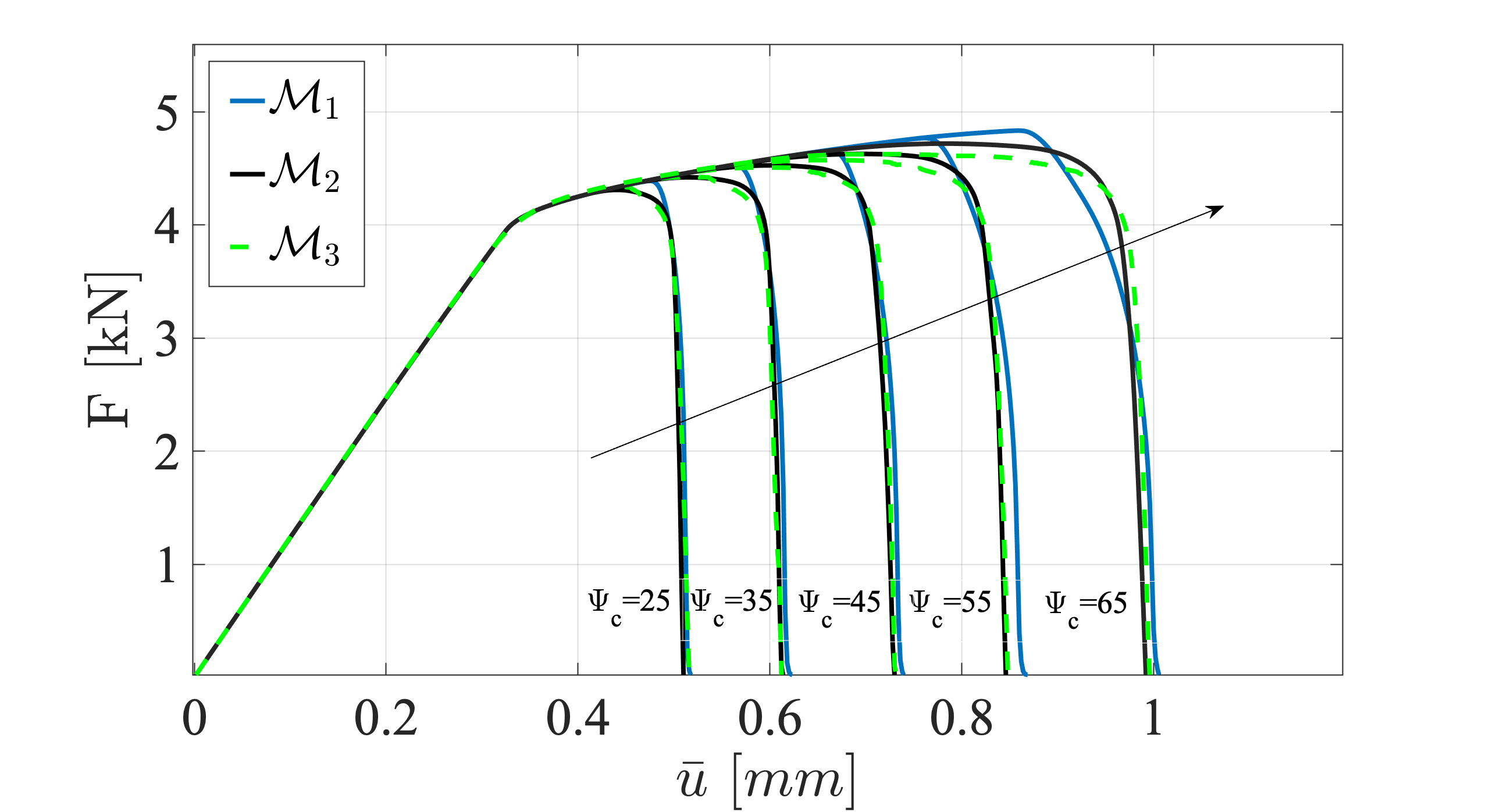}}  
	\caption{Example 1: the load-displacement curve for the inferred equivalent values $\psi_c$, with  $\psi_c$ varying between 25 and 65.  
	}
	\label{example1_psi}
\end{figure}

The resulting equivalent plastic strain ($\alpha$) and crack phase-field ($d$) at complete failure are shown in Figure \ref{Figure13}. The solutions are  based on the posterior density of the material parameters for different models, which are given in Table \ref{Example1_posterior}. Accordingly, the fracture path initiates within the maximum equivalent plastic regions, which appear near the notches. Next, the crack propagates in the plastic localization band, in which two cracks merge at the specimen center.  It can be observed that even though the load-displacement curves shown are practically identical in all models, the corresponding phase-field profiles, and thus, hardening profiles, are not; see Figure \ref{Figure13}. This can be explained, first of all, by the solution non-uniqueness of the phase-field fracture problem, and, secondly, by the fact that the different phase-field models in fact provide only the approximation of the fracture problem. Thus, the necessity of comparing the results with an experimental observation is crucial. Hereby, based on the experimental test provided in \cite{ambati2016} (second experiment), a sharp crack transition between two notches is expected. Thus $\calM_2$ and $\calM_3$ seem to yield a more accurate fracture pattern.
\begin{figure}[!]
	\centering
	{\includegraphics[clip,trim=0cm 1cm 12cm 6cm, width=13cm]{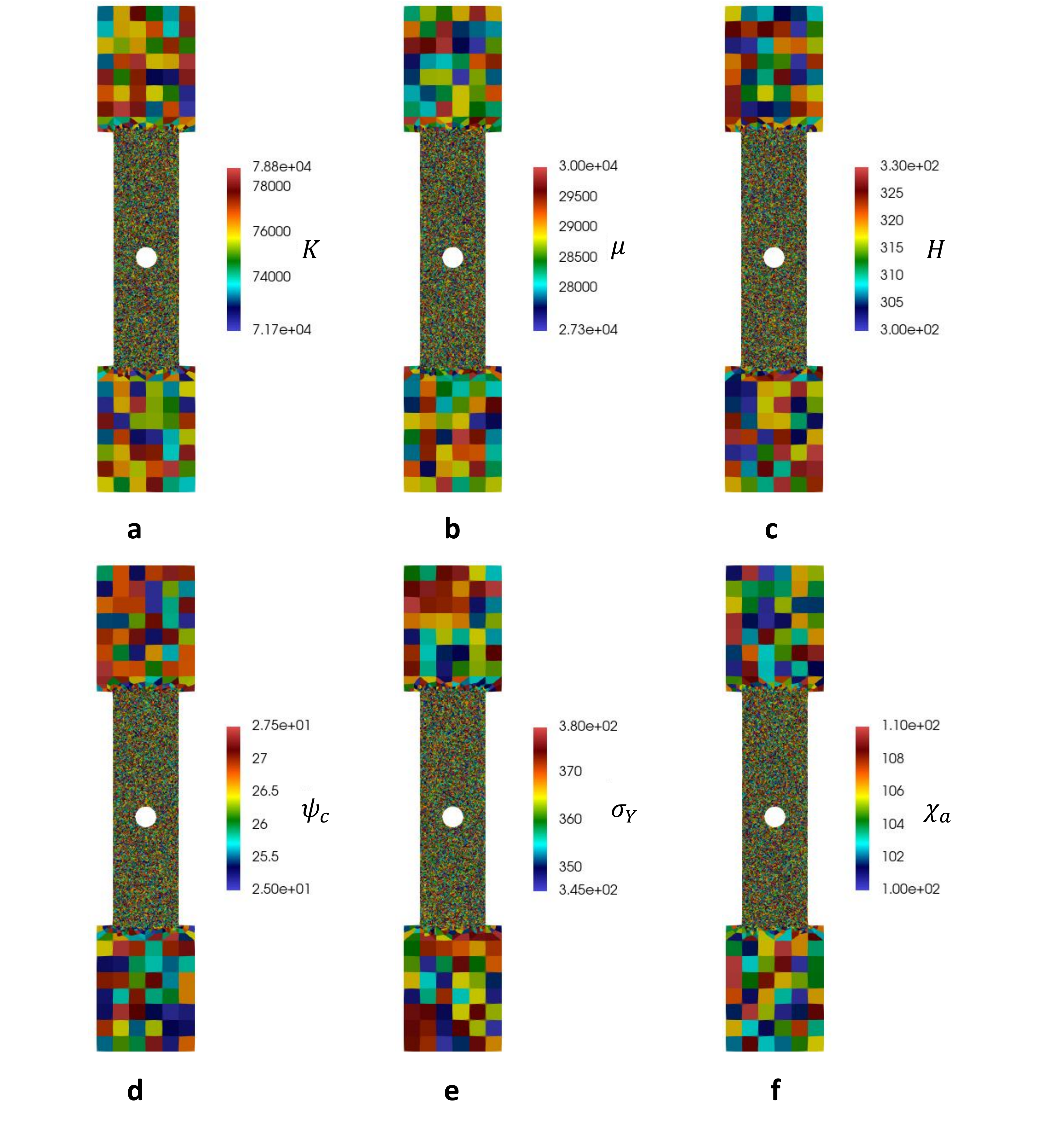}}  
	\caption{ Example 2: material parameters distribution (prior density) on each element of the domain $\calB$. Note that here, to avoid a strong deviation of the surrounding points due to the random normal distribution of the material properties, one could use an additional random distribution length-scale to achieve mesh objectivity, that is, the so-called \textit{heterogeneity length-scale}; see \cite{aquino2018coupled}. In this study, we have not used a heterogeneity length-scale since we have assumed that the material  distribution provides the synthetic observations.}
	\label{Figure21}
\end{figure}
%

%
\sectpb[Section52]{Example 2: I-shaped tensile specimen for anisotropic ductile fracture}

The main objective of this example is the adoption of Bayesian inversion for an  anisotropic ductile phase-field fracture process. The BVP depicted in Figure \ref{Figure1bvp}b consists of an I-shaped specimen with a circular void in the center of the domain. The geometrical dimensions in Figure \ref{Figure1bvp}b are set as $H_1=110$ mm, $H_2=28.6$ mm, $w_1=22$ mm, and $w_2=14.8$ mm, with the central void located in $(x,y)=(H_1/2,w_1/2)$, with a radius of $r=2.5$ mm.


Herein, we assume that the material constituents are \textit{not} distributed uniformly through the continuum domain, and thus, the material is divided into several phases. Hence, heterogeneity in strength from one area of the domain to another one is expected. Note, however, that by means of the Bayesian inversion framework, we aim to determine the \emph{effective} mechanical parameters. Here, we consider the parameters as a random field (with given mean and variation).
Figure \ref{Figure21} illustrates the fluctuation of different material parameters (on the element-wise basis) with spatial correlation where a  10\% variation is included. For instance, for the parameter $K$, the expectation is assumed as $K$=75\,000\,MPa,  with  a variation between 71\,700\,MPa and 78\,800\,MPa. This fluctuation will be used to provide the reference observation. Specifically, we will replicate 500 simulations (with a specific mesh size) to estimate the reference observation considering the mentioned variation. The distribution of the parameters on the geometry is shown in Figure \ref{Figure21}.

The numerical example is performed by applying a monotonic displacement increment  ${{\Delta \bar{u}}_y}=2\times10^{-3}$ mm in the vertical direction at the top boundary of the specimen (Figure \ref{Figure1bvp}). To remove the rigid body motion, the bottom edge is fixed in the $x-y$ directions. The minimum finite element size is $0.45$ mm. The two-dimensional I-shaped domain partition contains 17038 elements. 

 \begin{figure}[t!]
	\vspace{0.1cm}
	\hspace{-1cm}
	\subfloat{\includegraphics[width=0.37\textwidth]{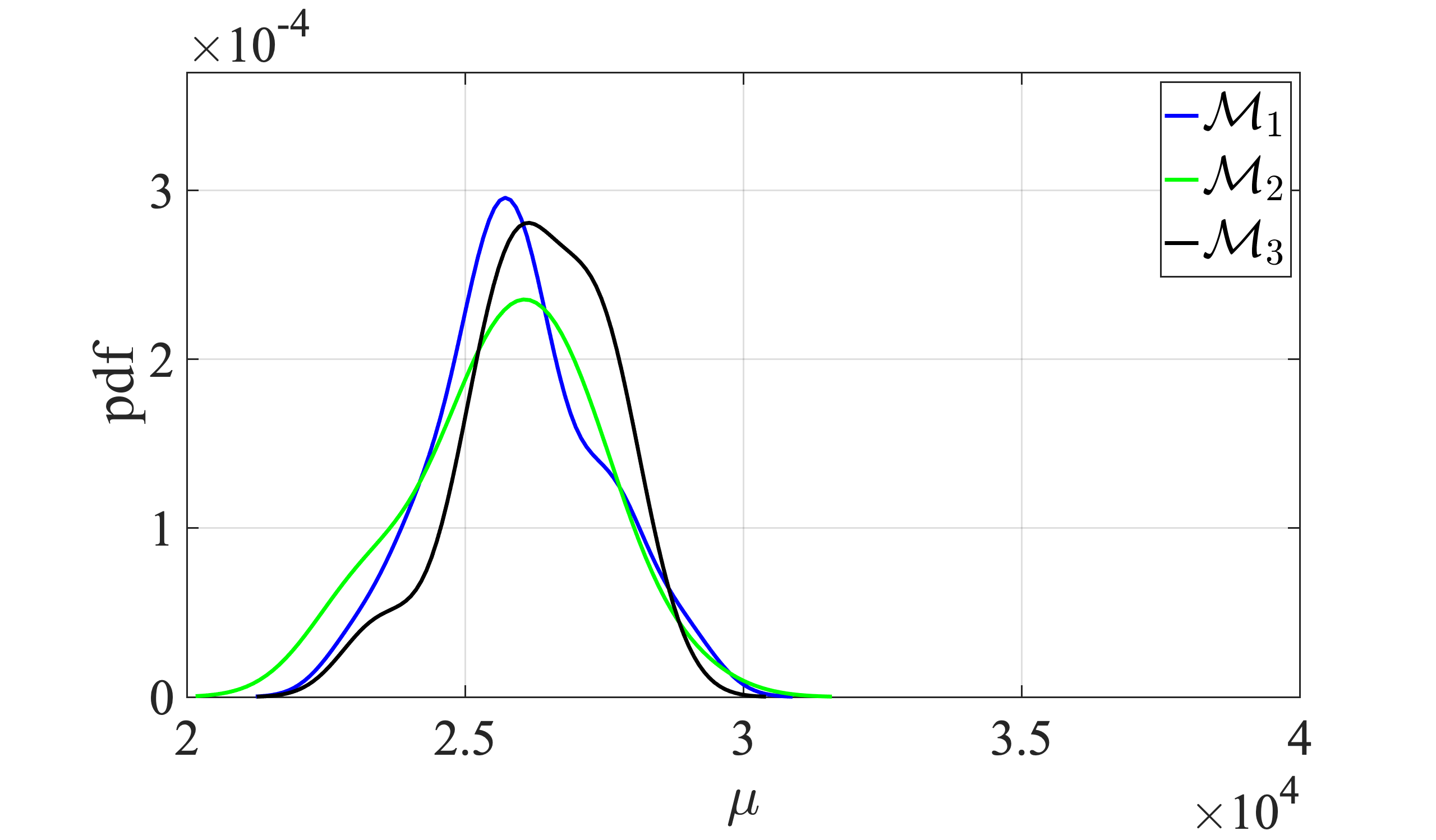}}   \subfloat{\includegraphics[width=0.37\textwidth]{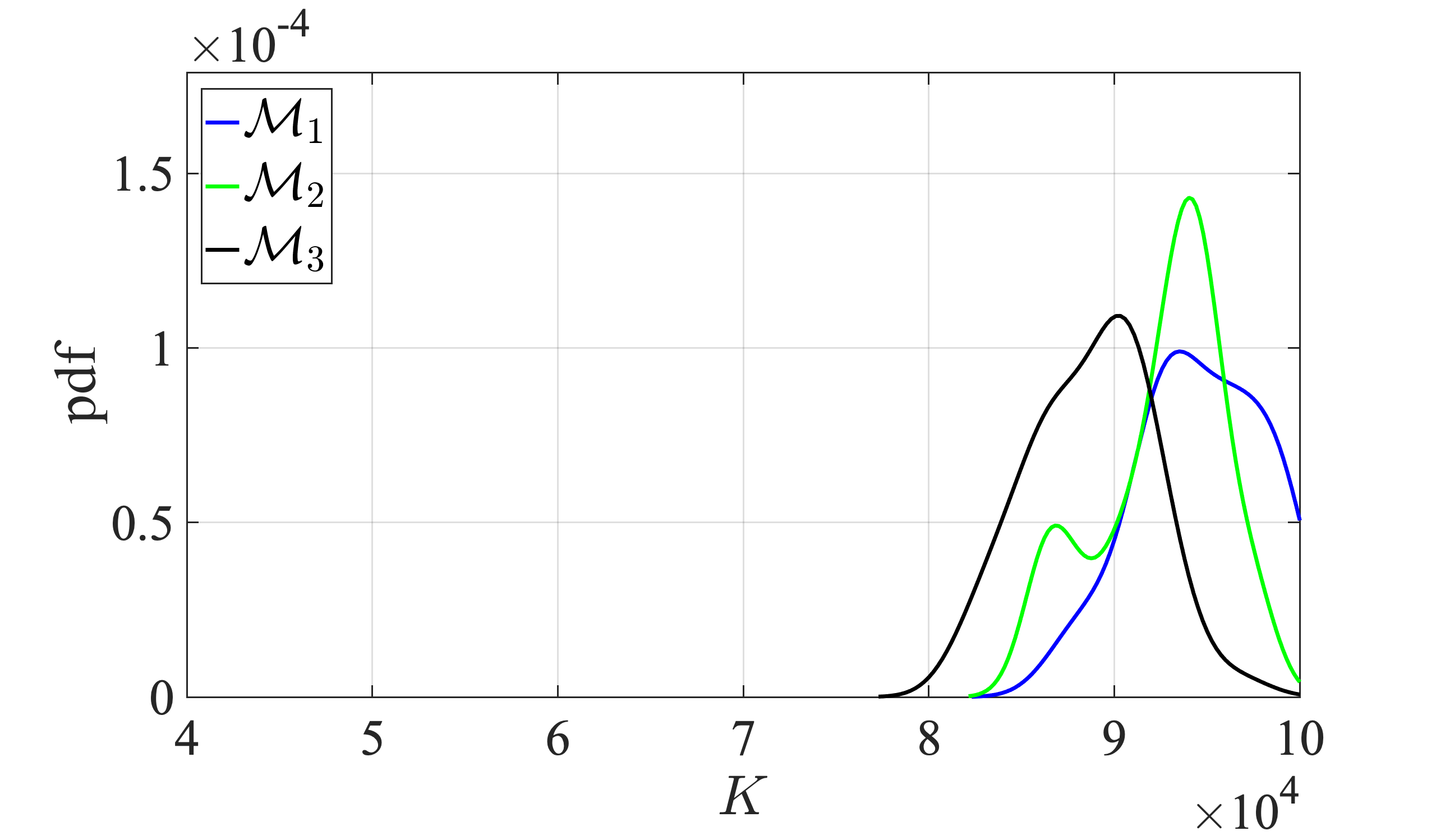}}  	\subfloat{\includegraphics[width=0.37\textwidth]{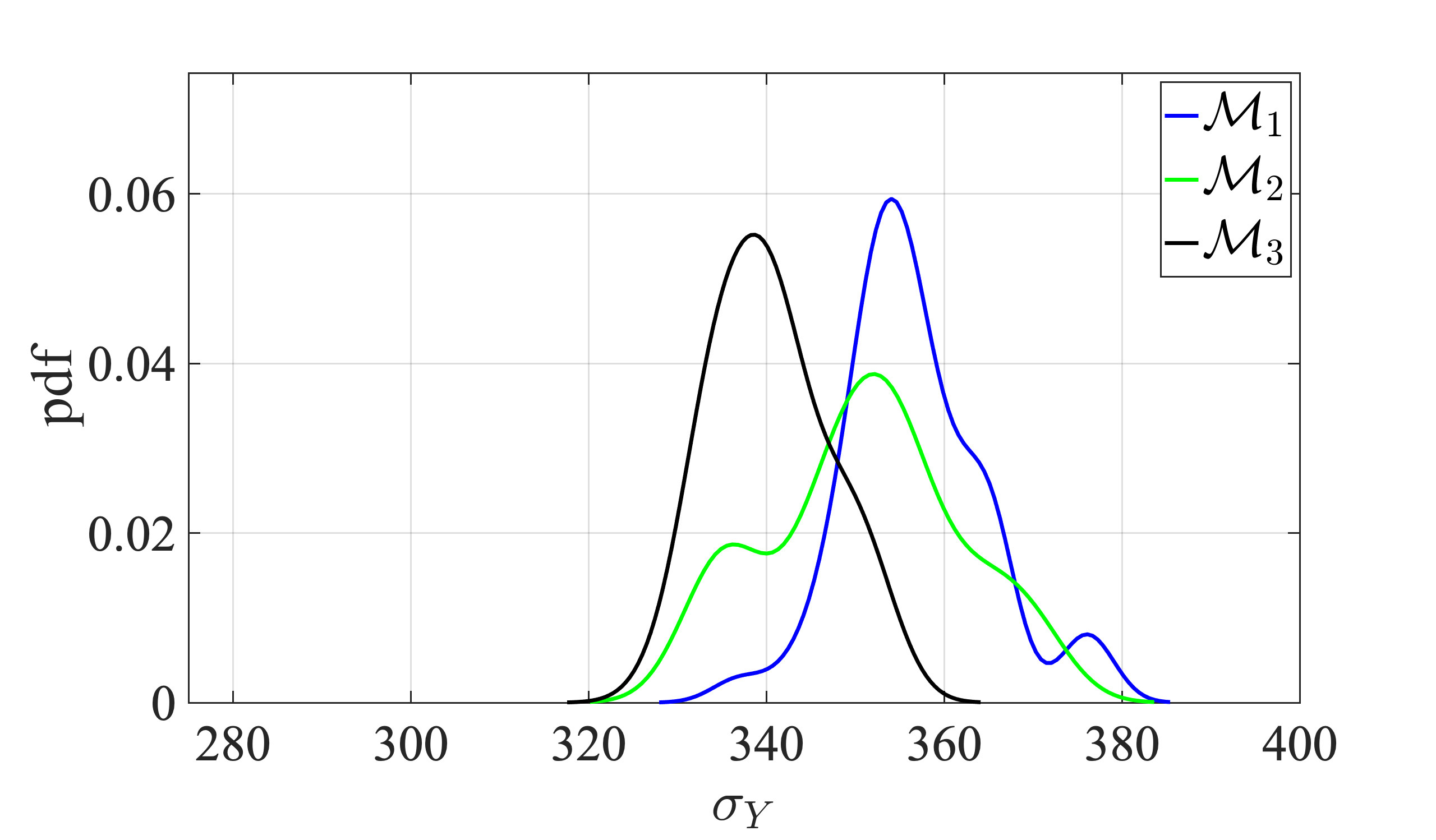}}
	\hfill 
	\vspace{0.1cm}
	\hspace{-1.1cm}
	\subfloat{\includegraphics[width=0.37\textwidth]{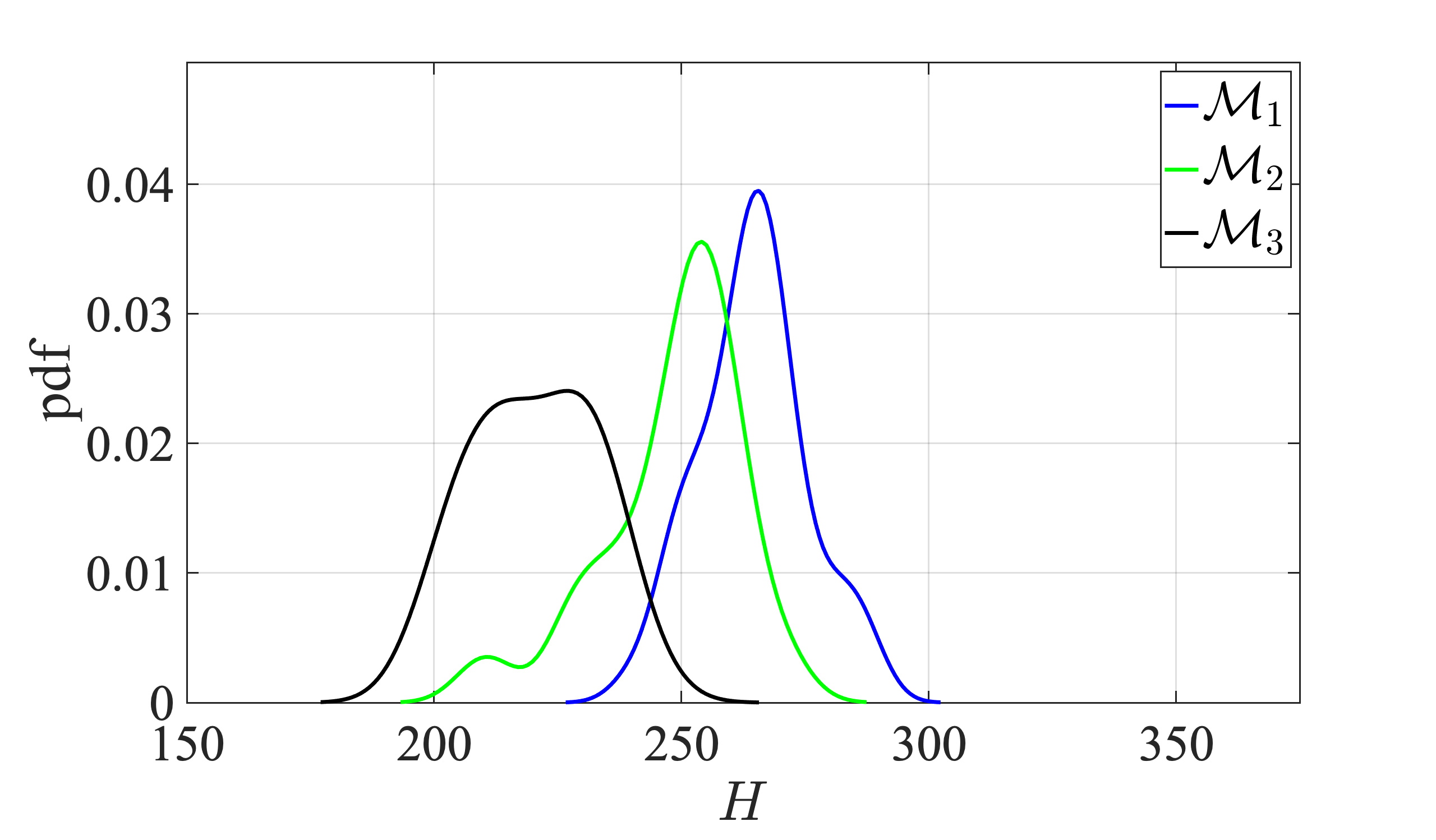}}   \subfloat{\includegraphics[width=0.37\textwidth]{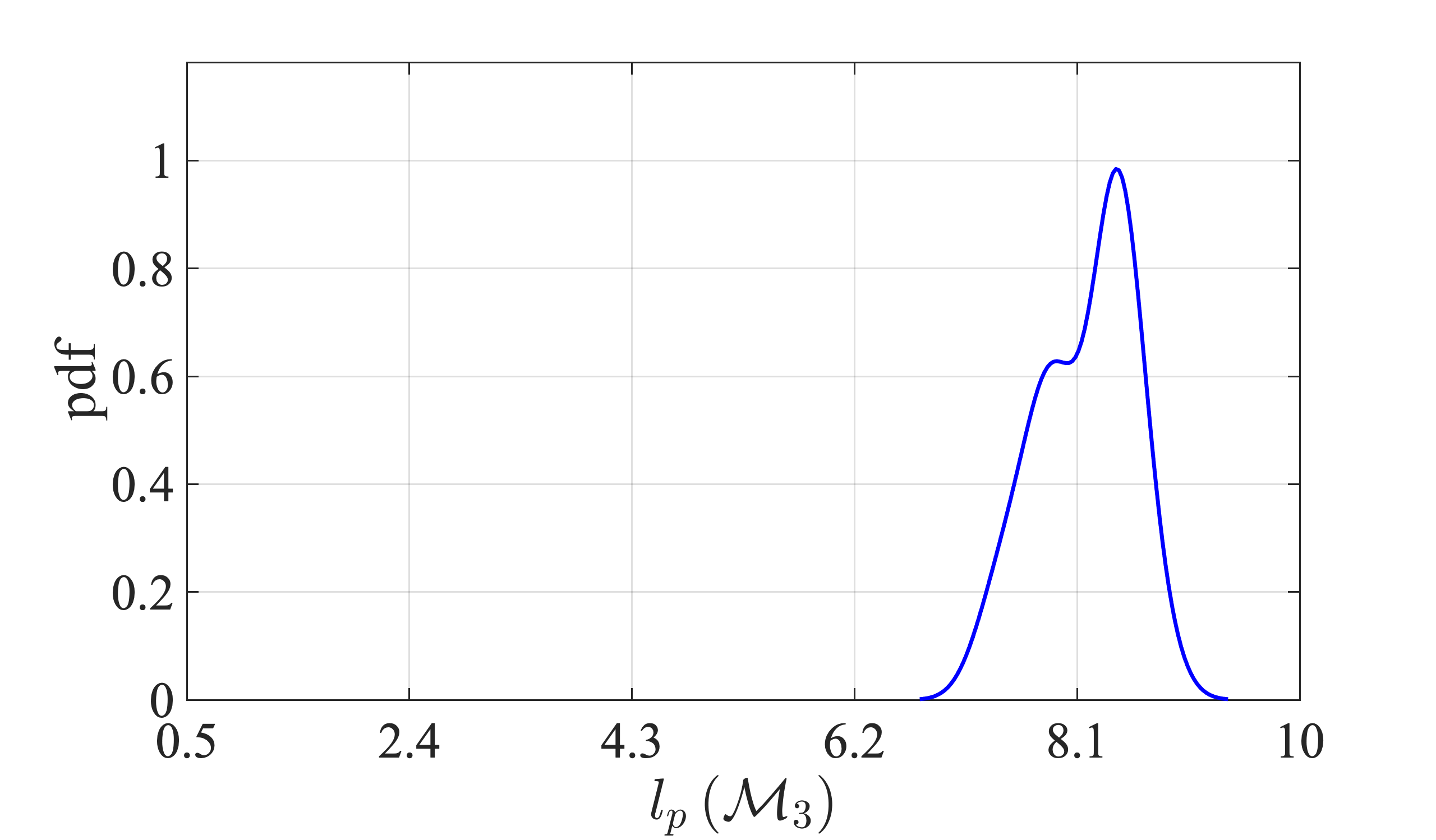}}  	\subfloat{\includegraphics[width=0.37\textwidth]{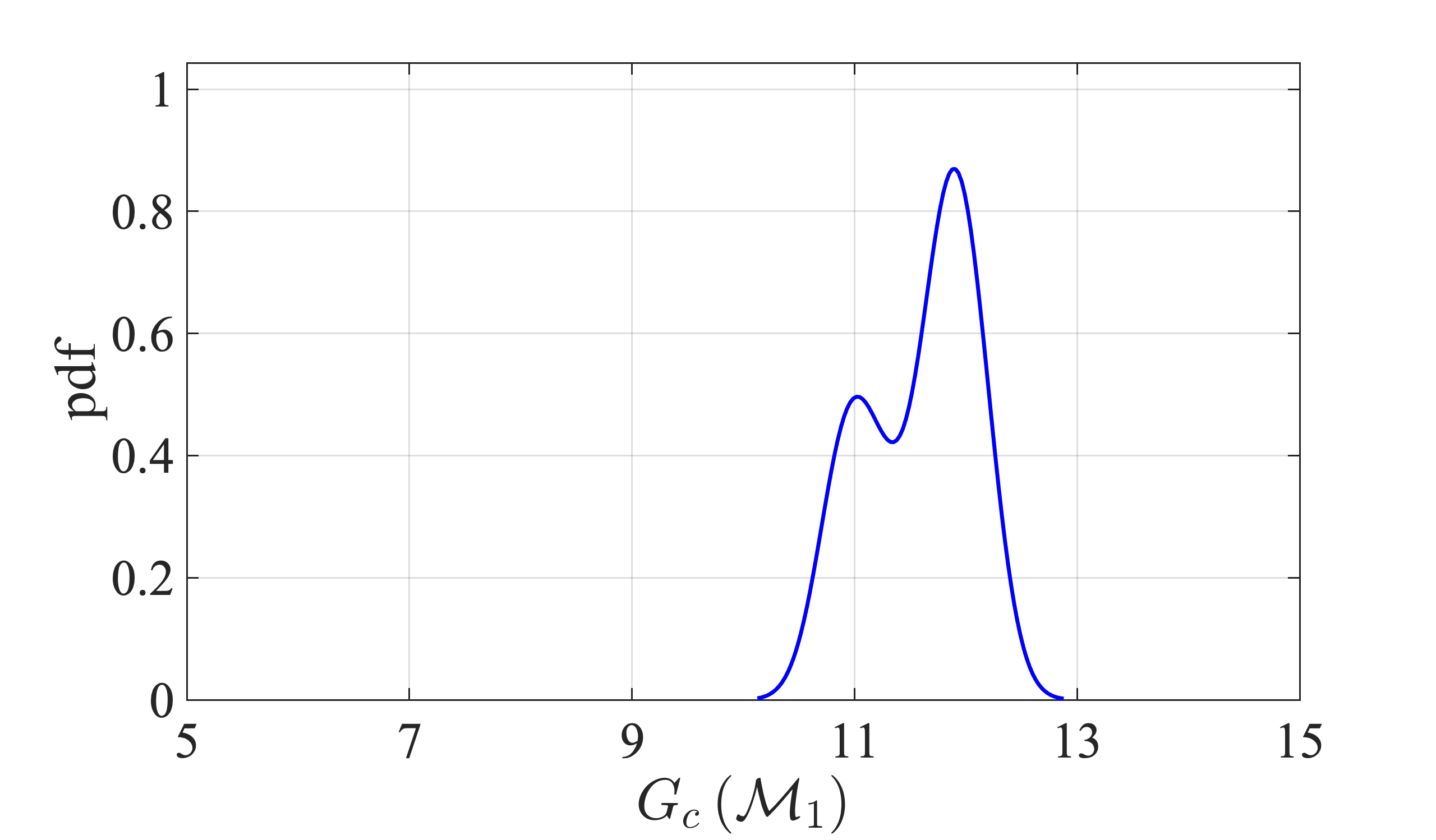}}
	\hfill 
	\vspace{0.5cm}
	\hspace{-1.1cm}
	\subfloat{\includegraphics[width=0.37\textwidth]{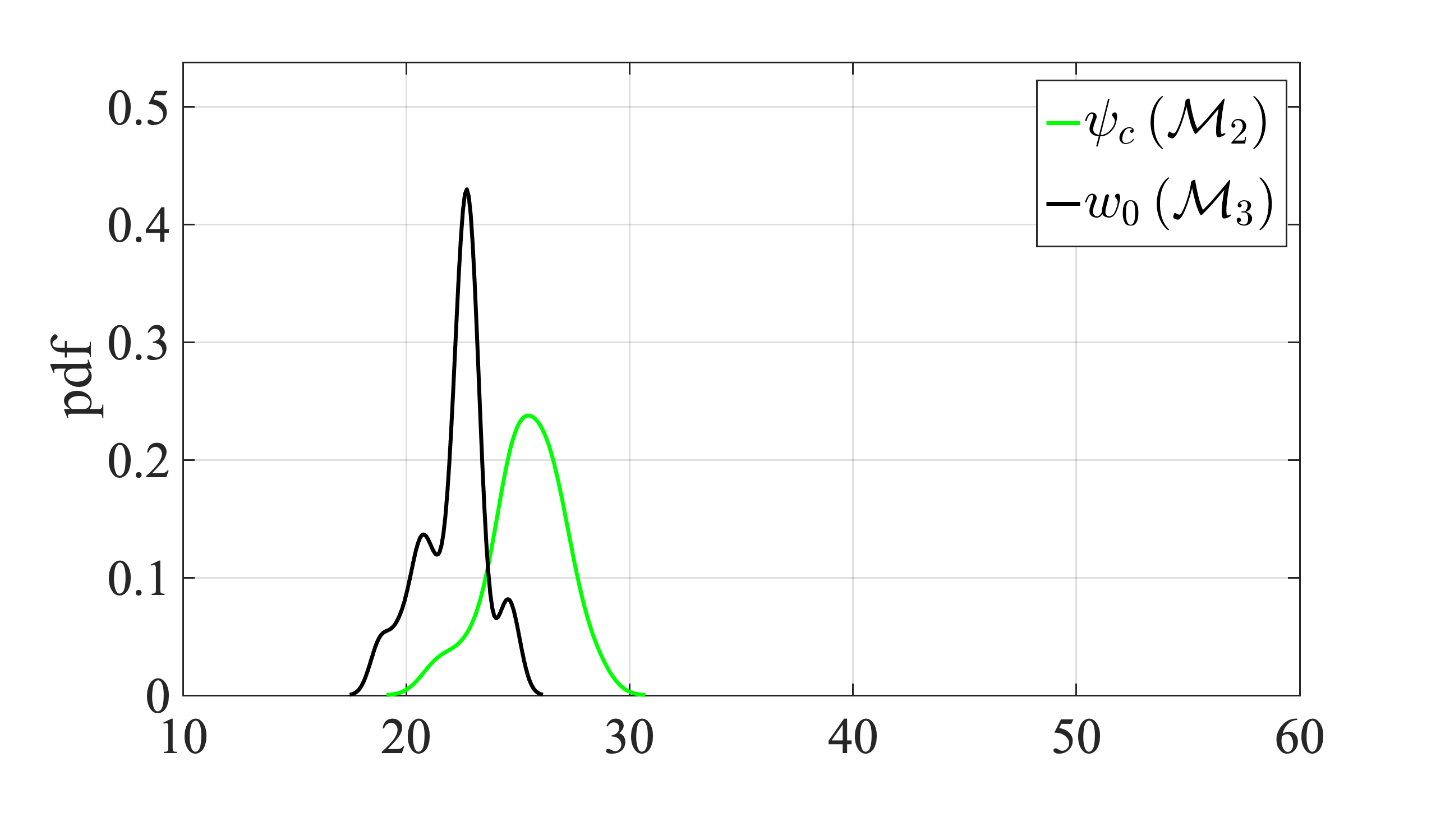}}   \subfloat{\includegraphics[width=0.37\textwidth]{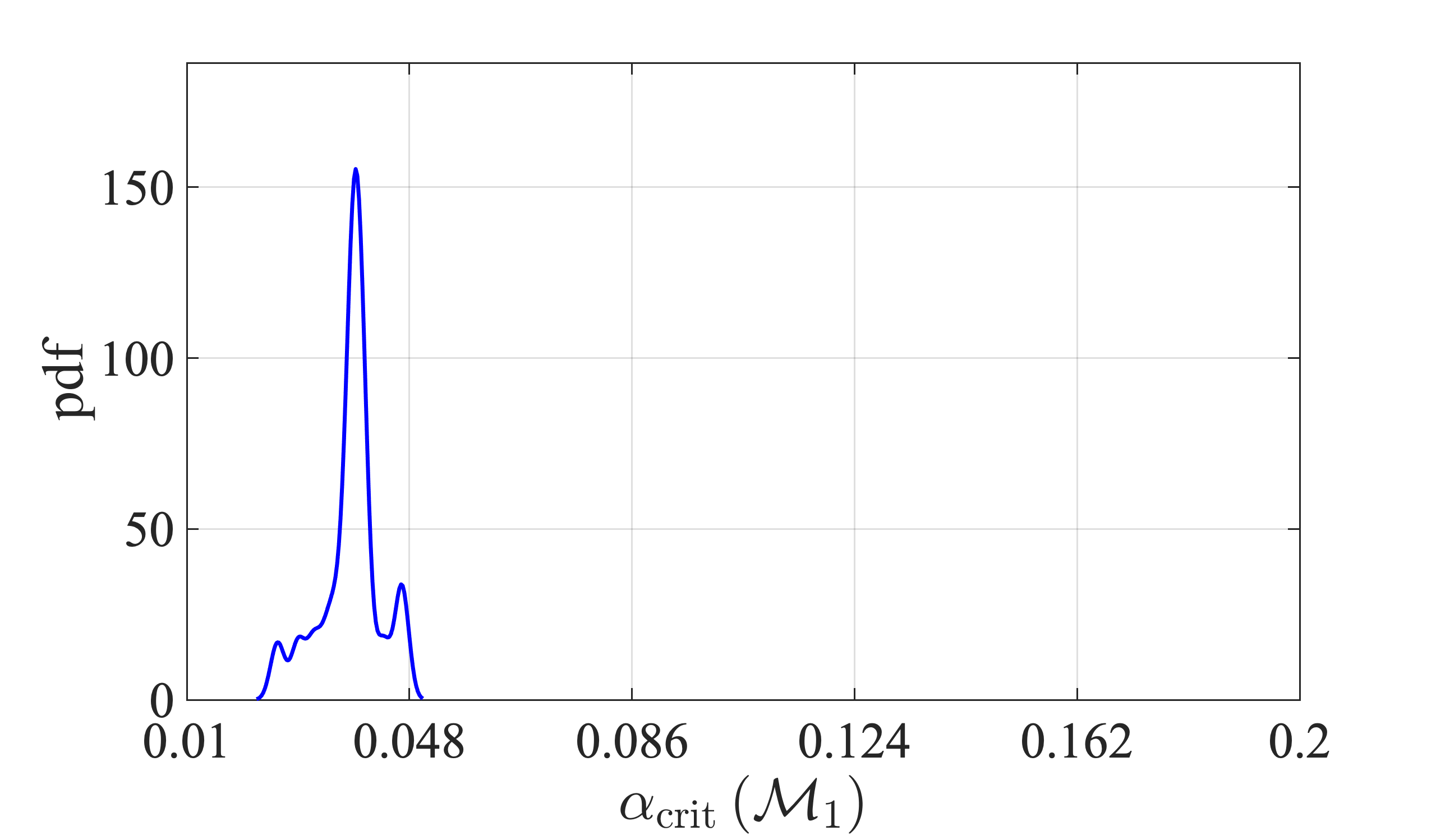}}  	\subfloat{\includegraphics[width=0.37\textwidth]{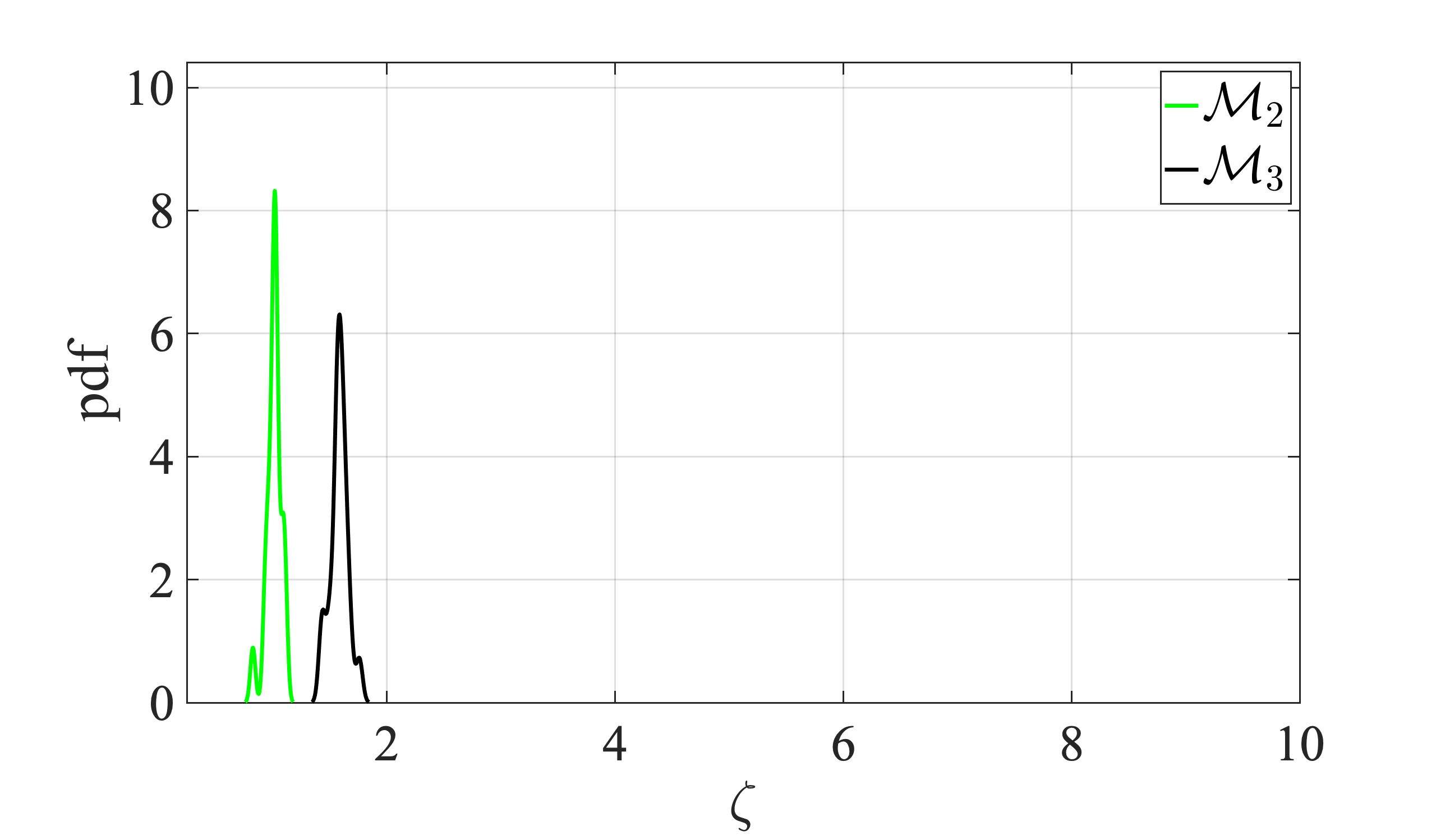}}
	\hfill 
	\vspace{0.1cm}
	\hspace{4.75cm}
	\subfloat{\includegraphics[width=0.37\textwidth]{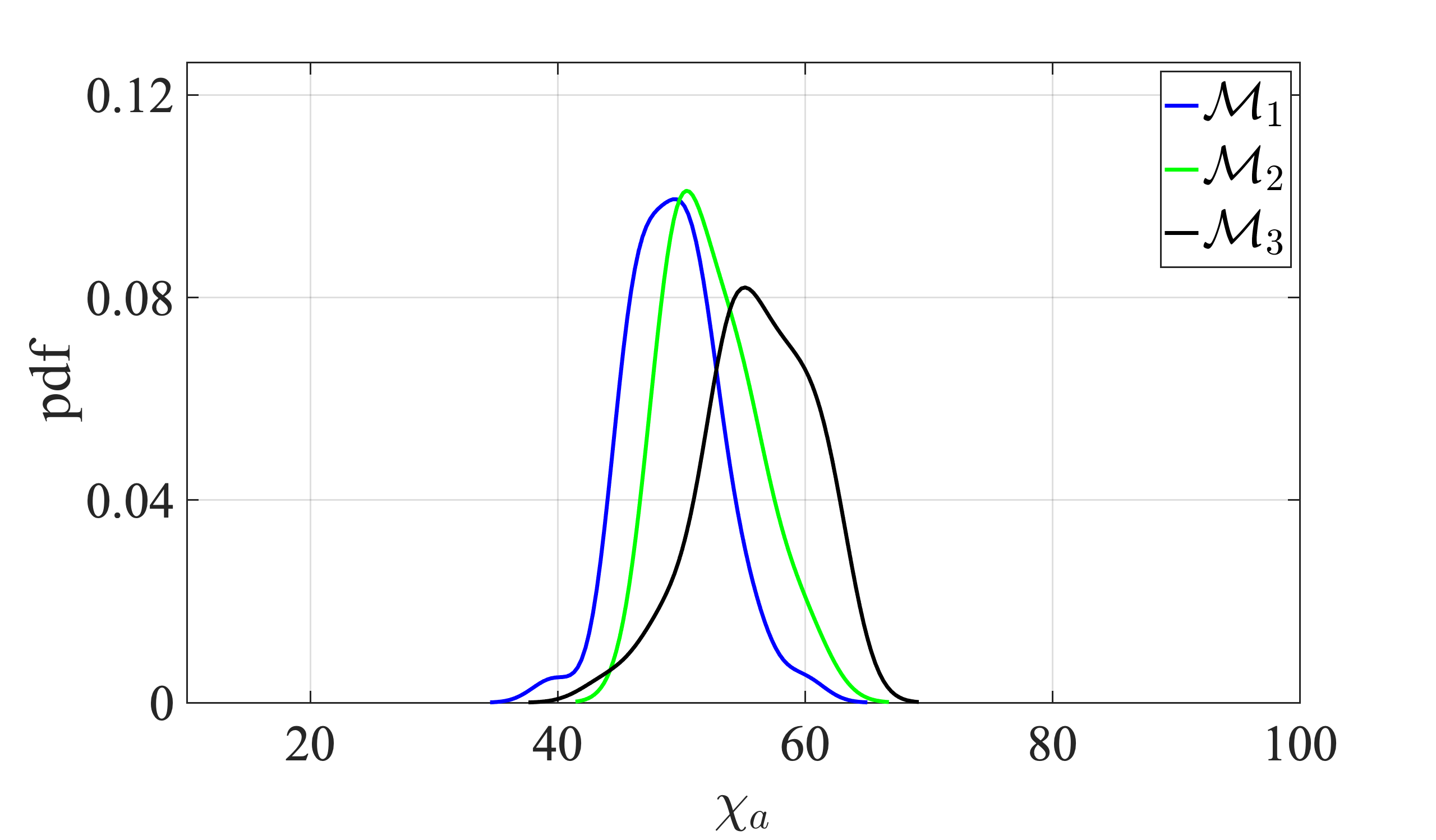}}  
	\caption{Example 2: the posterior distribution of the effective parameters using $\mathcal{M}_1$, $\mathcal{M}_2$, and $\mathcal{M}_3$.}
	\label{Example2_models}
\end{figure}

\begin{figure}[!]
	\centering
	{\includegraphics[width=0.85\textwidth]{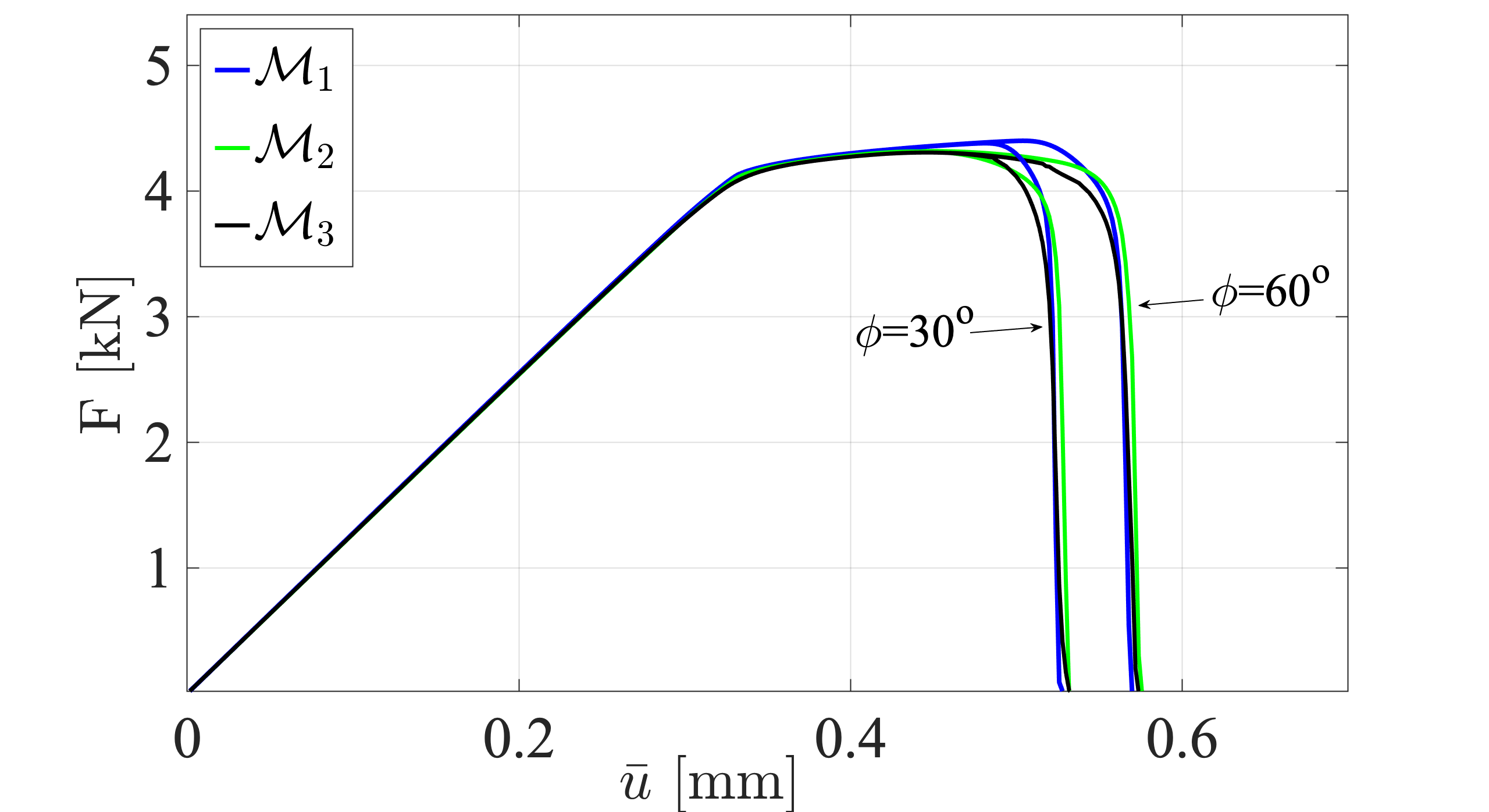}}  
	\caption{Example 2: the load-displacement curve obtained by the estimated values for $\phi=30^o$ and $\phi=60^o$ and all three models.
	}
	\label{example2-phi}
\end{figure}

In this example, to determine the effective mechanical parameters, the DRAM algorithm is used. Due to the anisotropic structure of the solid, in addition to the identified parameters in Example 1, the stiffness parameter $\chi_a$ must be estimated. Here, we select $\phi=45^o$ for the parameter identification and propose $N=10\,000$ candidates. A synthetic reference value using $\mathcal{N}_{\text{dof}}=26\,870$ is employed as the reference observation considering the already mentioned parameter variation. 
The prior (uniform) densities of the parameters are listed in Table \ref{range2}.  For each parameter, the inferred values by different models are relatively similar, which highlights the robustness of the Bayesian setting. The posterior distributions are depicted in Figure \ref{Example2_models}. The mean values of the posterior distributions are given in Table \ref{Example2_posterior} using $\mathcal{M}_1$, $\mathcal{M}_2$, and $\mathcal{M}_3$. Again, we solve the model equations employing the identified parameters to verify the effectiveness of the Bayesian framework, as shown in Figure \ref{example12} (right). All diagrams show that implementing the DRAM algorithm gives rise to a reasonable agreement between the models and also the reference observation.

\begin{figure}[!]
	\centering
	{\includegraphics[clip,trim=0cm 6cm 16cm 4.5cm, width=13cm]{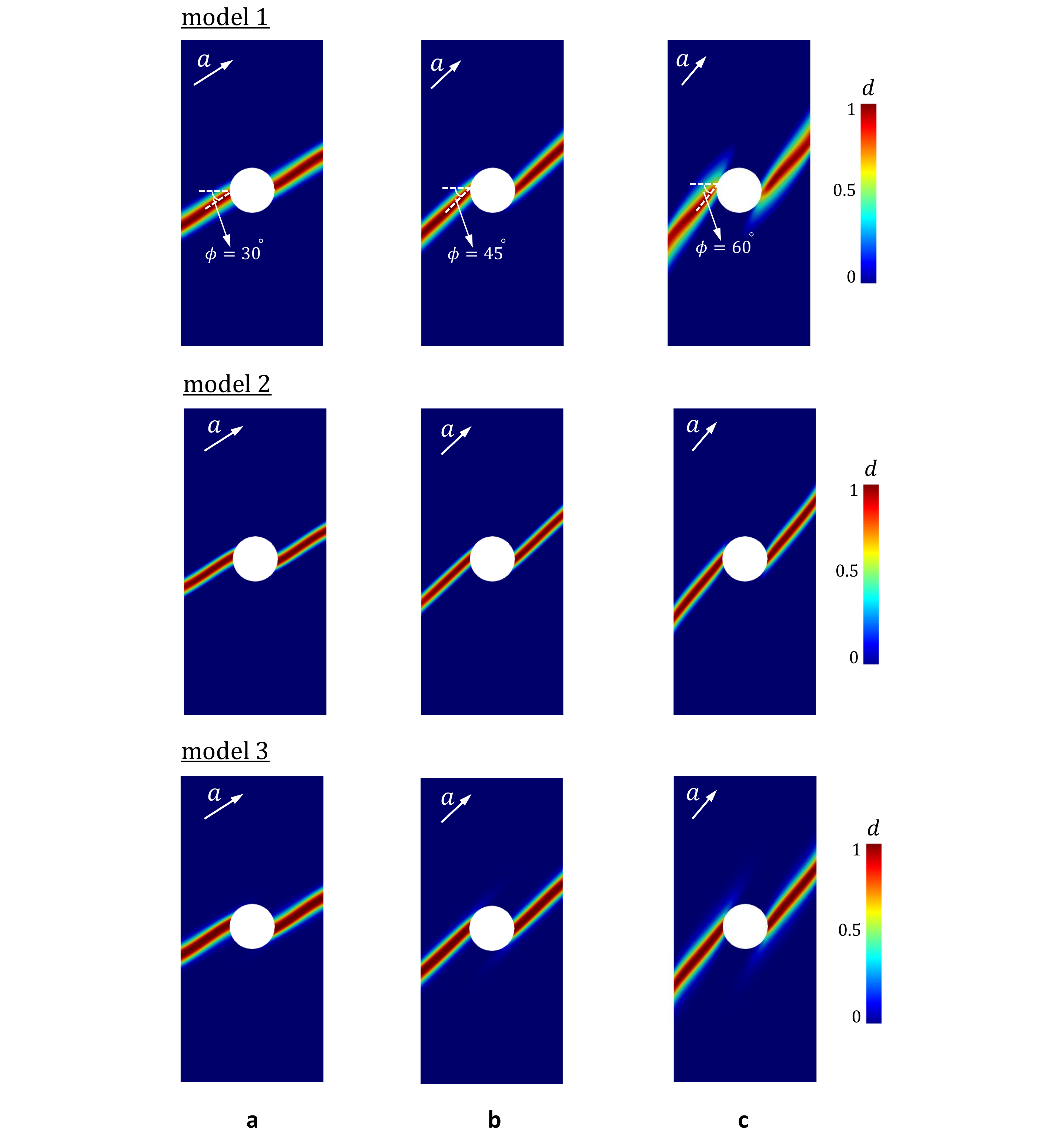}}  
	\caption{ Example 2: the crack phase-field at complete failure for all three models. Different preferential fiber directions are considered in the transversely isotropic setting: (a) $\phi=30^o$, (b) $\phi=45^o$, and (c) $\phi=60^o$.}
	\label{Figure22}
\end{figure}

For further investigation of the sensitivity of the inferred parameters obtained through the DRAM algorithm, two additional preferential fiber directions, namely $\phi=30^o$ and $\phi=60^o$, are employed. Again, Figure \ref{example2-phi} illustrates the robustness of the Bayesian setting, showing for all models and different orientations a consistent behavior in the in elastic, plastic, and fracture stages.  Figures \ref{Figure22} illustrates the crack phase-field solution at complete failure for $\phi=30^o$, $\phi=45^o$, and $\phi=60^o$. Note that the solutions are based on the posterior density of the material parameters, which are given in Table \ref{Example3_posterior}. An important observation is that both the equivalent plastic strain and the crack phase-field evolve in the direction of the preferred fiber orientation. 
\begin{table}[!]
	\caption{Example 2: The uniform prior distribution of the inferred parameters.}
	\vspace{1mm}
	\centering
	\begin{tabular}{llcccccccccccccc}
		&Parameter  &$H$        & $\mu$    &   $K$   &   $\sigma_Y$   &$G_c$&$\alpha_{\text{crit}}$     &$\psi_c$   &$w_0$ &$l_p$ &$\zeta$ &$\chi_a$  \\[2mm]\hline\\
		&   min   &150   & 20\,000  & 40\,000 &275 & 5 & 0.01 &10&10&0.5& 0.25& 10\\[4mm]
		& max &   375  &  40\,000  &  100\,000 & 400 &15 &0.2 & 60 &60&10& 10&100\\[4mm]
		\hline
		\label{range2}
	\end{tabular}
\end{table}
\begin{table}[!]
	\caption{Example 2: the mean value of posterior density of  the model parameters for the three models.}
	\vspace{1mm}
	\centering
	\begin{tabular}{ccccccccccccccc}
		&Model  &$H$        & $\mu$    &   $K$   &   $\sigma_Y$   &$G_c$&$\alpha_{\text{crit}}$     &$\psi_c$   &$w_0$ &$l_p$&$\zeta$ &$\chi_a$  \\[2mm]\hline\\
		&    $\mathcal{M}_1$   &265   & 26\,050 & 94\,010 &355 & 11.6 & 0.038  &--&--&--&-- &50 \\[4mm]
		& $\mathcal{M}_2$ &   245  & 26\,100  &  92\,100 & 354 &-- &-- & 25.25& &--&1.01&52\\[4mm]
		&   $\mathcal{M}_3$  &220  & 26\,300   & 88\,950 &340 & -- & -- &--& 22 & 8.29&1.6 &55 \\[4mm]
		\hline
		\label{Example2_posterior}
	\end{tabular}
\end{table}
\begin{table}[!]
	\caption{Example 3: The uniform prior distribution of the inferred parameters.}
	\vspace{1mm}
	\centering
	\begin{tabular}{llcccccccccccccc}
		&Parameter  &$H$        & $\mu$    &   $K$   &   $\sigma_Y$   &$G_c$&$\alpha_{\text{crit}}$     &$\psi_c$   &$w_0$ &$l_p$ &$\zeta$ & \\[2mm]\hline\\
		&   min   &10   & 20\,000  & 40\,000 &50 & 100 & 0.001 &5&5&0.001& 1\\[4mm]
		& max &   50  &  40\,000  &  100\,000 & 200 &300 &0.1 & 25 &25&10& 20\\[4mm]
		\hline
		\label{range33}
	\end{tabular}
\end{table}
\begin{table}[!]
	\caption{Example 3: the mean value of posterior density of  the model parameters for the three models.}
	\vspace{1mm}
	\centering
	\begin{tabular}{ccccccccccccccc}
		&Model  &$H$        & $\mu$    &   $K$   &   $\sigma_Y$   &$G_c$&$\alpha_{\text{crit}}$     &$\psi_c$   &$w_0$ &$l_p$&$\zeta$ & \\[2mm]\hline\\
		&    $\mathcal{M}_1$   &30   & 26\,500 & 73\,500 &115 & 248 & 0.0142  &--&--&--&--  \\[4mm]
		& $\mathcal{M}_2$ &   15  & 26\,200  &  73\,700 & 116 &-- &-- & 13.4& &--&2\\[4mm]
		&   $\mathcal{M}_3$  &15  & 30\,100   & 75\,050 &112 & -- & -- &--& 10.3 & 0.0018&15  \\[4mm]
		\hline
		\label{Example3_posterior}
	\end{tabular}
\end{table}


 \sectpb[Section53]{ Example 3: Flat I-shaped Al-5005 test under tensile loading}

\begin{figure}[!ht]
	\centering
	{\includegraphics[clip,trim=4cm 17cm 7cm 5cm, width=16cm]{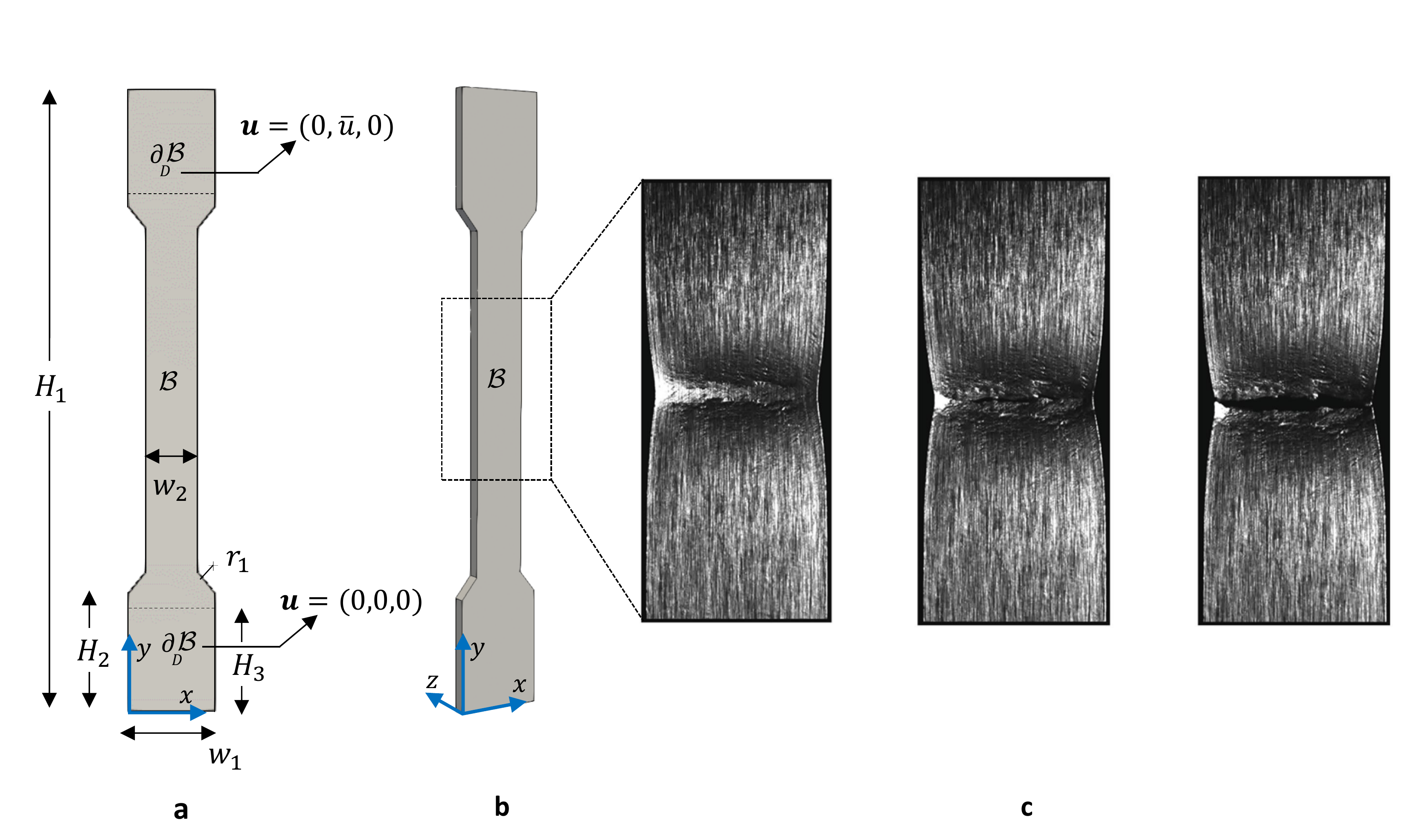}}  
	\caption{Example 3: flat I-shaped Al-5005 test under tensile loading. (a) Geometry of specimen with boundary conditions, (b) three-dimensional perspective, and (c) experimental observation taken from \cite{ambati2016}.}
	\label{exm3bvp}
\end{figure}

While synthetic observations have been used so far, the last two examples are concerned with experimental observations to estimate the posterior density of the material unknowns. The following example considers an I-shaped specimen for Al-5005 material under tensile loading. We aim at reproducing the experimentally observed ductile fracture process of the tensile test through the proposed Bayesian inversion framework. The BVP is shown in Figure \ref{exm3bvp}a. The experimental observations of necking and fracture are shown in Figure \ref{exm3bvp}c.  The geometrical dimensions are set as $H_1=144$ mm, $H_2=27$ mm, $H_3=22$ mm, $w_1=20$ mm, $w_2=12$ mm, and $r_1=14$ mm. The specimen domain has a $3$ mm thickness, as shown in Figure \ref{exm3bvp}b.

The numerical example is performed by applying a monotonic displacement increment  ${\Delta \bar{u}}_y=0.02$ mm in the vertical direction at the top boundary of the specimen for 300 time steps. The minimum finite element size is $1.5$ mm. The flat I-shaped domain partition contains 3230 hexahedron linear elements.

Both DRAM and ensemble-Kalman filter (EKF) are efficient MCMC techniques and have shown their computational performance reasonably. However, a fair comparison can determine which method will be more advantageous in ductile fracture.  
 
\sectpc{Convergence performance of the MCMC methods}
 In the already mentioned examples, the MCMC techniques have been used to identify the mechanical parameters. The main advantage of the MH algorithm is its ease of implementation. However, it suffers from slow convergence, and the starting value may affect the convergence status. The DRAM and EKF variants, as more advanced techniques, show more a productive performance. In this part, we strive to study their efficiency in the context of ductile fracture. 
 
Convergence diagnostics is essential in MCMC methods since it determines the accuracy of the parameter, and with how many iterations the chain converges to the target distribution. Here, we use multiple chains with different initialized values, expecting that a significantly large number of Markov chains gives rise to the same results. In other words, the candidate distribution from chains should be similar using multiple chains initial starting values.

 $\hat{R}$-convergence diagnostics \cite{brooks1998general,gelman1992inference} is an efficient tool to monitor the convergence of the MCMC by comparing the between and within chain estimates for model parameters and other univariate quantities of interest. Assuming $m$ parallel chains, we determine the variance between the chain means $B/N$ and calculate the average of the within chain variances $W$. The target variance is given by
\begin{align}
\mathcal{S}^2=(1-\frac{1}{N})W+\frac{B}{N},
\end{align}
 where $N$ is the length of the chain. Then, we calculate the potential scale reduction factor, or PSRF (also called $\hat{R}$-statistics) by
 \begin{align}
 \hat{R}=\frac{m+1}{m}\frac{\mathcal{S}^2}{W}-\frac{N-1}{m}.
 \end{align}
If the MCMC method converges appropriately, the chains are not affected by the starting point, and $\hat{R}$ reduces to 1. In other words, we can conclude that all chains are close to the target distribution \cite{brooks1998general}. In order to verify the method efficiency, a threshold can be defined, e.g., values less than 1.5 or 1.2 indicating a good convergence performance.
%
%
\begin{figure}[t!]
   	\subfloat{\includegraphics[width=0.51\textwidth]{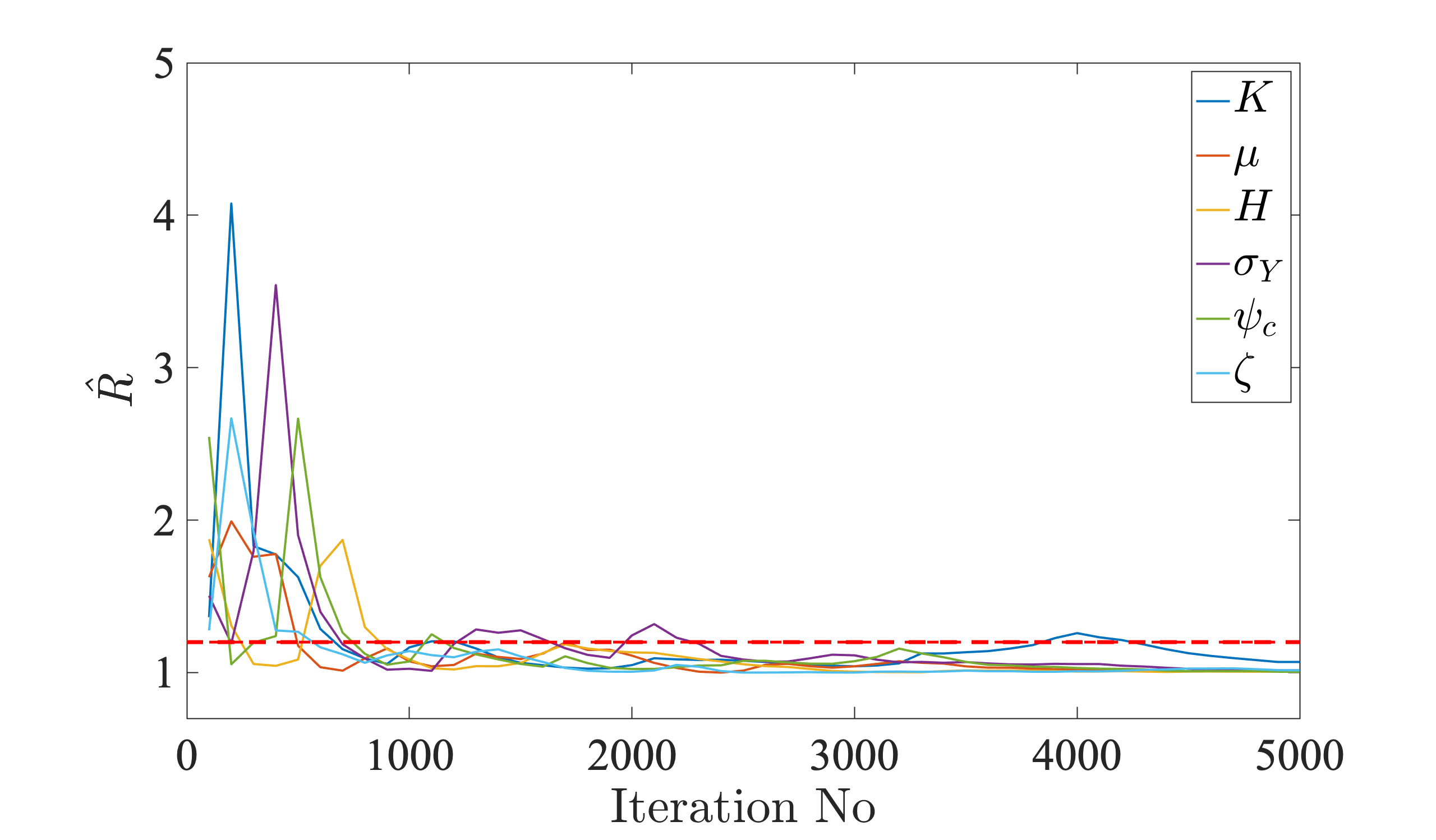}} 	\subfloat{\includegraphics[width=0.5\textwidth]{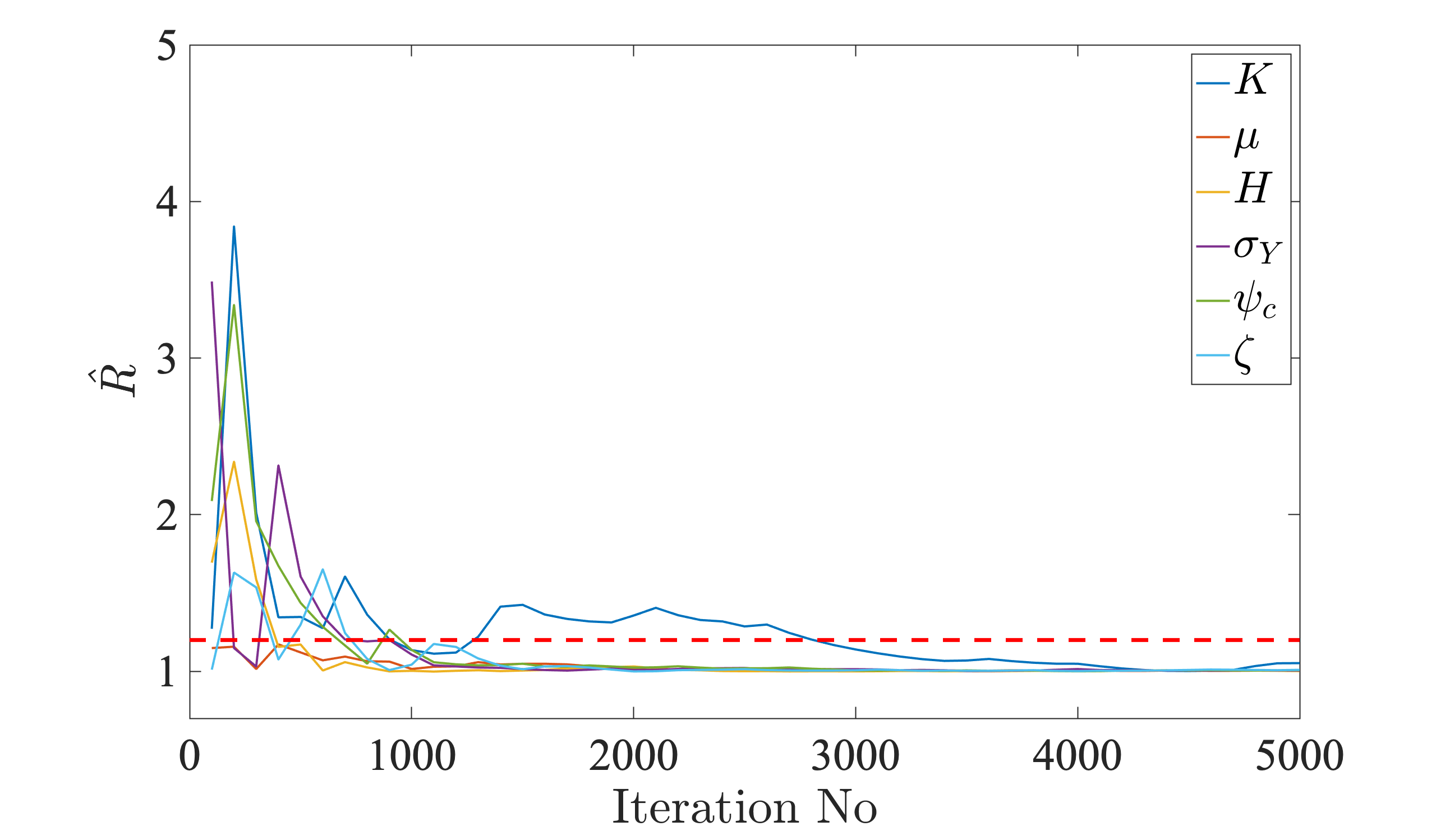}}
   	   	\vspace{0.4cm}  
   	\caption{The $\hat{R}$-statistics test for the convergence of all effective mechanical parameters of the I-shaped example. Here, the DRAM algorithm (left) is compared with EKF (right) using $\mathcal{M}_2$. The red dashed line shows a threshold of 1.2.}
   	\label{convergence2}
\end{figure}
%

In the I-shaped example, in order to draw a comparison between DRAM and EKF algorithms, we study their convergence to conclude which model shows a faster convergence taking all inferred parameters into account. We use $N=1\,000$ and five parallel MCMCs ($m=5$) with a uniform distribution indicated in Table \ref{range33}. Using $\mathcal{M}_2$, Figure \ref{convergence2} shows that by employing EKF in all parameters, fewer candidates are necessary to converge to the posterior density. Indeed, for all parameters excluding $K$, after $1\,500$ samples, $\hat{R}$-statistics converges to 1, showing a high level of accuracy. The performance of the DRAM is acceptable, since most of the variables after $2\,500$ samples are below the threshold, although again, the bulk modulus shows more variation (probably due to the large chosen prior density).
 \begin{figure}[t!]
 	\vspace{0.1cm}
 	\hspace{-1cm}
 	\subfloat{\includegraphics[width=0.37\textwidth]{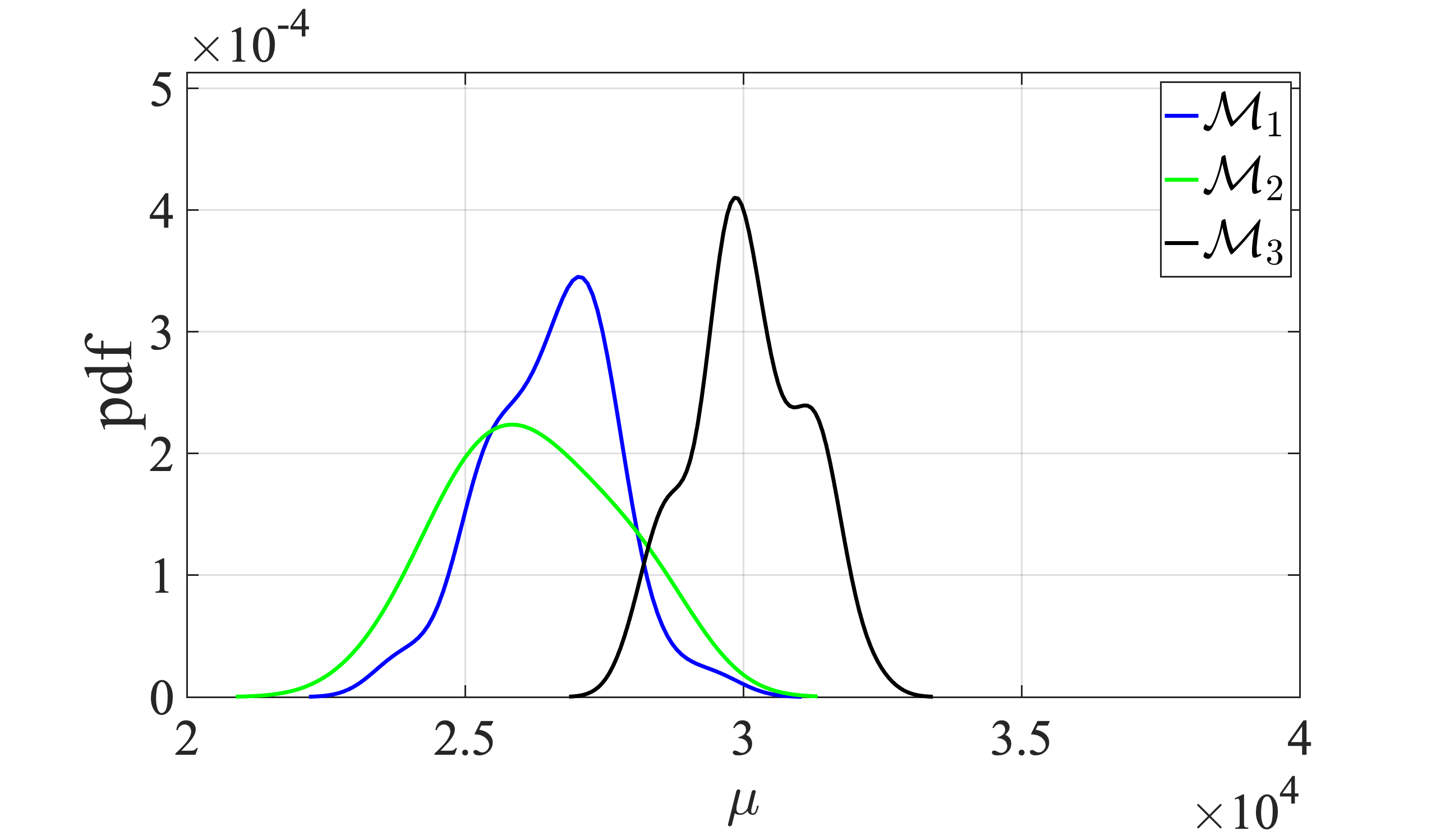}}   \subfloat{\includegraphics[width=0.37\textwidth]{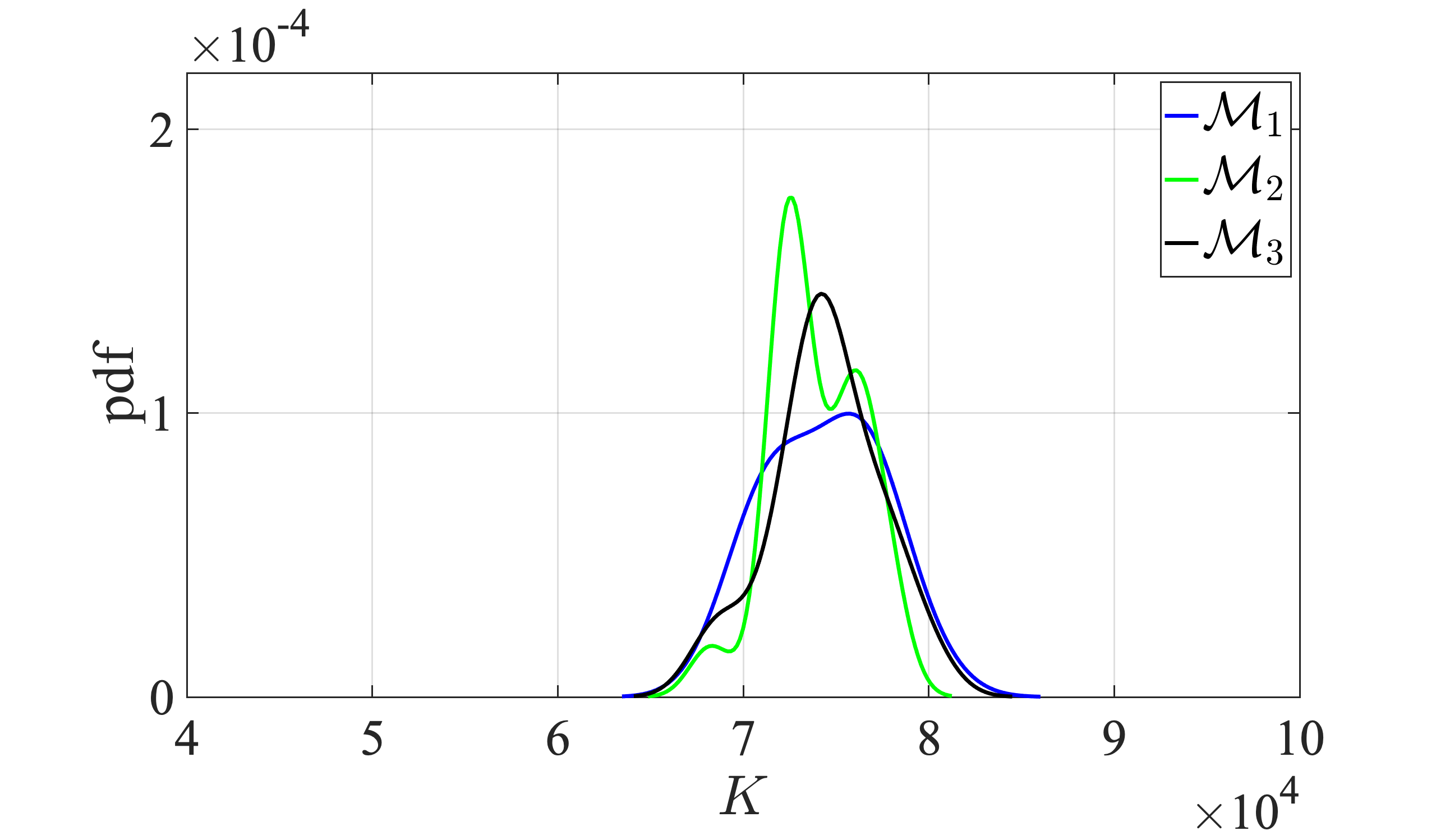}}  	\subfloat{\includegraphics[width=0.37\textwidth]{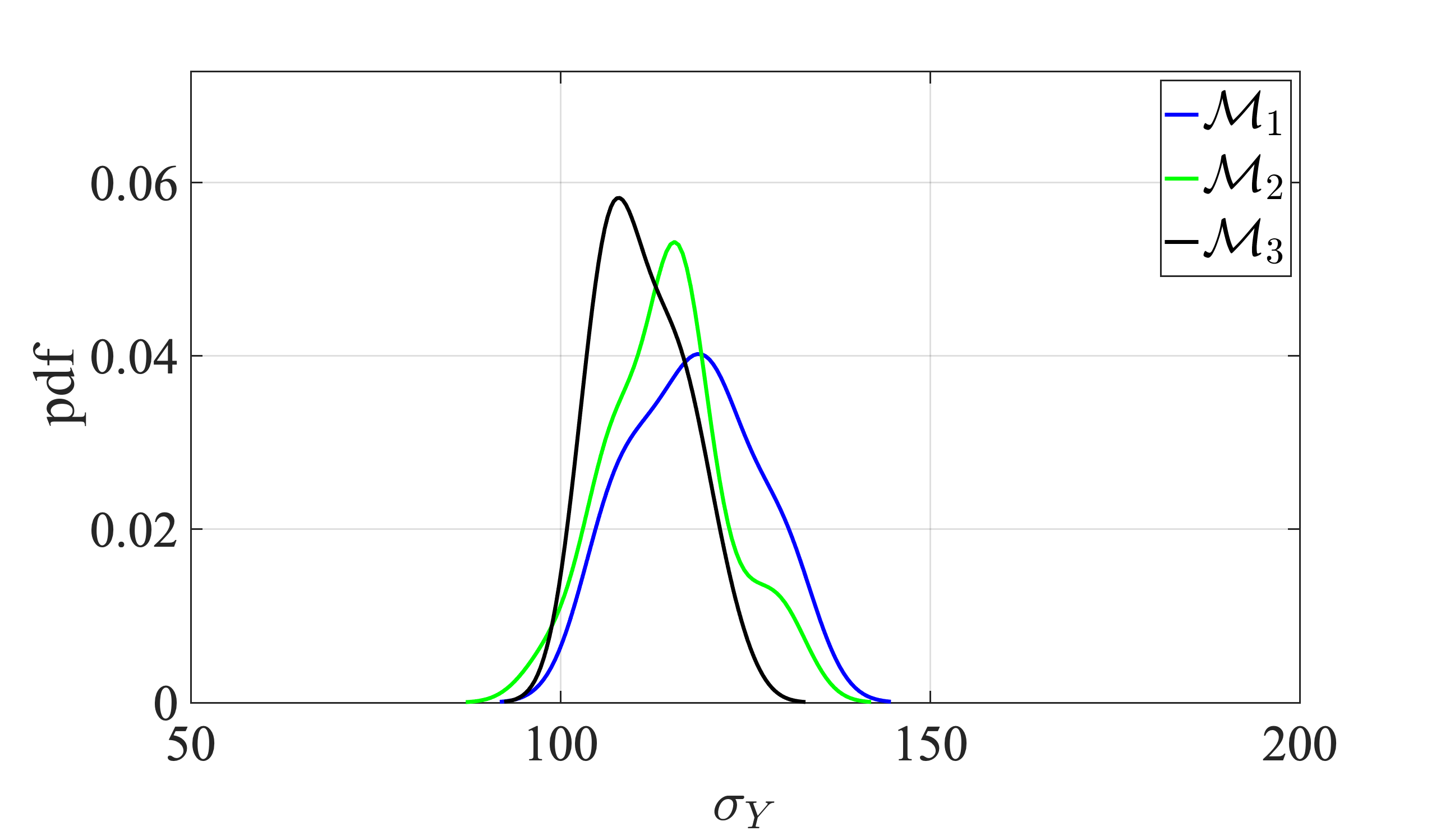}}
 	\hfill 
 	\vspace{0.1cm}
 	\hspace{-1.1cm}
 	\subfloat{\includegraphics[width=0.37\textwidth]{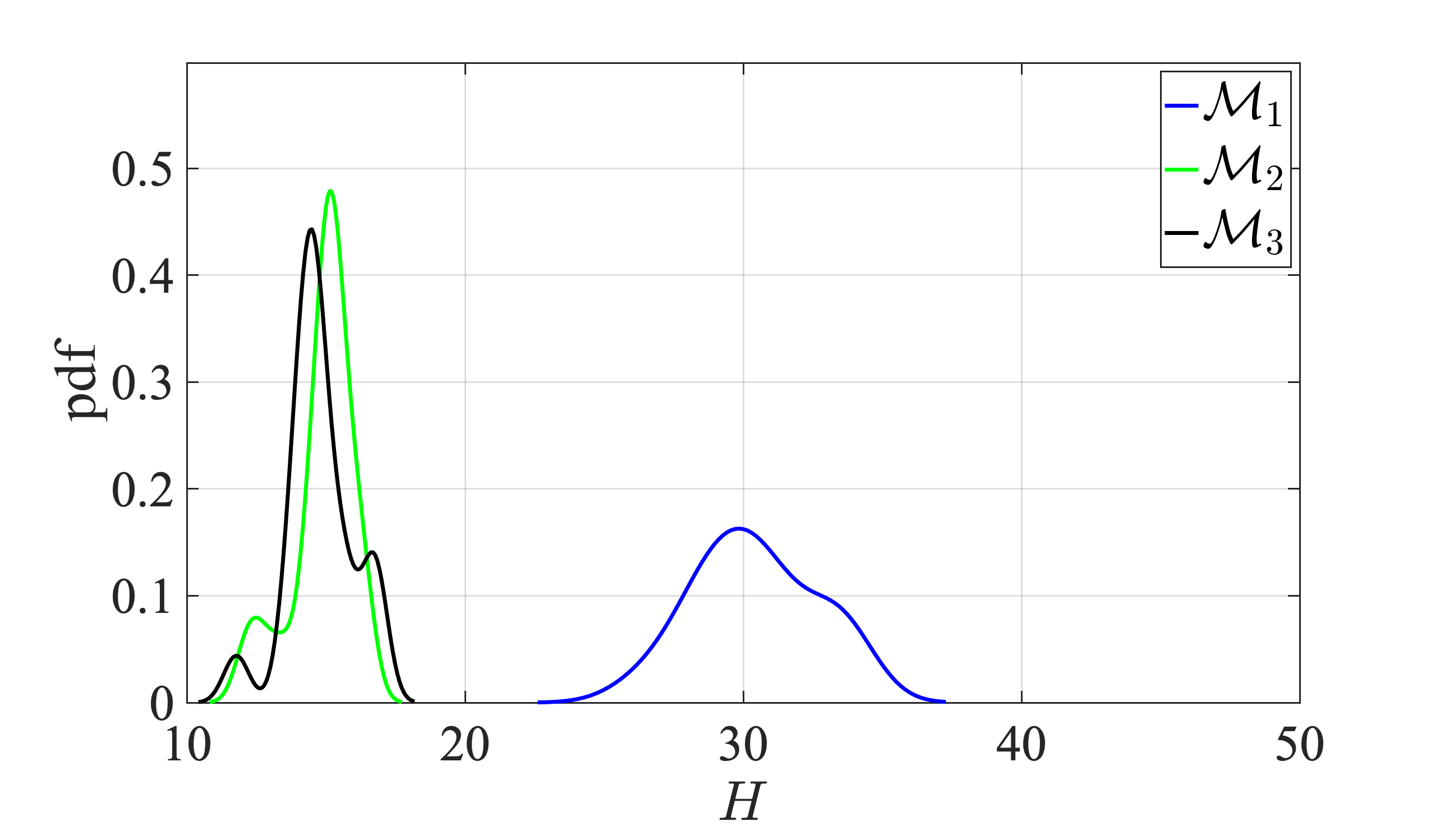}}   \subfloat{\includegraphics[width=0.37\textwidth]{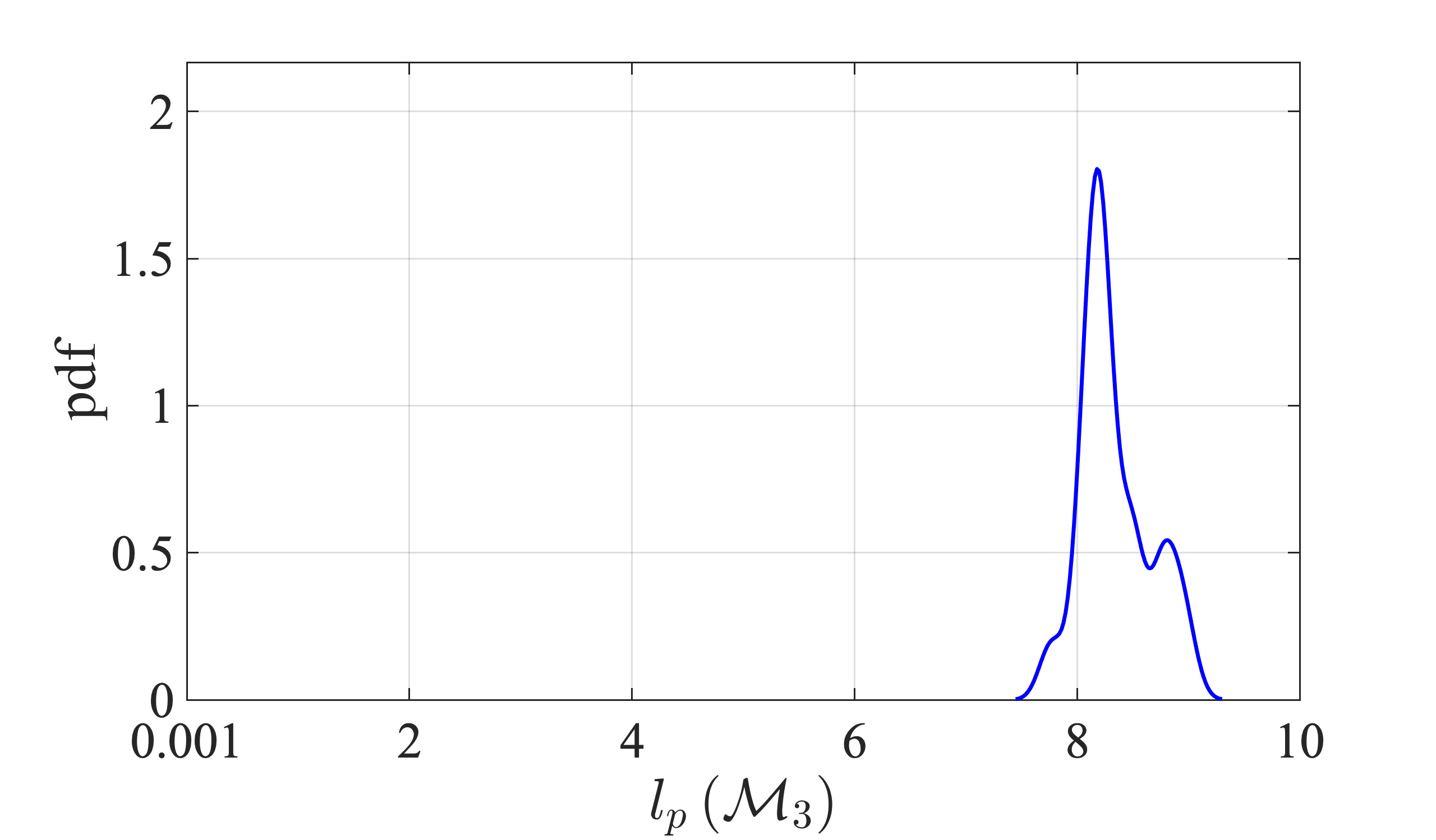}}  	\subfloat{\includegraphics[width=0.37\textwidth]{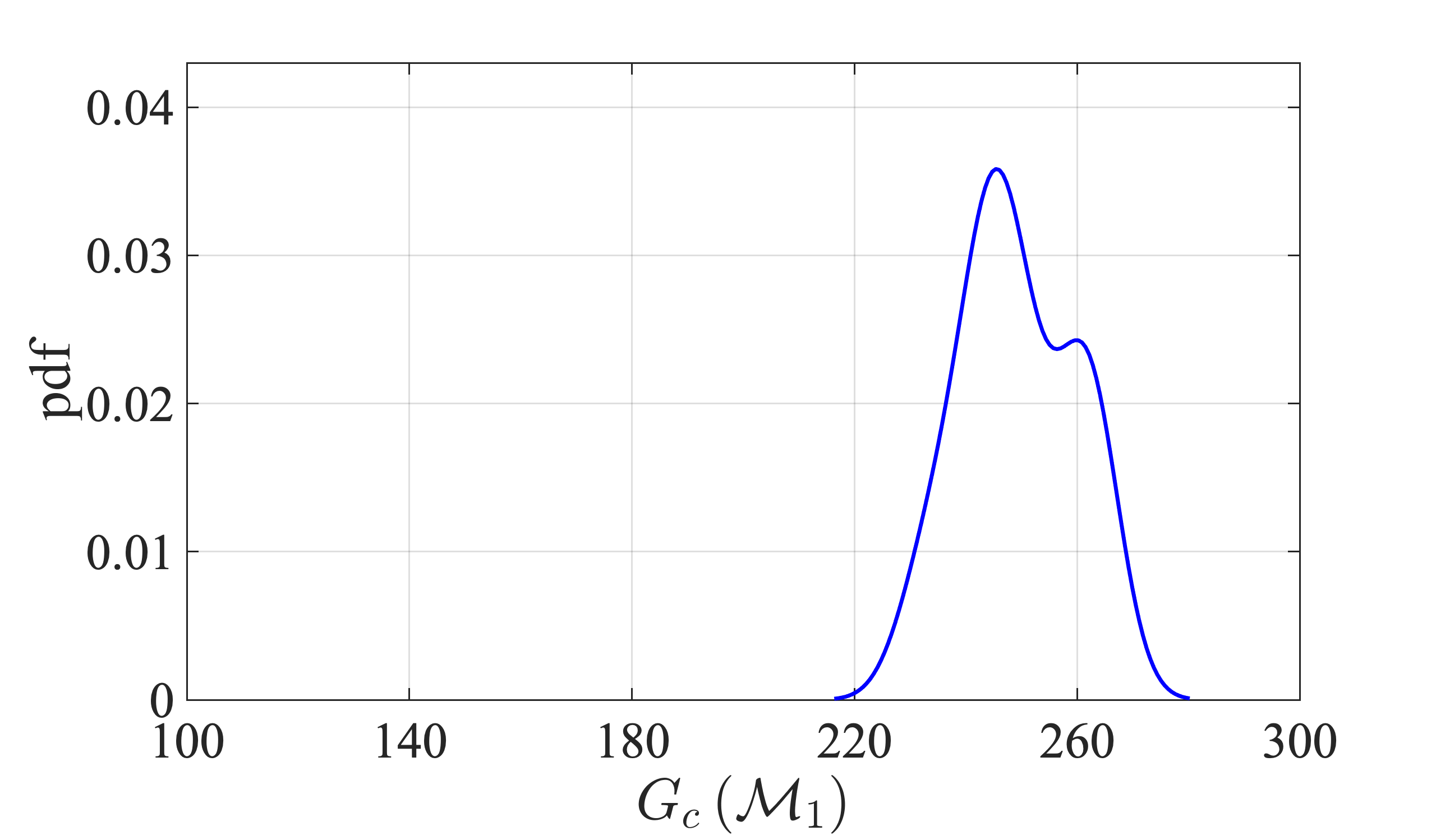}}
 	\hfill 
 	\vspace{0.5cm}
 	\hspace{-1.1cm}
 	\subfloat{\includegraphics[width=0.37\textwidth]{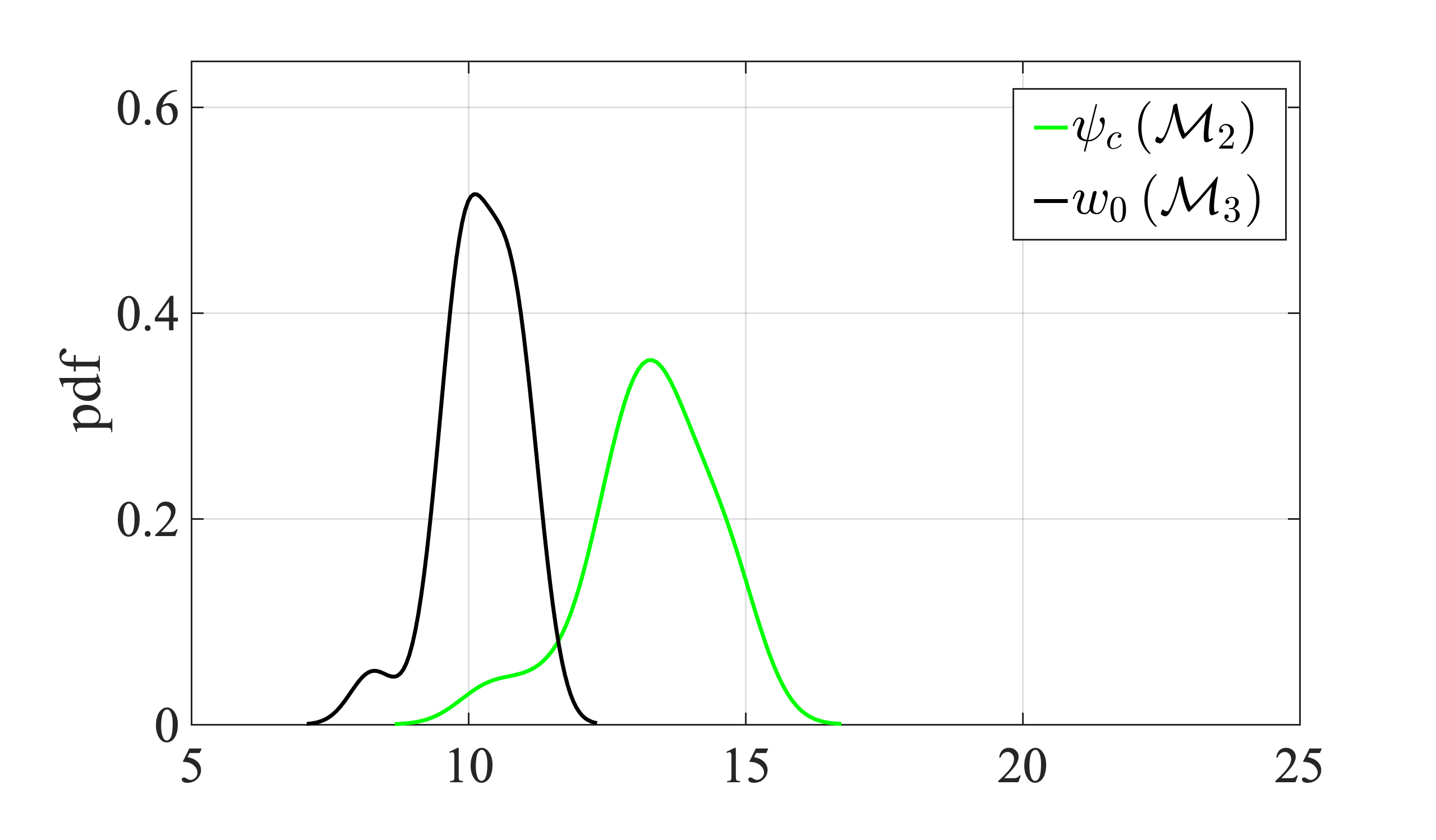}}   \subfloat{\includegraphics[width=0.37\textwidth]{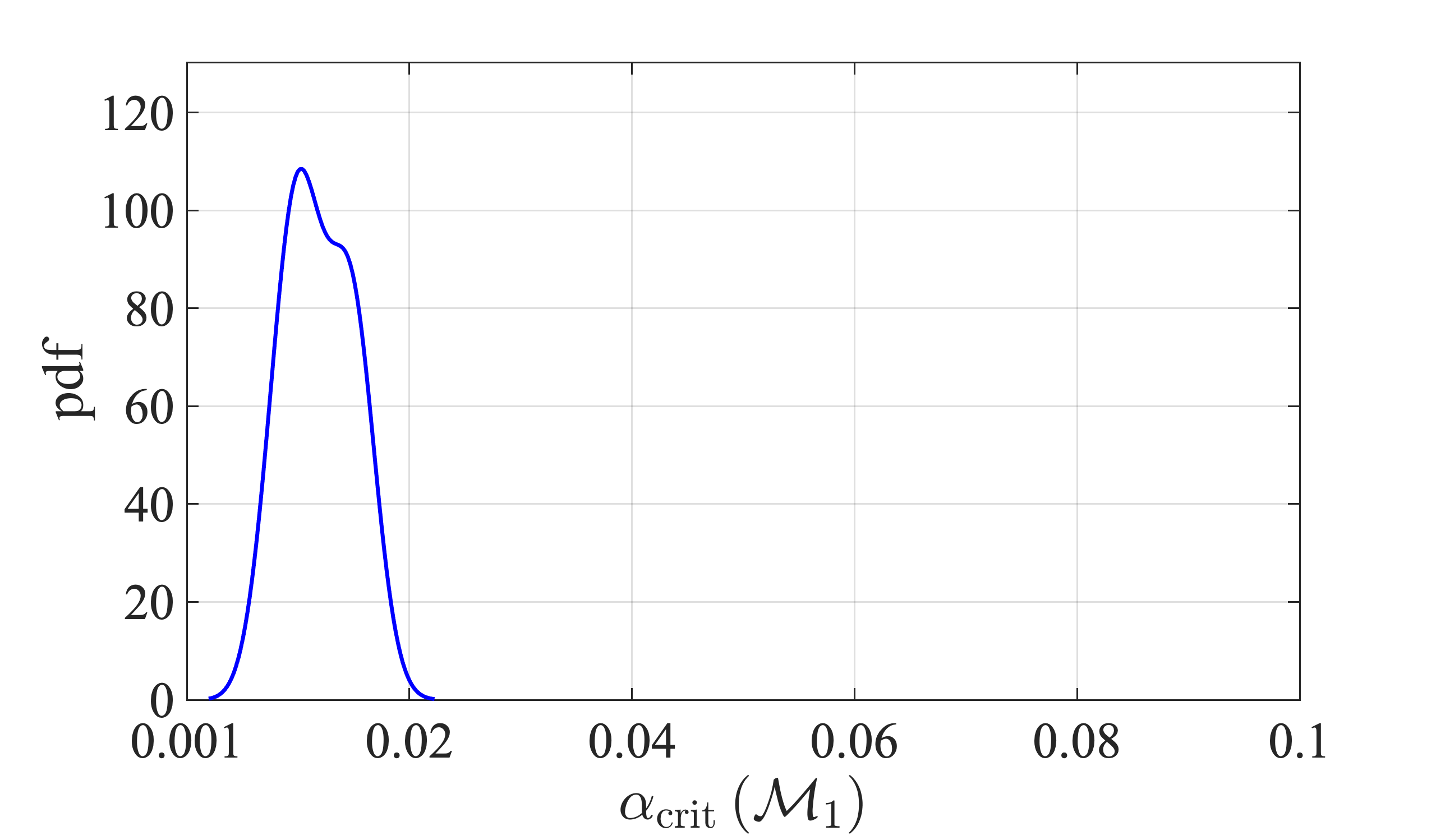}}  	\subfloat{\includegraphics[width=0.37\textwidth]{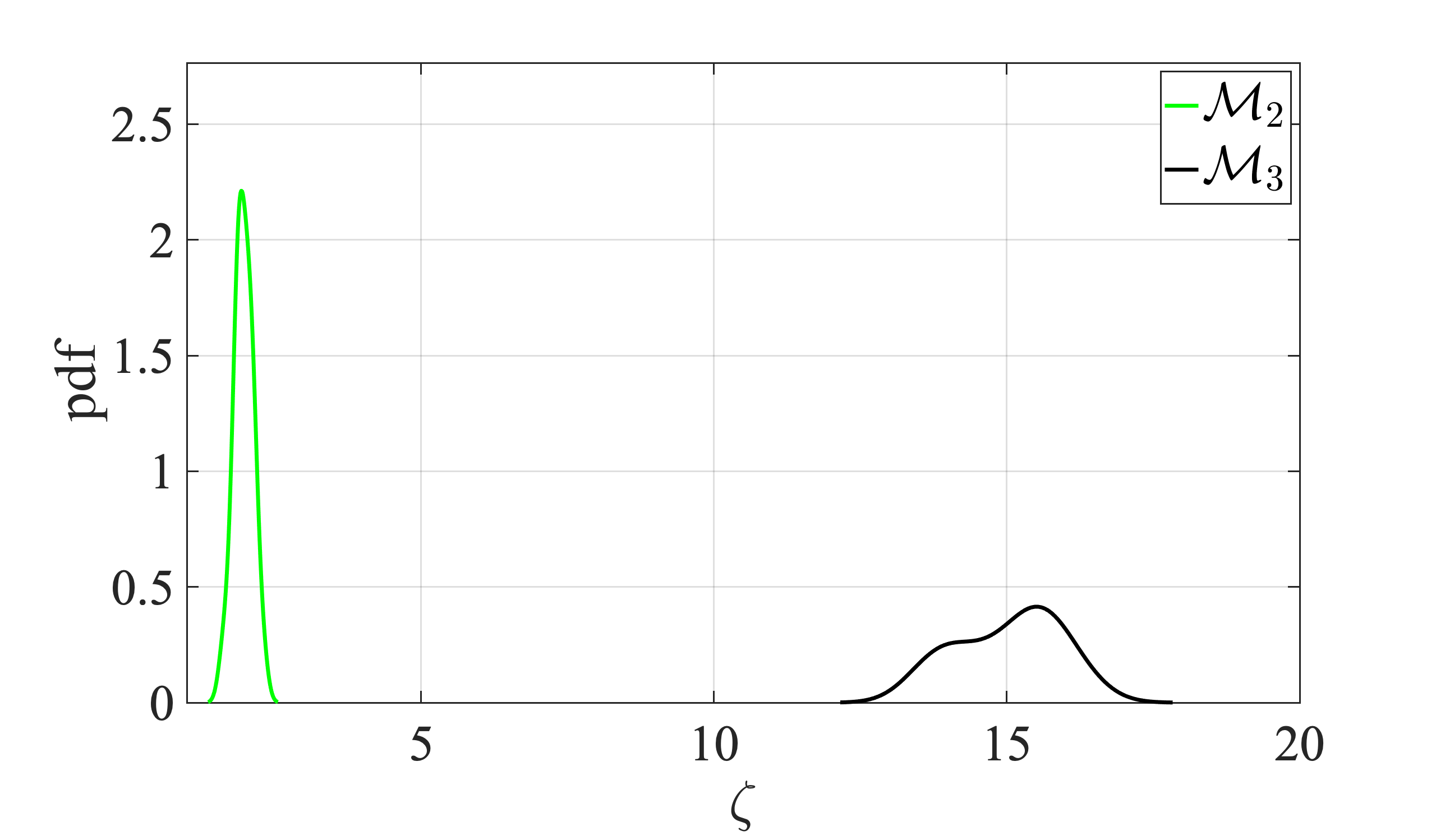}}
 	\caption{Example 3: the posterior distribution of the effective parameter using $\mathcal{M}_1$, $\mathcal{M}_2$, and $\mathcal{M}_3$.}
 	\label{Example3_models}
 \end{figure}
 \begin{figure}[t!]
 	\subfloat{\includegraphics[width=0.5\textwidth]{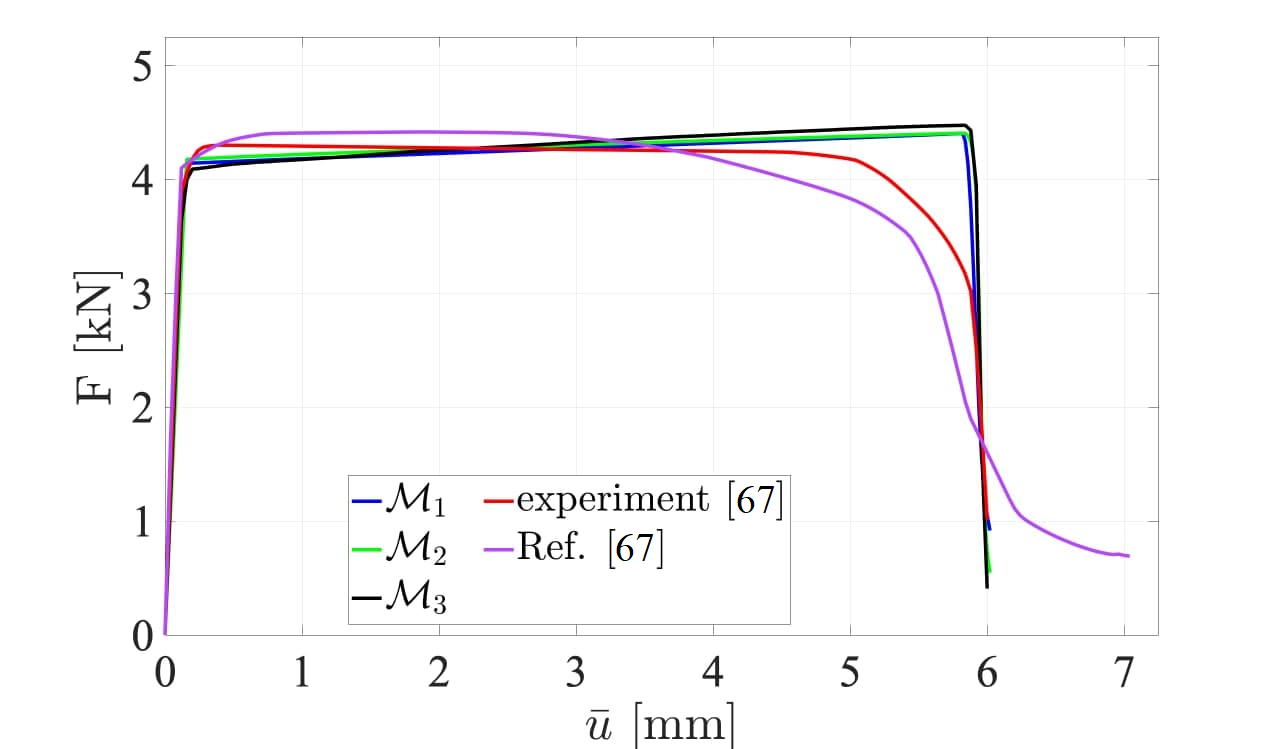}} 	\subfloat{\includegraphics[width=0.5\textwidth]{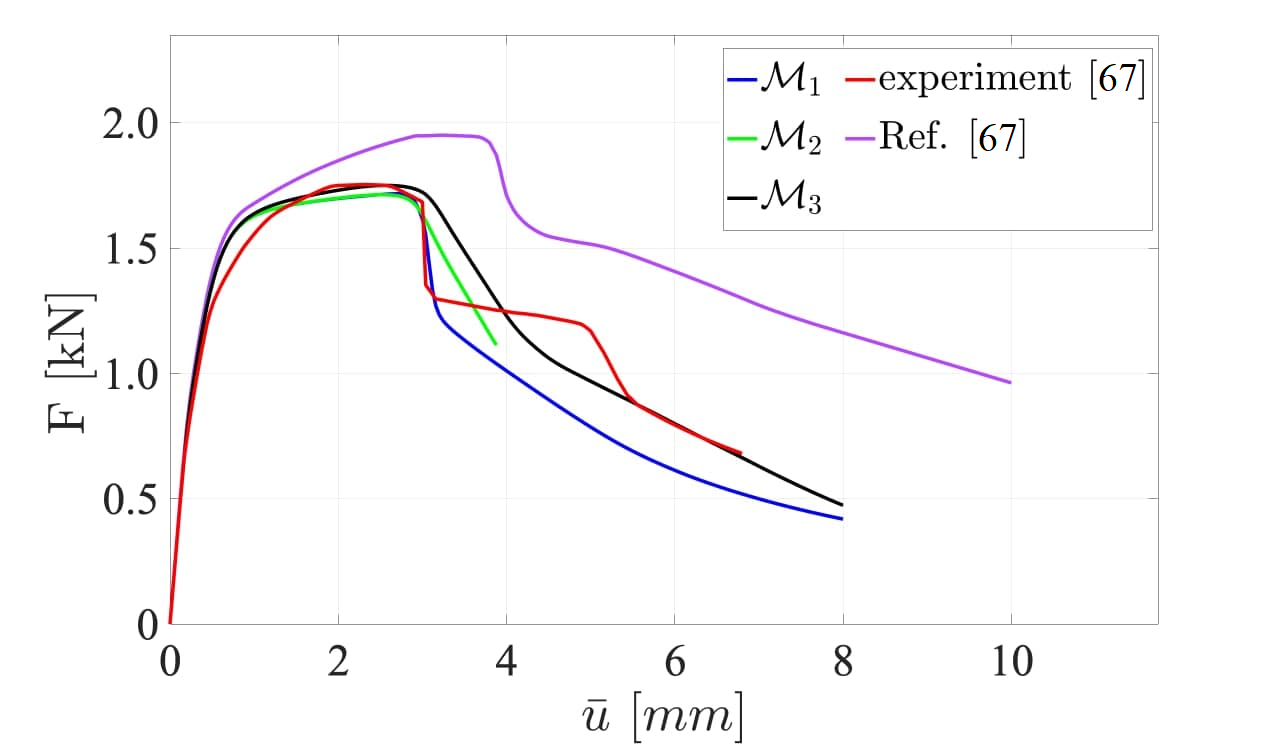}}  
 	\caption{Load-displacement curve computed using the inferred values, employing $\mathcal{M}_1$, $\mathcal{M}_2$, and $\mathcal{M}_3$.  The experimental data are depicted with red. Here, Table \ref{Example3_posterior} and Table \ref{Example4_posterior} are used for Example 3 (left) and Example 4 (right). The results are also compared with the simulation results of \cite{ambati2016}. }
 	\label{example34}
 \end{figure}
 
\begin{figure}[!]
	\centering
	\subfloat{{\includegraphics[clip,trim=2cm 21cm 8cm 4.5cm, width=15cm]{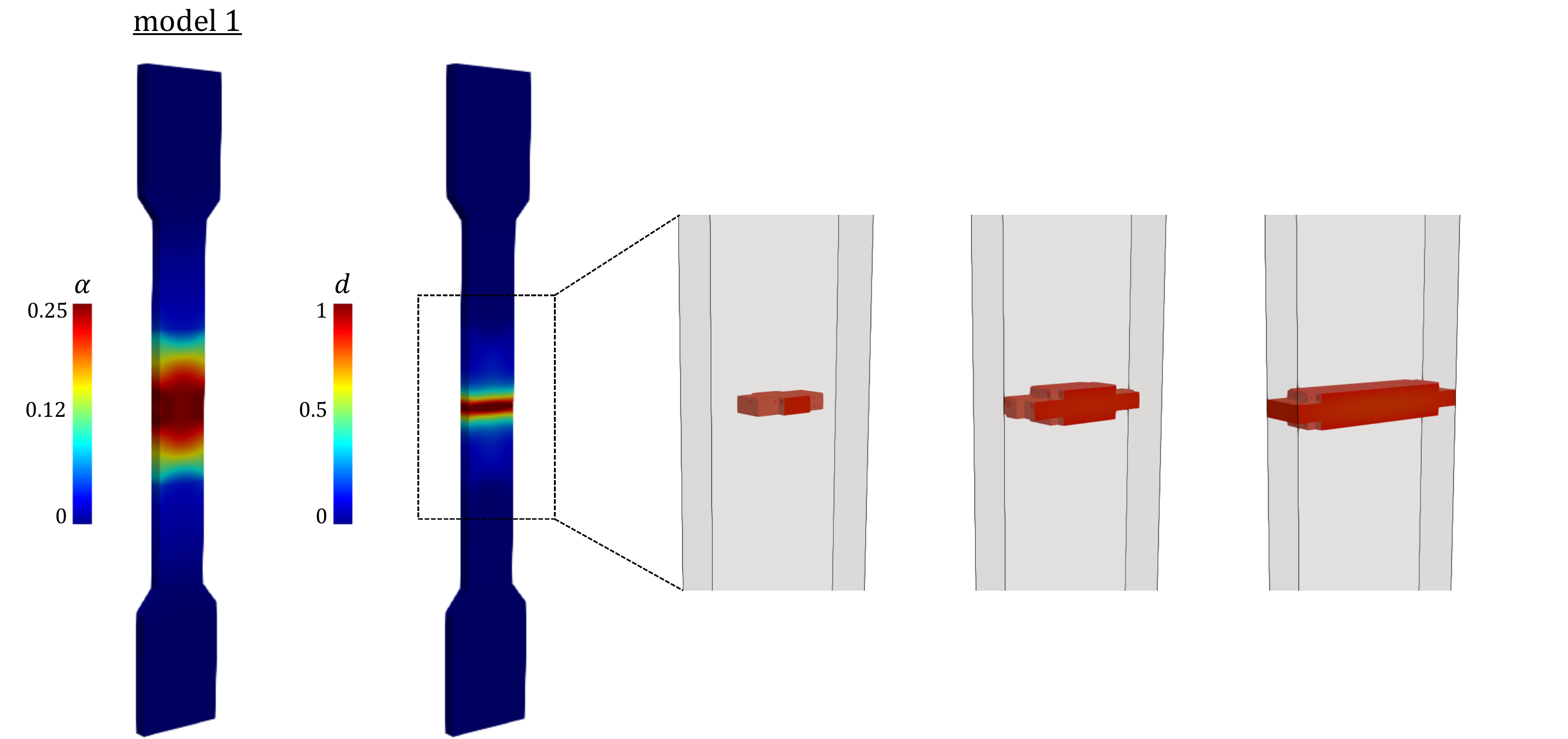}}}\\  
	\subfloat{{\includegraphics[clip,trim=2cm 21cm 8cm 4.5cm, width=15cm]{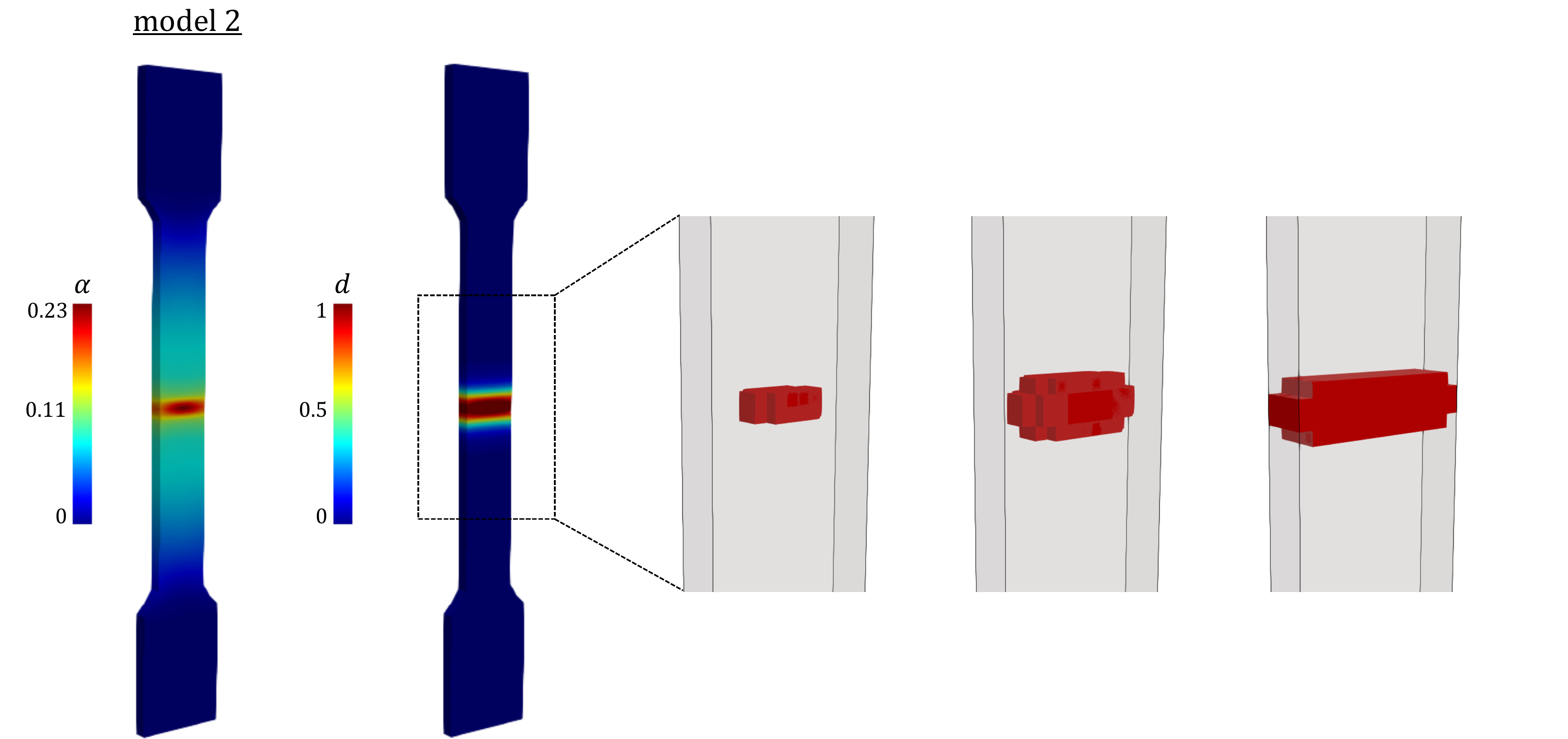}}}\\  
	\subfloat{{\includegraphics[clip,trim=2cm 19cm 8cm 4.5cm, width=15cm]{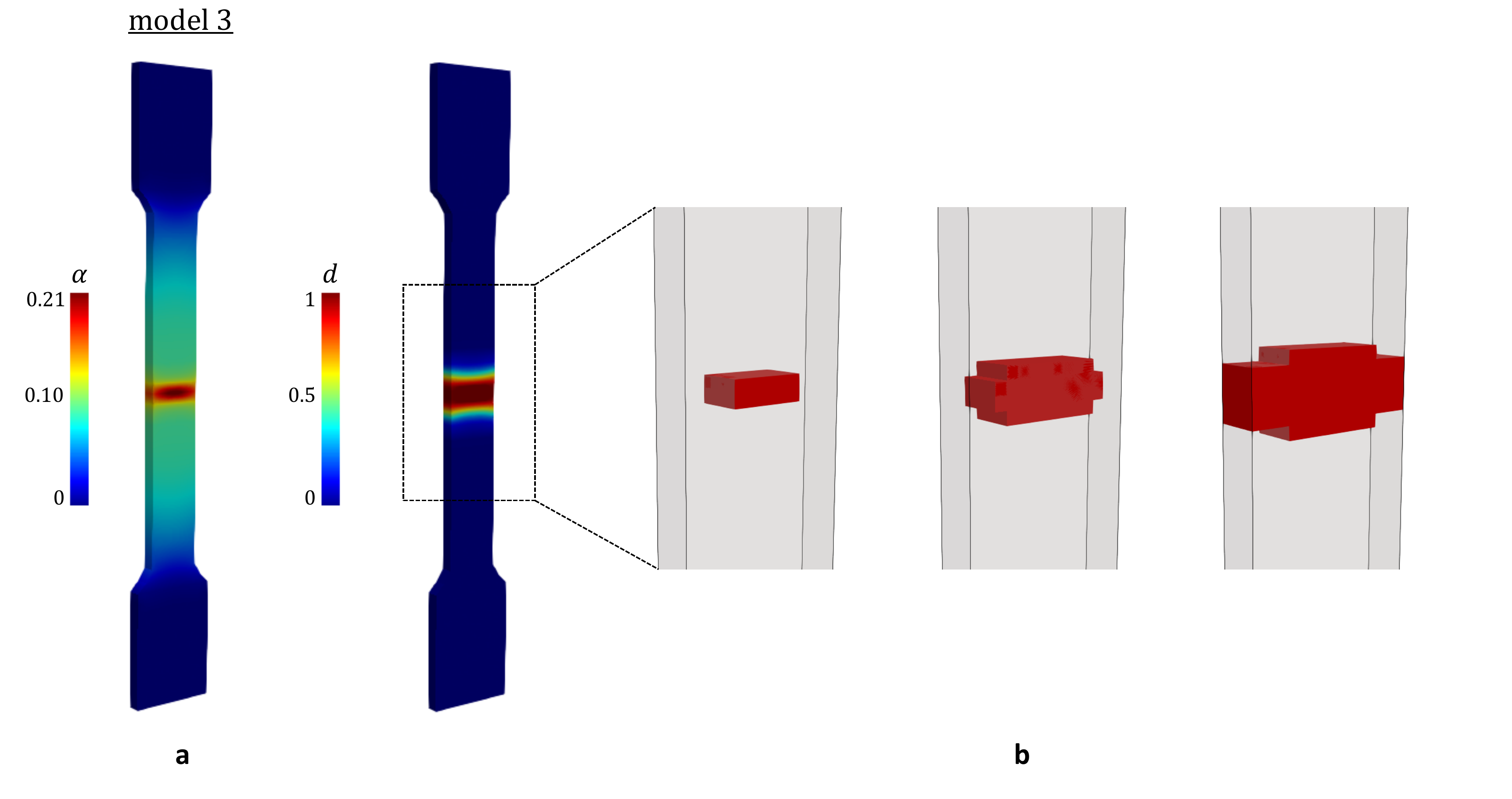}}}  
	\caption{ Example 3: (a) the equivalent plastic strain at complete failure, and (b) the evolution of the crack phase-field  for different deformation stages up to final failure at $\bar{u}_y=6$ mm.} 
	\label{Figure31}
\end{figure}

According to the above-mentioned discussion, we choose the EKF technique for the rest of the examples. We use five parallel MCMCs with $2\,000$ samples. For the observation, we use experimental data taken from \cite{ambati2016}. Figure \ref{Example3_models} shows the posterior density of estimated parameters, and the mean values are summarized in Table \ref{Example3_posterior}. We then employ the identified quantities in all three models. Figure \ref{Figure31} shows a comparison between the load-displacement curve obtained by the models and the experiments. Interestingly, by employing the Bayesian framework, all crack propagation stages (until fracture reaches the boundary) are modeled accurately.

%

Next, we investigate the ductile failure response employing the posterior density of the material parameters given in Table \ref{Example3_posterior}. The evolution of the crack phase-field $d$ is provided in Figure \ref{Figure31} at three deformation stages up to complete failure. Additionally, the equivalent plastic strain $\alpha$ at final failure is shown. It can be grasped that, regardless of the formulations, fracture initiation appears at the center of the specimen, and then evolves towards the two edges of the specimen until complete failure. The first important observation is that the simulation results are in well-agreement with the experimental failure pattern shown in Figure \ref{exm3bvp}c. Another important observation is that the crack phase-field in $\calM_2$ and $\calM_3$ is more diffuse than $\calM_1$. The main reason for this is that the crack driving force of $\calM_1$ is scaled by $\alpha$; thus, the phase-field diffusivity is strongly coupled to the ductile response; see Remark \ref{Cdriv_m1}. 

%
\sectpb[Section54]{Example 4:  Sandia fracture challenge}
%
The last example aims at estimating the posterior density of the unknown material properties for a specimen frequently used in the literature, namely, a Sandia fracture challenge \cite{boyce2014sandia}. The 2014 fracture challenge problem launched by the Sandia National Lab \cite{boyce2014sandia} has provided an ideal platform to assess the computational  capability and limitations of each participating team \cite{zhang2014modeling}. Specifically, this challenge aims to evaluate the computational ability to predict crack initiation and propagation of ductile fracture with respect to the experimental observation.
A recent comparative literature overview was conducted in \cite{DiLiWiTy21}.
 As reported in \cite{ambati2016}, standard phase-field formulations without estimating accurate ductile material properties quantitatively overestimate the post-yielding load-displacement response, which can be improved by performing an accurate calibration of the plasticity and phase-field parameters. Hence, we aim at reproducing the experimentally observed Sandia fracture challenge through the proposed Bayesian inversion framework.

\begin{figure}[!ht]
	\centering
	{\includegraphics[clip,trim=0cm 20cm 0cm 5cm, width=16cm]{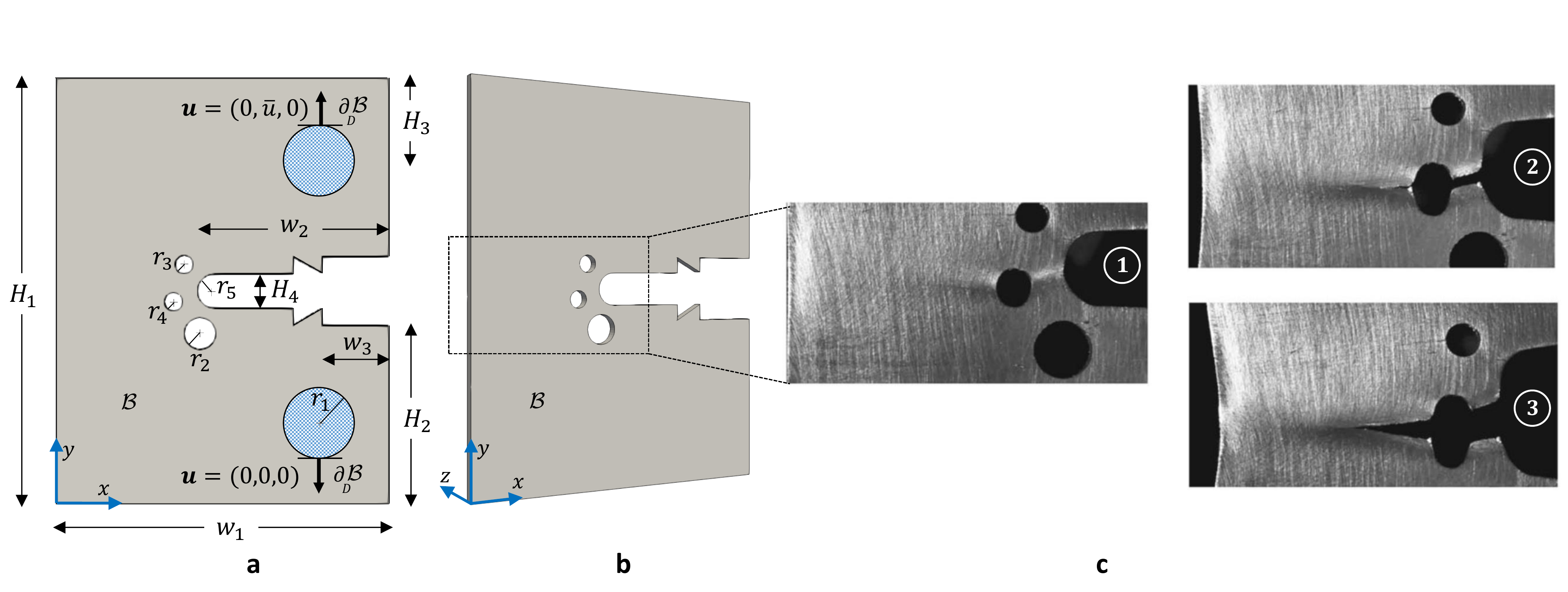}}  
	\caption{Example 4: Sandia fracture challenge. (a) Geometry of specimen with boundary conditions, (b) three-dimensional perspective, and (c) experimental observation taken from \cite{ambati2016}.}
	\label{exm4bvp}
\end{figure}

The experiments are based on the material Al-5052 H34, which experimentally induces a complex failure mode; see \cite{guo2013experimental, ambati2016}. The configuration is shown in Figure \ref{exm4bvp}a, while the experimental observations are shown in Figure \ref{exm4bvp}c.  The geometrical configuration includes two pins. The top pin is displaced vertically, while the lower pin is fixed in all directions. The two pins are considered to be rigid (here taken 10 times stiffer than the rest of the domain). The geometrical dimensions are set as $H_1=80$ mm, $H_2=35$ mm, $H_3=15$ mm, $H_4=6.5$ mm, $w_1=22$ mm, $w_2=36$ mm, and $w_3=12.5$ mm. The pins have an identical radius of $r_1=6$ mm, while the horizontal notch is rounded with a radius of $r_5=3.25$ mm. The specimen includes three voids with centers and radii $c_2=(x_2,y_2)=(27,32)$ mm and $r_2=3$ mm, $c_3=(x_3,y_3)=(24,45)$ mm and $r_3=1.75$ mm, and $c_4=(x_4,y_4)=(22,38)$ mm and $r_4=1.75$ mm, respectively. The specimen domain has a $2$ mm thickness, as shown in Figure~\ref{exm4bvp}b.

The numerical example is performed by applying a monotonic displacement increment  ${\Delta \bar{u}}_y=0.02$ mm in the vertical direction at the top pin for 400 time steps. The minimum finite element size is $0.32$ mm. Consequently, the Sandia specimen domain partition contains 17980 hexahedron linear elements.

\begin{figure}[b!]
 	\vspace{0.1cm}
 	\hspace{-1cm}
 	\subfloat{\includegraphics[width=0.37\textwidth]{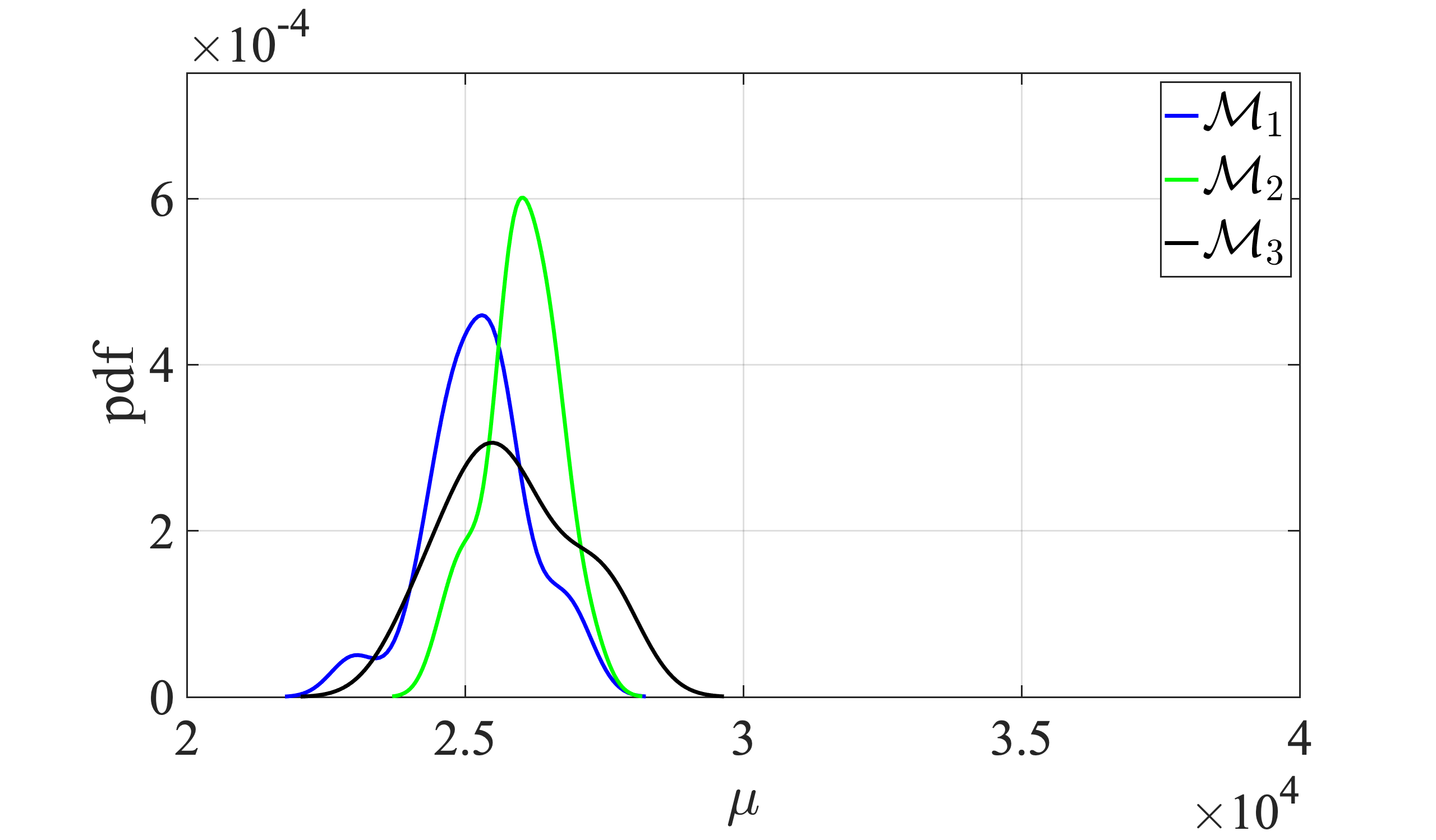}}   \subfloat{\includegraphics[width=0.37\textwidth]{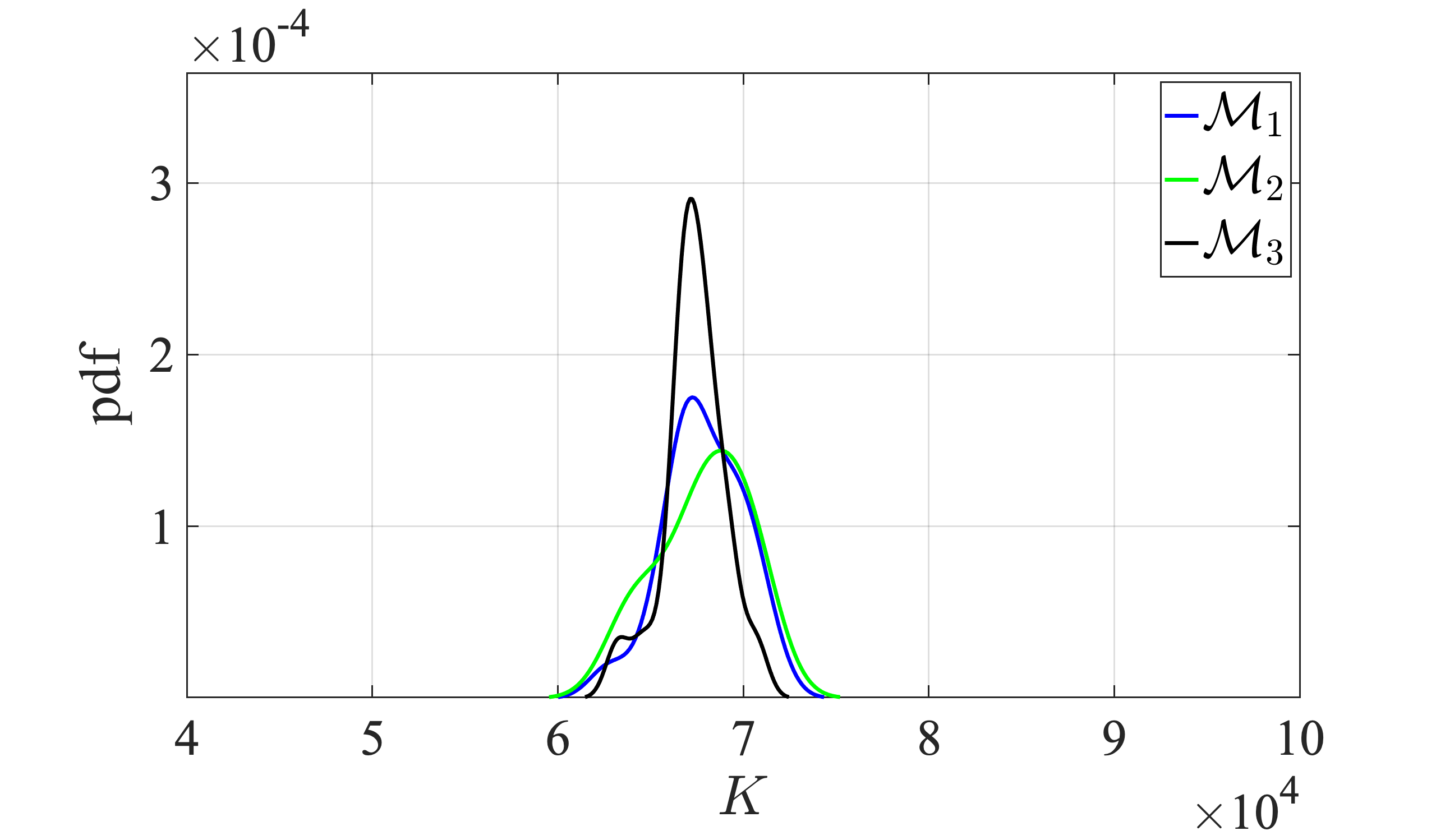}}  	\subfloat{\includegraphics[width=0.37\textwidth]{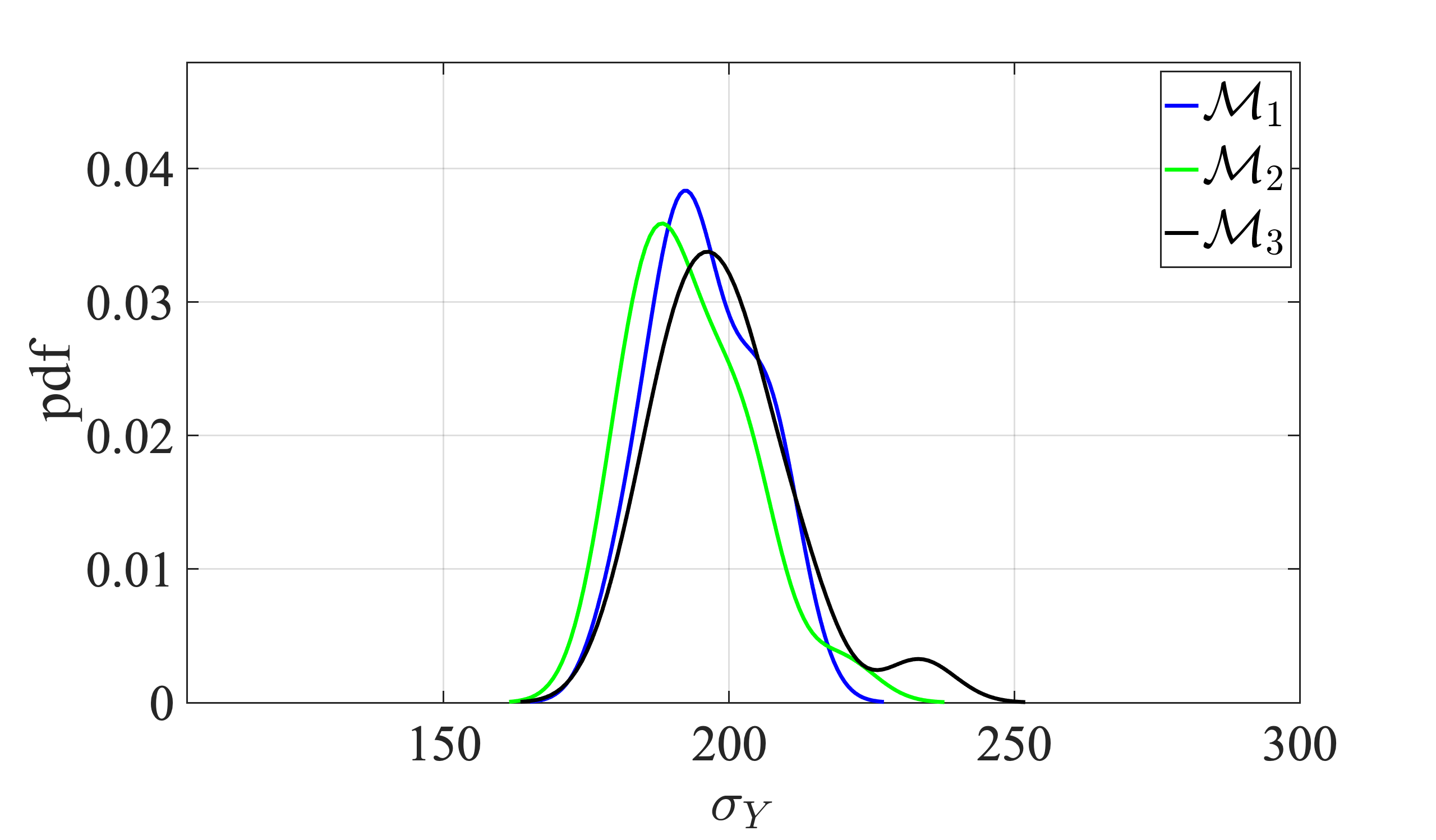}}
 	\hfill 
 	\vspace{0.1cm}
 	\hspace{-1.1cm}
 	\subfloat{\includegraphics[width=0.37\textwidth]{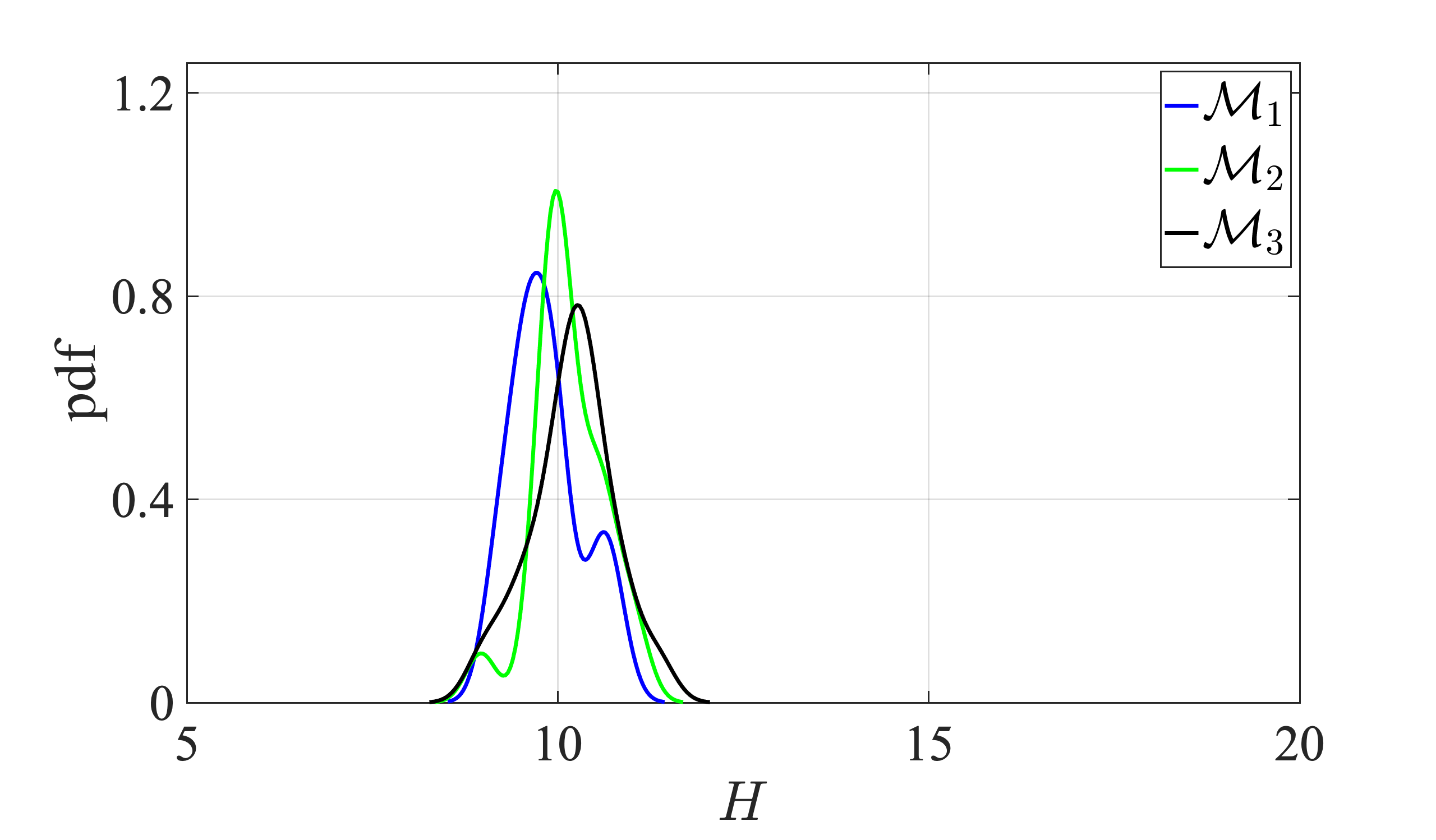}}   \subfloat{\includegraphics[width=0.37\textwidth]{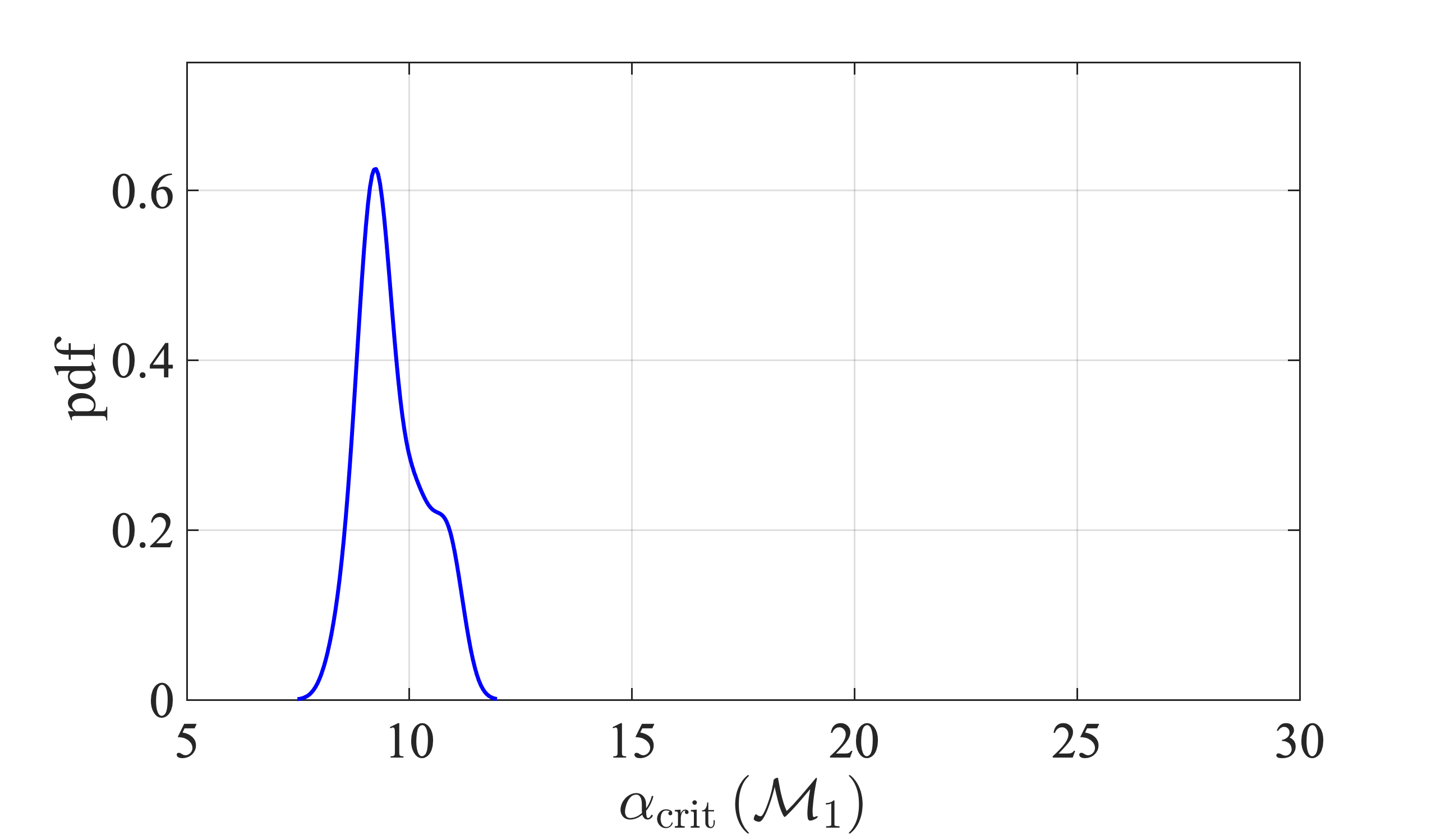}}  	\subfloat{\includegraphics[width=0.37\textwidth]{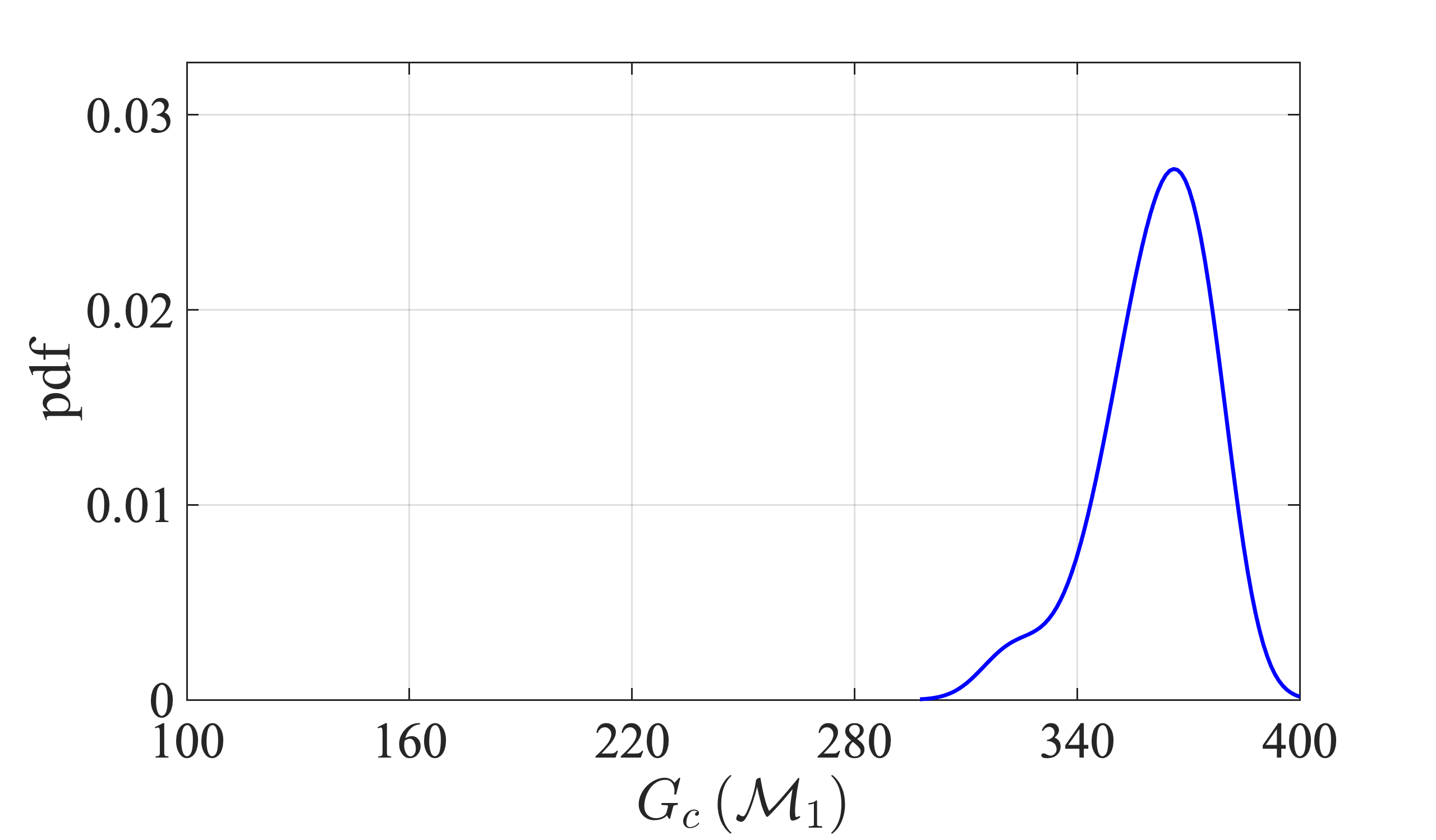}}
 	\hfill 
 	\vspace{0.5cm}
 	\hspace{-1.1cm}
 	\subfloat{\includegraphics[width=0.37\textwidth]{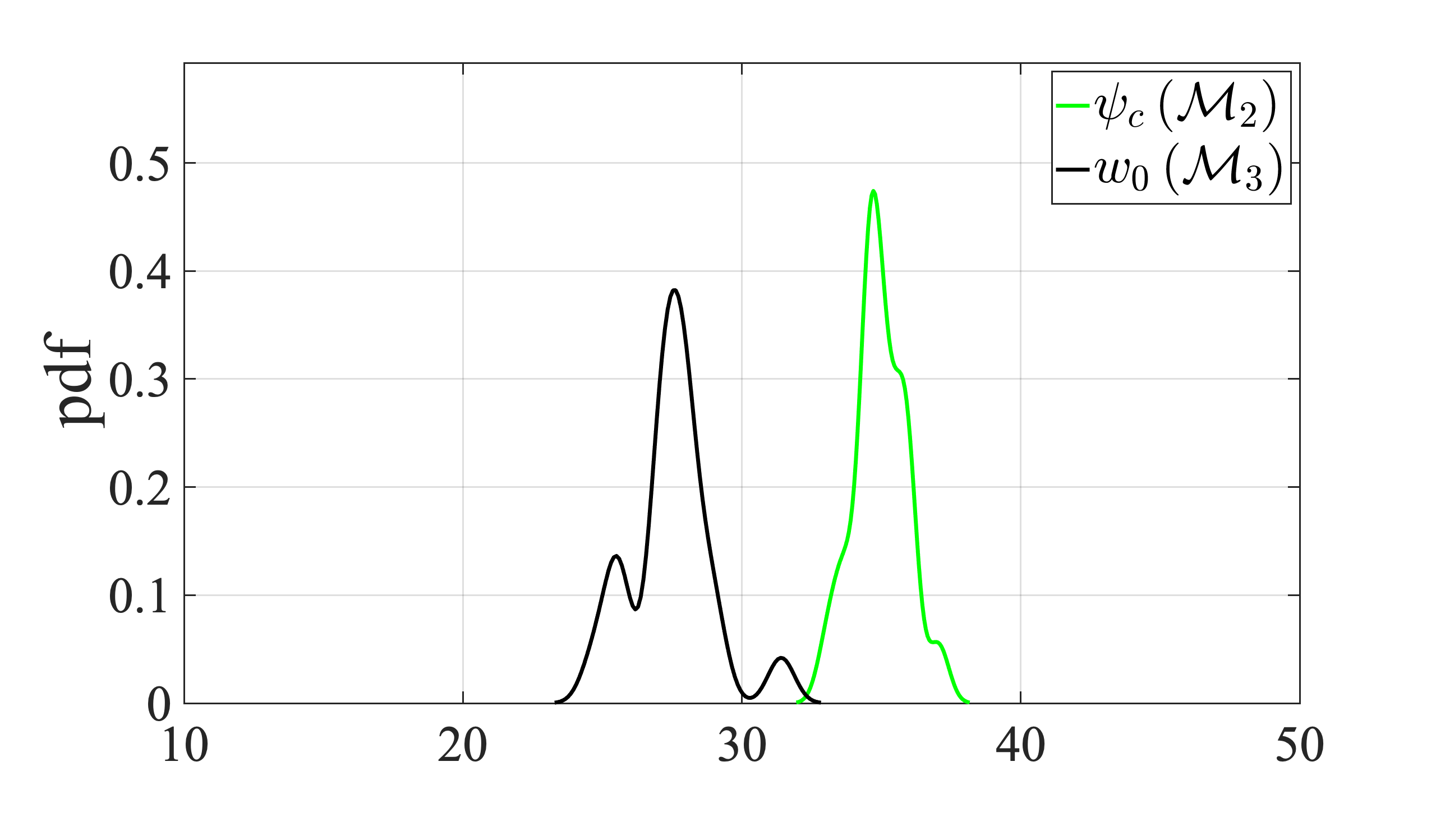}}   \subfloat{\includegraphics[width=0.37\textwidth]{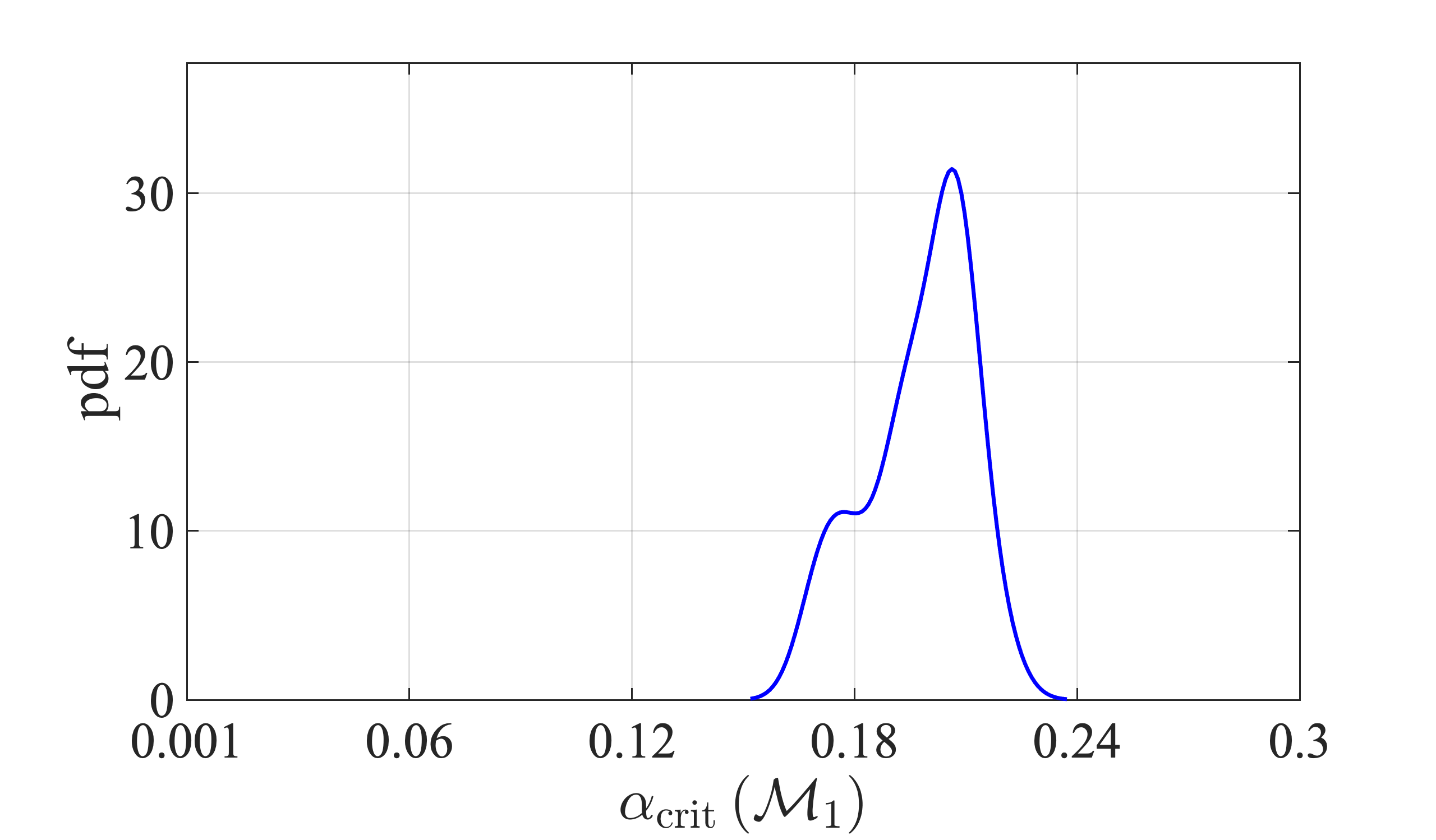}}  	\subfloat{\includegraphics[width=0.37\textwidth]{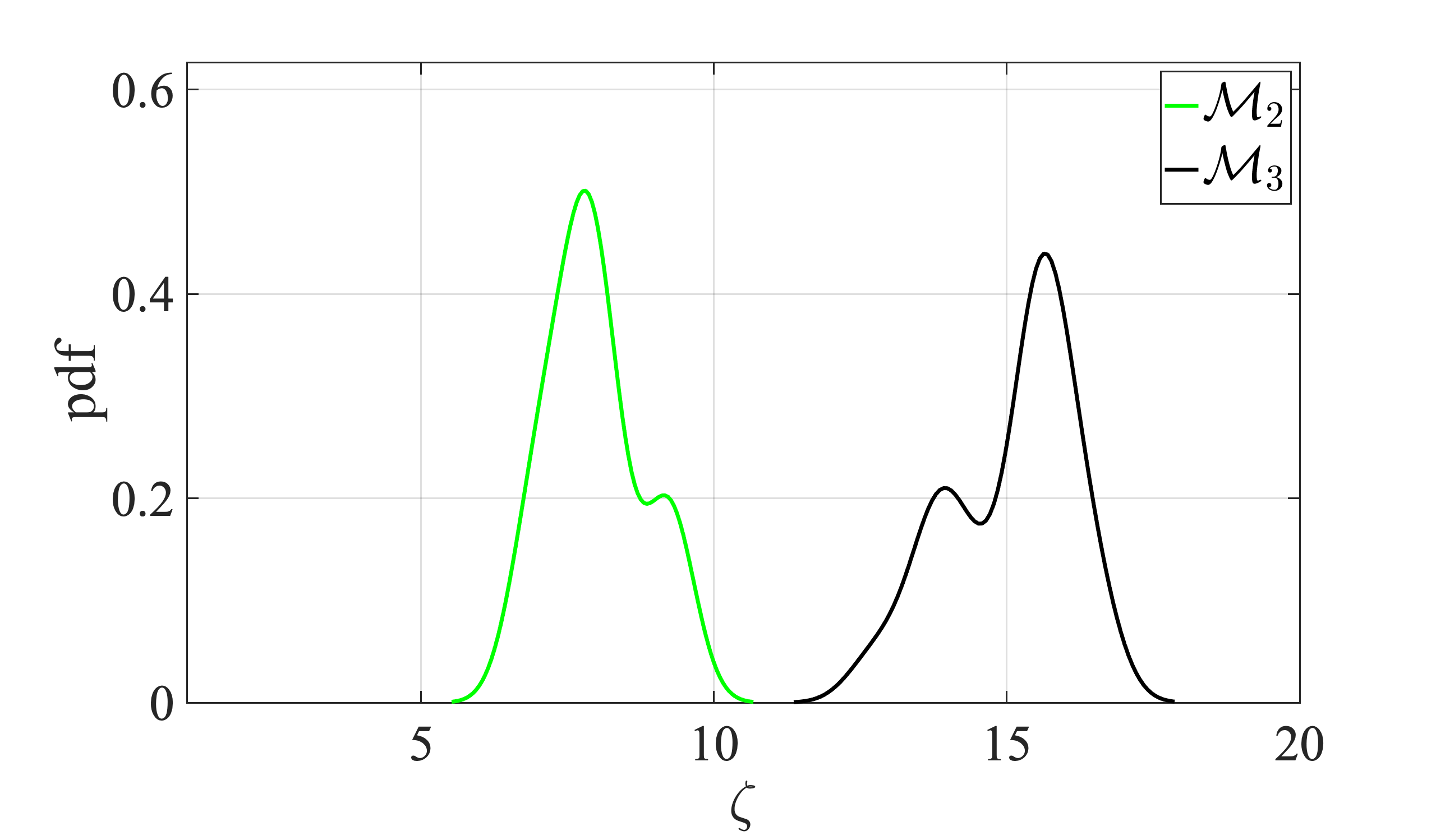}}
 	\caption{Example 4: the posterior distribution of the effective parameters using $\mathcal{M}_1$, $\mathcal{M}_2$, and $\mathcal{M}_3$.}
 	\label{Example4_models}
 \end{figure}

Due to its efficiency, the EKF technique is again considered for this example to identify the parameters. We use four parallel MCMCs with $2\,500$ samples. The experimental data (reference observation) is taken from \cite{ambati2016} and the prior densities are indicated in Table \ref{range3}. The posterior distributions are shown in Figure \ref{Example4_models}, while the mean values are summarized in Table \ref{Example4_posterior}. Finally, once again, we solve $\mathcal{M}_1$, $\mathcal{M}_2$, and $\mathcal{M}_3$ using the inferred values, and compare the results with those obtained in \cite{ambati2016}. As shown in Figure \ref{example34}, by employing the proposed Bayesian inversion framework, the parameters are estimated accurately, showing a very good agreement between the simulated data in all models and the experiments. This implementation also enhanced the computational capabilities of $\mathcal{M}_1$, showing an improved model accuracy compared to the results obtained in \cite{ambati2016}.
\begin{table}[!]
	\caption{Example 4: the uniform distribution of the inferred parameters in the Sandia example.}
	\vspace{1mm}
	\centering
	\begin{tabular}{llcccccccccccccc}
		&Parameter  &$H$        & $\mu$    &   $K$   &   $\sigma_Y$   &$G_c$&$\alpha_{\text{crit}}$     &$\psi_c$   &$w_0$ &$l_p$ &$\zeta$ & \\[2mm]\hline\\
		&   min   &5   & 20\,000  & 40\,000 &150 & 100 & 0.001 &10&10&5& 1\\[4mm]
		& max &   20  &  40\,000  &  100\,000 & 300 &400 &0.3 & 50 &50&30& 20\\[4mm]
		\hline
		\label{range3}
	\end{tabular}
\end{table}
\begin{table}[!]
	\caption{Example 4: the mean value of posterior density of  the model parameters for the three models.}
	\vspace{1mm}
	\centering
	\begin{tabular}{ccccccccccccccc}
		&Model  &$H$        & $\mu$    &   $K$   &   $\sigma_Y$   &$G_c$&$\alpha_{\text{crit}}$     &$\psi_c$   &$w_0$ &$l_p$&$\zeta$ & \\[2mm]\hline\\
		&    $\mathcal{M}_1$   &9.8   & 25\,500 & 67\,500 &195 & 360 & 0.2  &--&--&--&--  \\[4mm]
		& $\mathcal{M}_2$ &   10  & 26\,200  &  67\,200 & 194 &-- &-- & 35& &--&8\\[4mm]
		&   $\mathcal{M}_3$  &10.1  & 25\,900   & 67\,800 &195 & -- & -- &--& 28 & 9.64&15  \\[4mm]
		\hline
		\label{Example4_posterior}
	\end{tabular}
\end{table}

\begin{figure}[!t]
	\centering
	\subfloat{{\includegraphics[clip,trim=0cm 23cm 0cm 4.5cm, width=16cm]{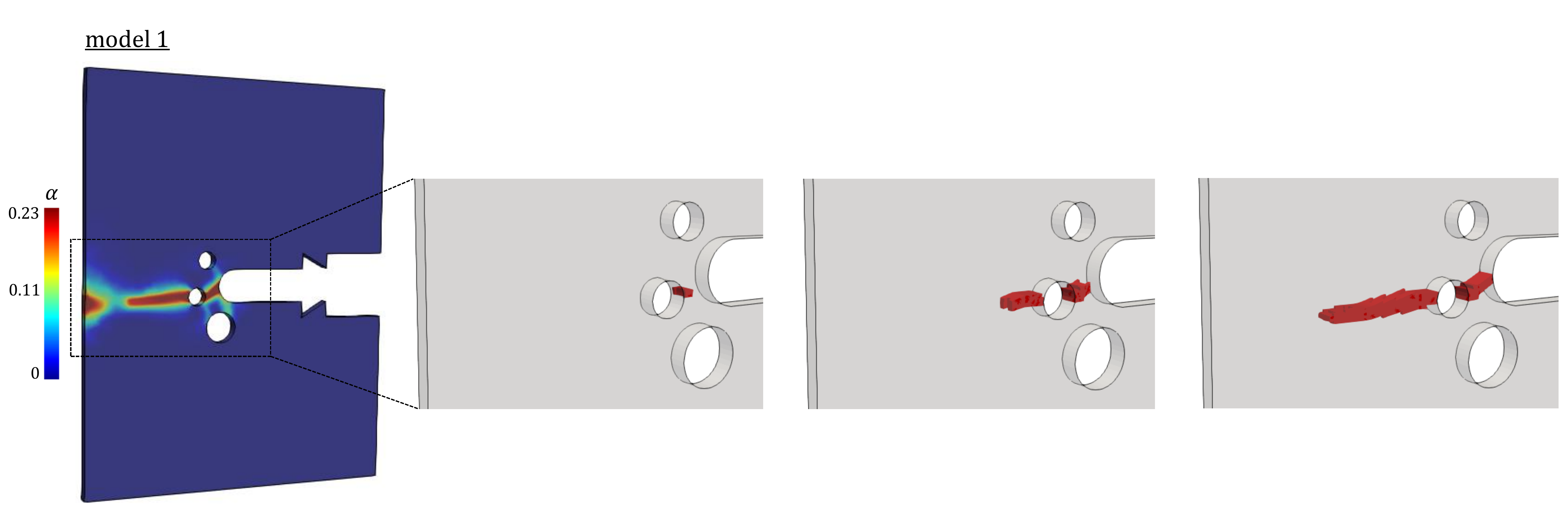}}}\\  
	\subfloat{{\includegraphics[clip,trim=0cm 23cm 0cm 4.5cm, width=16cm]{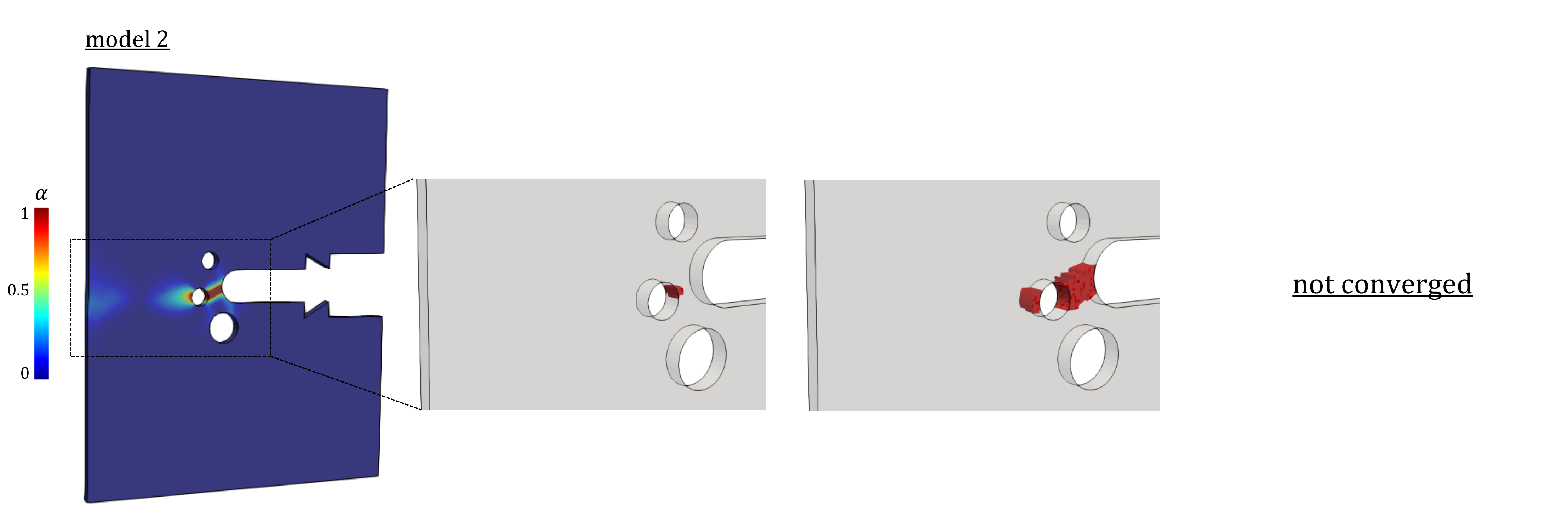}}}\\  
	\subfloat{{\includegraphics[clip,trim=0cm 21cm 0cm 4.5cm, width=16cm]{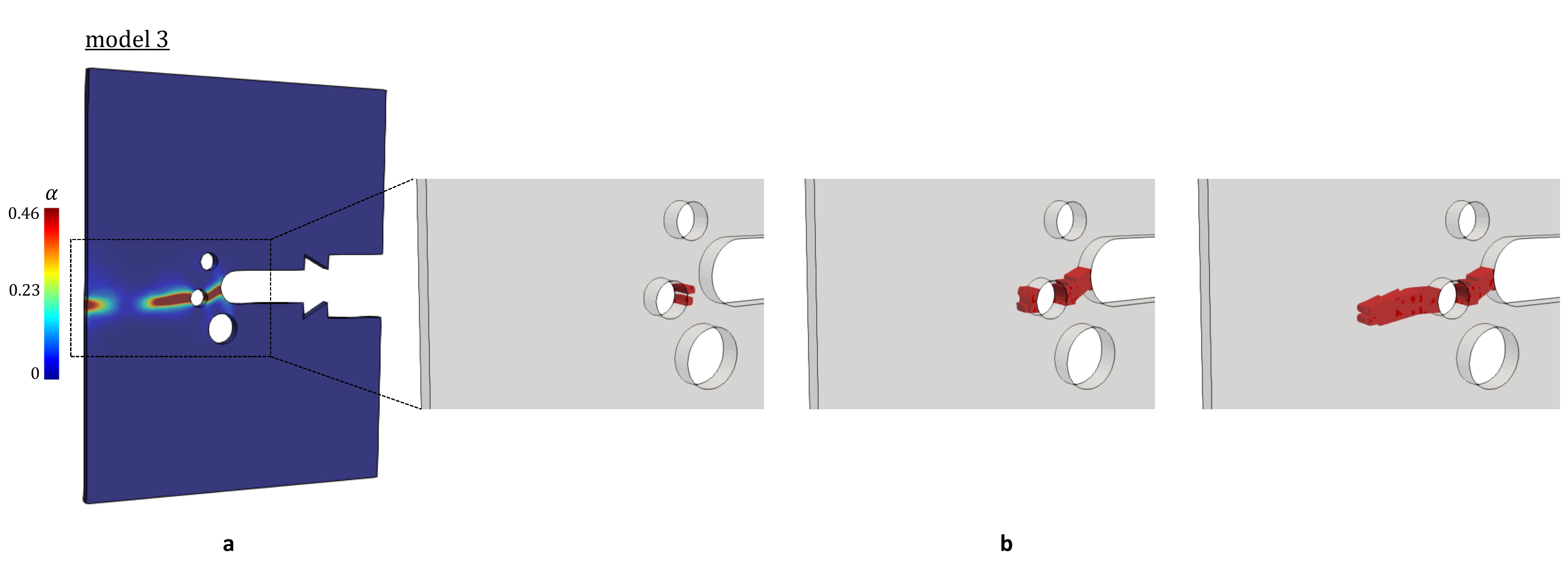}}}  
	\caption{ Example 4: (a) The equivalent plastic strain at complete failure, and (b) the evolution of the crack phase-field for different deformation stages up to final failure at $\bar{u}_x=8$ mm. The final fracture stage was not captured with model 2 due to lack of convergence.}
	\label{Figure41}
\end{figure}

For a better insight into the fracture process in all models, the evolution of the crack phase-field $d$ is provided in Figure \ref{Figure41} at three deformation stages up to complete failure. Additionally, the equivalent plastic strain $\alpha$ at final failure is shown. Note that the solutions are based on the posterior density of the material parameters, which are given in Table \ref{Example4_posterior}.

It can be grasped that, for all three models, the fracture path first initiates at the void located in the middle of the specimen, and afterwards, evolves towards the right edge of the domain. In addition, a secondary crack initiates from the central void, but this time from its left side, and then propagates towards the left edge of the specimen until complete failure. Note further that the simulation results are in well-agreement with the experimental failure pattern shown in Figure \ref{exm4bvp}c. 

Based on our numerical result, it is worth noting that the stage of the secondary crack was no longer predicted by model 2. To estimate the posterior density of the material parameters through Bayesian inversion, several candidates are required. Thus, a stable forward method is crucial. Otherwise, deviating material properties will result in unstable solutions. In this context, we  highlight that model 3 provides the most stable solutions, at the cost of an additional PDE that must be solved to obtain the plastic response (as opposed to local plasticity in model 1 and model 2).

\sectpa[Section7]{Conclusion}

In this work, we have proposed a robust and efficient step-wise Bayesian inversion method for ductile fracture problems using phase-field models. In particular, a Bayesian inversion framework (as a probabilistic technique) based on MCMC is developed to identify unknown ductile fracture parameters. Three common MCMC methods, namely the MH algorithm, DRAM algorithm, and EKF-MCMC have been used to estimate the effective parameters in ductile fracture. The posterior density results from the inverse problem are evaluated with synthetic measurements (for the first two examples) as well as experimental observations (for the last two examples). To approximate ductile failure, a phase-field fracture formulation is used for a ductile material exhibiting $J_2$-plasticity in a quasi-static kinematically linear regime. To do so, we have presented a unified formulation for phase-field modeling of ductile fracture, which is resolved through an incremental energy minimization approach. The overall formulation is revisited and extended to the case of anisotropic ductile fracture. Three different specific models are subsequently recovered by certain choices of parameters and constitutive functions.

In the first numerical example, the equivalence of the parameters in the different models is provided. In the second example (anisotropic ductile fracture), we investigate the evolution of the failure response through the posterior density function obtained by the proposed step-wise Bayesian inversion method. We have shown that the equivalent plastic strain $\alpha$ as well as the crack phase-field $d$ evolve in the direction of the preferred fiber orientation. The last two examples are concerned with the experimental observations to estimate the posterior density of the material unknowns. We observed that, although the MH algorithm can be implemented easily, it is sensitive to the initial guess, leading to slow convergence, and may depend on the prior density. As more advanced techniques, we compared the convergence of the DRAM and EKF-MCMC methods by employing a reliable convergence diagnostic tool, namely $\hat{R}$-convergence. Using the Kalman filter improves considerably the convergence of different MCMCs, i.e., fewer iterations are needed to obtain a high level of accuracy. We conclude that this method is more efficient compared to the DRAM algorithm.

Through our findings, the coupling scheme between the step-wise Bayesian inversion framework and ductile fracture simulations results in an accurate and a reliable information related to the model parameters. As a consequence, an excellent agreement was obtained between the results of all examined models and the experimental observations.

\subsection*{Acknowledgment}
N. Noii  and F. Aldakheel were founded by the Priority Program \texttt{DFG-SPP 2020} within its second funding phase. T. Wick and P. Wriggers were funded by the Deutsche Forschungsgemeinschaft (DFG, German Research Foundation) under Germany's Excellence Strategy within the Cluster of Excellence PhoenixD, \texttt{EXC 2122} (project number: 390833453). 

 
 {\normalsize
 	\begin{spacing}{0.8}
 \bibliographystyle{ieeetr}
\bibliography{./newlit}

\begin{thebibliography}{100}

\bibitem{miehe2011}
C.~Miehe, ``A multi-field incremental variational framework for
  gradient-extended standard dissipative solids,'' {\em Journal of the
  Mechanics and Physics of Solids}, vol.~59, no.~4, pp.~898--923, 2011.

\bibitem{miehe2013mixed}
C.~Miehe, F.~Aldakheel, and S.~Mauthe, ``Mixed variational principles and
  robust finite element implementations of gradient plasticity at small
  strains,'' {\em International Journal for Numerical Methods in Engineering},
  vol.~94, no.~11, pp.~1037--1074, 2013.

\bibitem{peerlings1998gradient}
R.~H. Peerlings, R.~de~Borst, W.~Brekelmans, and M.~G. Geers,
  ``Gradient-enhanced damage modelling of concrete fracture,'' {\em Mechanics
  of Cohesive-frictional Materials: An International Journal on Experiments,
  Modelling and Computation of Materials and Structures}, vol.~3, no.~4,
  pp.~323--342, 1998.

\bibitem{kiefer2018gradient}
B.~Kiefer, T.~Waffenschmidt, L.~Sprave, and A.~Menzel, ``A gradient-enhanced
  damage model coupled to plasticity-multi-surface formulation and algorithmic
  concepts,'' {\em International Journal of Damage Mechanics}, vol.~27, no.~2,
  pp.~253--295, 2018.

\bibitem{junker2021efficient}
P.~Junker, J.~Riesselmann, and D.~Balzani, ``Efficient and robust numerical
  treatment of a gradient-enhanced damage model at large deformations,'' {\em
  arXiv preprint arXiv:2102.08819}, 2021.

\bibitem{barfusz2021single}
O.~Barfusz, T.~Brepols, T.~van~der Velden, J.~Frischkorn, and S.~Reese, ``A
  single gauss point continuum finite element formulation for gradient-extended
  damage at large deformations,'' {\em Computer Methods in Applied Mechanics
  and Engineering}, vol.~373, p.~113440, 2021.

\bibitem{miehe+welschinger+hofacker10a}
C.~Miehe, F.~Welschinger, and M.~Hofacker, ``Thermodynamically consistent
  phase-field models of fracture: Variational principles and multi-field fe
  implementations,'' {\em International Journal for Numerical Methods in
  Engineering}, vol.~83, pp.~1273--1311, 2010.

\bibitem{BourFraMar08}
B.~Bourdin, G.~Francfort, and J.-J. Marigo, ``The variational approach to
  fracture,'' {\em Journal of Elasticity}, vol.~91, pp.~5--148, 2008.

\bibitem{kuhn2010continuum}
C.~Kuhn and R.~M{\"u}ller, ``A continuum phase field model for fracture,'' {\em
  Engineering Fracture Mechanics}, vol.~77, no.~18, pp.~3625--3634, 2010.

\bibitem{de1996some}
R.~De~Borst and J.~Pamin, ``Some novel developments in finite element
  procedures for gradient-dependent plasticity,'' {\em International Journal
  for Numerical Methods in Engineering}, vol.~39, no.~14, pp.~2477--2505, 1996.

\bibitem{polizzotto1998thermodynamics}
C.~Polizzotto and G.~Borino, ``A thermodynamics-based formulation of
  gradient-dependent plasticity,'' {\em European Journal of
  Mechanics-A/Solids}, vol.~17, no.~5, pp.~741--761, 1998.

\bibitem{liebe2001theory}
T.~Liebe and P.~Steinmann, ``Theory and numerics of a thermodynamically
  consistent framework for geometrically linear gradient plasticity,'' {\em
  International Journal for Numerical Methods in Engineering}, vol.~51, no.~12,
  pp.~1437--1467, 2001.

\bibitem{miehe2014variational}
C.~Miehe, F.~Welschinger, and F.~Aldakheel, ``Variational gradient plasticity
  at finite strains. part ii: Local--global updates and mixed finite elements
  for additive plasticity in the logarithmic strain space,'' {\em Computer
  Methods in Applied Mechanics and Engineering}, vol.~268, pp.~704--734, 2014.

\bibitem{aldakheel16}
F.~Aldakheel, {\em Mechanics of Nonlocal Dissipative Solids: Gradient
  Plasticity and Phase Field Modeling of Ductile Fracture}.
\newblock PhD thesis, Institute of Applied Mechanics (CE), Chair I, University
  of Stuttgart, 2016.
\newblock http://dx.doi.org/10.18419/opus-8803.

\bibitem{NoiiGL18}
T.~Gerasimov, N.~Noii, O.~Allix, and L.~De~Lorenzis, ``A non-intrusive
  global/local approach applied to phase-field modeling of brittle fracture,''
  {\em Advanced Modeling and Simulation in Engineering Sciences}, 2018.
\newblock https://doi.org/10.1186/s40323-018-0105-8.

\bibitem{KUMAR2020104027}
A.~Kumar, B.~Bourdin, G.~A. Francfort, and O.~Lopez-Pamies, ``Revisiting
  nucleation in the phase-field approach to brittle fracture,'' {\em Journal of
  the Mechanics and Physics of Solids}, p.~104027, 2020.

\bibitem{TaTiBouMaMau17}
E.~Tann{\'e}, T.~Li, B.~Bourdin, J.-J. Marigo, and C.~Maurini, ``Crack
  nucleation in variational phase-field models of brittle fracture,'' {\em
  Journal of the Mechanics and Physics of Solids}, vol.~110, pp.~80--99, 2018.

\bibitem{GoeNov10}
N.~Van~Goethem and A.~Novotny, ``Crack nucleation sensitivity analysis,'' {\em
  Mathematical Methods in the Applied Sciences}, vol.~33, no.~16,
  pp.~1978--1994, 2010.

\bibitem{NoiiWick2019}
N.~Noii and T.~Wick, ``A phase-field description for pressurized and
  non-isothermal propagating fractures,'' {\em Computer Methods in Applied
  Mechanics and Engineering}, vol.~351, pp.~860 -- 890, 2019.

\bibitem{Wi20_book}
T.~Wick, {\em Multiphysics Phase-Field Fracture: Modeling, Adaptive
  Discretizations, and Solvers}.
\newblock Berlin, Boston: De Gruyter, 2020.

\bibitem{alessi2014}
R.~Alessi, J.~Marigo, and S.~Vidoli, ``Gradient damage models coupled with
  plasticity and nucleation of cohesive cracks,'' {\em Archive for Rational
  Mechanics and Analysis}, vol.~214, no.~2, pp.~575--615, 2014.

\bibitem{du2015}
F.~Duda, A.~Ciarbonetti, P.~S{\'a}nchez, and A.~Huespe, ``A
  phase-field/gradient damage model for brittle fracture in elastic--plastic
  solids,'' {\em International Journal of Plasticity}, vol.~65, pp.~269--296,
  2015.

\bibitem{ambati2015}
M.~Ambati, T.~Gerasimov, and L.~De~Lorenzis, ``Phase-field modeling of ductile
  fracture,'' {\em Computational Mechanics}, vol.~55, no.~5, pp.~1017--1040,
  2015.

\bibitem{alessi2017}
R.~Alessi, J.~Marigo, C.~Maurini, and S.~Vidoli, ``Coupling damage and
  plasticity for a phase-field regularisation of brittle, cohesive and ductile
  fracture: one-dimensional examples,'' {\em International Journal of
  Mechanical Sciences}, vol.~149, pp.~559--576, 2018.

\bibitem{borden2016}
M.~Borden, T.~Hughes, C.~Landis, A.~Anvari, and I.~Lee, ``A phase-field
  formulation for fracture in ductile materials: Finite deformation balance law
  derivation, plastic degradation, and stress triaxiality effects,'' {\em
  Computer Methods in Applied Mechanics and Engineering}, vol.~312,
  pp.~130--166, 2016.

\bibitem{kuhn2016}
C.~Kuhn, T.~Noll, and R.~M{\"u}ller, ``On phase field modeling of ductile
  fracture,'' {\em GAMM-Mitteilungen}, vol.~39, no.~1, pp.~35--54, 2016.

\bibitem{ulloa2016}
J.~Ulloa, P.~Rodr{\'\i}guez, and E.~Samaniego, ``On the modeling of dissipative
  mechanisms in a ductile softening bar,'' {\em Journal of Mechanics of
  Materials and Structures}, vol.~11, no.~4, pp.~463--490, 2016.

\bibitem{alessi2018comparison}
R.~Alessi, M.~Ambati, T.~Gerasimov, S.~Vidoli, and L.~De~Lorenzis, ``Comparison
  of phase-field models of fracture coupled with plasticity,'' in {\em Advances
  in Computational Plasticity}, pp.~1--21, Springer, 2018.

\bibitem{choo2018}
J.~Choo and W.~Sun, ``Coupled phase-field and plasticity modeling of geological
  materials: From brittle fracture to ductile flow,'' {\em Computer Methods in
  Applied Mechanics and Engineering}, vol.~330, pp.~1--32, 2018.

\bibitem{kienle2019}
D.~Kienle, F.~Aldakheel, and M.-A. Keip, ``A finite-strain phase-field approach
  to ductile failure of frictional materials,'' {\em International Journal of
  Solids and Structures}, vol.~172, pp.~147--162, 2019.

\bibitem{aldakheel+wriggers+miehe18}
F.~Aldakheel, P.~Wriggers, and C.~Miehe, ``{A modified Gurson-type plasticity
  model at finite strains: Formulation, numerical analysis and phase-field
  coupling},'' {\em Computational Mechanics}, vol.~62, pp.~815--833, 2018.

\bibitem{dittmann2020}
M.~Dittmann, F.~Aldakheel, J.~Schulte, F.~Schmidt, M.~Kr{\"u}ger, P.~Wriggers,
  and C.~Hesch, ``Phase-field modeling of porous-ductile fracture in non-linear
  thermo-elasto-plastic solids,'' {\em Computer Methods in Applied Mechanics
  and Engineering}, vol.~361, p.~112730, 2020.

\bibitem{aldakheel2019virtual}
F.~Aldakheel, B.~Hudobivnik, and P.~Wriggers, ``Virtual element formulation for
  phase-field modeling of ductile fracture,'' {\em International Journal for
  Multiscale Computational Engineering}, vol.~17, no.~2, 2019.

\bibitem{storm2021comparative}
J.~Storm, M.~Pise, D.~Brands, J.~Schr{\"o}der, and M.~Kaliske, ``A comparative
  study of micro-mechanical models for fiber pullout behavior of reinforced
  high performance concrete,'' {\em Engineering Fracture Mechanics}, vol.~243,
  p.~107506, 2021.

\bibitem{aldakheel2020microscale}
F.~Aldakheel, ``A microscale model for concrete failure in poro-elasto-plastic
  media,'' {\em Theoretical and Applied Fracture Mechanics}, vol.~107,
  p.~102517, 2020.

\bibitem{heider2020phase}
Y.~Heider and W.~Sun, ``A phase field framework for capillary-induced fracture
  in unsaturated porous media: Drying-induced vs. hydraulic cracking,'' {\em
  Computer Methods in Applied Mechanics and Engineering}, vol.~359, p.~112647,
  2020.

\bibitem{aldakheel2020global}
F.~Aldakheel, N.~Noii, T.~Wick, and P.~Wriggers, ``A global--local approach for
  hydraulic phase-field fracture in poroelastic media,'' {\em Computers \&
  Mathematics with Applications}, 2020.

\bibitem{yin2020ductile}
B.~Yin and M.~Kaliske, ``A ductile phase-field model based on degrading the
  fracture toughness: Theory and implementation at small strain,'' {\em
  Computer Methods in Applied Mechanics and Engineering}, vol.~366, p.~113068,
  2020.

\bibitem{fang2019}
J.~Fang, C.~Wu, J.~Li, Q.~Liu, C.~Wu, G.~Sun, and L.~Qing, ``Phase field
  fracture in elasto-plastic solids: variational formulation for multi-surface
  plasticity and effects of plastic yield surfaces and hardening,'' {\em
  International Journal of Mechanical Sciences}, vol.~156, pp.~382--396, 2019.

\bibitem{ulloa2021}
J.~Ulloa, J.~Wambacq, R.~Alessi, G.~Degrande, and S.~Fran\c{c}ois,
  ``Phase-field modeling of fatigue coupled to cyclic plasticity in an
  energetic formulation,'' {\em Computer Methods in Applied Mechanics and
  Engineering}, vol.~373, p.~113473, 2021.

\bibitem{alessi2015}
R.~Alessi, J.~Marigo, and S.~Vidoli, ``Gradient damage models coupled with
  plasticity: variational formulation and main properties,'' {\em Mechanics of
  Materials}, vol.~80, pp.~351--367, 2015.

\bibitem{tanne2017}
E.~Tanne, {\em Variational phase-field models from brittle to ductile fracture:
  nucleation and propagation}.
\newblock PhD thesis, Universit{\'e} Paris-Saclay (ComUE), 2017.

\bibitem{miehe2017phase}
C.~Miehe, F.~Aldakheel, and S.~Teichtmeister, ``Phase-field modeling of ductile
  fracture at finite strains: A robust variational-based numerical
  implementation of a gradient-extended theory by micromorphic
  regularization,'' {\em International Journal for Numerical Methods in
  Engineering}, vol.~111, no.~9, pp.~816--863, 2017.

\bibitem{muhlhaus1991variational}
H.-B. M{\"u}hlhaus and E.~Alfantis, ``A variational principle for gradient
  plasticity,'' {\em International Journal of Solids and Structures}, vol.~28,
  no.~7, pp.~845--857, 1991.

\bibitem{rodriguez2018}
P.~Rodriguez, J.~Ulloa, C.~Samaniego, and E.~Samaniego, ``A variational
  approach to the phase field modeling of brittle and ductile fracture,'' {\em
  International Journal of Mechanical Sciences}, vol.~144, pp.~502--517, 2018.

\bibitem{dittmann2018}
M.~Dittmann, F.~Aldakheel, J.~Schulte, P.~Wriggers, and C.~Hesch, ``Variational
  phase-field formulation of non-linear ductile fracture,'' {\em Computer
  Methods in Applied Mechanics and Engineering}, vol.~342, pp.~71--94, 2018.

\bibitem{miehe2016bphase}
C.~Miehe, F.~Aldakheel, and A.~Raina, ``Phase field modeling of ductile
  fracture at finite strains: A variational gradient-extended plasticity-damage
  theory,'' {\em International Journal of Plasticity}, vol.~84, pp.~1--32,
  2016.

\bibitem{miehe+hofacker+schaenzel+aldakheel15}
C.~Miehe, M.~Hofacker, L.-M. Sch\"anzel, and F.~Aldakheel, ``Phase field
  modeling of fracture in multi-physics problems. {P}art {II}.
  brittle-to-ductile failure mode transition and crack propagation in
  thermo-elastic-plastic solids,'' {\em Computer Methods in Applied Mechanics
  and Engineering}, vol.~294, pp.~486--522, 2015.

\bibitem{noii2020adaptive}
N.~Noii, F.~Aldakheel, T.~Wick, and P.~Wriggers, ``An adaptive global--local
  approach for phase-field modeling of anisotropic brittle fracture,'' {\em
  Computer Methods in Applied Mechanics and Engineering}, vol.~361, p.~112744,
  2020.

\bibitem{smith2013uncertainty}
R.~C. Smith, {\em Uncertainty quantification: theory, implementation, and
  applications}, vol.~12.
\newblock SIAM, 2013.

\bibitem{khodadadian2020bayesian}
A.~Khodadadian, N.~Noii, M.~Parvizi, M.~Abbaszadeh, T.~Wick, and C.~Heitzinger,
  ``{A Bayesian estimation method for variational phase-field fracture
  problems},'' {\em Computational Mechanics}, vol.~66, pp.~827--849, 2020.

\bibitem{haario2006dram}
H.~Haario, M.~Laine, A.~Mira, and E.~Saksman, ``{DRAM: efficient adaptive
  MCMC},'' {\em Statistics and Computing}, vol.~16, no.~4, pp.~339--354, 2006.

\bibitem{laloy2012high}
E.~Laloy and J.~A. Vrugt, ``{High-dimensional posterior exploration of
  hydrologic models using multiple-try DREAM (ZS) and high-performance
  computing},'' {\em Water Resources Research}, vol.~48, no.~1, 2012.

\bibitem{emerick2012combining}
A.~A. Emerick, A.~C. Reynolds, {\em et~al.}, ``{Combining the ensemble Kalman
  filter with Markov-chain Monte Carlo for improved history matching and
  uncertainty characterization},'' {\em Spe Journal}, vol.~17, no.~02,
  pp.~418--440, 2012.

\bibitem{adeli2020effect}
E.~Adeli, B.~Rosi{\'c}, H.~G. Matthies, and S.~Reinst{\"a}dler, ``{Effect of
  Load Path on Parameter Identification for Plasticity Models using Bayesian
  Methods},'' in {\em Quantification of Uncertainty: Improving Efficiency and
  Technology}, pp.~1--13, Springer, 2020.

\bibitem{adeli2020comparison}
E.~Adeli, B.~Rosi{\'c}, H.~G. Matthies, S.~Reinst{\"a}dler, and D.~Dinkler,
  ``{Comparison of Bayesian methods on parameter identification for a
  viscoplastic model with damage},'' {\em Metals}, vol.~10, no.~7, p.~876,
  2020.

\bibitem{mirsian2019new}
S.~Mirsian, A.~Khodadadian, M.~Hedayati, A.~Manzour-ol Ajdad,
  R.~Kalantarinejad, and C.~Heitzinger, ``{A new method for selective
  functionalization of silicon nanowire sensors and Bayesian inversion for its
  parameters},'' {\em Biosensors and Bioelectronics}, vol.~142, p.~111527,
  2019.

\bibitem{khodadadian2020bayesiann}
A.~Khodadadian, B.~Stadlbauer, and C.~Heitzinger, ``{Bayesian inversion for
  nanowire field-effect sensors},'' {\em Journal of Computational Electronics},
  vol.~19, no.~1, pp.~147--159, 2020.

\bibitem{noii2020bayesian}
N.~Noii, A.~Khodadadian, and T.~Wick, ``Bayesian inversion for anisotropic
  hydraulic phase-field fracture,'' {\em arXiv preprint arXiv:2007.16038},
  2020.

\bibitem{kuryaeva1997influence}
R.~Kuryaeva and V.~Kirkinskii, ``Influence of high pressure on the refractive
  index and density of tholeiite basalt glass,'' {\em Physics and chemistry of
  minerals}, vol.~25, no.~1, pp.~48--54, 1997.

\bibitem{pariseau2017design}
W.~G. Pariseau, {\em Design analysis in rock mechanics}.
\newblock CRC Press, 2017.

\bibitem{riggleman2010antiplasticization}
R.~A. Riggleman, J.~F. Douglas, and J.~J. de~Pablo, ``Antiplasticization and
  the elastic properties of glass-forming polymer liquids,'' {\em Soft Matter},
  vol.~6, no.~2, pp.~292--304, 2010.

\bibitem{giancoli2016physics}
D.~C. Giancoli, {\em Physics: principles with applications}.
\newblock Boston: Pearson, 2~ed., 2016.

\bibitem{callister2018materials}
W.~D. Callister and D.~G. Rethwisch, {\em Materials science and engineering: an
  introduction}.
\newblock Wiley New York, 9~ed., 2014.

\bibitem{guo2013experimental}
J.~Guo, {\em An experimental and numerical investigation on damage evolution
  and ductile fracture mechanism of aluminum alloy}.
\newblock PhD thesis, PhD dissertation, The University of Tokushima, 2013.

\bibitem{grif1}
D.~Roylance, ``Introduction to fracture mechanics.''
  \url{https://web.mit.edu/course/3/3.11/www/modules/frac.pdf}, 2001.

\bibitem{ambati2016}
M.~Ambati, R.~Kruse, and L.~De~Lorenzis, ``A phase-field model for ductile
  fracture at finite strains and its experimental verification,'' {\em
  Computational Mechanics}, vol.~57, no.~1, pp.~149--167, 2016.

\bibitem{eller2014plasticity}
T.~Eller, L.~Greve, M.~Andres, M.~Medricky, A.~Hatscher, V.~T. Meinders, and
  A.~H. van~den Boogaard, ``Plasticity and fracture modeling of
  quench-hardenable boron steel with tailored properties,'' {\em Journal of
  Materials Processing Technology}, vol.~214, no.~6, pp.~1211--1227, 2014.

\bibitem{li2017high}
Z.~Li, S.~Zhao, H.~Diao, P.~Liaw, and M.~Meyers, ``{High-velocity deformation
  of Al 0.3 CoCrFeNi high-entropy alloy: Remarkable resistance to shear
  failure},'' {\em Scientific reports}, vol.~7, no.~1, pp.~1--8, 2017.

\bibitem{gomatam2006comprehensive}
R.~R. Gomatam and E.~Sancaktar, ``A comprehensive fatigue life predictive model
  for electronically conductive adhesive joints under constant-cycle loading,''
  {\em Journal of Adhesion Science and Technology}, vol.~20, no.~1,
  pp.~87--104, 2006.

\bibitem{marigo2016}
J.-J. Marigo, C.~Maurini, and K.~Pham, ``An overview of the modelling of
  fracture by gradient damage models,'' {\em Meccanica}, vol.~51, no.~12,
  pp.~3107--3128, 2016.

\bibitem{chen2017flaw}
C.~Chen, Z.~Wang, and Z.~Suo, ``Flaw sensitivity of highly stretchable
  materials,'' {\em Extreme Mechanics Letters}, vol.~10, pp.~50--57, 2017.

\bibitem{reese2021using}
S.~Reese, T.~Brepols, M.~Fassin, L.~Poggenpohl, and S.~Wulfinghoff, ``Using
  structural tensors for inelastic material modeling in the finite strain
  regime--a novel approach to anisotropic damage,'' {\em Journal of the
  Mechanics and Physics of Solids}, vol.~146, p.~104174, 2021.

\bibitem{maugin1990}
G.~Maugin, ``Infernal variables and dissipative structures,'' {\em Journal of
  Non-Equilibrium Thermodynamics}, vol.~15, no.~2, pp.~173--192, 1990.

\bibitem{fremond1996}
M.~Fr{\'e}mond and B.~Nedjar, ``Damage, gradient of damage and principle of
  virtual power,'' {\em International Journal of Solids and Structures},
  vol.~33, no.~8, pp.~1083--1103, 1996.

\bibitem{mielke2006}
A.~Mielke, ``A mathematical framework for generalized standard materials in the
  rate-independent case,'' in {\em Multifield Problems in Solid and Fluid
  Mechanics}, pp.~399--428, Springer, 2006.

\bibitem{mielke2015}
A.~Mielke and T.~Roub{\'\i}{\v{c}}ek, ``Rate-independent systems,'' {\em Theory
  and Application (in preparation)}, 2015.

\bibitem{pham2010}
K.~Pham and J.-J. Marigo, ``Approche variationnelle de l'endommagement: I. les
  concepts fondamentaux,'' {\em Comptes Rendus M{\'e}canique}, vol.~338, no.~4,
  pp.~191--198, 2010.

\bibitem{bourdin2000numerical}
B.~Bourdin, G.~Francfort, and J.-J. Marigo, ``Numerical experiments in
  revisited brittle fracture,'' {\em Journal of the Mechanics and Physics of
  Solids}, vol.~48, no.~4, pp.~797--826, 2000.

\bibitem{francfort1998revisiting}
G.~Francfort and J.-J. Marigo, ``Revisiting brittle fracture as an energy
  minimization problem,'' {\em Journal of the Mechanics and Physics of Solids},
  vol.~46, no.~8, pp.~1319--1342, 1998.

\bibitem{teichtmeister2017phase}
S.~Teichtmeister, D.~Kienle, F.~Aldakheel, and M.-A. Keip, ``Phase field
  modeling of fracture in anisotropic brittle solids,'' {\em International
  Journal of Non-Linear Mechanics}, vol.~97, pp.~1--21, 2017.

\bibitem{kuhn2015}
C.~Kuhn, A.~Schl{\"u}ter, and R.~M{\"u}ller, ``On degradation functions in
  phase field fracture models,'' {\em Computational Materials Science},
  vol.~108, pp.~374--384, 2015.

\bibitem{wu2017}
J.-Y. Wu, ``A unified phase-field theory for the mechanics of damage and
  quasi-brittle failure,'' {\em Journal of the Mechanics and Physics of
  Solids}, vol.~103, pp.~72--99, 2017.

\bibitem{wu2018}
J.-Y. Wu, V.~Nguyen, C.~Nguyen, D.~Sutula, S.~Bordas, and S.~Sinaie, ``Phase
  field modeling of fracture,'' {\em Advances in Applied Mechancis: Multi-Scale
  Theory and Computation}, vol.~52, 2018.

\bibitem{gerasimov2019}
T.~Gerasimov and L.~De~Lorenzis, ``On penalization in variational phase-field
  models of brittle fracture,'' {\em Computer Methods in Applied Mechanics and
  Engineering}, vol.~354, pp.~990--1026, 2019.

\bibitem{wheeler2014}
M.~Wheeler, T.~Wick, and W.~Wollner, ``An augmented-lagrangian method for the
  phase-field approach for pressurized fractures,'' {\em Computer Methods in
  Applied Mechanics and Engineering}, vol.~271, pp.~69--85, 2014.

\bibitem{heister2015}
T.~Heister, M.~Wheeler, and T.~Wick, ``A primal-dual active set method and
  predictor-corrector mesh adaptivity for computing fracture propagation using
  a phase-field approach,'' {\em Computer Methods in Applied Mechanics and
  Engineering}, vol.~290, pp.~466--495, 2015.

\bibitem{wambacq2021}
J.~Wambacq, J.~Ulloa, G.~Lombaert, and S.~Fran{\c{c}}ois, ``Interior-point
  methods for the phase-field approach to brittle and ductile fracture,'' {\em
  Computer Methods in Applied Mechanics and Engineering}, vol.~375, p.~113612,
  2021.

\bibitem{mang2020phase}
K.~Mang, T.~Wick, and W.~Wollner, ``A phase-field model for fractures in nearly
  incompressible solids,'' {\em Computational Mechanics}, vol.~65, no.~1,
  pp.~61--78, 2020.

\bibitem{MieWelHof10b}
C.~Miehe, M.~Hofacker, and F.~Welschinger, ``A phase field model for
  rate-independent crack propagation: {R}obust algorithmic implementation based
  on operator splits,'' {\em Computer Methods in Applied Mechanics and
  Engineering}, vol.~199, pp.~2765--2778, 2010.

\bibitem{han1999}
W.~Han and B.~Reddy, {\em Plasticity: mathematical theory and numerical
  analysis}, vol.~9.
\newblock Springer Science \& Business Media, 1999.

\bibitem{ulloa2020}
J.~Ulloa, R.~Alessi, J.~Wambacq, G.~Degrande, and S.~Fran{\c{c}}ois, ``On the
  variational modeling of non-associative plasticity,'' {\em International
  Journal of Solids and Structures}, vol.~217-218, pp.~272--296, 2021.

\bibitem{rockafellar1970}
R.~Rockafellar, {\em Convex analysis}.
\newblock Princeton University Press, 1970.

\bibitem{rappel2020tutorial}
H.~Rappel, L.~A. Beex, J.~S. Hale, L.~Noels, and S.~Bordas, ``{A tutorial on
  Bayesian inference to identify material parameters in solid mechanics},''
  {\em Archives of Computational Methods in Engineering}, vol.~27, no.~2,
  pp.~361--385, 2020.

\bibitem{wang2020uncertainty}
Y.~Wang and D.~L. McDowell, {\em Uncertainty quantification in multiscale
  materials modeling}.
\newblock Woodhead Publishing, 2020.

\bibitem{metropolis1953equation}
N.~Metropolis, A.~W. Rosenbluth, M.~N. Rosenbluth, A.~H. Teller, and E.~Teller,
  ``Equation of state calculations by fast computing machines,'' {\em The
  Journal of Chemical Physics}, vol.~21, no.~6, pp.~1087--1092, 1953.

\bibitem{hastings1970monte}
W.~K. Hastings, ``{Monte Carlo sampling methods using Markov chains and their
  applications},'' {\em Biometrika}, vol.~57, no.~1, pp.~97--109, 1970.

\bibitem{haario2001adaptive}
H.~Haario, E.~Saksman, J.~Tamminen, {\em et~al.}, ``{An adaptive Metropolis
  algorithm},'' {\em Bernoulli}, vol.~7, no.~2, pp.~223--242, 2001.

\bibitem{green2001delayed}
P.~J. Green and A.~Mira, ``{Delayed rejection in reversible jump
  Metropolis--Hastings},'' {\em Biometrika}, vol.~88, no.~4, pp.~1035--1053,
  2001.

\bibitem{zhang2020improving}
J.~Zhang, J.~A. Vrugt, X.~Shi, G.~Lin, L.~Wu, and L.~Zeng, ``{Improving
  Simulation Efficiency of MCMC for Inverse Modeling of Hydrologic Systems with
  a Kalman-Inspired Proposal Distribution},'' {\em Water Resources Research},
  vol.~56, no.~3, pp.~1--24, 2020.

\bibitem{buljac2018calibration}
A.~Buljac, V.-M.~T. Navas, M.~Shakoor, A.~Bouterf, J.~Neggers, M.~Bernacki,
  P.-O. Bouchard, T.~F. Morgeneyer, and F.~Hild, ``{On the calibration of
  elastoplastic parameters at the microscale via X-ray microtomography and
  digital volume correlation for the simulation of ductile damage},'' {\em
  European Journal of Mechanics-A/Solids}, vol.~72, pp.~287--297, 2018.

\bibitem{abbaszadeh2021reduced}
M.~Abbaszadeh, M.~Dehghan, A.~Khodadadian, N.~Noii, C.~Heitzinger, and T.~Wick,
  ``{A reduced-order variational multiscale interpolating element free Galerkin
  technique based on proper orthogonal decomposition for solving Navier--Stokes
  equations coupled with a heat transfer equation: Nonstationary incompressible
  Boussinesq equations},'' {\em Journal of Computational Physics}, vol.~426,
  p.~109875, 2021.

\bibitem{noii2021glductile}
F.~Aldakheel, N.~Noii, T.~Wick, O.~Allix, and P.~Wriggers, ``Multilevel
  global-local techniques for adaptive ductile phase-field fracture,'' {\em
  arXiv preprint arXiv:2103.02377}, 2021.

\bibitem{ZiShen16}
V.~Ziaei-Rad and Y.~Shen, ``Massive parallelization of the phase field
  formulation for crack propagation with time adaptivity,'' {\em Computer
  Methods in Applied Mechanics and Engineering}, vol.~312, pp.~224--253, 2016.

\bibitem{FaMau17}
P.~Farrell and C.~Maurini, ``Linear and nonlinear solvers for variational
  phase-field models of brittle fracture,'' {\em International Journal for
  Numerical Methods in Engineering}, vol.~109, no.~5, pp.~648--667, 2017.

\bibitem{HeiWi18_pamm}
T.~Heister and T.~Wick, ``Parallel solution, adaptivity, computational
  convergence, and open-source code of 2d and 3d pressurized phase-field
  fracture problems,'' {\em PAMM}, vol.~18, no.~1, p.~e201800353, 2018.

\bibitem{KoKr20}
A.~Kopani{\v{c}}{\'a}kov{\'a} and R.~Krause, ``A recursive multilevel trust
  region method with application to fully monolithic phase-field models of
  brittle fracture,'' {\em Computer Methods in Applied Mechanics and
  Engineering}, vol.~360, p.~112720, 2020.

\bibitem{JoLaWi20}
D.~Jodlbauer, U.~Langer, and T.~Wick, ``Matrix-free multigrid solvers for
  phase-field fracture problems,'' {\em Computer Methods in Applied Mechanics
  and Engineering}, vol.~372, p.~113431, 2020.

\bibitem{JoLaWi20_parallel}
D.~Jodlbauer, U.~Langer, and T.~Wick, ``Parallel matrix-free higher-order
  finite element solvers for phase-field fracture problems,'' {\em Mathematical
  and Computational Applications}, vol.~25, no.~3, p.~40, 2020.

\bibitem{graeser2021truncated}
C.~Gr{\"a}ser, D.~Kienle, and O.~Sander, ``Truncated nonsmooth newton multigrid
  for phase-field brittle-fracture problems,'' {\em arXiv preprint
  arXiv:2007.12290}, 2020.

\bibitem{aquino2018coupled}
T.~Aquino and M.~Dentz, ``A coupled time domain random walk approach for
  transport in media characterized by broadly-distributed heterogeneity length
  scales,'' {\em Advances in Water Resources}, vol.~119, pp.~60--69, 2018.

\bibitem{brooks1998general}
S.~P. Brooks and A.~Gelman, ``General methods for monitoring convergence of
  iterative simulations,'' {\em Journal of Computational and Graphical
  Statistics}, vol.~7, no.~4, pp.~434--455, 1998.

\bibitem{gelman1992inference}
A.~Gelman, D.~B. Rubin, {\em et~al.}, ``Inference from iterative simulation
  using multiple sequences,'' {\em Statistical Science}, vol.~7, no.~4,
  pp.~457--472, 1992.

\bibitem{boyce2014sandia}
B.~L. Boyce, S.~L. Kramer, H.~E. Fang, T.~E. Cordova, M.~K. Neilsen, K.~Dion,
  A.~K. Kaczmarowski, E.~Karasz, L.~Xue, A.~J. Gross, {\em et~al.}, ``The
  sandia fracture challenge: blind round robin predictions of ductile
  tearing,'' {\em International Journal of Fracture}, vol.~186, no.~1-2,
  pp.~5--68, 2014.

\bibitem{zhang2014modeling}
T.~Zhang, E.~Fang, P.~Liu, and J.~Lua, ``Modeling and simulation of 2012 sandia
  fracture challenge problem: phantom paired shell for abaqus and plane strain
  core approach,'' {\em International Journal of Fracture}, vol.~186, no.~1-2,
  pp.~117--139, 2014.

\bibitem{DiLiWiTy21}
P.~Diehl, R.~Lipton, T.~Wick, and T.~Mayank, ``A comparative review of
  peridynamics and phase-field models for engineering fracture mechanics,'' Mar
  2021.
\newblock https://engrxiv.org/gty2b/.

\end{thebibliography}
\end{spacing}}
\end{document}